\documentclass[11pt,a4paper]{article}
\usepackage[utf8]{inputenc}

\usepackage{mathrsfs,amsthm,xcolor,verbatim,bbm,amsmath,amsfonts,amssymb,nicefrac,enumitem,hyperref,bm,mathtools,xparse,etoolbox}
\usepackage{algorithm,algpseudocode} 
\usepackage{textcomp} 
\usepackage[margin=1in]{geometry}
\usepackage[capitalise]{cleveref} 

\crefname{enumi}{item}{items}
\crefname{equation}{}{}
\crefname{subsection}{Subsection}{Subsections}
\crefname{figure}{Figure}{Figures}
\crefname{algorithm}{Algorithm}{Algorithms}

\newcommand{\Va}{V}
\newcommand{\Vb}{W}

\newcommand{\supertau}{\tau}
\newcommand{\param}{c}
\newcommand{\sparam}{\varXi}

\newcommand{\dd}{d}
\newcommand{\spp}{\operatorname{Supp}}
\newcommand{\Ysum}{\chi}
\newcommand{\X}{Z}
\newcommand{\olr}{\operatorname{OLR}}

\hypersetup{
    colorlinks,
    linkcolor={blue!80!black},
    citecolor={green},
    urlcolor={blue!80!black}
}

\theoremstyle{plain}
\newtheorem{theorem}{Theorem}[section]
\newtheorem{lemma}[theorem]{Lemma}
\newtheorem{prop}[theorem]{Proposition}
\newtheorem{cor}[theorem]{Corollary}
\newtheorem{framework}[theorem]{Framework}

\newtheorem{remark}[theorem]{Remark}
\newtheorem{definition}[theorem]{Definition}

\DeclareMathAlphabet{\mathpzc}{OT1}{pzc}{m}{it}

\DeclareFontEncoding{LS1}{}{}
\DeclareFontSubstitution{LS1}{stix}{m}{n}
\DeclareMathAlphabet{\mathscr}{LS1}{stixscr}{m}{n}

\newcommand{\E}{\mathbb{E}}
\renewcommand{\P}{\mathbb{P}}
\newcommand{\Q}{\mathbb{Q}}
\newcommand{\R}{\mathbb{R}}
\newcommand{\N}{\mathbb{N}}
\newcommand{\Z}{\mathbb{Z}}

\newcommand{\ssum}{\textstyle\sum}
\newcommand{\ssuml}{\textstyle\sum\limits}

\newcommand{\with}{\curvearrowleft}

\newcommand{\cA}{\mathcal{A}}
\newcommand{\cB}{\mathcal{B}}

\newcommand{\cD}{\mathcal{D}}
\newcommand{\cE}{\mathcal{E}}
\newcommand{\cF}{\mathcal{F}}
\newcommand{\cG}{\mathcal{G}}
\newcommand{\cH}{\mathcal{H}}

\newcommand{\cM}{\mathcal{M}}
\newcommand{\cN}{\mathcal{N}}

\newcommand{\cT}{\mathcal{T}}

\newcommand{\cX}{\mathcal{X}}

\newcommand{\scrD}{\mathscr{D}}

\newcommand{\scrG}{\mathscr{G}}

\newcommand{\scrN}{\mathscr{N}}

\newcommand{\fL}{\mathfrak{L}}
\newcommand{\fM}{\mathbf{M}}
\newcommand{\fN}{\mathfrak{N}}

\newcommand{\fd}{d}
\newcommand{\los}{\ell}
\newcommand{\xdim}{\mathfrak{d}}

\newcommand{\fe}{\mathfrak{e}}

\newcommand{\fl}{\ell}

\newcommand{\ft}{\mathfrak{t}}

\newcommand{\kkk}{\mathscr{k}}
\newcommand{\sss}{\mathscr{s}}
\newcommand{\ttt}{\mathscr{t}}

\renewcommand{\emptyset}{\varnothing}

\DeclarePairedDelimiter{\norm}{\lVert}{\rVert}

\DeclarePairedDelimiter{\rbr}{(}{)}
\DeclarePairedDelimiter{\br}{[}{]}
\DeclarePairedDelimiter{\cu}{\{}{\}}
\DeclarePairedDelimiter{\spro}{\langle}{\rangle}

\renewcommand{\d}{ \mathrm{d}}

\newcommand{\qandq}{\quad \text{and} \quad }
\newcommand{\qqandqq}{\qquad\text{and}\qquad}

\newcommand{\indicator}[1]{\mathbbm{1}_{\smash{#1}}}

\ExplSyntaxOn

\bool_new:N \g_noteobserve

\NewDocumentCommand{\setnote}{}{
  \bool_gset_true:N \g_noteobserve
}

\NewDocumentCommand{\setobserve}{}{
  \bool_gset_false:N \g_noteobserve
}

\NewDocumentCommand{\nobs}{ o }{
  \IfValueT{#1}{
    \str_if_eq:noTF {note} {#1} {
      \bool_gset_true:N \g_noteobserve
    } {
      \str_if_eq:noTF {Note} {#1} {
        \bool_gset_true:N \g_noteobserve
      } {
        \bool_gset_false:N \g_noteobserve
      }
    }
  }
  \bool_if:nTF { \g_noteobserve } {
    \bool_gset_false:N \g_noteobserve
    note
  } {
    \bool_gset_true:N \g_noteobserve
    observe
  }
  \IfValueF{#1}{~}
}

\NewDocumentCommand{\Nobs}{ o }{
  \IfValueT{#1}{
    \str_if_eq:noTF {note} {#1} {
      \bool_gset_true:N \g_noteobserve
    } {
      \str_if_eq:noTF {Note} {#1} {
        \bool_gset_true:N \g_noteobserve
      } {
        \bool_gset_false:N \g_noteobserve
      }
    }
  }
  \bool_if:nTF { \g_noteobserve } {
    \bool_gset_false:N \g_noteobserve
    Note
  } {
    \bool_gset_true:N \g_noteobserve
    Observe
  }
  \IfValueF{#1}{~}
}

\int_new:N \g_furthermore

\NewDocumentCommand{\Moreover}{ o o }{
  \IfValueT{#1}{
    \str_case:nn {#1} {
      {Furthermore} {\int_set:Nn {\g_furthermore} {0}}
      {Moreover} {\int_set:Nn {\g_furthermore} {1}}
      {In~addition} {\int_set:Nn {\g_furthermore} {2}}
      {note} {\bool_gset_true:N \g_noteobserve}
      {observe} {\bool_gset_false:N \g_noteobserve}
    }
    \IfValueT{#2}{
      \str_case:nn {#2} {
        {Furthermore} {\int_set:Nn {\g_furthermore} {0}}
        {Moreover} {\int_set:Nn {\g_furthermore} {1}}
        {In~addition} {\int_set:Nn {\g_furthermore} {2}}
        {note} {\bool_gset_true:N \g_noteobserve}
        {observe} {\bool_gset_false:N \g_noteobserve}
      }
    }
  }
  \int_case:nn { \int_mod:nn {\g_furthermore} {3} } {
    { 0 } { Furthermore,~\nobs that}
    { 1 } { Moreover,~\nobs that}
    { 2 } { In~addition,~\nobs that}
  }
  \int_incr:N \g_furthermore
  \IfValueF{#1}{~}
}

\bool_new:N \g_hencetherefore

\NewDocumentCommand{\hence}{}{
  \bool_if:nTF { \g_hencetherefore } {
    \bool_gset_false:N \g_hencetherefore
    hence~
  } {
    \bool_gset_true:N \g_hencetherefore
    therefore~
  }
}

\NewDocumentCommand{\Hence}{}{
  \bool_if:nTF { \g_hencetherefore } {
    \bool_gset_false:N \g_hencetherefore
    Hence,~we~obtain~
  } {
    \bool_gset_true:N \g_hencetherefore
    Therefore,~we~obtain~
  }
}

\seq_new:N \g_cflist_loaded
\seq_new:N \g_cflist_pending

\NewDocumentCommand{\cfadd}{ m }
{
  \seq_if_in:NnF \g_cflist_loaded { #1 } {
    \seq_if_in:NnF \g_cflist_pending { #1 } {
      \seq_gput_right:Nn \g_cflist_pending { #1 }
    }
  }
}

\NewDocumentCommand{\cfconsiderloaded}{ m }{
  \seq_gput_right:Nn \g_cflist_loaded {#1}
}

\NewDocumentCommand{\cfremove}{ m }
{
  \seq_gremove_all:Nn \g_cflist_pending { #1 }
}

\NewDocumentCommand{\cfload}{ o }
{
  \seq_if_empty:NTF \g_cflist_pending {\unskip} {
    (cf.\ \cref{\seq_use:Nn \g_cflist_pending {,}})\IfValueTF{#1}{#1~}{\unskip}
    \seq_gconcat:NNN \g_cflist_loaded \g_cflist_loaded \g_cflist_pending
    \seq_gclear:N \g_cflist_pending
  }
}

\NewDocumentCommand{\cfclear} {} {
  \seq_gclear:N \g_cflist_loaded
  \seq_gclear:N \g_cflist_pending
}

\NewDocumentCommand{\cfout}{ o }
{
  \seq_if_empty:NTF \g_cflist_pending {\unskip} {
    (cf.\ \cref{\seq_use:Nn \g_cflist_pending {,}})\IfValueTF{#1}{#1~}{\unskip}
    \seq_gclear:N \g_cflist_pending
  }
}

\NewDocumentCommand{\ifnocf} { m } {
  \seq_if_empty:NT \g_cflist_pending { #1 }
}

\ExplSyntaxOff

\NewDocumentEnvironment{cproof}{m}
{\begin{proof}[Proof of \cref{#1}]}%
{\noindent The proof of \cref{#1} is thus complete.
\end{proof}}

\NewDocumentEnvironment{cproof2}{m}
{\begin{proof}[Proof of \cref{#1}]}%
{\noindent This completes the proof of \cref{#1}.
\end{proof}}

\hypersetup{colorlinks=true}

\usepackage{tcolorbox}
\tcbuselibrary{skins, breakable, listings}  
\newcounter{algorithmcounter}  
\renewcommand{\thealgorithmcounter}{\arabic{algorithmcounter}}  
\newtcolorbox[use counter=algorithmcounter]{myalgorithm}[3][]{
    enhanced,
    breakable,
    fonttitle=\bfseries,
    title=Algorithm~\thealgorithmcounter: #2,
    label={#3},  
    label type=algorithmcounter,
    #1,
    colframe=black,
    colback=white,
    coltitle=black,
    colbacktitle=white,
    sharp corners,
    boxrule=0.5pt,
    boxsep=1mm,
    top=1mm,
    bottom=1mm,
    left=0mm,
    right=0mm
}

\newcommand{\mycomment}[1]{\hfill \textcolor{gray}{\textit{\# #1}}}

\crefname{algorithmcounter}{Algorithm}{Algorithms}
\Crefname{algorithmcounter}{Algorithm}{Algorithms}
\crefname{line}{line}{lines}
\Crefname{line}{Line}{Lines}

\makeatletter
\ExplSyntaxOn
\seq_new:N \g_abbrs
\prop_new:N \g_abbr_counts
\tl_new:N \l_abbr_count_tl

\NewDocumentCommand{\abbr}{m m O{#1} m m O{#4}}{
	\expandafter\newcommand\csname#3\endcsname[1][]{
		\seq_if_in:NnTF \g_abbrs {#1} {
			\prop_get:NnN \g_abbr_counts {#1} \l_abbr_count_tl
			\prop_gput:Nnx \g_abbr_counts {#1} {\int_eval:n {\l_abbr_count_tl + 1}}
			\hyperref[#1]{#1}
		} {
			\seq_gput_left:Nn \g_abbrs {#1}
			\prop_gput:Nnn \g_abbr_counts {#1} {1}
			\expandafter\gdef\csname#1@def\endcsname{#2}
			\phantomsection\label{#1}
			\str_if_eq:nnTF{##1}{}{\emph{#2}}{##1}~(\hyperref[#1]{#1})
		}
	}
	\expandafter\newcommand\csname#6\endcsname[1][]{
		\seq_if_in:NnTF \g_abbrs {#1} {
			\prop_get:NnN \g_abbr_counts {#1} \l_abbr_count_tl
			\prop_gput:Nnx \g_abbr_counts {#1} {\int_eval:n {\l_abbr_count_tl + 1}}
			\hyperref[#1]{#4}
		} {
			\expandafter\gdef\csname#1@def\endcsname{#5}
			\seq_gput_left:Nn \g_abbrs {#1}
			\prop_gput:Nnn \g_abbr_counts {#1} {1}
			\phantomsection\label{#1}
			\str_if_eq:nnTF{##1}{}{\emph{#5}}{##1}~(\hyperref[#1]{#4})
		}
	}
}

\ExplSyntaxOff
\makeatother

\abbr{ANN}{artificial neural network}{ANNs}{artificial neural networks}
\abbr{SGD}{stochastic gradient descent}{SGDs}{stochastic gradient descent}
\abbr{GD}{gradient descent}{GDs}{gradient descent}
\abbr{PDE}{partial differential equation}{PDEs}{partial differential equations}
\abbr{DKM}{deep Kolmogorov method}{DKMs}{deep Kolmogorov methods}
\abbr{DGM}{deep Galerkin method}{DGMs}{deep Galerkin methods}
\abbr{PINN}{physics-informed neural network}{PINNs}{physics-informed neural networks}
\abbr{DRM}{deep Ritz method}{DRMs}{deep Ritz methods}
\abbr{GELU}{Gaussian error linear unit}{GELUs}{Gaussian error linear units}
\abbr{BSDE}{backward stochastic differential equation}{BSDEs}{backward stochastic differential equations}
\abbr{SDE}{stochastic differential equation}{SDEs}{stochastic differential equations}
\abbr{ReLU}{rectified linear unit}{ReLUs}{rectified linear units}

\title{Learning rate adaptive stochastic gradient descent optimization methods: 
numerical simulations 
for deep learning methods for partial differential equations 
and convergence analyses}

\author{Steffen Dereich$^{1}$, Arnulf Jentzen$^{2,3}$, and Adrian Riekert$^{4}$\bigskip\\
\small{$^1$ Institute for Mathematical Stochastics, University of M\"unster,}\vspace{-0.1cm}\\
\small{Germany; e-mail: \texttt{steffen.dereich}\textcircled{\texttt{a}}\texttt{uni-muenster.de}}\smallskip\\
\small{$^2$ School of Data Science and Shenzhen Research Institute of Big Data,}\vspace{-0.1cm}\\
\small{The Chinese University of Hong Kong, Shenzhen (CUHK-Shenzhen),}\vspace{-0.1cm}\\
\small{China; e-mail: \texttt{ajentzen}\textcircled{\texttt{a}}\texttt{cuhk.edu.cn}}\smallskip\\
\small{$^3$ Applied Mathematics: Institute for Analysis and Numerics,}\vspace{-0.1cm}\\
\small{University of M\"unster, Germany; e-mail: \texttt{ajentzen}\textcircled{\texttt{a}}\texttt{uni-muenster.de}}\smallskip\\
\small{$^4$ Applied Mathematics: Institute for Analysis and Numerics,}\vspace{-0.1cm}\\
\small{University of M\"unster, Germany; e-mail: \texttt{ariekert}\textcircled{\texttt{a}}\texttt{uni-muenster.de}}}

\date{\today}

\begin{document}

\maketitle

\begin{abstract}
The standard stochastic gradient descent (SGD) optimization method, as well as adaptive methods such as the Adam optimizer fail to converge if the learning rates do not converge to zero (particularly, in the situation of constant learning rates). In practice, human-tuned deterministic learning rate schedules or small constant learning rates are often used, and implementations in machine learning frameworks like {\sc Tensorflow} and {\sc Pytorch} typically employ constant learning rates. We propose a learning-rate-adaptive approach for SGD methods, adjusting the learning rate based on empirical estimates for the objective function values. Specifically, we propose a learning-rate-adaptive variant of the Adam optimizer and implement it for several machine learning problems, including deep learning methods for partial differential equations such as deep Kolmogorov methods, physics-informed neural networks, and deep Ritz methods. We refer to \url{https://github.com/deeplearningmethods/adaptive-learning-rate} for the Python source codes for the numerical simulations in this work. Our results show that the proposed adaptive Adam variant achieves faster reductions of the objective function value compared to Adam with default learning rates. For certain quadratic minimization problems, we rigorously prove that an adaptive SGD variant converges to the global minimizer. This proof uses properties of invariant measures of the SGD dynamics and a generalized convergence analysis for SGD with random predictable learning rates which we develop in this work. 
\end{abstract}

\pagebreak 

\tableofcontents

\pagebreak

\section{Introduction}

It is known that the plain vanilla standard \SGD\ optimization method, 
accelerated \SGD\ methods such as the momentum \SGD\ optimization method, 
adaptive \SGD\ methods such the Adagrad and the RMSprop optimization methods, 
as well as accelerated and adaptive \SGD\ optimization methods such as the Adam optimizer 
fail to converge if the learning rates do not converge to zero 
as, for instance, in the situation of constant learning rates; cf., for example, 
\cite[Section~7.2.2.2]{JentzenKuckuckvonWurstemberger2023arXiv}. 
In numerical simulations SGD optimization methods are sometimes applied with learning rate schedules that are determined through human based tuning procedures, are sometimes applied with fixed deterministic decaying learning rate schedules, and are also sometimes applied with small but constant learning rates. In standard machine learning implementation frameworks such as {\sc Pytorch} (cf. \cite{Paszke2017}) and {\sc Tensorflow} (cf. \cite{Abadi2016Tensorflow}) the default learning rate schedules for SGD optimization methods are constant learning rates (e.g., $\tfrac{1}{1000}$ in the case of SGD in {\sc Pytorch}, $\frac{1}{100}$ in the case of SGD in {\sc Tensorflow}, $\frac{1}{1000}$ in the case of Adam in {\sc Pytorch}, and $\frac{1}{1000}$ in the case of Adam in {\sc Tensorflow}).

In this work we propose (cf.\ \cref{sec:algorithmic_description}), 
study numerically (cf.\ \cref{sec:numerics}), and 
study analytically (cf.\ \cref{sec:random_learning_rates,sec:invariant_measure,sec:test_loss_strictly_increases,sec:convergence_analysis}) 
a \emph{learning-rate-adaptive approach for \SGD\ optimization methods} 
in which the learning rate is adjusted based 
on empirical estimates for the values of the objective function of the optimization problem 
(the function that one intends to minimize in the considered optimization problem). 
In particular, we propose a learning-rate-adaptive variant of the Adam optimizer 
and implement it in {\sc Pytorch} in case of several \ANN\ learning problems, 
particularly, in the context of deep learning approximation methods 
for \PDEs\ such as 
\DKMs\ (cf.\ \cref{sec:DKM}), \PINNs\ (cf.\ \cref{sec:PINNs}), 
and \DRMs\ (cf.\ \cref{sec:DRM}).
In each of the presented learning problems 
the proposed learning-rate-adaptive variant of the Adam optimizer 
reduces the value of the objective function 
significantly faster than the Adam optimizer with the default learning rate. 
If one uses suitable deterministic decreasing learning rate schedules the observed convergence speed can in some cases be faster. 
Nevertheless, these schedules usually require additional effort for hyperparameter tuning for each machine learning problem to be solved.

For a simple class of quadratic minimization problems 
we also rigorously prove in this work that a learning-rate-adaptive variant 
of the standard \SGD\ optimization method converges to the minimizer 
of the considered minimization problem; see \cref{statement_adaptive_convergence} in \cref{sec:convergence_analysis} below.  
A special case of this convergence result is the subject 
of \cref{thm:intro} in this introductory section. 

In \cref{thm:intro} the random variables $X_{n, m }$ describe i.i.d.~data that are used for the \SGD\ process.
Here the mini-batch used in the $n$-th training step consists of the $X_{n, m }$, $m \in \cu{1, 2, \ldots, M }$, 
while the $X_{n, -m}$, $m \in \cu{1, 2, \ldots, \fM }$, are used to compute an independent test error.
We intend to minimize the expected loss function $\E [ \ell ( \theta , X_{1, 1 } )]$, the global minimum of which is $\E [ X_{1, 1 } ]$.

The sequence $(\nu_n)_{n \in \N}$ describes the decreasing pre-determined step-sizes.
The stochastic process $\gamma \colon \N \times \Omega \to \R$ describes the random adaptive learning rates used during the training. 
Our learning rate selection strategy is to compute the average loss with respect to the data points $X_{n, -m}$, $m \in \cu{1, 2, \ldots, \fM }$, in the $n$-th step. If this loss has increased, we choose as the next learning rate $\gamma_{n+1}$ the next element in the sequence $(\nu_k)$, while otherwise we choose $\gamma_{n+1} = \gamma_n$.
This ensures that the learning rates $\gamma_n$ are always decreasing (see \cref{thm_intro_adaptive_learning_rates} for details).
It should be pointed out that this is a simplified version of the learning rate description strategy we consider in our algorithmic descriptions in \cref{sec:algorithmic_description}
and our numerical experiments in \cref{sec:numerics}, respectively.

The conclusion \cref{thm_intro_eq_conclusion} of \cref{thm:intro} reveals that the \SGD\ process $\Theta$, at the random times $N_t$, converges in the $L^1$-sense to the global optimum $\E [ X_{1, 1 } ]$ as $t \to \infty $. 
We now present the precise statement of \cref{thm:intro}.

\begin{theorem}[Convergence of \SGD\ with adaptive learning rates]
\label{thm:intro}
Let 
$ \dd \in \N $, $ a \in \R $, $ b \in (a,\infty) $, 
let 
$ ( \Omega, \mathcal{F}, \P ) $ be a probability space, 
let $ X_{ n, m } \colon \Omega \to [a,b]^{ \fd } $, $ (n, m) \in \Z^2 $, be i.i.d.\ random variables, 
assume\footnote{Note that for all $ n \in \N $, $ x = ( x_1, \dots, x_n ) \in \R^n $ 
it holds that $ \| x \| = ( \sum_{ m = 1 }^n | x_m |^2 )^{ \nicefrac{ 1 }{ 2 } } $.} 
for all $ x \in [a,b]^d $, $ \varepsilon \in (0,1) $ that 
$
  \P( \| X_{ 0, 0 } - x \| < \varepsilon ) > 0
$,
let $ \fl = ( \fl( \theta, x ) )_{ ( \theta, x ) \in \R^{ \fd } \times \R^{ \fd } } \colon \R^{ \fd } \times \R^{ \fd } \to \R $ 
satisfy for all $ \theta, x \in \R^{ \fd } $ that 
\begin{equation}
  \fl ( \theta, x ) 
  =
  \| \theta - x \|^2
  ,
\end{equation}
let $ M, \fM \in \N $, 
let 
$ \Theta \colon \N_0 \times \Omega \to \R^{ \fd } $
and 
$ \gamma \colon \N \times \Omega \to \R $
be stochastic processes 
which satisfy for all $ n \in \N $ that 
\begin{equation}
  \Theta_n 
  = \Theta_{ n - 1 } 
  - 
  \frac{ \gamma_n }{ M }
  \left[ 
    \sum_{ m = 1 }^M
    ( \nabla_{ \theta } \fl )( \Theta_{ n - 1 }, X_{ n, m } ) 
  \right]
  ,
\end{equation}
assume that $ ( X_{ n, m } )_{ ( n, m ) \in \Z^2 } $ and $ \Theta_0 $ are independent, 
assume $ \E[ \| \Theta_0 \|^2 ] < \infty $, 
let $ ( \nu_n )_{ n \in \N } \subseteq (0,1) $ be strictly decreasing,
assume 
$
  \sum_{ n = 1 }^{ \infty } \nu_n = \infty 
$
and 
$
  \lim_{ n \to \infty } \nu_n = 0
$, 
assume for all $ n \in \N $ that 
$
  \gamma_1 = \nu_1
$
and 
\begin{equation}
\label{thm_intro_adaptive_learning_rates}
  \gamma_{ n + 1 } = 
  \begin{cases}
    \gamma_n
  &
    \colon
      \sum_{ m = 1 }^{ \fM }
      \fl ( \Theta_n, X_{ n, - m } ) 
    \leq 
      \sum_{ m = 1 }^{ \fM }
      \fl ( \Theta_{ n - 1 }, X_{ n, - m } ) 
  \\
    \max\bigl(
      ( 0, \gamma_n )
      \cap 
      [
        \cup_{
          k = 1
        }^{ \infty }
        \{ \nu_k \} 
      ]
    \bigr)
  &
    \colon
      \sum_{ m = 1 }^{ \fM }
      \fl( \Theta_n, X_{ n, - m } ) 
    >
      \sum_{ m = 1 }^{ \fM }
      \fl( \Theta_{ n - 1 }, X_{ n, - m } ) 
      ,
  \end{cases}
\end{equation}
and for every $ t \in [0,\infty) $ let $ N_t \colon \Omega \to \R $ satisfy 
\begin{equation}
\textstyle 
  N_t = 
  \inf\bigl\{ 
    n \in \N_0 \colon 
    \sum_{ k = 1 }^n \gamma_k \geq t
  \bigr\}
  .
\end{equation}
Then 
\begin{equation}
	\label{thm_intro_eq_conclusion}
  \lim_{ t \to \infty }
  \E\bigl[
    \|  
      \Theta_{ N_t } - \E[ X_{ 1, 1 } ]
    \|
  \bigr]
  = 0
  .
\end{equation}
\end{theorem}

\cref{thm:intro} is an immediate consequence of \cref{statement_adaptive_convergence}.
Our convergence proofs of \cref{thm:intro} and \cref{statement_adaptive_convergence}, respectively, 
are based on an analysis of the laws of invariant measures 
of the \SGD\ optimization method (cf.\ \cref{sec:invariant_measure} below)
as well as on a more general convergence analysis 
for the \SGD\ optimization method with random but predictable learning rates 
which we develop in this work (cf.\ \cref{sec:random_learning_rates}).


We next provide a brief overview of related works in the literature regarding \GD\ and \SGD\ optimization methods with adaptive learning rates.
A classical idea to find a good learning rate for \SGD\ is to use a line-search strategy which chooses learning rates that lead to a sufficient decrease of the loss in each training step \cite{VaswaniMishkin2019,YuChen1995}.
Another approach is to determine a learning rate schedule via an additional optimization procedure using gradient methods or methods from reinforcement learning \cite{Andrychowicz2016,LiMalik2016,Song2023,XuDai2019}.
For further ideas regarding deterministic learning rate schedules with non-constant step sizes we refer, for example, to \cite{Loshchilov2016,Smith2017}.
Zeroth-order adaptive gradient methods, where the learning rate is increased or decreased during training depending of the past behavior of the loss, have been proposed, for instance, in \cite{LiRen2021,Mukherjee2019,Shim2023}.
In this work we pursue a similar approach for more general gradient methods than the standard \SGD\ optimizer.
An automatic \GD\ method for \ANN\ training where the learning rate is computed depending on the network architecture and parameters has been proposed in \cite{Bernstein2023}.
In \cite{Horvath2022,KoehneKreis2023,Wang2021,Wojcik2019,YedidaSaha2020} additional approaches which estimate a good learning rate for \SGD\ adaptively during training have been introduced.
A different approach which computes the learning rate schedule as the solution of an optimal control problem has been introduced in \cite{LiTaiE2017}.

A popular methodology relies on computing dynamic learning rates based on an adaptive estimation of the second moments of the past gradients, and algorithms of this type include in particular the RMSProp (see \cite{Hinton_RMSProp}), Adagrad (see \cite{DuchiHazan2011}), and Adam (see \cite{KingmaBa2014}) optimizers.
Modifications of the Adam optimizer have been proposed, for example, in \cite{LoshchilovHutter2017,Luo2019,ZaheerReddi2018}.
For further related methods we refer, for instance, to \cite{SchaulSixin2013,Tong2022,WuWard2018,ZhuangTang2020}.

Regarding theoretical convergence results for gradient methods with adaptive learning rates we refer, in particular, to \cite{Crawshaw2022,DeMukherjee2018,Defossez2020,FawTziotis2022,ReddiKale2019,ZhangChen2022,ZouShen2019} for convergence analyses for Adam and similar optimization methods for general non-convex objective functions under a Lipschitz continuity assumption on the gradients. 
Such results usually show convergence of the gradient norms to zero as the training time increases to infinity.
More generally, convergence results for Adam and related methods under weaker local smoothness conditions have been obtained in \cite{Barakat2021,HongLin2024,Li2023,WangZhang2023}.
Under similar assumptions, in \cite{ZhangFang2020} the convergence behavior of \SGD\ with clipped learning rates has been analyzed.
Similar results for non-smooth objective functions can be found in \cite{XiaoHu2024}.
In a general setting, which is most similar to this work, stochastic gradient methods with random adaptive learning rates have been investigated, for example, in \cite{ChenLiuSunHong2018,GodichonBaggioni2023,LiuJiang2019}.

Finally, for further overviews on gradient-based optimization methods and theoretical convergence results 
we refer the interested reader, for example, to 
the survey articles and monographs \cite{JentzenKuckuckvonWurstemberger2023arXiv,Nesterov2014,Ruder2017overview} 
and the references mentioned therein.

The remainder of this work is organized as follows. 
In \cref{sec:algorithmic_description} we describe 
in details the learning-rate-adaptive algorithmic procedures 
that we propose, first, in the situation where the underlying 
SGD method is the plain vanilla SGD method (see \cref{sec:SGD})
and, thereafter, in the more general situation 
where the underlying SGD method belongs 
to a class of general SGD methods which includes, for example, 
the Adam optimizer (see \cref{sec:Adam}). 
In \cref{sec:numerics} we test the proposed strategy 
numerically in the case of 
several \ANN\ learning problems, 
particularly, in the context of deep learning approximation methods 
for \PDEs\ such as \DKMs\ (cf.\ \cref{sec:DKM}), \PINNs\ (cf.\ \cref{sec:PINNs}), 
and \DRMs\ (cf.\ \cref{sec:DRM}). 
In \cref{sec:random_learning_rates} 
we present a general convergence analysis 
for the \SGD\ optimization method with random but predictable learning rates 
which converge almost surely to zero. 
In \cref{sec:test_loss_strictly_increases} 
we prove for a class of simple quadratic optimization problems, 
loosely speaking, that the adaptive learning rates of the SGD method from \cref{thm:intro} 
(see \cref{thm_intro_adaptive_learning_rates} above) 
do converge almost surely to zero under the assumption that the test loss strictly decreases 
with positive probability under the associated invariant measure of the SGD method. 
In \cref{sec:invariant_measure} we show, roughly speaking, 
that this assumption that the test loss strictly descreases 
with positive probability under the invariant measure 
is fulfilled (cf.\ \cref{prop:positive:prob} for details) 
and thus that the adaptive learning rates of the SGD method 
from \cref{thm:intro} (cf.\ \cref{thm_intro_adaptive_learning_rates}) 
do converge almost surely to zero (cf.\ \cref{cor:SGD_constant_nondeg} for details). 
In \cref{sec:convergence_analysis} we combine 
\cref{cor:SGD_constant_nondeg} from \cref{sec:invariant_measure} 
with the findings from \cref{sec:random_learning_rates}
to establish \cref{statement_adaptive_convergence}. 
\cref{thm:intro} above is an immediate consequence of \cref{statement_adaptive_convergence}.

\section{Algorithmic descriptions}
\label{sec:algorithmic_description}

In this section we describe our approach of learning-rate-adaptive \SGD\ methods.  
In \cref{sec:SGD} we specify learning-rate-adaptive modifications of 
the plain vanilla \SGD\ optimization method and in \cref{sec:Adam} we combine 
the learning-rate-adaptive approach with more sophisticated \SGD\ optimization methods 
such as the RMSProp (see \cite{Hinton_RMSProp}), 
the Adagrad (see \cite{DuchiHazan2011}), and 
the Adam (see \cite{KingmaBa2014}) optimizers.

\subsection{Learning-rate-adaptive plain-vanilla SGD methods}
\label{sec:SGD}

\begin{framework}[Algorithmic description of a learning-rate-adaptive \SGD\ method]
\label{framework:adaptive_SGD}
Let 
$ \fd, \xdim \in \N $, $ a \in \R $, $ b \in (a,\infty) $, 
let $ ( \Omega, \mathcal{F}, \P ) $ be a probability space, 
let $X_{n,m}^{p, \lambda} \colon \Omega \to [a,b]^{ \xdim } $, 
$ (n, m, p, \lambda) \in \Z^3 \times \R $, be i.i.d.\ random variables, 
let 
$ \los = ( \los( \theta, x ) )_{ ( \theta, x ) \in \R^{ \fd } \times \R^{ \xdim } } \colon \R^{ \fd } \times \R^{ \xdim } \to \R $ 
be measurable, 
assume for all $ x \in \R^{ \xdim } $ 
that $ \R^{ \fd } \ni \theta \mapsto \los( \theta, x ) \in \R $ is differentiable,
let $ M, N, \sss, \ttt \in \N $, 
for every random variable $ \xi \colon \Omega \to \R^\fd $ 
and every $ p \in \Z $, $ \lambda \in \R $  
let $ \Phi^{ \xi, p, \lambda } \colon \N_0 \times \Omega \to \R^{ \fd } $
satisfy for all $ n \in \N $ that 
\begin{equation}
  \Phi^{ \xi, p, \lambda }_0 = \xi 
\qquad 
\text{and}
\qquad 
  \Phi^{\xi, p, \lambda}_n 
  = 
  \Phi ^{\xi, p, \lambda}_{n-1}
  - 
  \frac{ \lambda }{ M } 
  \left[ 
    \sum_{ m = 1 }^M 
    ( \nabla_{ \theta } \los )( 
      \Phi^{ \xi, p, \lambda }_{ n - 1 }, X_{ n, m }^{ p, \lambda } 
    ) 
  \right] 
  ,
\end{equation}
for every $ \lambda \in \R $ 
let $ h( \lambda ) \subseteq \R $ 
be a finite set, 
let 
$ \Theta \colon \N_0 \times \Omega \to \R^{ \fd } $, 
$ \gamma \colon \N_0 \times \Omega \to \R $, 
$ \alpha \colon \N_0 \times \Omega \to \R $, 
and 
$ \tau \colon \N_0 \times \Omega \to \N_0 $
satisfy for all $ n \in \N $ that 
\begin{equation}
  \alpha_0 = \tau_0 = 0 , 
  \qquad 
  \Theta_n 
  = \Theta_{ n - 1 } 
  - 
  \frac{ \gamma_n }{ M }
  \left[ 
    \sum_{ m = 1 }^M
    ( \nabla_{ \theta } \los )( \Theta_{ n - 1 }, X_{ n, m }^{0 , 0} ) 
  \right]
  ,
\end{equation}
\begin{equation}
  \gamma_n 
  \in 
  \begin{cases}
    \operatorname{argmin}_{ 
      \lambda \in h( \gamma_{ n - 1 } ) 
    }
    \bigl[
      \sum_{ m = 1 }^N 
      \los( 
        \Phi_{ \sss }^{ \Theta_{ n - 1 }, n, \lambda }, X_{ 0, m }^{ n, \lambda } 
      )
    \bigr]
  &
    \colon
    \tau_{ n - 1 } = \alpha_{ n - 1 } = 0
  \\
    \{ \gamma_{ n - 1 } \} 
  &
    \colon
    \text{else}
  ,
  \end{cases}
\end{equation}
\begin{equation}
  \alpha_n = 
  \begin{cases}
    1 
  & 
    \colon \ssum_{ m = 1 }^M 
    \los( \Theta_{ n - 1 }, X_{ n, m }^{ 0, 0 } ) 
    = \min_{ k \in \cu{ 1, 2, \dots, n } } 
    \bigl[
      \sum_{ m = 1 }^M 
      \los( \Theta_{ k - 1 }, X_{ k, m }^{ 0, 0 } ) 
    \bigr]
\\
    0 
  & 
    \colon \text{else}
    ,
  \end{cases}
\end{equation}
and 
\begin{equation}
  \tau_n = 
  \begin{cases}
    0 
  & 
    \colon \alpha_n = 1
  \\
    \tau_{n-1} + 1 
  & 
    \colon ( \alpha_n = 0 ) \wedge ( \tau_{ n - 1 } < \ttt - 1 ) 
  \\
    0 
  & 
    \colon ( \alpha_n = 0 ) \wedge ( \tau_{ n - 1 } = \ttt - 1 ) .
  \end{cases}
\end{equation}
\end{framework}

\begin{remark}
In the pseudocodes in \cref{alg:find_lr} 
and \cref{alg:sgd_dyn} below and 
in our implementations in \cref{sec:numerics} 
we have that 
the sets $ h( \lambda ) $, $ \lambda \in \R $, 
in \cref{framework:adaptive_SGD} 
satisfy that 
there exist 
$ \kkk \in \N $, 
$ \eta \in (1, \infty) $ 
such that for all $ \gamma \in (0, \infty) $ 
it holds that 
\begin{equation}
  h( \gamma ) 
  = 
  \bigl(
  \cup_{ j = 1 }^{ \kkk }
  \bigl\{ 
    \eta^{j - (\kkk+1)/2} 
    \gamma 
  \bigr\} 
  \bigr)
  .
\end{equation}
\end{remark}

In the descriptions in \cref{alg:find_lr,alg:sgd_dyn,alg:find_lr_adam,alg:adam_dyn} we consider the problem of minimizing 
a stochastic loss function 
$ L \colon \R^\fd \to \R $
defined by
\begin{equation}
\label{eq:optimization_problem_for_algorithms}
  L( \theta ) = \E[ \los(\theta , X )] 
\end{equation}
for all $ \theta \in \R^{ \fd } $
where $\fd, \xdim \in \N$ are natural numbers, 
where 
$ \los \colon \R^\fd \times \R^\xdim \to \R$ is assumed to be sufficiently regular, 
and where 
$ X \colon \Omega \to \R^\xdim $ is a random variable with distribution $ \mu = \P_X $ 
on a probability space $ ( \Omega, \cF, \P ) $.

In the first pseudocode in \cref{alg:find_lr} below 
we describe how we find a good learning rate by doing a few training steps 
for different values of the learning rate and selecting the best one from a geometric grid.
In \cref{alg:find_lr} 
\begin{itemize}
\item 
we note that $ M \in \N $ is the batch size for the \SGD\ method, 
\item 
we note that $ N \in \N $ is the size of the independent datasets used for testing,
\item 
we note that $ \sss \in \N $ is the number of testing \SGD\ steps, 
\item 
we note that $ \kkk \in \N $ is the number of trial learning rate values among which we select the best,
\item 
we note that $ \gamma_0 \in (0 , \infty) $ is the initial/default value of the learning rate which we set as the center of the search grid,
\item 
we note that $ \eta \in (1 , \infty )$ is a factor we use for the grid search, 
and 
\item 
we note that $ \theta_0 \in \R^{ \dd } $ is the initial/current state of the \SGD\ method. 
\end{itemize}
Typical values we tested 
in our numerical simulations in \cref{sec:numerics} 
are 
$ \sss \in \cu{ 30, 50} $, 
$ \kkk \in \cu{ 3, 5 } $, 
$ \eta \in \cu{ 3, 4, 5} $.

\begin{myalgorithm}{Find optimal learning rate for \SGD}{alg:find_lr}
	
	\begin{algorithmic}[1] 
		\Statex\textbf{Setting:} The optimization problem in \cref{eq:optimization_problem_for_algorithms}
		\Statex\textbf{Input:} 
		 $ M, N, \sss, \kkk \in \N $,
		 $ \gamma_0 \in (0 , \infty ) $,
		 $ \eta \in ( 1, \infty ) $,
		 $ \theta_0 \in \R^{ \dd } $

		\Statex\textbf{Output:} 
		Learning rate $\gamma _* = \olr( M, N, \sss, \kkk, \gamma_0, \eta, \theta_0 ) \in (0, \infty) $

		\vspace{-2mm} 
		\noindent\hspace*{-8.1mm}\rule{\dimexpr\linewidth+9.24mm}{0.4pt}  
		\vspace{-3mm} 

		\State $ l_{opt} \gets \infty $
		\State $ \gamma_* \gets \gamma_0 $
		\For{$ \gamma \in \cu{ \gamma_0 \cdot \eta^{ j - ( \kkk +1) / 2 } \colon j \in \cu{ 1, 2, \dots, \kkk } } $}  \mycomment{Loop over geometric grid}
		\State $ \theta \gets \theta_0 $
		\For{$ i \in \cu{ 1, 2, \ldots, \sss} $} \mycomment{Do $\sss$ \SGD\ steps with learning rate $\gamma$}
		\State $X_1, X_2, \ldots, X_M  \gets$ i.i.d.~samples of $\mu$
		\State $\theta \gets \theta - \gamma M^{-1} \sum_{m=1}^M (\nabla_\theta \los ) ( \theta , X_m)$
		\EndFor
		\State $Y_1, Y_2, \ldots, Y_N  \gets$ i.i.d.~samples of $\mu$
		\State $l \gets N^{-1} \sum_{n=1}^N \los(\theta , Y_n)$ \mycomment{Compute average loss over $N$ independent data points}
		\If{$l < l_{opt}$} \mycomment{Update learning rate if new loss is smaller}
		\State $\gamma_* \gets \gamma$
		\State $l_{opt} \gets l$
		\EndIf
		\EndFor
		
		\noindent
		\Return $\gamma _*$ \mycomment{Return learning rate $\gamma_*$ with smallest loss}
	\end{algorithmic}
\end{myalgorithm}

In the pseudocode in \cref{alg:sgd_dyn} below we describe how we employ the learning rate selection method 
in \cref{alg:find_lr} above for the \SGD\ training procedure.

In addition to the parameters introduced in \cref{alg:find_lr},
the natural number $T \in \N$ in \cref{alg:sgd_dyn} is the total number of \SGD\ training steps.
The natural number $\ttt \in \N$ is a tolerance parameter in the sense that we select a new learning rate whenever the training loss has not improved for the last $\ttt $ steps.
In particular, the loss is allowed to increase for more than one step.
This is a more general setting compared to our theoretical results, where we assume for simplicity that $\ttt = 1$, i.e., a new learning rate is chosen in every step where the loss increases.

\begin{myalgorithm}{\SGD\ with dynamic learning rates}{alg:sgd_dyn}
\begin{algorithmic}[1]
    \Statex\textbf{Setting:} The optimization problem in \cref{eq:optimization_problem_for_algorithms}
	\Statex\textbf{Input:} 
	$ M, N, \sss, \ttt, \kkk, T \in \N$, 
	$\gamma_0 \in (0 , \infty )$,
	$\eta \in (1 , \infty )$,
	$\theta_0 \in \R^\fd$
	
	\Statex\textbf{Output:} 
	Approximate critical point $\theta \in \R^\fd $
	
	\vspace{-2mm} 
	\noindent\hspace*{-8.1mm}\rule{\dimexpr\linewidth+9.24mm}{0.4pt}  
	\vspace{-3mm} 
\State $\theta \gets \theta_0$
\State $l_{opt} \gets \infty$
\State $\gamma \gets \olr( M, N, \sss, \kkk, \gamma_0, \eta, \theta)$ \mycomment{Determine initial learning rate}
\State $c \gets 0$
\For{$i \in \cu{1, 2, \ldots, T}$}
			\State  $X_1, X_2, \ldots, X_M \gets$ i.i.d.~samples of $\mu$
			\State $l \gets M^{-1} \sum_{m=1}^M \los ( \theta , X_m)$
			\State $\theta \gets \theta - \gamma M^{-1} \sum_{m=1}^M (\nabla_\theta \los ) ( \theta , X_m)$ \mycomment{\SGD\ step}
			\If{$l < l_{opt}$}
				\State $l_{opt} \gets l $
				\State $c \gets 0$
			\Else
				\State $c \gets c+1$ \mycomment{Counts steps since last improvement}
			\EndIf
			\If{$c = \ttt $} \mycomment{If no improvement occurred for $\ttt $ steps}
				\State $\gamma \gets \olr (k, s, M, N, \gamma, \eta, \theta)$ \mycomment{Determine new learning rate}
				\State $c \gets 0$
			\EndIf
		\EndFor
		
		\noindent 
		\Return $\theta$
\end{algorithmic}
\end{myalgorithm}

\subsection{Learning-rate-adaptive general SGD methods}
\label{sec:Adam}

We next introduce in \cref{framework:adaptive_Adam} 
a generalization of \cref{framework:adaptive_SGD} 
which includes more general optimization methods 
than \SGD\ including, in particular, the Adam optimizer.
Here the functions $\psi_n$, $n \in \N$, determine the descent direction in the $n$-th step depending on all the previous mini-batch gradients with respect to the data $X_{i, m}$, $i \in \cu{1, 2, \ldots, n-1}$.
Essentially all commonly used \SGD\ type optimization methods such as momentum \SGD, RMSprop, Adam, etc.~can be described this way.

\begin{framework}[Algorithmic description of a learning-rate-adaptive general \SGD\ method]
\label{framework:adaptive_Adam}
Let 
$ \fd, \xdim \in \N $, $ a \in \R $, $ b \in (a,\infty) $, 
let $ ( \Omega, \mathcal{F}, \P ) $ be a probability space, 
let $X_{n,m}^{p, \lambda} \colon \Omega \to [a,b]^{ \xdim } $, 
$ (n, m, p, \lambda) \in \Z^3 \times \R $, be i.i.d.\ random variables, 
let 
$ \los = ( \los( \theta, x ) )_{ ( \theta, x ) \in \R^{ \fd } \times \R^{ \xdim } } \colon \R^{ \fd } \times \R^{ \xdim } \to \R $ 
be measurable, 
assume for all $ x \in \R^{ \xdim } $ 
that $ \R^{ \fd } \ni \theta \mapsto \los( \theta, x ) \in \R $ is differentiable,
let $ M, N, \sss, \ttt \in \N $, 
for every $ n \in \N $ 
let 
$
  \psi_n \colon ( \R^{ \dd } )^n \to \R^{ \dd }
$
be measurable,
for every random variable $ \xi \colon \Omega \to \R^\fd $ 
and every $ p \in \Z $, $ \lambda \in \R $  
let $ \Phi^{ \xi, p, \lambda } \colon \N_0 \times \Omega \to \R^{ \fd } $
satisfy for all $ n \in \N $ that 
$
  \Phi^{ \xi, p, \lambda }_0 = \xi 
$
and 
\begin{equation}
  \Phi^{ \xi, p, \lambda }_n 
  = 
  \Phi^{ \xi, p, \lambda }_{n-1}
  - 
  \lambda 
  \,
  \psi_n\biggl( 
    \frac{ 1 }{ M } 
    \left[ 
\textstyle 
      \sum\limits_{ m = 1 }^M 
\displaystyle 
      ( \nabla_{ \theta } \los )( 
        \Phi^{ \xi, p, \lambda }_0, X_{ 1, m }^{ p, \lambda } 
      ) 
    \right]
    ,
    \dots 
    ,
    \frac{ 1 }{ M } 
    \left[ 
\textstyle 
      \sum\limits_{ m = 1 }^M 
      ( \nabla_{ \theta } \los )( 
        \Phi^{ \xi, p, \lambda }_{ n - 1 }, X_{ n, m }^{ p, \lambda } 
      ) 
    \right]
  \biggr)
  ,
\end{equation}
for every $ \lambda \in \R $ 
let $ h( \lambda ) \subseteq \R $ 
be a finite set, 
let
$ \Theta \colon \N_0 \times \Omega \to \R^{ \fd } $, 
$ \gamma \colon \N_0 \times \Omega \to \R $, 
$ \alpha \colon \N_0 \times \Omega \to \R $, 
and 
$ \tau \colon \N_0 \times \Omega \to \N_0 $
satisfy for all $ n \in \N $ that 
$\delta_0 = \gamma_0$,
\begin{equation}
  \Theta_n 
  = 
  \Theta_{ n - 1 } 
  - 
  \gamma_n 
  \,
  \psi_n\biggl(
    \frac{ 1 }{ M }
    \left[ 
      \textstyle 
      \sum\limits_{ m = 1 }^M
      \displaystyle 
      ( \nabla_{ \theta } \los )( \Theta_0, X_{ 1, m }^{0 , 0} ) 
    \right]
    ,
    \dots 
    ,
    \frac{ 1 }{ M }
    \left[ 
      \textstyle 
      \sum\limits_{ m = 1 }^M
      \displaystyle 
      ( \nabla_{ \theta } \los )( \Theta_{ n - 1 }, X_{ n, m }^{0 , 0} ) 
    \right]
  \biggr)
  ,
\end{equation}
\begin{equation}
  \gamma_n 
  \in 
  \begin{cases}
    \operatorname{argmin}_{ 
      \lambda \in h( \gamma_{ n - 1 } ) 
    }
    \bigl[
      \sum_{ m = 1 }^N 
      \los( 
        \Phi_{ \sss }^{ \Theta_{ n - 1 }, n, \lambda }, X_{ 0, m }^{ n, \lambda } 
      )
    \bigr]
  &
    \colon
    \tau_{ n - 1 } = \alpha_{ n - 1 } = 0 \\
    \{ \gamma_{ n - 1 } \} 
  &
    \colon
    \text{else}
  ,
  \end{cases}
\end{equation}
\begin{equation}
  \alpha_n = 
  \begin{cases}
    1 
  & 
    \colon \ssum_{ m = 1 }^M 
    \los( \Theta_{ n - 1 }, X_{ n, m }^{ 0, 0 } ) 
    = \min_{ k \in \cu{ 1, 2, \dots, n } } 
    \bigl[
      \sum_{ m = 1 }^M 
      \los( \Theta_{ k - 1 }, X_{ k, m }^{ 0, 0 } ) 
    \bigr]
\\
    0 
  & 
    \colon \text{else}
    ,
  \end{cases}
\end{equation}
\begin{equation}
  \alpha_0 = \tau_0 = 0 , 
  \qqandqq 
  \tau_n = 
  \begin{cases}
    0 
  & 
    \colon \alpha_n = 1
  \\
    \tau_{n-1} + 1 
  & 
    \colon ( \alpha_n = 0 ) \wedge ( \tau_{ n - 1 } < \ttt - 1 ) 
  \\
    0 
  & 
    \colon ( \alpha_n = 0 ) \wedge ( \tau_{ n - 1 } = \ttt - 1 ) .
  \end{cases}
\end{equation}
\end{framework}

We now describe in \cref{alg:find_lr_adam} a variant of the learning rate selection method in \cref{alg:find_lr} which relies on the Adam optimizer instead of the standard \SGD\ method.
In \cref{alg:find_lr_adam} the parameters $\beta_1, \beta_2 \in [0, 1 ]$ describe the first and second order momentum parameters of the Adam algorithm while $\epsilon \in (0 , \infty )$ is the regularization parameter.
The remaining hyperparameters are analogous to \cref{alg:find_lr}.

\begin{myalgorithm}{Find optimal learning rate for Adam}{alg:find_lr_adam}
	\begin{algorithmic}[1]
		\Statex\textbf{Setting:} The optimization problem in \cref{eq:optimization_problem_for_algorithms}
		\Statex\textbf{Input:} 
		$ M, N, \sss, \kkk \in \N $,
		$ \gamma_0, \epsilon \in (0 , \infty ) $,
		$ \eta \in ( 1, \infty ) $,
		$\beta_1, \beta_2 \in [0 , 1 ]$,
		$ \theta_0 \in \R^{ \dd } $
		
		\Statex\textbf{Output:} 
		Learning rate $\gamma _* = \olr _A ( M, N, \sss, \kkk, \gamma_0, \eta, \theta_0 ) \in (0, \infty) $
		\vspace{-2mm} 
		\noindent\hspace*{-8.1mm}\rule{\dimexpr\linewidth+9.24mm}{0.4pt}  
		\vspace{-3mm} 
		
		\State $ l_{opt} \gets \infty $
		\State $ \gamma_* \gets \gamma_0 $
		\State $m \gets 0$
		\State $v \gets 0 $
		\For{$ \gamma \in \cu{ \gamma_0 \cdot \eta^{ j - ( \kkk +1) / 2 } \colon j \in \cu{ 1, 2, \dots, \kkk } } $} 
		\State $ \theta \gets \theta_0 $
		\For{$ i \in \cu{ 1, 2, \ldots, \sss} $} \mycomment{Do $\sss$ Adam steps with learning rate $\gamma$}
		\State $X_1, X_2, \ldots, X_M  \gets$ i.i.d.~samples of $\mu$
		\State $g \gets M^{-1} \sum_{j=1}^M (\nabla_\theta \los ) ( \theta , X_j )$
		\State $m \gets \beta_1 m + (1 - \beta_1) g$
		\State $v \gets \beta_2 v + (1 - \beta_2 ) g^2$ \mycomment{Square $g^2$ is understood componentwise}
		\State $m \gets m / (1 - \beta_1 ^i )$
		\State $v \gets v / (1 - \beta_2 ^i )$
		\State $\theta \gets \theta - \gamma m / (\sqrt{v} + \epsilon)$ \mycomment{Root $\sqrt{v}$ is understood componentwise}
		\EndFor
		\State $Y_1, Y_2, \ldots, Y_N  \gets$ i.i.d.~samples of $\mu$
		\State $l \gets N^{-1} \sum_{n=1}^N \los(\theta , Y_n)$ \mycomment{Compute average loss over $N$ independent data points}
		\If{$l < l_{opt}$} \mycomment{Update learning rate if new loss is smaller}
		\State $\gamma_* \gets \gamma$
		\State $l_{opt} \gets l$
		\EndIf
		\EndFor
		
		\noindent
		\Return $\gamma _*$ \mycomment{Return learning rate $\gamma_*$ with smallest loss}
	\end{algorithmic}
\end{myalgorithm}

In the pseudocode in \cref{alg:adam_dyn} below we describe how we employ the learning rate selection method in \cref{alg:find_lr_adam} above for the Adam training procedure.

\begin{myalgorithm}{Adam with dynamic learning rates}{alg:adam_dyn}
	\begin{algorithmic}[1]
		\Statex\textbf{Setting:} The optimization problem in \cref{eq:optimization_problem_for_algorithms}
		\Statex\textbf{Input:} 
		$ M, N, \sss, \ttt, \kkk, T \in \N$,
		$\gamma _0 , \epsilon \in (0 , \infty )$,
		$\eta \in (1 , \infty )$,
		$\beta_1, \beta_2 \in [0 , 1 ]$,
		$\theta_0 \in \R^\fd$
		
		\Statex\textbf{Output:} 
		Approximate critical point $\theta \in \R^\fd $
		
		\vspace{-2mm} 
		\noindent\hspace*{-8.1mm}\rule{\dimexpr\linewidth+9.24mm}{0.4pt}  
		\vspace{-3mm} 
		\State $\theta \gets \theta_0$
		\State $l_{opt} \gets \infty$
		\State $m \gets 0$
		\State $v \gets 0 $
		\State $\gamma \gets \olr _A ( M, N, \sss, \kkk, \gamma_0, \eta, \theta)$ \mycomment{Determine initial learning rate}
		\State $c \gets 0$
		\For{$i \in \cu{1, 2, \ldots, T}$}
		\State $X_1, X_2, \ldots, X_M  \gets$ i.i.d.~samples of $\mu$
		\State $l \gets M^{-1} \sum_{m=1}^M \los ( \theta , X_m)$
		\State $g \gets M^{-1} \sum_{j=1}^M (\nabla_\theta \los ) ( \theta , X_j )$
		\State $m \gets \beta_1 m + (1 - \beta_1) g$
		\State $v \gets \beta_2 v + (1 - \beta_2 ) g^2$ \mycomment{Square $g^2$ is understood componentwise}
		\State $m \gets m / (1 - \beta_1 ^i )$
		\State $v \gets v / (1 - \beta_2 ^i )$
		\State $\theta \gets \theta - \gamma m / (\sqrt{v} + \epsilon)$ \mycomment{Root $\sqrt{v}$ is understood componentwise}
		\If{$l < l_{opt}$}
			\State $l_{opt} \gets l $
			\State $c \gets 0$
		\Else
			\State $c \gets c+1$ \mycomment{Counts steps since last improvement}
		\EndIf
		\If{$c = \ttt $} \mycomment{If no improvement occurred for $\ttt $ steps}
			\State $\gamma \gets \olr_A (k, s, M, N, \gamma, \eta, \theta)$ \mycomment{Determine new learning rate}
			\State $c \gets 0$
		\EndIf
		\EndFor
		
		\noindent 
		\Return $\theta$
	\end{algorithmic}
\end{myalgorithm}

\section{Numerical experiments}
\label{sec:numerics}
In this section we test the proposed learning-rate-adaptive strategy 
numerically in the case of several \ANN\ learning problems. 
In \cref{subsec:simple_supervised} we apply the 
proposed learning-rate-adaptive approach 
in the situation of a simple $ 6 $-dimensional supervised learning problem 
with an explicitly given polynomial target function 
(see \cref{eq:simple_supervised}) 
that we want to 
approximately learn through 
fully-connected feedforward \ANNs\ with one hidden layer. 
Thereafter, in \cref{sec:DKM,sec:DRM,sec:PINNs} we test the proposed 
learning-rate-adaptive approach in the context of 
deep learning approximation methods for \PDEs. 
In particular, we apply the proposed 
learning-rate-adaptive approach 
in the situation of a $ 25 $-dimensional heat equation 
(see \cref{sec:Heat_PDE})
and a $ 10 $-dimensional Black--Scholes equation
(see \cref{sec:BS}). 

It should be pointed out that in all our numerical experiments the gradient steps for our adaptive method include both the training steps and those used to determine a new learning rate.
In this sense, the graphs really present a fair comparison in terms of computational complexity.

In the recent years there has been a large interest 
in the scientific literature in the design and the study of 
deep learning based numerical approximation methods 
for \PDEs, especially, in the situation of high-dimensional \PDE\ problems. 
For instance, 
based on suitable \BSDE\ representations 
for \PDEs\ (generalized Feynman--Kac formulas) 
in 
\cite{MR3736669,HanJentzenE2018}
the so-called 
deep \BSDE\ method has been proposed 
for approximating \PDEs, 
based on classical formulations for \PDEs\ in \cite{SirignanoSpiliopoulos2018} 
the so-called \DGM\ and 
in \cite{MR3881695}
so-called \PINNs\ 
have been proposed for 
approximating \PDEs\ (cf.\ \cref{sec:PINNs} below), 
in \cite{EYu2018}
the so-called \DRM\ for 
approximating elliptic \PDE\ problems 
based on variational formulations 
of \PDEs\ has been proposed (cf.\ \cref{sec:DRM} below), 
and 
in \cite{Beck2021Kolmogorov}
a suitable \DKM\ for approximating 
\PDEs\ based on forward \SDE\ representations 
for \PDEs\ (Feynman--Kac formulas) 
has been proposed (cf.\ \cref{sec:DKM} below). 
The above named references for deep learning based approximation methods 
for (high-dimensional) \PDE\ problems just provide a small selection 
of a nowadays extensive literature in this area of research 
and for additional references and details we refer, for example, 
to the overview and review articles \cite{MR4270459,MR4356985,MR4547087,Karniadakisetal2021PIML,Germainetal2021Overview}, 
the monograph \cite[Chapters~16, 17, and 18]{JentzenKuckuckvonWurstemberger2023arXiv}, 
and the references mentioned therein.

Each of the below presented numerical experiments 
has been implemented in {\sc Python} using the 
machine learning library {\sc Pytorch} (cf.~\cite{Paszke2017}). 
The {\sc Python} source codes for each of the below presented 
numerical experiments
can be downloaded on {\sc GitHub} 
from \url{https://github.com/deeplearningmethods/adaptive-learning-rate}
and 
can also be downloaded from the arXiv webpage of this article 
as a .tar.gz file containing separate .py files 
by clicking on the arXiv webpage of this article, first, 
at \emph{Other Formats} and, thereafter, at 
\emph{Download source}.

\subsection{Supervised learning problem}
\label{subsec:simple_supervised}

In this subsection we train standard fully-connected feedforward \ANNs\ to approximate the target function
\begin{equation}
	\label{eq:simple_supervised}
	\begin{split}
		[-1, 1 ] ^d \ni x \mapsto 1 + \ssum_{i=1}^d (d + 1 - 2 i ) ( x_i )^3 \in \R
	\end{split}
\end{equation}
for $d =6$.
We use \ANNs\ with the \ReLU\ activation and 
one hidden layer with $ 128 $ neurons on the hidden layer. 
For the probability distribution of the input data 
we choose the uniform distribution on $ [-1, 1]^d $.
We use a batch size of 256 and 
we evaluate the models on independently generated test datasets of size 2000. 
As further parameters we choose $\sss = 50$, $\kkk = 5$, $\eta = 4$ 
and the tolerance $\ttt = 400$.

\begin{figure}
    \includegraphics[scale=0.4]{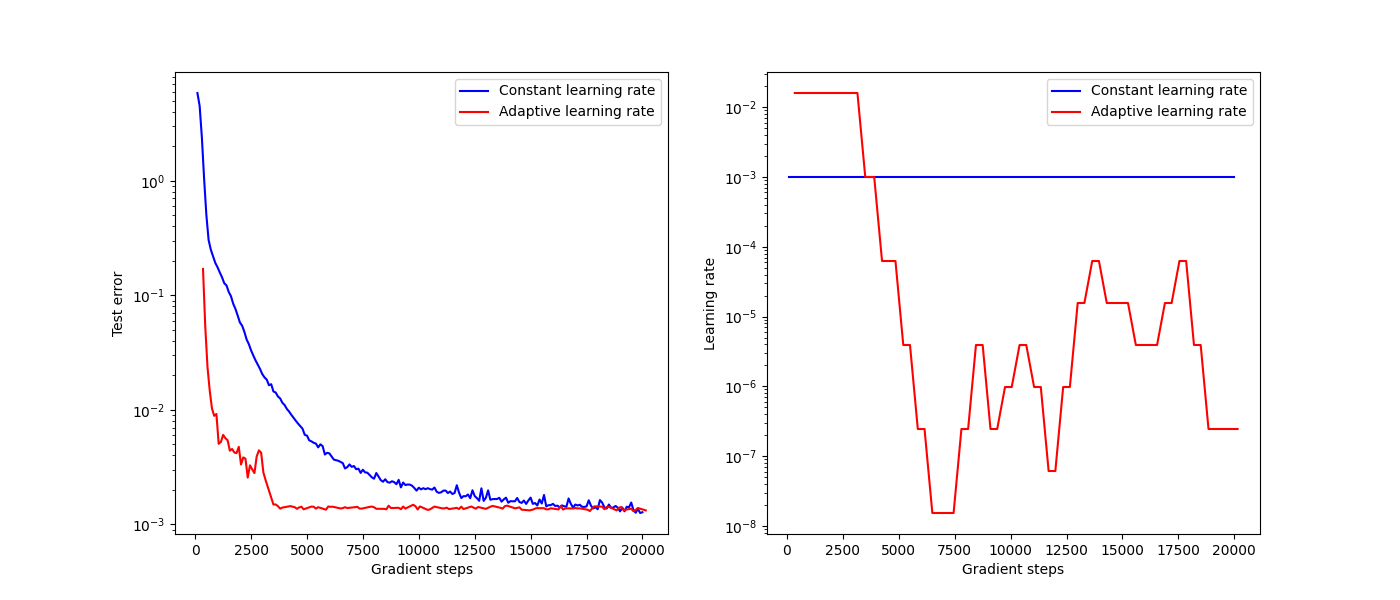}
	\caption{\label{fig_supervised_ann}Numerical results for the 
	supervised learning problem 
	from \cref{subsec:simple_supervised} 
	with the target function in  \cref{eq:simple_supervised}.}
\end{figure}
The results are visualized in \cref{fig_supervised_ann}.
As one can see from the graphics, while the standard Adam optimizer with constant learning rate eventually reaches 
almost the same loss value as the proposed learning-rate-adaptive variant of it, 
it takes significantly more steps.
It seems that the reason for this is that the proposed learning-rate-adaptive variant of Adam first chooses a larger learning rate 
compared to the default value, which leads to a steeper decrease of the loss curve, and later choose a smaller learning rate, which refines the parameters and reduces oscillations.

\subsection{Deep Kolmogorov method (DKM)}
\label{sec:DKM}

In this subsection we combine the proposed learning-rate-adaptive approach 
with the \DKM\ described in Beck et al.~\cite{Beck2021Kolmogorov} to approximately solve 
some linear \PDEs\ of the Kolmogorov type. 

\subsubsection{DKM for a heat PDE}
\label{sec:Heat_PDE}

We first consider the heat \PDE\ on $\R^d$ for $d = 25$.
Specifically, we consider $ d = 25 $ 
and we attempt to approximate the solution $u \colon [0 , \infty ) \times \R^d \to \R$ of the 
initial value heat \PDE\ problem
\begin{equation}
	\label{eq:kolmogorov_heat}
	\begin{split}
		\tfrac{ \partial }{ \partial t } u & = \Delta_x u ,
		\qquad 
		u(0, x ) 
		= \norm{x} ^2
	\end{split}
\end{equation}
for $ t \in [0,\infty) $, $ x \in \R^d $
at time $ T = 2 $ on the set $ [ -1, 1]^d \subseteq \R^d $.
The initial value \PDE\ problem in \cref{eq:kolmogorov_heat} 
can be reformulated as a stochastic minimization problem (cf.~Beck et al.~\cite{Beck2021Kolmogorov}) and thus the Adam optimizer can be used to compute an approximate minimizer.
We employ fully-connected feedforward \ANNs\ with the smooth \GELU\ activation, 
which seems to be more suitable for approximating smooth \PDE\ problems.
We use \ANNs\ with three hidden layers with $ 50 $ neurons on the first hidden layer, $100$ neurons on the second hidden layer, and $50$ neurons on the third hidden layer 
and we use a batch size of 4096. 
As further parameters we choose $\sss = 50$, $\kkk = 5$, $\eta = 4$ and the tolerance $\ttt = 400$.

\begin{figure}
  \includegraphics[scale=0.45]{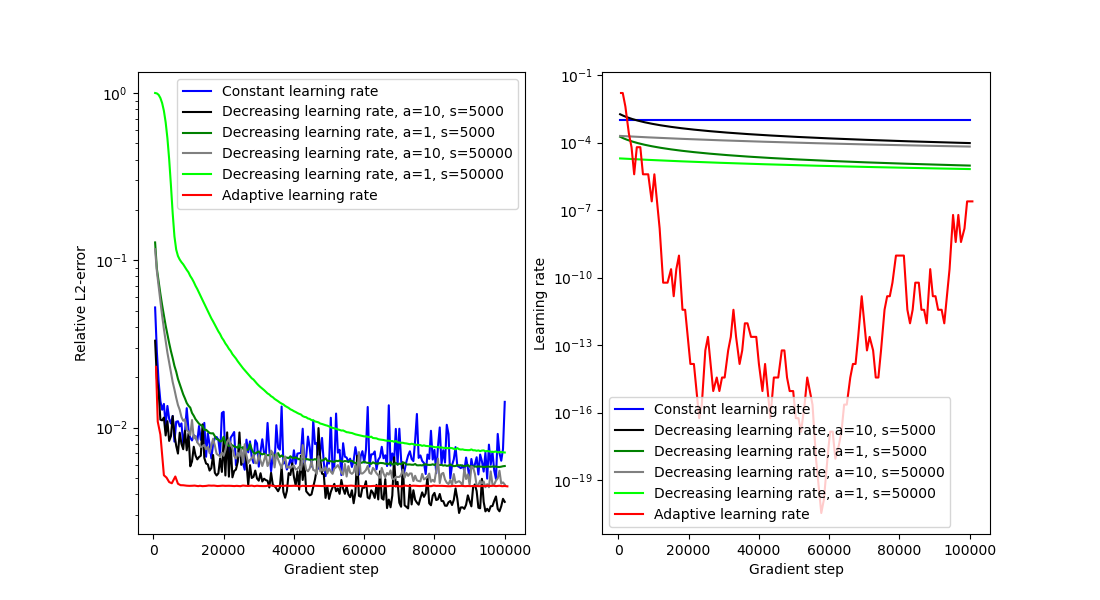}
	\caption{\label{fig_kolmogorov}Numerical results for the heat \PDE\ in \cref{eq:kolmogorov_heat} using the \DKM.}
\end{figure}

The results are visualized in \cref{fig_kolmogorov}. 
To approximately compute the relative $ L^2 $-error 
in the left graph of \cref{fig_kolmogorov} 
we compare the output of the \DKM\ with the exact solution 
$ u \colon [0,\infty) \times \R^d \to \R $
of \cref{eq:kolmogorov_heat}, 
which satisies for all $ t \in [0,\infty) $, $ x \in \R^d $ that 
$ u(t,x) = \norm{x}^2 + 2 d t $.
In addition, we also compare the results when using deterministic decreasing learning rate schedules of the form $\gamma_n = \frac{a}{n+s}$ for different values $a, s \in \N$.

\subsubsection{DKM for a Black--Scholes PDE}
\label{sec:BS}
Next we consider a Black-Scholes \PDE\ with $ d = 10 $ space dimensions.
Specifically, 
we consider $ d = 10 $, 
for every $ i \in \{ 1, 2, \dots, d \} $ 
we consider $ \sigma_i = \frac{ i + 1 }{ 2 d } $, 
we consider 
$ \mu = \frac{ 1 }{ 10 } $, 
$ r = \frac{ 5 }{ 100 } $, 
$ K = 100 $, 
and 
$ T = 1 $, 
and 
we attempt to approximate the solution $ u \colon [0, T] \times \R^d \to \R $ 
of the initial value Black--Scholes \PDE\ problem 
\begin{equation}
\label{eq:kolmogorov_black_scholes}
\textstyle 
		\frac{ \partial }{ \partial t } u 
		= \frac{ 1 }{ 2 } 
		\biggl[ 
		\sum\limits_{ i = 1 }^d 
		(
		  \sigma_i x_i
        )^2 
        \frac{ \partial^2 }{ 
          \partial x_i^2 
        }
        u 
        \biggr]
		+ 
		\mu 
		\left[ 
		\sum\limits_{i=1}^d x_i \frac{ \partial }{ \partial x_i } u 
		\right] 
		, 
		\qquad 
		u( 0 , x ) 
		= 
		\frac{
          \max \cu*{ \max \cu{x_1, x_2, \dots, x_d } - K , 0 }
        }{
          \exp( r T ) 
        }
\end{equation}
for $ t \in [0,T] $, $ x \in \R^d $ 
at the final time $ T = 1 $
on the set $ [90, 110]^d \subseteq \R^d $ 
using the \DKM. 
We employ fully-connected feedforward \ANNs\ with the \GELU\ activation 
and three hidden layers 
with $ 50 $ neurons on the first hidden layer, 
$ 100 $ neurons on the second hidden layer, 
and $ 50 $ neurons on the third hidden layer. 
In addition, we use batch normalization after the input layer 
and we use a batch size of $ 2048 $.
As further parameters we choose $\sss = 50$, $\kkk = 5$, $\eta = 4$ and the tolerance $ \ttt = 500 $.

\begin{figure}
	\includegraphics[scale=0.45]{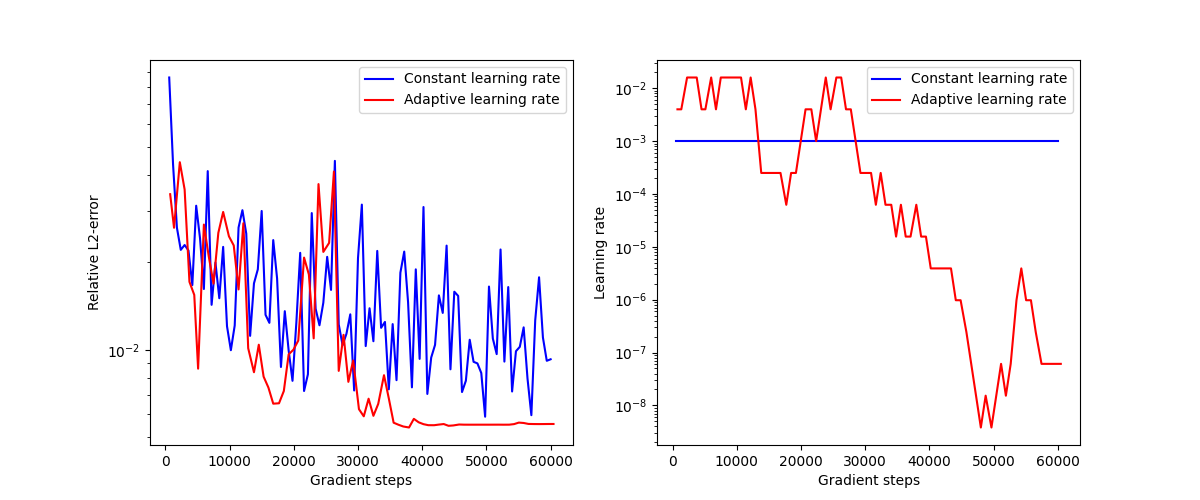}
	\caption{\label{fig_blackscholes}Numerical results for the Black-Scholes \PDE\ in \cref{eq:kolmogorov_black_scholes} using the \DKM.}
\end{figure}

The results are visualized in \cref{fig_blackscholes}. 
To approximately compute the relative $ L^2 $-error 
in the left graph of \cref{fig_blackscholes} 
we compare the output of the \DKM\ with an approximation of the exact solution 
$ u \colon [0,T] \times \R^d \to \R $ of \cref{eq:kolmogorov_black_scholes}
computed with a Feynman--Kac formula based Monte Carlo method with $ 409'600 $ 
Monte Carlo samples.

\subsubsection{DKM for a stochastic Lorenz equation}
As the next example we consider a stochastic Lorenz initial value \PDE\ problem. 
Specifically, 
we consider $ T = 1 $, 
$ \alpha = ( \alpha_1, \alpha_2, \alpha_3 ) = (10, 14, \nicefrac{8}{3} ) \in \R^3 $,
$ \beta = \nicefrac{3}{20} $
and we aim to approximate the solution 
$
  u \colon [0,T] \times \R^d \to \R 
$
of the stochastic Lorenz initial value \PDE\ problem 
\begin{equation}
	\label{eq:kolmogorov_lorenz}
	\begin{split}
		\tfrac{ \partial }{ \partial t } u (t, x) 
		&= \tfrac{\beta^2}{2} \Delta_x u (t, x) 
		+ \alpha_1(x_2 - x_1) \tfrac{ \partial }{ \partial x_1 } u (t, x)
		+ (\alpha_2 x_1 - x_2 - x_1 x_3) \tfrac{ \partial }{ \partial x_2 } u (t, x) \\
		& \qquad + (x_1 x_2 - \alpha_3 x_3) \tfrac{ \partial }{ \partial x_3 } u (t, x) , 
\\
		u(0 , x ) &= \norm{x} ^2
	\end{split}
\end{equation}
for $ t \in [0,T] $, $ x \in \R^d $ 
at time $ T = 1 $ on the set 
$ [ \nicefrac{ 1 }{ 2 }, \nicefrac{ 3 }{ 2 } ] \times [8,10] \times [10, 12] \subseteq \R^3 $
using the \DKM. 
To approximate the solution of the \PDE\ in \cref{eq:kolmogorov_lorenz} 
we employ the \DKM\ with 
fully-connected feedforward \ANNs\ with the \GELU\ activation 
and two hidden layers with $ 100 $ neurons on the first hidden layer and $ 100 $ neurons 
on the second hidden layer. 
We use a batch size of $ 512 $ and we approximate the solutions of 
the \SDE\ associated to \cref{eq:kolmogorov_lorenz} 
(which is often referred to as stochastic Lorenz system with additive noise) 
using $ 100 $ Euler--Maruyama time steps, as described in \cite{Beck2021Kolmogorov}. 
As further parameters we chose $\sss = 50$, $\kkk = 5$, $\eta = 3$ and the tolerance $\ttt = 1000$.

\begin{figure}
	\includegraphics[scale=0.5]{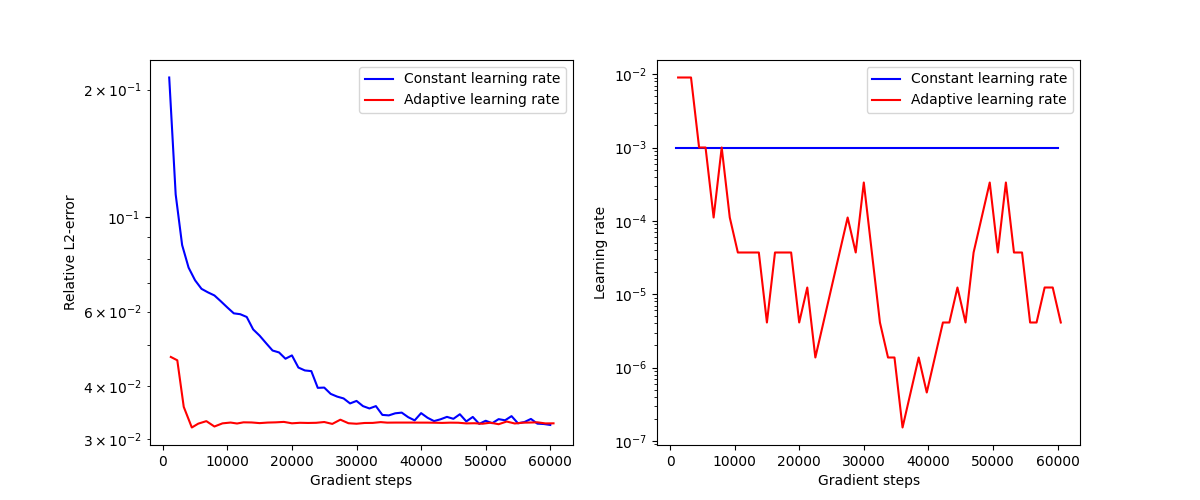}
	\caption{\label{fig_lorenz}Numerical results for the Lorenz \PDE\ in \cref{eq:kolmogorov_lorenz} using the \DKM.}
\end{figure}

The results are visualized in \cref{fig_lorenz}. 
To approximately compute the relative $ L^2 $-error in the left graph side of \cref{fig_lorenz} 
we compare the output of the \DKM\ with an approximation of the exact solution 
$ u \colon [0,T] \times \R^d \to \R $ of \cref{eq:kolmogorov_lorenz} 
computed with a Feynman--Kac formula based Monte Carlo method 
using $ 800'000 $ Monte Carlo samples.

\subsection{Physics-informed neural networks (PINNs)}
\label{sec:PINNs}

In this subsection we employ \PINNs\ (cf.\ \cite{MR3881695}) and 
the \DGM\ (cf.\ \cite{SirignanoSpiliopoulos2018}), respectively, 
to approximately compute solutions 
of a $ 2 $-dimensional sine-Gordon type \PDE\ (cf.\ \cref{sec:sine_Gordon}) and 
a $ 2 $-dimensional Allen-Cahn \PDE\ (cf.\ \cref{sec:Allen_Cahn}).

\subsubsection{PINNs for a sine-Gordon type PDE}
\label{sec:sine_Gordon}

In this example we use \PINNs\ to approximately solve a semilinear heat \PDE\ 
on the domain $ D = [0,2] \times [0,1] $ with time horizon $T=1$. 
Specifically, 
we consider $ D = [0,2] \times [0,1] $ 
and $ T = 1 $ 
and we attempt to approximate 
the solution $ u \colon [0, T] \times D \to \R $ of the initial value sine-Gordon type \PDE\ problem 
\begin{equation}
	\label{eq:pinn_sine}
		\tfrac{ \partial }{ \partial t } u(t,x) = \tfrac{ 1 }{ 20 } \Delta_x u(t,x) + \sin( u(t,x) ) ,
\qquad 
		u(0, x ) = \tfrac{3}{2} \left| \sin( \pi x_1 ) \sin( \pi x_2 ) \right|^2
\end{equation}
for $ t \in [0,T] $, $ x \in (0,2) \times (0,1) \subseteq \R^2 $ 
equipped with Dirichlet boundary conditions.
We employ fully-connected feedforward \ANNs\ with the \GELU\ activation and $ 3 $ hidden layers 
with $ 32 $ neurons on the first hidden layer, $ 64 $ neurons on the second hidden layer, 
and $ 32 $ neurons on the third hidden layer.
We use a batch size of $256$.
As further parameters we chose $\sss = 50$, $\kkk = 5$, $\eta = 4$ and the tolerance $\ttt = 400$.

\begin{figure}
	\includegraphics[scale=0.5]{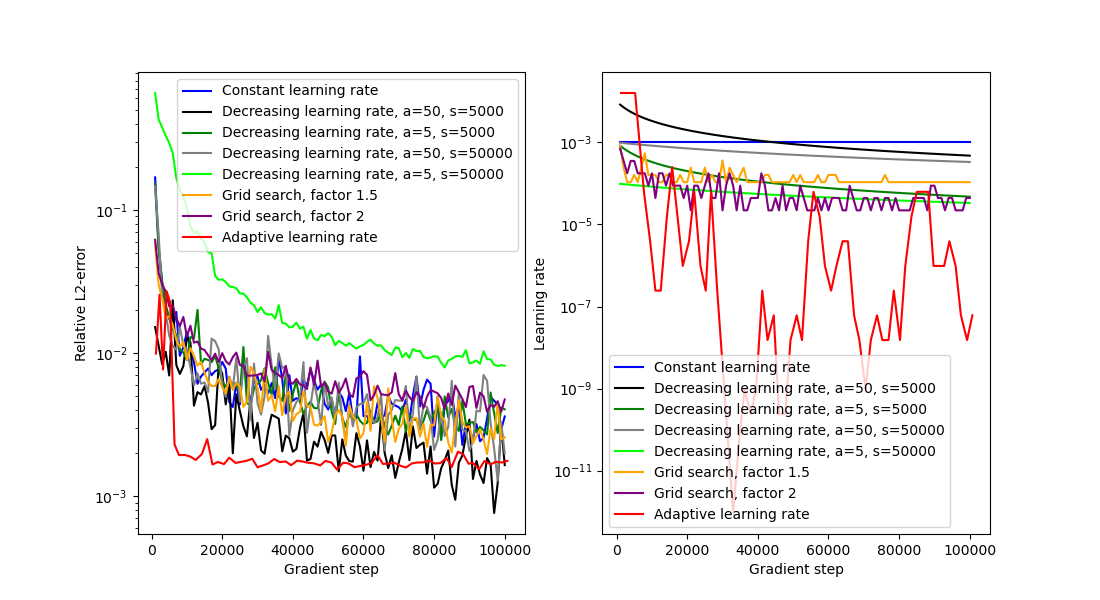}
	\caption{\label{fig_pinn_sine}Numerical results for the sine-Gordon type \PDE\ in \cref{eq:pinn_sine} using \PINNs.}
	\end{figure}

The results are visualized in \cref{fig_pinn_sine}.
To approximately compute the relative $ L^2 $-error 
in the left graph of \cref{fig_pinn_sine} 
we compare the output with an approximation of the exact solution 
$
  u \colon [0,T] \times D \to \R
$
of \cref{eq:pinn_sine} computed with a linear-implicit Runge--Kutta finite element method 
using $101^2$ degrees of freedom in space and $500$ linear-implicit Runge--Kutta time steps.
As above, we also compare the results for deterministic decreasing learning rate schedules of the form $\gamma_n = \frac{a}{n+s}$.
Furtermore, as a comparison we attempt to determine the optimal learning rate in each step by using a grid search over $11$ grid points with the geometric grid factors $1.5$ and $2$, respectively.
It can be seen that this method does not work as well as our adaptive learning rate method, even though it has significantly higher computational cost due to evaluating the loss function at multiple points in each gradient step.

\subsubsection{PINNs for an Allen-Cahn PDE}
\label{sec:Allen_Cahn}
In this example we use \PINNs\ to approximately solve an Allen-Cahn \PDE\ on 
the $ 2 $-dimensional domain 
$ (0,2) \times (0,1) \subseteq \R^2 $ 
with time horizon $ T = 4 $. 
Specifically, 
we consider 
$ D = [0,2] \times [0,1] $, 
we consider the time horizon $ T  = 4 $, 
and we attempt to approximate the solution $u \colon [0, T] \times D \to \R$ of the 
Allen-Cahn initial value \PDE\ problem 
\begin{equation}
\label{eq:pinn_allencahn}
\textstyle 
  \frac{ \partial }{ \partial t } u(t,x) 
  = \frac{ 1 }{ 100 } \Delta_x u(t,x) 
  + u(t,x) - ( u(t,x) )^3 
  ,
\qquad 
  u(0, x ) = 
  \sin( \pi x_1 ) \sin( \pi x_2 )
\end{equation}
for $ t \in [0,T] $, $ x = ( x_1, x_2 ) \in (0,2) \times (0,1) $
with Dirichlet boundary conditions.
We employ fully-connected feedforward \ANNs\ with 
the \GELU\ activation 
and 2 hidden layers with $ 64 $ neurons on the first hidden layer and $ 64 $ neurons 
on the second hidden layer. 
We use a batch size of $256$. 
As further parameters we chose $\sss = 50$, $\kkk = 5$, $\eta = 4$ 
and the tolerance $\ttt = 400$.

\begin{figure}
	\centering
	\includegraphics[scale=0.5]{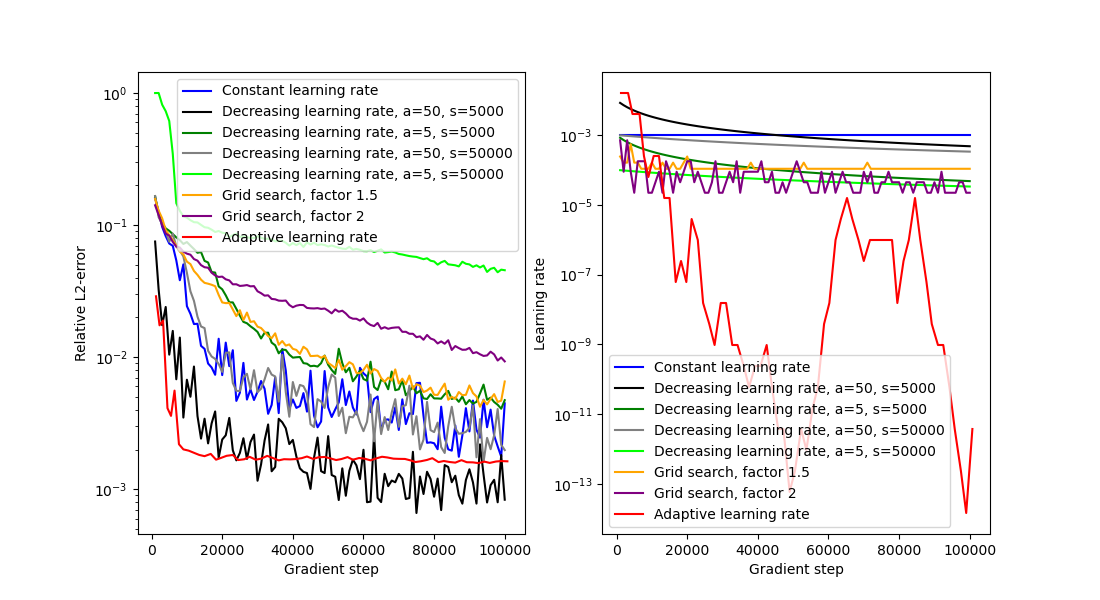}
	\caption{\label{fig_pinn_ac}Numerical results for the Allen-Cahn \PDE\ in \cref{eq:pinn_allencahn} using \PINNs.}
\end{figure}

\begin{figure}
	\centering 
	\includegraphics[scale=0.5]{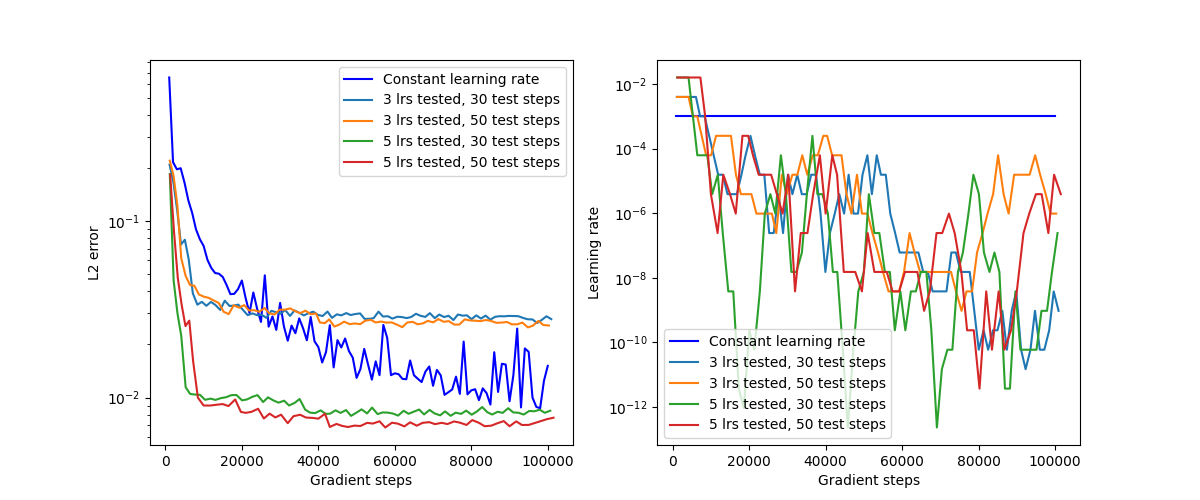}
	\caption{\label{fig:pinn-multi}Comparison of parameter choices for the Allen-Cahn \PDE\ solved with \PINNs.}
\end{figure}

The results are visualized in \cref{fig_pinn_ac}.
To compute the $ L^2 $-error 
in the left graph of \cref{fig_pinn_ac} 
we compare the output of the \PINNs\ with an approximation 
of the exact solution 
$ u \colon [0,T] \times D \to \R $ of 
\cref{eq:pinn_allencahn} computed 
with a linear-implicit Runge--Kutta finite element method 
using $ 101^2 $ degrees of freedom in space 
and $ 500 $ linear-implicit Runge--Kutta time steps. 
As before, we also compare deterministic decreasing learning rate schedules and optimized learning rates in each step using a line search method.
Moreover, in \cref{fig:pinn-multi} we compare multiple choices 
of the parameters $\kkk \in \cu{3, 5}$ and $\sss \in \cu{30, 50 }$ 
for the \PINN\ training problem for the Allen-Cahn \PDE\ in \cref{eq:pinn_allencahn}.

\subsection{Deep Ritz method for a Poisson equation}
\label{sec:DRM}

In this subsection we employ the \DRM\ from E \& Yu~\cite{EYu2018} 
to approximately solve a $ 10 $-dimensional Poission equation. 
Specifically, we consider $ d = 10 $ 
and we attempt to approximation the solution 
$ 
  u \colon [-1, 1 ]^d \to \R
$
of the Poisson equation 
\begin{equation}
\label{eq:poisson_ritz}
\begin{split}
  \Delta u( x ) = 0 , 
  \qquad 
  u(y) = \ssum_{ i = 1 }^5 
  y_{ 2 i - 1 } y_{ 2 i } 
\end{split}
\end{equation}
for $ x = ( x_1, \dots, x_d ) \in ( -1, 1 )^d $, 
$ 
  y = ( y_1, \dots, y_d ) \in (  [ -1, 1 ]^d \backslash (-1,1)^d )
$.
We employ a residual \ANNs\ architecture with 4 blocks, 
each consisting of two layers with $ 16 $ neurons each and 
a residual connection as described in \cite{EYu2018}, and the \GELU\ activation. 
We use a batch size of $ 4096 $. 
As further parameters we chose $ \sss = 50 $, $ \kkk = 5 $, $ \eta = 4 $ 
and the tolerance $ \ttt = 400 $.

\begin{figure}
	\caption{Results for the Poisson equation in \cref{eq:poisson_ritz} using the \DRM.}
	\label{fig:ritz}
	\centering 
	\includegraphics[scale=0.5]{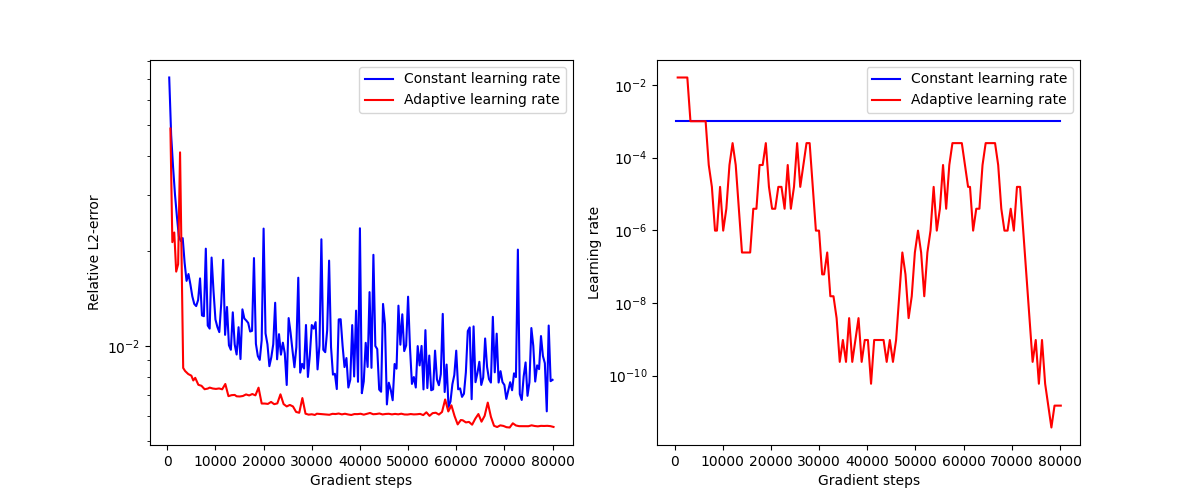}
\end{figure}

The results are visualized in \cref{fig:ritz}. 
To approximately compute the relative $ L^2 $-error 
in the left graph of \cref{fig:ritz} we compare the output of the \DRM\ with 
the exact solution $ u \colon [-1,1]^d \times \R^d \to \R $ 
of \cref{eq:poisson_ritz}, 
which satisfies for all $ x = ( x_1, \dots, x_d ) \in [-1,1]^d $ 
that 
$
  u( x ) = \ssum_{ i = 1 }^5 x_{ 2 i - 1 } x_{ 2 i } 
$.

\section{Convergence analysis of SGD with random learning rates}
\label{sec:random_learning_rates}

The main goal of this section is to prove in \cref{cor:convergence_in_probability_ex}
an abstract convergence result for \SGD\ processes with random decreasing learning rates in the case of an objective function satisfying a certain coercivity condition.
\cref{cor:convergence_in_probability_ex} is a consequence of \cref{prop:convergence_in_probability}, one of the main results of this article.

\subsection{Bounds for linear implicit Euler approximations}

In this section we recall in \cref{lem:exp_bound} an elementary estimate for the exponential function which will be employed later in the convergence proofs.
\cref{cor:exp_bound} is a simple consequence of \cref{lem:exp_bound}.

\begin{lemma}[Bounds for the linear implicit Euler approximation]
\label{lem:exp_bound}
Let $ c \in [1,\infty) $. 
Then it holds for all 
$ x \in [ 0, c^{ - 1 } \ln( c ) ] $ 
that
\begin{equation}
\label{eq:lem_exp_bound}
  \exp(
    - c x
  )
  \leq 1 - x \leq \exp( - x )
  .
\end{equation}
\end{lemma}

\begin{cproof}{lem:exp_bound}
Throughout this proof let $ \varphi \colon \R \to \R $
satisfy for all $ x \in \R $ that
\begin{equation}
\label{eq:definition_varphi}
  \varphi(x) = 1 - x - \exp( - c x )
  .
\end{equation}
\Nobs that \cref{eq:definition_varphi} ensures 
for all $ x \in \R $ that 
\begin{equation}
  \varphi'(x) = - 1 + c \exp( - c x )
  .
\end{equation}
\Hence that for all $ x \in \R $ it holds that 
$
  \varphi'(x) \geq 0 
$
if and only if 
\begin{equation}
  \exp( - c x ) \geq c^{ - 1 }
  .
\end{equation}
\Hence that for all $ x \in \R $ it holds that 
$
  \varphi'(x) \geq 0 
$
if and only if 
\begin{equation}
  - c x \geq \ln( c^{ - 1 } ) = - \ln( c )
  .
\end{equation}
\Hence that for all $ x \in \R $ it holds that 
$
  \varphi'(x) \geq 0 
$
if and only if 
\begin{equation}
  x \leq \frac{ \ln( c ) }{ c }
  .
\end{equation}
\end{cproof}

\begin{cor}[Bounds for the linear implicit Euler approximation]
\label{cor:exp_bound}
Let $ x \in [ 0, \exp( - 1 ) ] $. 
Then 
\begin{equation}
\label{eq:cor_exp_bound}
  \exp( x )
  \leq 
  ( 1 - x )^{ - 1 } 
  \leq 
  \exp\bigl(
    \exp( 1 ) x
  \bigr)  
  .
\end{equation}
\end{cor}

\begin{cproof}{cor:exp_bound}
\Nobs that \cref{lem:exp_bound} establishes \cref{eq:cor_exp_bound}. 
\end{cproof}

\subsection{Properties for convergence in probability}

In thus subsection we present in \cref{lem:almost_sure_monotone,lem:almost_sure,lem:convergence_in_probability_charac}
a few well-known facts on convergence in probability and 
almost sure convergence of random variables. 
A detailed proof for \cref{lem:almost_sure_monotone} can be found in the \href{https://arxiv.org/abs/2406.14340}{arXiv-version} of this article. 

\begin{lemma}[Monotone convergence in probability]
\label{lem:almost_sure_monotone}
Let $ ( \Omega, \mathcal{F}, \P ) $ be a probability space, 
for every $ n \in \N $ 
let $ X_n \colon \Omega \to [0,\infty] $
be a random variable, 
assume 
for all $ \varepsilon \in (0,\infty) $ 
that
$
  \lim_{ n \to \infty }
  \P( 
    X_n > \varepsilon 
  ) = 0
$, 
and assume for all $ n \in \N $ that 
$
  \P(
    X_{ n + 1 } \leq X_n
  ) = 1
$. 
Then 
\begin{equation}
\textstyle 
  \P\bigl(
    \lim_{ n \to \infty }
    X_n = 0
  \bigr)
  = 1 .
\end{equation}
\end{lemma}

The next result is well-known; see, e.g., \cite[Remark 6.4]{Klenke2014}.

\begin{lemma}[Almost sure convergence implies convergence in probability]
\label{lem:almost_sure}
Let $ ( \Omega, \mathcal{F}, \P ) $ be a probability space, 
for every $ n \in \N $ 
let $ X_n \colon \Omega \to [0,\infty] $
be a random variable, 
let $ \varepsilon \in (0,\infty) $, 
and assume 
$
  \P\bigl( 
    \lim_{ n \to \infty }
    \allowbreak 
    X_n 
    = 0
  \bigr) = 1
$. 
Then 
\begin{equation}
  \lim_{ n \to \infty }
  \P\bigl(
    X_n \geq \varepsilon
  \bigr)
  = 0 .
\end{equation}
\end{lemma}

The next lemma is also well-known and is a consequence of, e.g., \cite[Theorem 6.7]{Klenke2014}.

\begin{lemma}[Characterization of convergence in probability]
\label{lem:convergence_in_probability_charac}
Let $ ( \Omega, \mathcal{F}, \P ) $ be a probability space 
and for every $ n \in \N $ 
let $ X_n \colon \Omega \to [0,\infty] $
be a random variable. 
Then the following two statements are equivalent:
\begin{enumerate}[label = (\roman*)]
\item 
It holds for all $ \varepsilon \in (0,\infty) $ that 
$
  \lim_{ n \to \infty }
  \P(
    X_n \geq \varepsilon
  )
  = 0
$.
\item 
It holds that
$
  \lim_{ n \to \infty }
  \E[ \min\{ 1, X_n \} ] = 0
$. 
\end{enumerate}
\end{lemma}

\subsection{Convergence in probability to zero}

We now establish in \cref{prop:convergence_in_probability} our first convergence result for \SGD\ processes with random learning rates.
In \cref{prop:convergence_in_probability} we think of $g \colon \R^\fd \to \R^\fd$ as the negative gradient of the objective function and as of the random variables $D_n \colon \Omega \to \R^\fd$, $n \in \N$, as random error terms.
\Nobs that we assume in \cref{eq:coercivity_for_zero} a certain coercivity property of $g$, but no further conditions such as Lipschitz continuity, which is often used when analyzing \SGD\ processes.
It should also be pointed out that an objective function whose negative gradient satisfies \cref{eq:coercivity_for_zero} is not necessarily convex.

The random learning rates $\gamma_n \colon \Omega \to [0 , \infty )$, $n \in \N$,
are assumed to be previsible, decreasing, convergent to zero, and not summable.
It should be noted that we do not require a summability condition like $\sum_{n=1} ^\infty \gamma_n^2 < \infty$, which is common in the analysis of \SGD\ methods.

Regarding the error terms $D_n$, we assume that these form a martingale difference sequence and satisfy the boundedness condition $\E [ \sup_{n \le \tau_R } \norm{D_n } ^2 ] < \infty $ for every $R > 0$,
where the stopping time $\tau_R$ is the first time where the process exits a ball of radius $\sqrt{R}$ (see \cref{eq:definition_of_tau_R} below).
In other words, as long as the \SGD\ process remains in a bounded domain, the second moment of the error terms can be uniformly bounded from above.
This condition is usually satisfied for \SGD\ processes under mild boundedness assumptions on the considered data.

The conclusion of \cref{prop:convergence_in_probability} reveals that the iterates $\Theta_{N_t}$ at the stopping times $N_t$ introduced in \cref{eq:N_n_definition} converge to zero in probability as $t \to \infty$ when restricting to the event $\cE \subseteq \Omega$ where the process remains bounded.

\begin{prop}[Convergence in probability]
\label{prop:convergence_in_probability}
Let $ ( \Omega, \mathcal{F}, ( \mathbb{F}_n )_{ n \in \N_0 }, \P ) $ be a filtered probability space, 
for every $ n \in \N $ 
let $ \gamma_n \colon \Omega \to [0,\infty) $ be $ \mathbb{F}_{ n - 1 } $-measurable, 
assume 
$
  \{ 
    \lim_{ n \to \infty }
    \gamma_n 
    =
    0
  \} 
  =
  \{ 
    \sum_{ n = 1 }^{ \infty } \gamma_n 
    =
    \infty
  \} 
  =
  \cap_{ n = 1 }^{ \infty }
  \{ 
    \gamma_{ n + 1 } \leq \gamma_n 
  \} 
  = \Omega 
$, 
let $ \fd \in \N $, 
let $ g \colon \R^{ \fd } \to \R ^\fd$ be measurable and locally bounded, 
let $ c \in (0,\infty) $ satisfy for all $ \theta \in \R^{ \fd } $ that 
\begin{equation}
\label{eq:coercivity_for_zero}
  \left< \theta, g( \theta ) \right> 
  \leq 
  - c \| \theta \|^2 
  ,
\end{equation}
let $ D \colon \N \times \Omega \to \R^{ \fd } $
be a stochastic process,
assume for all $ n \in \N $ that $\E [ \norm{D_n} ] < \infty $ and, $ \P $-a.s.,
$
  \E[ 
    D_n | \mathbb{F}_{ n - 1 }
  ] 
  = 0
$,
let 
$ \Theta \colon \N_0 \times \Omega \to \R^{ \fd } $
be an $ ( \mathbb{F}_n )_{ n \in \N_0 } $-adapted stochastic process, 
assume for all $ n \in \N $ that 
\begin{equation}
\label{eq:Theta_n_recursion_in_prop}
  \Theta_n 
  = \Theta_{ n - 1 } 
  + 
  \gamma_n 
  g( \Theta_{ n - 1 } )
  +
  \gamma_n 
  D_n 
  ,
\end{equation}
for every $ R \in \R $ let $ \tau_R \colon \Omega \to ( \N_0 \cup \{ \infty \} ) $ satisfy
\begin{equation}
	\label{eq:definition_of_tau_R}
	\tau_R = 
	\inf\bigl(
	\bigl\{ 
	n \in \N_0 \colon 
	\| \Theta_n \|^2 
	\geq R
	\bigr\} 
	\cup 
	\{ \infty \} 
	\bigr)
	,
\end{equation}
assume for all $R \in \R$ that $\E [ \sup_{n \le \tau_R } \norm{D_n } ^2 ] < \infty $,
for every $ t \in [0,\infty) $ 
let 
$
  N_t \colon \Omega \to \N_0
$
satisfy 
\begin{equation}
\label{eq:N_n_definition}
\textstyle 
  N_t
  =
  \inf\bigl(
    \bigl\{
      n \in \N_0
      \colon 
      \sum_{ m = 1 }^n
      \gamma_m
      \geq 
      t
    \bigr\}
  \bigr)
  ,
\end{equation}
and let $ \cE \subseteq \Omega $ satisfy 
$
  \cE = 
  \{ 
    \sup_{ n \in \N }
    \| \Theta_n \|  
    < \infty 
  \} 
$.
Then
\begin{enumerate}[label = (\roman*)]
\item
\label{item:i}
it holds for all $ \delta \in (0,\infty) $ that 
\begin{equation}
  \lim_{ n \to \infty }
  \E\bigl[
    \min\bigl\{ 
      1 ,
      \| 
        \Theta_{ N_{ n \delta } }
      \|
      \mathbbm{1}_{ \cE }
    \bigr\} 
  \bigr]
  = 0
  ,
\end{equation}
\item
\label{item:ii}
it holds for all $ \delta \in (0,\infty) $ that 
\begin{equation}
  \lim_{ n \to \infty }
  \E\bigl[
    \min\bigl\{ 
      1 ,
      \bigl(
        \max\nolimits_{ 
          N_{ n \delta } \leq m \leq N_{ ( n + 1 ) \delta }
        }
        \| 
          \Theta_m
        \|
      \bigr)
      \mathbbm{1}_{ \cE }
    \bigr\} 
  \bigr]
  = 0
  ,
\end{equation}
and 
\item
\label{item:iii}
it holds that 
\begin{equation}
  \lim_{ t \to \infty }
  \E\bigl[ 
    \min\bigl\{ 
      1 , 
      \| \Theta_{ N_t } \|
      \mathbbm{1}_{ \cE }
    \bigr\}
  \bigr]
  = 0
  .
\end{equation}
\end{enumerate}
\end{prop}

To establish the first statement of \cref{prop:convergence_in_probability},
we need to show convergence in probability, i.e., that for every $\rho > 0$
one has $\lim_{n \to \infty } \P ( \norm{\Theta_{N_{n \delta } } } \indicator{\cE} > \rho ) = 0$.
To do this, we use the estimate \cref{eq:fundamental_estimate_for_item_i} and show that each of the three summands on the right-hand side vanishes.
\begin{itemize}
	\item For the first term, we use the coercivity assumption \cref{eq:coercivity_for_zero} to derive the sufficient decrease property \cref{prop:convergence:eq_decrease} for the norm $\norm{\Theta_n}^2$.
	Roughly speaking, restricting to the event 
	that $\cu{\gamma_{ \scrN^R_{ m \delta } + 1 } \leq \varepsilon  } \cap \{ \tau_R > N_{ n \delta } \}$, the error terms in \cref{prop:convergence:eq_decrease} will be small due to the boundedness of the moments of $D_n$.
	Furthermore, by definition of the stopping times $N_t$ we know that the sum of learning rates $\ft_{ \scrN^R_{ n \delta } }$ until the considered time is at least $n \delta$,
	and thus the first term on the right-hand side of \cref{prop:convergence:eq_decrease} will tend to zero as $n \to \infty $.
	\item The second term on the right-hand side of \cref{eq:fundamental_estimate_for_item_i} vanishes due to the assumption that $\gamma_n \to 0$.
	\item Finally, the third term on the right-hand side of \cref{eq:fundamental_estimate_for_item_i} is zero by definition of the event $\cE$.
\end{itemize}
The second statement is then proved in a similar way, by using an analogous error decomposition (see \cref{eq:fundamental_estimate_for_item_ii} for details).
Finally, \cref{item:iii} is a relatively simple consequence of \cref{item:ii}.
We now present the detailed proof.

\begin{cproof}{prop:convergence_in_probability}
Throughout this proof for every $ n \in \N_0 $ let 
$ \ft_n \colon \Omega \to \R $ satisfy 
\begin{equation}
\label{eq:in_proof_definition_of_ft}
  \ft_n = 
  \sum_{ m = 1 }^n \gamma_m
  ,
\end{equation}
for every 
$ R \in \R $, $ t \in [0,\infty) $ let $ \scrN_t^R \colon \Omega \to \N_0 $ satisfy 
$
  \scrN_t^R = 
  \min\{ \tau_R , N_t \}
$, 
for every $ R \in [0,\infty) $ let $ \scrG_R \in \R $ and $ \scrD_R \colon \Omega \to \R $ satisfy 
\begin{equation}
	\label{eq:def_GR_DR}
  \scrG_R 
  =
  \sup\bigl(
    \cup_{ 
      \theta \in 
      \{ 
        \vartheta \in \R^{ \fd } \colon \| \vartheta \| \leq R 
      \}
    }
    \bigl\{ 
      \| g( \theta ) \|^2
    \bigr\} 
  \bigr)
  \qqandqq \scrD_R = \sup\nolimits_{ n \le \tau_R } \| D_n \|^2 ,
\end{equation}
and for every $ n \in \N $ 
let 
$
  \cD_n \colon \Omega \to \R
$
satisfy 
\begin{equation}
\label{eq:definition_cD_n}
  \cD_n
=
  ( \gamma_n )^2
  \| g( \Theta_{ n - 1 } ) \|^2
  +
  ( \gamma_n )^2 
  \| D_n \|^2
  +
  2 ( \gamma_n )^2
  \langle 
    g( \Theta_{ n - 1 } )
    ,
    D_n
  \rangle
  +
  2 
  \gamma_n 
  \langle 
    \Theta_{ n - 1 } 
    ,
    D_n
  \rangle
  .
\end{equation}
We first establish \cref{item:i}. 
We show \cref{item:i} through 
an application of \cref{lem:convergence_in_probability_charac}. 
We thus need to verify that for all 
for all $ \rho, \delta \in (0,\infty) $
it holds that 
$
  \lim_{ n \to \infty }
  \P(
    \| \Theta_{ N_{ n \delta } } \| \mathbbm{1}_{ \cE } > \rho 
  )
  = 0
$
(see \cref{eq:convergence_in_probability_item_i} below). 
For this \nobs that 
the fact that for all $ A, B \in \cF $ 
it holds that 
$
  \P( A ) 
  = 
  \P( A \cap B ) + \P( A \cap ( \Omega \backslash B ) )
$
shows that for all $ \rho, R, \delta, \varepsilon \in (0,\infty) $, 
$ m, n \in \N_0 $
with $ m \leq n $
it holds that
\begin{equation}
\begin{split}
&
  \P\bigl(
    \| \Theta_{ N_{ n \delta } } \| \mathbbm{1}_{ \cE }
    > \rho
  \bigr)
=
  \P\bigl(
    \{ 
      \| \Theta_{ N_{ n \delta } } \| 
      > \rho 
    \}
    \cap 
    \cE
  \bigr)
\\ &
\leq 
  \P\bigl(
    \{ 
      \| \Theta_{ N_{ n \delta } } \| 
      > \rho 
    \}
    \cap 
    \cE
    \cap 
    \{ 
      \tau_R > N_{ n \delta }
    \}
  \bigr)
  +
  \P\bigl(
    \cE
    \cap 
    \{ 
      \tau_R \leq N_{ n \delta }
    \} 
  \bigr)
\\ &
=
  \P\bigl(
    \{ 
      \| \Theta_{ \min\{ \tau_R, N_{ n \delta } \} } \| 
      > \rho 
    \}
    \cap 
    \cE
    \cap 
    \{ 
      \tau_R > N_{ n \delta }
    \}
  \bigr)
  +
  \P\bigl(
    \cE
    \cap 
    \{ 
      \tau_R \leq N_{ n \delta }
    \} 
  \bigr)
\\ &
=
  \P\bigl(
    \{ 
      \| \Theta_{ \scrN^R_{ n \delta } } \| 
      > \rho 
    \}
    \cap 
    \cE
    \cap 
    \{ 
      \tau_R > N_{ n \delta }
    \}
  \bigr)
  +
  \P\bigl(
    \cE
    \cap 
    \{ 
      \tau_R \leq N_{ n \delta }
    \} 
  \bigr)
\\ &
\leq 
  \P\bigl(
    \bigl\{ 
      \| \Theta_{ \scrN^R_{ n \delta } } \| 
      > \rho 
    \bigr\}
    \cap 
    \{ 
      \tau_R > N_{ n \delta }
    \}
  \bigr)
  +
  \P\bigl(
    \cE
    \cap 
    \{ 
      \tau_R < \infty
    \} 
  \bigr)
\\ &
=
  \P\bigl(
    \bigl\{ 
      \| \Theta_{ \scrN^R_{ n \delta } } \| 
      > \rho 
    \bigr\}
    \cap 
    \{ 
      \tau_R > N_{ n \delta }
    \}
    \cap 
    \{ 
      \gamma_{ N_{ m \delta } + 1 } \leq \varepsilon  
    \} 
  \bigr)
\\ &
\quad 
  +
  \P\bigl(
    \{ 
      \| \Theta_{ \scrN^R_{ n \delta } } \| 
      > \rho 
    \}
    \cap 
    \{ 
      \tau_R > N_{ n \delta }
    \}
    \cap 
    \{ 
      \gamma_{ N_{ m \delta } + 1 } > \varepsilon  
    \} 
  \bigr)
  +
  \P\bigl(
    \cE
    \cap 
    \{ 
      \tau_R < \infty
    \} 
  \bigr)
\\ &
\leq 
  \P\bigl(
    \bigl\{ 
      \| \Theta_{ \scrN^R_{ n \delta } } \| 
      > \rho 
    \bigr\}
    \cap 
    \{ 
      \tau_R > N_{ n \delta }
    \}
    \cap 
    \{ 
      \gamma_{ \scrN^R_{ m \delta } + 1 } \leq \varepsilon  
    \} 
  \bigr)
  +
  \P\bigl(
    \gamma_{ N_{ m \delta } + 1 } > \varepsilon  
  \bigr)
  +
  \P\bigl(
    \cE
    \cap 
    \{ 
      \tau_R < \infty
    \} 
  \bigr)
\\ &
=
  \P\bigl(
    \bigl\{ 
      \| \Theta_{ \scrN^R_{ n \delta } } \| 
      \mathbbm{1}_{
        \{ 
          \gamma_{ \scrN^R_{ m \delta } + 1 } \leq \varepsilon  
        \} 
      }
      > \rho 
    \bigr\}
    \cap 
    \{ 
      \tau_R > N_{ n \delta }
    \}
  \bigr)
  +
  \P\bigl(
    \gamma_{ N_{ m \delta } + 1 } > \varepsilon  
  \bigr)
  +
  \P\bigl(
    \cE
    \cap 
    \{ 
      \tau_R < \infty
    \} 
  \bigr)
  .
\end{split}
\end{equation}
The fact that for all $ s, t \in [0,\infty) $
with $ s \leq t $ it holds that 
$
  N_s \leq N_t 
$
and the Markov inequality \hence demonstrate 
that for all $ \rho, R, \delta, \varepsilon \in (0,\infty) $, 
$ m, n \in \N_0 $
with $ m \leq n $
it holds that
\begin{equation}
\begin{split}
&
  \P\bigl(
    \| \Theta_{ N_{ n \delta } } \| \mathbbm{1}_{ \cE }
    > \rho
  \bigr)
\\ &
\leq 
  \P\bigl(
    \bigl\{ 
      \| \Theta_{ \scrN^R_{ n \delta } } \|^2 
      \mathbbm{1}_{
        \{ 
          \gamma_{ \scrN^R_{ m \delta } + 1 } \leq \varepsilon  
        \} 
      }
      > \rho^2 
    \bigr\}
    \cap 
    \{ 
      \tau_R > N_{ n \delta }
    \} 
  \bigr)
  +
  \P\bigl(
    \gamma_{ N_{ m \delta } + 1 } > \varepsilon  
  \bigr)
  +
  \P\bigl(
    \cE
    \cap 
    \{ 
      \tau_R < \infty
    \} 
  \bigr)
\\ &
=
  \P\bigl(
      \| \Theta_{ \scrN^R_{ n \delta } } \|^2 
      \mathbbm{1}_{
        \{ 
          \gamma_{ \scrN^R_{ m \delta } + 1 } \leq \varepsilon  
        \} 
        \cap 
        \{ 
          \tau_R > N_{ n \delta }
        \} 
      }
      > \rho^2 
  \bigr)
  +
  \P\bigl(
    \gamma_{ N_{ m \delta } + 1 } > \varepsilon  
  \bigr)
  +
  \P\bigl(
    \cE
    \cap 
    \{ 
      \tau_R < \infty
    \} 
  \bigr)
\\ &
\leq
  \rho^{ - 2 }
  \E\bigl[
    \| \Theta_{ \scrN^R_{ n \delta } } \|^2 
    \mathbbm{1}_{
      \{ 
        \gamma_{ \scrN^R_{ m \delta } + 1 } \leq \varepsilon  
      \} 
      \cap 
      \{ 
        \tau_R > N_{ n \delta }
      \} 
    }
  \bigr]
  +
  \P\bigl(
    \gamma_{ N_{ m \delta } + 1 } > \varepsilon  
  \bigr)
  +
  \P\bigl(
    \cE
    \cap 
    \{ 
      \tau_R < \infty 
    \} 
  \bigr)
  .
\end{split}
\end{equation}
\Hence
for all $ \rho, \delta, \varepsilon \in (0,\infty) $ 
that
\begin{equation}
\label{eq:fundamental_estimate_for_item_i}
\begin{split}
&
  \limsup_{ n \to \infty }
  \P\bigl(
    \| \Theta_{ N_{ n \delta } } \| \mathbbm{1}_{ \cE }
    > \rho
  \bigr)
=
  \limsup_{ R \to \infty }
  \limsup_{ m \to \infty }
  \limsup_{ n \to \infty }
  \P\bigl(
    \| \Theta_{ N_{ n \delta } } \| \mathbbm{1}_{ \cE }
    > \rho
  \bigr)
\\ &
\leq
  \rho^{ - 2 }
  \biggl[
    \limsup_{ R \to \infty }
    \limsup_{ m \to \infty }
    \limsup_{ n \to \infty }
    \E\bigl[
      \| \Theta_{ \scrN^R_{ n \delta } } \|^2 
      \mathbbm{1}_{
        \{ 
          \gamma_{ \scrN^R_{ m \delta } + 1 } \leq \varepsilon  
        \}
        \cap 
        \{ \tau_R > N_{ n \delta } \} 
      }
    \bigr]
  \biggr]
\\ &
  +
  \biggl[
    \limsup_{ m \to \infty }
    \P\bigl(
      \gamma_{ N_{ m \delta } + 1 } > \varepsilon  
    \bigr)
  \biggr]
  +
  \biggl[
    \limsup_{ R \to \infty }
    \P\bigl(
      \cE 
      \cap 
      \{ 
        \tau_R < \infty 
      \} 
    \bigr)
  \biggr]
  .
\end{split}
\end{equation}
To establish \cref{item:i}, it thus remains to show that 
each of the three summands on the right hand side of 
\cref{eq:fundamental_estimate_for_item_i} vanishes 
and, thereafter, to combine this with \cref{eq:fundamental_estimate_for_item_i} 
and \cref{lem:convergence_in_probability_charac} 
(see \cref{eq:convergence_in_probability_item_i} below). 
We now first show that the second summand on 
the right hand side of 
\cref{eq:fundamental_estimate_for_item_i} vanishes. 
\Nobs that \cref{eq:N_n_definition} 
and the assumption that 
$
  \{ 
    \lim_{ n \to \infty } \gamma_n = 0
  \} = \Omega
$
show that 
for all $ \omega \in \Omega $ 
it holds that 
$
  \lim_{ n \to \infty }
  (
    N_{ n \delta }( \omega ) 
  ) 
  = \infty
$.
\Hence for all $ \omega \in \Omega $ that 
\begin{equation}
\label{eq:Nm_limit}
  \lim_{ m \to \infty }
  \bigl(
    N_{ m \delta }( \omega ) + 1
  \bigr) 
  = \infty
\end{equation}
\Moreover 
the assumption that 
$
  \{ 
    \lim_{ n \to \infty } \gamma_n = 0
  \} 
  = \Omega
$
ensures that 
for all $ \omega \in \Omega $ 
it holds that 
$
  \lim_{ n \to \infty }
  \gamma_n( \omega ) = 0
$. 
This and \cref{eq:Nm_limit} show 
for all $ \omega \in \Omega $ that
$
  \lim_{ m \to \infty }
  \gamma_{ N_{ m \delta }( \omega ) + 1 }( \omega ) = 0
$. 
\Hence 
\begin{equation}
\label{eq:almost_sure_convergence_gamma_N_m_delta}
\textstyle 
  \P\bigl(
    \lim_{ m \to \infty }
    \gamma_{ N_{ m \delta } + 1 } 
    = 0
  \bigr) 
  = 1
  .
\end{equation}
\cref{lem:almost_sure} 
\hence shows for all $ \varepsilon \in (0,\infty) $ that  
\begin{equation}
\label{eq:item_i_second_summand}
  \lim_{ m \to \infty }
  \P\bigl(
    \gamma_{ N_{ m \delta } + 1 } > \varepsilon
  \bigr) 
  = 0
  .
\end{equation}
In the next step we show that the third summand 
on the right hand side of \cref{eq:fundamental_estimate_for_item_i} vanishes. 
\Nobs that \cref{eq:definition_of_tau_R} ensures that 
for all $ R \in (0,\infty) $ it holds that 
\begin{equation}
  \{ \tau_R < \infty \} 
  =
  \bigl\{ 
    \exists \, n \in \N_0 \colon \| \Theta_n \|^2 \geq R
  \bigr\}
  \subseteq 
  \bigl\{ 
  \textstyle 
    \sup_{ n \in \N } \| \Theta_n \|^2 \geq R
  \big\} 
  .
\end{equation}
The fact that the measure $ \P $ is continuous from above 
\hence shows that 
\begin{equation}
\label{eq:item_i_third_summand}
\begin{split}
&
  \limsup_{ R \to \infty }
  \P\bigl(
    \cE 
    \cap 
    \{ 
      \tau_R < \infty 
    \} 
  \bigr)
\leq 
  \limsup_{ R \to \infty }
  \P\biggl(
    \cE 
    \cap 
    \biggl\{ 
      \sup_{ n \in \N } \| \Theta_n \|^2 \geq R
    \biggr\} 
  \biggr)
\\ & =
  \P\Biggl(
    \textstyle 
    \bigcap\limits_{ R \in (0,\infty) }
    \displaystyle
    \biggl(
      \cE 
      \cap 
      \biggl\{ 
        \sup_{ n \in \N } \| \Theta_n \|^2 \geq R
      \biggr\} 
    \biggr)
  \Biggr)
=
  \P\Biggl(
    \cE 
    \cap 
    \biggl(
      \textstyle 
      \bigcap\limits_{ R \in (0,\infty) }
      \displaystyle
      \biggl\{ 
        \sup_{ n \in \N } \| \Theta_n \|^2 \geq R
      \biggr\} 
    \biggr)
  \Biggr)
\\ & 
=
  \P\Bigl(
    \cE 
    \cap 
    \bigl\{ 
      \forall \, R \in (0,\infty) \colon 
\textstyle 
      \sup_{ n \in \N } \| \Theta_n \|^2 \geq R
    \bigr\}
  \Bigr)
=
  \P\bigl(
    \cE 
    \cap 
    \bigl\{ 
      \sup_{ n \in \N } \| \Theta_n \|^2 = \infty
    \bigr\}
  \bigr)
\\ & 
=
\textstyle 
  \P\bigl(
    \cE 
    \cap 
    \bigl\{ 
      \sup_{ n \in \N } \| \Theta_n \| = \infty
    \bigr\}
  \bigr)
  =
  \P\bigl(
    \bigl\{ 
      \sup_{ n \in \N } \| \Theta_n \| < \infty
    \bigr\}
    \cap 
    \bigl\{ 
      \sup_{ n \in \N } \| \Theta_n \| = \infty
    \bigr\}
  \bigr)
\\ & 
=
  \P( \emptyset ) = 0
  .
\end{split}
\end{equation}
In the next step we show that the third summand 
on the right hand side of \cref{eq:fundamental_estimate_for_item_i} vanishes. 
\Nobs that \cref{eq:Theta_n_recursion_in_prop}
shows for all $ n \in \N $ that
\begin{equation}
\begin{split}
  \| \Theta_n \|^2
& 
=
  \left\|
    \Theta_{ n - 1 } 
    +
    \gamma_n
    g( \Theta_{ n - 1 } )
    +
    \gamma_n
    D_n
  \right\|^2
\\ &
=
  \left< 
    \Theta_{ n - 1 } 
    +
    \gamma_n
    g( \Theta_{ n - 1 } )
    +
    \gamma_n
    D_n
    ,
    \Theta_{ n - 1 } 
    +
    \gamma_n
    g( \Theta_{ n - 1 } )
    +
    \gamma_n
    D_n
  \right>
\\ & =
  \left\| \Theta_{ n - 1 } \right\|^2
  +
  ( \gamma_n )^2
  \| g( \Theta_{ n - 1 } ) \|^2
  +
  ( \gamma_n )^2 
  \| D_n \|^2
\\ &
\quad 
  +
  2 \gamma_n
  \langle 
    \Theta_{ n - 1 }
    ,
    g( \Theta_{ n - 1 } )
  \rangle
  +
  2 ( \gamma_n )^2
  \langle 
    g( \Theta_{ n - 1 } )
    ,
    D_n
  \rangle
  +
  2 
  \gamma_n 
  \langle 
    \Theta_{ n - 1 } 
    ,
    D_n
  \rangle
  .
\end{split}
\end{equation}
Combining this with \cref{eq:coercivity_for_zero,eq:definition_cD_n} proves for all $ n \in \N $ that
\begin{equation}
\begin{split}
  \| \Theta_n \|^2
& 
\leq 
  ( 1 - 2 c \gamma_n )
  \left\| \Theta_{ n - 1 } \right\|^2
  +
  ( \gamma_n )^2
  \| g( \Theta_{ n - 1 } ) \|^2
  +
  ( \gamma_n )^2 
  \| D_n \|^2
\\ &
\quad
  +
  2 ( \gamma_n )^2
  \langle 
    g( \Theta_{ n - 1 } )
    ,
    D_n
  \rangle
  +
  2 
  \gamma_n 
  \langle 
    \Theta_{ n - 1 } 
    ,
    D_n
  \rangle
\\ & =
  ( 1 - 2 c \gamma_n )
  \left\| \Theta_{ n - 1 } \right\|^2
  +
  \cD_n
  .
\end{split}
\end{equation}
\Hence 
for all $ m, n \in \N_0 $ with $ m < n $ that
\begin{equation}
\begin{split}
  \| \Theta_n \|^2
& \leq 
  \left( 1 - 2 c \gamma_n \right)
  \left\| \Theta_{ n - 1 } \right\|^2
  +
  \cD_n
\\ & \leq 
  \left( 1 - 2 c \gamma_n \right)
  \left( 1 - 2 c \gamma_{ n - 1 } \right)
  \left\| \Theta_{ n - 2 } \right\|^2
  +
  \left( 1 - 2 c \gamma_n \right)
  \cD_{ n - 1 }
  +
  \cD_n
\\ & 
\leq 
  \dots 
\\ & \leq 
  \left[ 
    \prod_{ k = m + 1 }^n
    \left( 1 - 2 c \gamma_k \right)
  \right] 
  \left\| \Theta_m \right\|^2
  +
  \sum_{ k = m + 1 }^n
  \left[ 
    \prod_{ l = k + 1 }^n
    \left( 1 - 2 c \gamma_l \right)
  \right] 
  \cD_k
  .
\end{split}
\end{equation}
\Hence 
for all $ m, n \in \N $ with $ m \leq n $ that
\begin{equation}
\begin{split}
  \| \Theta_n \|^2
& \leq 
  \left[ 
    \prod_{ k = m }^n
    \left( 1 - 2 c \gamma_k \right)
  \right] 
  \left\| \Theta_{ m - 1 } \right\|^2
  +
  \sum_{ k = m }^n
  \left[ 
    \prod_{ l = k + 1 }^n
    \left( 1 - 2 c \gamma_l \right)
  \right] 
  \cD_k
  .
\end{split}
\end{equation}
Combining \cref{eq:in_proof_definition_of_ft}, 
\cref{eq:definition_cD_n}, and the fact that for all $ x \in \R $ it holds that 
$ 1 + x \leq \exp( x ) $ 
\hence proves that for all $ m, n \in \N $ with 
$ m \leq n $ it holds that 
\begin{equation}
	\label{prop:convergence:eq_decrease}
\begin{split}
  \| \Theta_n \|^2
& \leq
  \exp\!\left(
    - 2 c \left( \ft_n - \ft_{ m - 1 } \right)
  \right)
  \left\| \Theta_{ m - 1 } \right\|^2
\\ &
  +
  \left[ 
  \sum_{ k = m }^n
  \exp\!\left(
    - 2 c \left( \ft_n - \ft_k \right)
  \right)
  \left( \gamma_k \right)^2
    \bigl[  
      \| g( \Theta_{ k - 1 } ) \|^2
      +
      \| D_k \|^2
      +
      \| g( \Theta_{ k - 1 } ) \|
      \| D_k \|
    \bigr] 
  \right]
\\ &
  +
  \left[
    \sum_{ k = m }^{ n }
    \left[ 
      \prod_{ l = k + 1 }^n
      \left( 1 - 2 c \gamma_l \right)
    \right] 
    2
    \gamma_k 
    \langle 
      \Theta_{ k - 1 } 
      ,
      D_k
    \rangle
  \right] 
  .
\end{split}
\end{equation}
\Hence 
for all $ m, n \in \N_0 $ with 
$ m < n $ that 
\begin{equation}
\begin{split}
&
  \| \Theta_n \|^2
  \mathbbm{1}_{ 
    \{ 2 c \gamma_{ m + 1 } < 1 \}
  }
\leq
  \exp\!\left(
    - 2 c \left( \ft_n - \ft_m \right)
  \right)
  \left\| \Theta_m \right\|^2
  \mathbbm{1}_{ 
    \{ 2 c \gamma_{ m + 1 } < 1 \}
  }
\\ &
  +
  \left[ 
  \sum_{ k = m + 1 }^n
  \exp\!\left(
    - 2 c \left( \ft_n - \ft_k \right)
  \right)
  \left( \gamma_k \right)^2
    \bigl[  
      \| g( \Theta_{ k - 1 } ) \|^2
      +
      \| D_k \|^2
      +
      \| g( \Theta_{ k - 1 } ) \|
      \| D_k \|
    \bigr] 
    \mathbbm{1}_{ 
      \{ 2 c \gamma_{ m + 1 } < 1 \}
    }
  \right]
\\ &
  +
  \left[
    \sum_{ k = m + 1 }^n
    \left[ 
      \prod_{ l = k + 1 }^n
      \left( 1 - 2 c \gamma_l \right)
    \right]
    2 \gamma_k 
    \left< 
      \Theta_{ k - 1 } 
      ,
      D_k
    \right>
    \mathbbm{1}_{ 
      \{ 2 c \gamma_{ m + 1 } < 1 \}
    }
  \right] 
  .
\end{split}
\end{equation}
\Hence 
for all $ m, n \in \N_0 $ with 
$ m < n $ that 
\begin{equation}
\begin{split}
&
  \| \Theta_n \|^2
  \mathbbm{1}_{ 
    \{ 2 c \gamma_{ m + 1 } < 1 \}
  }
\leq
  \exp\!\left(
    - 2 c \left( \ft_n - \ft_m \right)
  \right)
  \left\| \Theta_m \right\|^2
  \mathbbm{1}_{ 
    \{ 2 c \gamma_{ m + 1 } < 1 \}
  }
\\ &
  +
  \left[ 
  \sum_{ k = m + 1 }^n
  \exp\!\left(
    - 2 c \left( \ft_n - \ft_k \right)
  \right)
  \left( \gamma_k \right)^2
    \bigl[  
      \| g( \Theta_{ k - 1 } ) \|^2
      +
      \| D_k \|^2
      +
      \| g( \Theta_{ k - 1 } ) \|
      \| D_k \|
    \bigr] 
    \mathbbm{1}_{ 
      \{ 2 c \gamma_{ m + 1 } < 1 \}
    }
  \right]
\\ &
  +
  \left[ 
    \prod_{ l = m + 1 }^n
    \left( 1 - 2 c \gamma_l \right)
  \right]
\\ &
\cdot 
  \left[
    \sum_{ k = m + 1 }^n
    \left[ 
      \prod_{ l = m + 1 }^k
      \left(         
        1 - 2 c \gamma_l \mathbbm{1}_{ \{ 2 c \gamma_{ m + 1 } < 1 \} } 
      \right)
    \right]^{ - 1 }
    2 \gamma_k 
    \left< 
      \Theta_{ k - 1 } 
      ,
      D_k
    \right>
    \mathbbm{1}_{ 
      \{ 2 c \gamma_{ m + 1 } < 1 \}
    }
  \right] 
  .
\end{split}
\end{equation}
The fact that for all $ x, y \in \R $ it holds that 
$
  x y \leq \frac{ 1 }{ 2 } ( x^2 + y^2 )
$
and the fact that for all $ x \in \R $ it holds that 
$ 1 + x \leq \exp( x ) $ 
\hence shows that for all $ m, n \in \N_0 $ with $ m < n $ 
it holds that 
\begin{equation}
\begin{split}
&
  \| \Theta_n \|^2
  \mathbbm{1}_{ 
    \{ 2 c \gamma_{ m + 1 } < 1 \}
  }
\leq
  \exp\!\left(
    - 2 c \left( \ft_n - \ft_m \right)
  \right)
  \left\| \Theta_m \right\|^2
  \mathbbm{1}_{ 
    \{ 2 c \gamma_{ m + 1 } < 1 \}
  }
\\ &
  +
  \frac{ 3 }{ 2 }
  \left[ 
  \sum_{ k = m + 1 }^n
  \exp\!\left(
    - 2 c \left( \ft_n - \ft_k \right)
  \right)
  \left( \gamma_k \right)^2
    \bigl[  
      \| g( \Theta_{ k - 1 } ) \|^2
      +
      \| D_k \|^2
    \bigr] 
    \mathbbm{1}_{ 
      \{ 2 c \gamma_{ m + 1 } < 1 \}
    }
  \right]
\\ &
  +
  2
  \exp\!\left(
    - 2 c ( \ft_n - \ft_m )
  \right)
\\ &
  \cdot 
  \left|
  \sum_{ k = m + 1 }^n
  \left[ 
    \prod_{ l = m + 1 }^k
    \left( 
      1 
      - 
      2 c \gamma_l 
      \mathbbm{1}_{ 
        \{ 2 c \gamma_{ m + 1 } < 1 \}
      }
    \right)^{ - 1 }
  \right]
  \gamma_k 
  \left< 
    \Theta_{ k - 1 } 
    ,
    D_k
  \right>
  \mathbbm{1}_{ 
    \{ 2 c \gamma_{ m + 1 } < 1 \}
  }
  \right|
  .
\end{split}
\end{equation}
The fact that for all $ n \in \N $ it holds that $ \gamma_{ n + 1 } \leq \gamma_n $ 
\hence shows 
that for all 
$ v, m, n \in \N_0 $ with $ v \leq m \leq n $ 
it holds that 
\begin{equation}
\label{eq:Theta_n_estimate_recursion}
\begin{split}
&
  \| \Theta_n \|^2
  \mathbbm{1}_{
    [ 0, ( 2 c )^{ - 1 } )
  }(
    \gamma_{ v + 1 }
  )
\leq
  \exp\!\left(
    - 2 c \left( \ft_n - \ft_m \right)
  \right)
  \left\| \Theta_m \right\|^2
  \mathbbm{1}_{
    [ 0, ( 2 c )^{ - 1 } )
  }(
    \gamma_{ v + 1 }
  )
\\ &
  +
  \frac{ 3 }{ 2 }
  \left[ 
  \sum_{ k = m + 1 }^n
  \exp\!\left(
    - 2 c \left( \ft_n - \ft_k \right)
  \right)
  \left( \gamma_k \right)^2
    \bigl[  
      \| g( \Theta_{ k - 1 } ) \|^2
      +
      \| D_k \|^2
    \bigr] 
    \mathbbm{1}_{
      [ 0, ( 2 c )^{ - 1 } )
    }(
      \gamma_{ v + 1 }
    )
  \right]
\\ &
  +
  2
  \exp\!\left(
    - 2 c ( \ft_n - \ft_m )
  \right)
\\ &
  \cdot 
  \left|
  \sum_{ k = m + 1 }^n
  \left[ 
    \prod_{ l = m + 1 }^k
    \left( 
      1 
      - 
      2 c \gamma_l 
      \mathbbm{1}_{
        [ 0, ( 2 c )^{ - 1 } )
      }(
        \gamma_{ v + 1 }
      )
    \right)^{ - 1 }
  \right]
  \gamma_k 
  \left< 
    \Theta_{ k - 1 } 
    ,
    D_k
  \right>
  \mathbbm{1}_{
    [ 0, ( 2 c )^{ - 1 } )
  }(
    \gamma_{ v + 1 }
  )
  \right|
  .
\end{split}
\end{equation}
\Hence 
for all 
$ R, \delta \in (0,\infty) $, 
$ m, n \in \N_0 $ 
with $ m < n $ that 
\begin{equation}
\begin{split}
&
  \| \Theta_{ \scrN^R_{ n \delta } } \|^2
  \mathbbm{1}_{
    [ 0, ( 2 c )^{ - 1 } )
  }(
    \gamma_{ \scrN^R_{ m \delta } + 1 }
  )
\leq
  \exp\bigl(
    - 2 c ( \ft_{ \scrN^R_{ n \delta } } - \ft_{ \scrN^R_{ ( n - 1 ) \delta } } )
  \bigr)
  \| \Theta_{ \scrN^R_{ ( n - 1 ) \delta } } \|^2
  \mathbbm{1}_{
    [ 0, ( 2 c )^{ - 1 } )
  }(
    \gamma_{ \scrN^R_{ m \delta } + 1 }
  )
\\ &
  +
  \frac{ 3 }{ 2 }
  \left[ 
  \sum_{ k = \scrN^R_{ ( n - 1 ) \delta } + 1 }^{ \scrN^R_{ n \delta } }
  \exp\bigl(
    - 2 c ( \ft_{ \scrN^R_{ n \delta } } - \ft_k )
  \bigr)
  \left( \gamma_k \right)^2
    \bigl[  
      \| g( \Theta_{ k - 1 } ) \|^2
      +
      \| D_k \|^2
    \bigr] 
    \mathbbm{1}_{ [ 0, ( 2 c )^{ - 1 } ) }(
      \gamma_{ \scrN^R_{ m \delta } + 1 }
    )
  \right]
\\ &
  +
  2
  \exp\bigl(
    - 2 c ( \ft_{ \scrN^R_{ n \delta } } - \ft_{ \scrN^R_{ ( n - 1 ) \delta } } )
  \bigr)
  \Biggl|
\textstyle 
  \sum\limits_{ k = \scrN^R_{ ( n - 1 ) \delta } + 1 }^{ \scrN^R_{ n \delta } }
  \left[ 
    \prod\limits_{ l = \scrN^R_{ ( n - 1 ) \delta } + 1 }^k
    \left( 
      1 
      - 
      2 c \gamma_l 
      \mathbbm{1}_{ [ 0, ( 2 c )^{ - 1 } ) }(
        \gamma_{ \scrN^R_{ m \delta } + 1 }
      )
    \right)^{ - 1 }
  \right]
\\ &
\cdot 
  \gamma_k 
  \left< 
    \Theta_{ k - 1 } 
    ,
    D_k
  \right>
    \mathbbm{1}_{ [ 0, ( 2 c )^{ - 1 } ) }(
      \gamma_{ \scrN^R_{ m \delta } + 1 }
    )
  \Biggr|
  .
\end{split}
\end{equation}
Induction \hence shows for all 
$ R, \delta \in (0,\infty) $, 
$ m, n \in \N_0 $ 
with $ m < n $ that 
\begin{equation}
\begin{split}
&
  \| \Theta_{ \scrN^R_{ n \delta } } \|^2
  \mathbbm{1}_{
    [ 0, ( 2 c )^{ - 1 } )
  }(
    \gamma_{ \scrN^R_{ m \delta } + 1 }
  )
\leq
  \exp\bigl(
    - 2 c ( \ft_{ \scrN^R_{ n \delta } } - \ft_{ \scrN^R_{ m \delta } } )
  \bigr)
  \| \Theta_{ \scrN^R_{ m \delta } } \|^2
  \mathbbm{1}_{
    [ 0, ( 2 c )^{ - 1 } )
  }(
    \gamma_{ \scrN^R_{ m \delta } + 1 }
  )
\\ &
  +
  \frac{ 3 }{ 2 }
  \left[ 
  \sum_{ k = \scrN^R_{ m \delta } + 1 }^{ \scrN^R_{ n \delta } }
  \exp\bigl(
    - 2 c ( \ft_{ \scrN^R_{ n \delta } } - \ft_k )
  \bigr)
  \left( \gamma_k \right)^2
    \bigl[  
      \| g( \Theta_{ k - 1 } ) \|^2
      +
      \| D_k \|^2
    \bigr] 
    \mathbbm{1}_{ [ 0, ( 2 c )^{ - 1 } ) }(
      \gamma_{ \scrN^R_{ m \delta } + 1 }
    )
  \right]
\\ &
  +
  2
  \Biggl[ 
  \sum_{ v = m }^{ n - 1 }
  \exp\bigl(
    - 2 c ( \ft_{ \scrN^R_{ n \delta } } - \ft_{ \scrN^R_{ v \delta } } )
  \bigr)
  \biggl|
\textstyle 
  \sum\limits_{ k = \scrN^R_{ v \delta } + 1 }^{ \scrN^R_{ (v + 1) \delta } }
  \biggl[ 
    \prod\limits_{ l = \scrN^R_{ v \delta } + 1 }^k
    \left( 
      1 
      - 
      2 c \gamma_l 
      \mathbbm{1}_{ [ 0, ( 2 c )^{ - 1 } ) }(
        \gamma_{ \scrN^R_{ m \delta } + 1 }
      )
    \right)^{ - 1 }
  \biggr]
\\ &
\cdot 
  \gamma_k
  \left< 
    \Theta_{ k - 1 } 
    ,
    D_k
  \right>
    \mathbbm{1}_{ [ 0, ( 2 c )^{ - 1 } ) }(
      \gamma_{ \scrN^R_{ m \delta } + 1 }
    )
  \biggr|
  \Biggr]
  .
\end{split}
\end{equation}
\Hence for all 
$ R, \delta \in (0,\infty) $, 
$ \varepsilon \in (0, ( 2 c )^{ - 1 } ) $, $ m, n \in \N_0 $ 
with $ m < n $ that 
\begin{equation}
\begin{split}
&
  \| \Theta_{ \scrN^R_{ n \delta } } \|^2
  \mathbbm{1}_{
    [ 0, \varepsilon ]
  }(
    \gamma_{ \scrN^R_{ m \delta } + 1 }
  )
\leq
  \exp\bigl(
    - 2 c ( \ft_{ \scrN^R_{ n \delta } } - \ft_{ \scrN^R_{ m \delta } } )
  \bigr)
  \| \Theta_{ \scrN^R_{ m \delta } } \|^2
  \mathbbm{1}_{ 
    [ 0, \varepsilon ]
  }(
    \gamma_{ \scrN^R_{ m \delta } + 1 }
  )
\\ &
  +
  \frac{ 3 }{ 2 }
  \left[ 
  \sum_{ k = \scrN^R_{ m \delta } + 1 }^{ \scrN^R_{ n \delta } }
  \exp\bigl(
    - 2 c ( \ft_{ \scrN^R_{ n \delta } } - \ft_k )
  \bigr)
  \left( \gamma_k \right)^2
    \bigl[  
      \scrG_R
      +
      \scrD _R
    \bigr] 
    \mathbbm{1}_{ 
      [ 0, \varepsilon ]
    }(
      \gamma_{ \scrN^R_{ m \delta } + 1 }
    )
  \right]
\\ &
  +
  2
  \Biggl[ 
  \sum_{ v = m }^{ n - 1 }
  \exp\bigl(
    - 2 c ( \ft_{ \scrN^R_{ n \delta } } - \ft_{ \scrN^R_{ v \delta } } )
  \bigr)
  \biggl|
\textstyle 
  \sum\limits_{ k = \scrN^R_{ v \delta } + 1 }^{ \scrN^R_{ (v + 1) \delta } }
  \biggl[ 
    \prod\limits_{ l = \scrN^R_{ v \delta } + 1 }^k
    \left( 
      1 
      - 
      2 c \gamma_l 
      \mathbbm{1}_{ 
        [ 0, \varepsilon ]
      }(
        \gamma_{ \scrN^R_{ m \delta } + 1 }
      )
    \right)^{ - 1 }
  \biggr]
\\ &
\cdot 
  \gamma_k 
  \left< 
    \Theta_{ k - 1 } 
    ,
    D_k
  \right>
    \mathbbm{1}_{ 
      [ 0, \varepsilon ]
    }(
      \gamma_{ \scrN^R_{ m \delta } + 1 }
    )
  \biggr|
  \Biggr]
  .
\end{split}
\end{equation}
\Hence for all 
$ R, \delta \in (0,\infty) $, 
$ \varepsilon \in (0, ( 2 c )^{ - 1 } ) $, $ m, n \in \N_0 $ 
with $ m < n $ that 
\begin{equation}
\begin{split}
&
  \| \Theta_{ \scrN^R_{ n \delta } } \|^2
  \mathbbm{1}_{
    \{ 
      \gamma_{ \scrN^R_{ m \delta } + 1 } \leq \varepsilon
    \}
    \cap 
    \{ \tau_R > N_{ n \delta } \}
  }
\leq
  \exp\bigl(
    - 2 c ( \ft_{ \scrN^R_{ n \delta } } - \ft_{ \scrN^R_{ m \delta } } )
  \bigr)
  \| \Theta_{ \scrN^R_{ m \delta } } \|^2
  \mathbbm{1}_{
    \{ 
      \gamma_{ \scrN^R_{ m \delta } + 1 } \leq \varepsilon
    \}
    \cap 
    \{ \tau_R > N_{ n \delta } \}
  }
\\ &
  +
  \frac{ 3 }{ 2 }
  \left[ 
  \sum_{ k = \scrN^R_{ m \delta } + 1 }^{ \scrN^R_{ n \delta } }
  \exp\bigl(
    - 2 c ( \ft_{ \scrN^R_{ n \delta } } - \ft_k )
  \bigr)
  \left( \gamma_k \right)^2
    \bigl[  
      \scrG_R
      +
      \scrD _R
    \bigr] 
    \mathbbm{1}_{
      \{ 
        \gamma_{ \scrN^R_{ m \delta } + 1 } \leq \varepsilon
      \}
      \cap 
      \{ \tau_R > N_{ n \delta } \}
    }
  \right]
\\ &
  +
  2
  \Biggl[ 
  \sum_{ v = m }^{ n - 1 }
  \exp\bigl(
    - 2 c ( \ft_{ \scrN^R_{ n \delta } } - \ft_{ \scrN^R_{ v \delta } } )
  \bigr)
  \biggl|
\textstyle 
  \sum\limits_{ k = \scrN^R_{ v \delta } + 1 }^{ \scrN^R_{ (v + 1) \delta } }
  \biggl[ 
    \prod\limits_{ l = \scrN^R_{ v \delta } + 1 }^k
    \left( 
      1 
      - 
      2 c \gamma_l 
      \mathbbm{1}_{ 
        [ 0, \varepsilon ]
      }(
        \gamma_{ \scrN^R_{ m \delta } + 1 }
      )
    \right)^{ - 1 }
  \biggr]
\\ &
\cdot 
  \gamma_k 
  \left< 
    \Theta_{ k - 1 } 
    ,
    D_k
  \right>
  \biggr|
  \mathbbm{1}_{
    \{ 
      \gamma_{ \scrN^R_{ m \delta } + 1 } \leq \varepsilon
    \}
    \cap 
    \{ \tau_R > N_{ n \delta } \}
  }
  \Biggr]
  .
\end{split}
\end{equation}
\Hence 
for all 
$ R, \delta \in (0,\infty) $, 
$ \varepsilon \in (0, ( 2 c )^{ - 1 } ) $, $ m, n \in \N_0 $ 
with $ m < n $ that 
\begin{equation}
\begin{split}
&
  \| \Theta_{ \scrN^R_{ n \delta } } \|^2
  \mathbbm{1}_{
    \{ 
      \gamma_{ \scrN^R_{ m \delta } + 1 } \leq \varepsilon
    \}
    \cap 
    \{ \tau_R > N_{ n \delta } \}
  }
\leq
  \exp\bigl(
    - 2 c ( \ft_{ \scrN^R_{ n \delta } } - \ft_{ \scrN^R_{ m \delta } } )
  \bigr)
  \| \Theta_{ \scrN^R_{ m \delta } } \|^2
  \mathbbm{1}_{
    \{ 
      \gamma_{ \scrN^R_{ m \delta } + 1 } \leq \varepsilon
    \}
    \cap 
    \{ \tau_R > N_{ n \delta } \}
  }
\\ &
  +
  \frac{ 3 }{ 2 }
  \left[ 
  \sum_{ k = \scrN^R_{ m \delta } + 1 }^{ \scrN^R_{ n \delta } }
  \exp\bigl(
    - 2 c ( \ft_{ \scrN^R_{ n \delta } } - \ft_k )
  \bigr)
  \left( \gamma_k \right)^2
    \bigl[  
      \scrG_R
      +
      \scrD _R
    \bigr] 
  \mathbbm{1}_{
    \{ 
      \gamma_{ \scrN^R_{ m \delta } + 1 } \leq \varepsilon
    \}
    \cap 
    \{ \tau_R > N_{ n \delta } \}
  }
  \right]
\\ &
  +
  2
  \Biggl[ 
  \sum_{ v = m }^{ n - 1 }
    \exp\bigl(
      - 2 c ( n \delta - [ v \delta + \varepsilon ] )  
    \bigr)
  \biggl|
\textstyle 
  \sum\limits_{ k = \scrN^R_{ v \delta } + 1 }^{ \scrN^R_{ (v + 1) \delta } }
  \biggl[ 
    \prod\limits_{ l = \scrN^R_{ v \delta } + 1 }^k
    \left( 
      1 
      - 
      2 c \gamma_l 
      \mathbbm{1}_{ 
        [ 0, \varepsilon ]
      }(
        \gamma_{ \scrN^R_{ m \delta } + 1 }
      )
    \right)^{ - 1 }
  \biggr]
\\ &
\cdot 
  \gamma_k
  \left< 
    \Theta_{ k - 1 } 
    ,
    D_k
  \right>
  \biggr|
  \mathbbm{1}_{
    \{ 
      \gamma_{ \scrN^R_{ m \delta } + 1 } \leq \varepsilon
    \}
    \cap 
    \{ \tau_R > N_{ n \delta } \}
  }
  \Biggr]
  .
\end{split}
\end{equation}
\Hence 
for all 
$ R, \delta \in (0,\infty) $, 
$ \varepsilon \in (0, ( 2 c )^{ - 1 } ) $, $ m, n \in \N_0 $ 
with $ m < n $ that 
\begin{equation}
\begin{split}
&
  \| \Theta_{ \scrN^R_{ n \delta } } \|^2
  \mathbbm{1}_{
    \{ 
      \gamma_{ \scrN^R_{ m \delta } + 1 } \leq \varepsilon
    \}
    \cap 
    \{ \tau_R > N_{ n \delta } \}
  }
\leq
  \exp\bigl(
    - 2 c ( \ft_{ \scrN^R_{ n \delta } } - \ft_{ \scrN^R_{ m \delta } } )
  \bigr)
  \| \Theta_{ \scrN^R_{ m \delta } } \|^2
  \mathbbm{1}_{
    \{ 
      \gamma_{ \scrN^R_{ m \delta } + 1 } \leq \varepsilon
    \}
    \cap 
    \{ \tau_R > N_{ n \delta } \}
  }
\\ &
  +
  \frac{ 3 }{ 2 }
  \left[ 
  \sum_{ k = \scrN^R_{ m \delta } + 1 }^{ \scrN^R_{ n \delta } }
  \exp\bigl(
    - 2 c ( \ft_{ \scrN^R_{ n \delta } } - \ft_k )
  \bigr)
  \left( \gamma_k \right)^2
    \bigl[  
      \scrG_R
      +
      \scrD _R
    \bigr] 
    \mathbbm{1}_{
      \{ 
        \gamma_{ \scrN^R_{ m \delta } + 1 } \leq \varepsilon
      \}
      \cap 
      \{ \tau_R > N_{ n \delta } \}
    }
  \right]
\\ &
  +
  2
  \Biggl[ 
  \sum_{ v = m }^{ n - 1 }
    \exp\bigl(
      - 2 c \bigl[ ( n - v ) \delta - \varepsilon \bigr]
    \bigr)
  \biggl|
\textstyle 
  \sum\limits_{ k = \scrN^R_{ v \delta } + 1 }^{ \scrN^R_{ (v + 1) \delta } }
  \biggl[ 
    \prod\limits_{ l = \scrN^R_{ v \delta } + 1 }^k
    \bigl( 
      1 
      - 
      2 c \gamma_l 
      \mathbbm{1}_{
        \{ 
          \gamma_{ \scrN^R_{ m \delta } + 1 } \leq \varepsilon
        \}
      }
    \bigr)^{ - 1 }
  \biggr]
\\ &
\cdot 
  \gamma_k 
  \left< 
    \Theta_{ k - 1 } 
    ,
    D_k
  \right>
      \mathbbm{1}_{
        \{ 
          \gamma_{ \scrN^R_{ m \delta } + 1 } \leq \varepsilon
        \}
        \cap 
        \{ \tau_R > N_{ m \delta } \}
      }
    \biggr|
  \Biggr]
  .
\end{split}
\end{equation}
\Hence 
for all 
$ R, \delta \in (0,\infty) $, $ \varepsilon \in (0, ( 2 c )^{ - 1 } ) $, $ m, n \in \N_0 $ 
with $ m < n $ that 
\begin{equation}
\begin{split}
&
  \E\bigl[ 
    \| \Theta_{ \scrN^R_{ n \delta } } \|^2
    \mathbbm{1}_{
      \{ 
        \gamma_{ \scrN^R_{ m \delta } + 1 } \leq \varepsilon
      \}
      \cap 
      \{ \tau_R > N_{ n \delta } \}
    }
  \bigr]
\leq
  \E\Bigl[ 
    \exp\bigl(
      - 2 c ( \ft_{ \scrN^R_{ n \delta } } - \ft_{ \scrN^R_{ m \delta } } )
    \bigr)
    \| \Theta_{ \scrN^R_{ m \delta } } \|^2
    \mathbbm{1}_{
      \{ 
        \gamma_{ \scrN^R_{ m \delta } + 1 } \leq \varepsilon
      \}
      \cap 
      \{ \tau_R > N_{ n \delta } \}
    }
  \Bigr]
\\ &
  +
  \frac{ 3 }{ 2 }
  \,
  \E\!\!\left[ 
  \sum_{ k = \scrN^R_{ m \delta } + 1 }^{ \scrN^R_{ n \delta } }
  \exp\bigl(
    - 2 c ( \ft_{ \scrN^R_{ n \delta } } - \ft_k )
  \bigr)
  \left( \gamma_k \right)^2
    \bigl[  
      \scrG_R
      +
      \scrD _R
    \bigr] 
    \mathbbm{1}_{ 
      \{ \gamma_{ \scrN^R_{ m \delta } + 1 } \leq \varepsilon \} 
      \cap 
      \{ \tau_R > N_{ n \delta } \}
    }
  \right]
\\ &
  +
  2
  \,
  \E\!\Biggl[ 
  \sum_{ v = m }^{ n - 1 }
    \exp\bigl(
      2 c \varepsilon - 2 c ( n - v ) \delta 
    \bigr)
  \biggl|
\textstyle 
  \sum\limits_{ k = \scrN^R_{ v \delta } + 1 }^{ \scrN^R_{ (v + 1) \delta } }
  \biggl[ 
    \prod\limits_{ l = \scrN^R_{ v \delta } + 1 }^k
    \bigl( 
      1 
      - 
      2 c \gamma_l 
      \mathbbm{1}_{ 
        \{ \gamma_{ \scrN^R_{ m \delta } + 1 } \leq \varepsilon \} 
      }
    \bigr)^{ - 1 }
  \biggr]
\\ &
\cdot 
    \gamma_k
    \left< 
      \Theta_{ k - 1 } 
      ,
      D_k
    \right>
      \mathbbm{1}_{ 
        \{ \gamma_{ \scrN^R_{ m \delta } + 1 } \leq \varepsilon \} 
        \cap 
        \{ \tau_R > N_{ m \delta } \}
      }
    \biggr|
  \Biggr]
  .
\end{split}
\end{equation}
\Hence 
for all 
$ R, \delta \in (0,\infty) $, $ \varepsilon \in (0, ( 2 c )^{ - 1 } ) $, $ m, n \in \N_0 $ 
with $ m < n $ that 
\begin{equation}
\begin{split}
&
  \E\bigl[ 
    \| \Theta_{ \scrN^R_{ n \delta } } \|^2
    \mathbbm{1}_{
      \{ 
        \gamma_{ \scrN^R_{ m \delta } + 1 } \leq \varepsilon
      \}
      \cap 
      \{ \tau_R > N_{ n \delta } \}
    }
  \bigr]
\\ &
\leq
  \E\bigl[ 
    \exp\bigl(
      - 2 c \bigl[ ( n - m ) \delta - \varepsilon \bigr]
    \bigr)
    \| \Theta_{ N_{ m \delta } } \|^2
    \mathbbm{1}_{
      \{ 
        \gamma_{ N_{ m \delta } + 1 } \leq \varepsilon
      \}
      \cap 
      \{ \tau_R > N_{ n \delta } \}
    }
  \bigr]
\\ &
  +
  2 \,
  \E\!\!\left[ 
  \sum_{ k = \scrN^R_{ m \delta } + 1 }^{ \scrN^R_{ n \delta } }
  \exp\bigl(
    - 2 c ( \ft_{ \scrN^R_{ n \delta } } - \ft_k )
  \bigr)
  \left( \gamma_k \right)^2
    \bigl[  
      \scrG_R
      +
      \scrD _R
    \bigr] 
    \mathbbm{1}_{ 
      \{ \gamma_{ \scrN^R_{ m \delta } + 1 } \leq \varepsilon \} 
      \cap 
      \{ \tau_R > N_{ n \delta } \}
    }
  \right]
\\ &
  +
  2
  \,
  \E\!\Biggl[ 
  \sum_{ v = m }^{ n - 1 }
    \exp\bigl(
      1 - 2 c ( n - v ) \delta 
    \bigr)
  \biggl|
\textstyle 
  \sum\limits_{ k = \scrN^R_{ v \delta } + 1 }^{ \scrN^R_{ (v + 1) \delta } }
  \biggl[ 
    \prod\limits_{ l = \scrN^R_{ v \delta } + 1 }^k
    \bigl( 
      1 
      - 
      2 c \gamma_l 
      \mathbbm{1}_{ 
        \{ \gamma_{ \scrN^R_{ m \delta } + 1 } \leq \varepsilon \} 
      }
    \bigr)^{ - 1 }
  \biggr]
\\ &
\cdot 
    \gamma_k
    \left< 
      \Theta_{ k - 1 } 
      ,
      D_k
    \right>
      \mathbbm{1}_{ 
        \{ \gamma_{ \scrN^R_{ m \delta } + 1 } \leq \varepsilon \} 
        \cap 
        \{ \tau_R > N_{ m \delta } \}
      }
    \biggr|
  \Biggr]
  .
\end{split}
\end{equation}
\Hence 
for all 
$ R, \delta \in (0,\infty) $, 
$ \varepsilon \in (0, ( 2 c )^{ - 1 } ) $, $ m, n \in \N_0 $ 
with $ m < n $ that 
\begin{equation}
\begin{split}
&
  \E\bigl[ 
    \| \Theta_{ \scrN^R_{ n \delta } } \|^2
    \mathbbm{1}_{
      \{ 
        \gamma_{ \scrN^R_{ m \delta } + 1 } \leq \varepsilon
      \}
      \cap 
      \{ \tau_R > N_{ n \delta } \}
    }
  \bigr]
\\ & \leq
  3
  \,
  \E\bigl[ 
    \exp\bigl(
      - 2 c ( n - m ) \delta 
    \bigr)
    \| \Theta_{ N_{ m \delta } } \|^2
    \mathbbm{1}_{
      \{ 
        \gamma_{ N_{ m \delta } + 1 } \leq \varepsilon
      \}
      \cap 
      \{ \tau_R > N_{ n \delta } \}
    }
  \bigr]
\\ &
  +
  \frac{ 3 }{ 2 }
  \,
  \E\!\!\left[ 
  \sum_{ k = \scrN^R_{ m \delta } + 1 }^{ \scrN^R_{ n \delta } }
  \exp\bigl(
    2 c \gamma_k
    - 2 c ( \ft_{ \scrN^R_{ n \delta } } - \ft_{ k - 1 } )
  \bigr)
  \left( \gamma_k \right)^2
    \bigl[  
      \scrG_R
      +
      \scrD _R
    \bigr] 
    \mathbbm{1}_{ 
      \{ \gamma_{ \scrN^R_{ m \delta } + 1 } \leq \varepsilon \} 
      \cap 
      \{ \tau_R > N_{ n \delta } \}
    }
  \right]
\\ &
  +
  6
  \,
  \E\!\Biggl[ 
  \sum_{ v = m }^{ n - 1 }
    \exp\bigl(
      - 2 c ( n - v ) \delta 
    \bigr)
  \biggl|
\textstyle 
  \sum\limits_{ k = \scrN^R_{ v \delta } + 1 }^{ \scrN^R_{ (v + 1) \delta } }
  \biggl[ 
    \prod\limits_{ l = \scrN^R_{ v \delta } + 1 }^k
    \bigl( 
      1 
      - 
      2 c \gamma_l 
      \mathbbm{1}_{ 
        \{ \gamma_{ \scrN^R_{ m \delta } + 1 } \leq \varepsilon \} 
      }
    \bigr)^{ - 1 }
  \biggr]
\\ &
\cdot 
    \gamma_k
    \left< 
      \Theta_{ k - 1 } 
      ,
      D_k
    \right>
      \mathbbm{1}_{ 
        \{ \gamma_{ \scrN^R_{ m \delta } + 1 } \leq \varepsilon \} 
        \cap 
        \{ \tau_R > N_{ m \delta } \}
      }
    \biggr|
  \Biggr]
  .
\end{split}
\end{equation}
The H\"{o}lder inequality \hence shows 
for all 
$ R, \delta \in (0,\infty) $, 
$ \varepsilon \in (0, ( 2 c )^{ - 1 } ) $, $ m, n \in \N_0 $ 
with $ m < n $ that 
\begin{equation}
\label{eq:item_i_first_summand_before_martingale_estimate}
\begin{split}
&
  \E\bigl[ 
    \| \Theta_{ \scrN^R_{ n \delta } } \|^2
    \mathbbm{1}_{
      \{ 
        \gamma_{ \scrN^R_{ m \delta } + 1 } \leq \varepsilon
      \}
      \cap 
      \{ \tau_R > N_{ n \delta } \}
    }
  \bigr]
\\ &
\leq
  3
  \,
  \E\bigl[ 
    \exp\bigl(
      - 2 c ( n - m ) \delta 
    \bigr)
    \| \Theta_{ N_{ m \delta } } \|^2
    \mathbbm{1}_{
      \{ 
        \gamma_{ N_{ m \delta } + 1 } \leq \varepsilon
      \}
      \cap 
      \{ \tau_R > N_{ n \delta } \}
    }
  \bigr]
\\ &
  +
  \frac{ 3 }{ 2 }
  \,
  \E\!\!\left[ 
  \sum_{ k = \scrN^R_{ m \delta } + 1 }^{ \scrN^R_{ n \delta } }
  \int_{ \ft_{ k - 1 } }^{ \ft_k }
  \exp\bigl(
    2 c \varepsilon
    - 2 c ( \ft_{ \scrN^R_{ n \delta } } - \ft_{ k - 1 } )
  \bigr)
  \gamma_k 
    \bigl[  
      \scrG_R
      +
      \scrD _R
    \bigr] 
    \mathbbm{1}_{ 
      \{ \gamma_{ \scrN^R_{ m \delta } + 1 } \leq \varepsilon \} 
      \cap 
      \{ \tau_R > N_{ n \delta } \}
    }
    \, \d s
  \right]
\\ &
  +
  6
  \sum_{ v = m }^{ n - 1 }
    \exp\bigl(
      - 2 c ( n - v ) \delta 
    \bigr)
  \Biggl(
  \E\!\Biggl[ 
  \biggl|
\textstyle 
  \sum\limits_{ k = \scrN^R_{ v \delta } + 1 }^{ \scrN^R_{ (v + 1) \delta } }
  \biggl[ 
    \prod\limits_{ l = \scrN^R_{ v \delta } + 1 }^k
    \bigl( 
      1 
      - 
      2 c \gamma_l 
      \mathbbm{1}_{ 
        \{ \gamma_{ \scrN^R_{ m \delta } + 1 } \leq \varepsilon \} 
      }
    \bigr)^{ - 1 }
  \biggr]
\\ &
\cdot 
    \gamma_k
    \left< 
      \Theta_{ k - 1 } 
      ,
      D_k
    \right>
      \mathbbm{1}_{ 
        \{ \gamma_{ \scrN^R_{ m \delta } + 1 } \leq \varepsilon \} 
        \cap 
        \{ \tau_R > N_{ m \delta } \}
      }
    \biggr|^2
  \Biggr]
  \Biggr)^{ \!\! \nicefrac{ 1 }{ 2 } }
  .
\end{split}
\end{equation}
\Moreover the fact that for all $ R \in \R $, $ k \in \N_0 $ it holds that 
$
  \{ \tau_R \leq k \} 
  =
  \{ 
    \exists \, n \in \{ 0, 1, \dots, k \} 
    \colon
    \| \Theta_n \|^2 \geq R
  \} 
  =
  (
    \cup_{ n = 0 }^k \{ 
      \| \Theta_n \|^2 \geq R
    \}
  )
  \in \mathbb{F}_k 
$
ensures that for all $ R \in \R $, $ m, k \in \N_0 $ it holds that
\begin{equation}
\begin{split}
&
  \{ \tau_R > N_{ m \delta } \}
  \cap 
  \{ 
    \scrN^R_{ m \delta } = k
  \} 
=
  \bigl\{ 
    \tau_R > \min\{ \tau_R, N_{ m \delta } \}
  \bigr\}
  \cap 
  \{ 
    \scrN^R_{ m \delta } = k
  \} 
\\ & 
=
  \{ \tau_R > \scrN^R_{ m \delta } \}
  \cap 
  \{ 
    \scrN^R_{ m \delta } = k
  \} 
=
  \{ \tau_R > k \}
  \cap 
  \{ 
    \scrN^R_{ m \delta } = k
  \} 
\\ & 
=
  \bigl(
  (
    \Omega \backslash 
    \{ \tau_R \leq k \}
  )
  \cap 
  \{ 
    \scrN^R_{ m \delta } = k
  \} 
  \bigr)
  \in \mathbb{F}_k
  . 
\end{split}
\end{equation}
\Hence for all $ R \in \R $, $ m \in \N_0 $ that
\begin{equation}
\label{eq:item_i_first_summand_measurability_0}
\begin{split}
  \{ \tau_R > N_{ m \delta } \}
  \in 
  \bigl\{ 
    A \in \cF \colon
    \bigl(
      \forall \, k \in \N_0 \colon 
      (
        A \cap \{ \scrN^R_{ m \delta } = k \} 
      )
      \in \mathbb{F}_k
    \bigr)
  \bigr\} 
  =
  \mathbb{F}_{ \scrN^R_{ m \delta } }
  .
\end{split}
\end{equation}
The fact that for all 
$ R, \delta \in (0,\infty) $, $ m, v \in \N_0 $
with $ m \leq v $ 
it holds that 
$
  \scrN^R_{ m \delta } 
  =
  \min\{ \tau_R, N_{ m \delta } \} 
  \leq 
  \min\{ \tau_R, N_{ v \delta } \} 
  =
  \scrN^R_{ v \delta } 
$
\hence 
shows that 
for all $ R, \delta \in (0,\infty) $, $ m, v \in \N_0 $
with $ m \leq v $ 
it holds that
\begin{equation}
\label{eq:item_i_first_summand_measurability_1}
  \{ \tau_R > N_{ m \delta  }\}
  \in 
  \mathbb{F}_{ \scrN^R_{ m \delta } }
  \subseteq 
  \mathbb{F}_{ \scrN^R_{ v \delta } }
  .
\end{equation}
\Moreover the assumption that for all $ k \in \N_0 $
it holds that 
$ \gamma_{ k + 1 } $ is $ \mathbb{F}_k $-measurable 
implies that 
for all $ R, \delta, \varepsilon \in (0,\infty) $, 
$ m, k \in \N_0 $
it holds that 
\begin{equation}
  \{ 
    \gamma_{ \scrN^R_{ m \delta } + 1 }
    \leq \varepsilon
  \}
  \cap 
  \{ 
    \scrN^R_{ m \delta } = k
  \} 
  =
  \bigl(
    \{ 
      \gamma_{ k + 1 }
      \leq \varepsilon
    \} 
    \cap 
    \{ 
      \scrN^R_{ m \delta } = k
    \} 
  \bigr)
  \in 
  \mathbb{F}_k 
  .
\end{equation}
This, 
\cref{eq:item_i_first_summand_measurability_0}, 
and 
\cref{eq:item_i_first_summand_measurability_1} 
show for all 
$ R, \delta, \varepsilon \in (0,\infty) $, $ m, v \in \N_0 $ 
with $ m \leq v $ that 
\begin{equation}
  \{ 
    \gamma_{ \scrN^R_{ m \delta } + 1 }
    \leq \varepsilon
  \}
  \in 
  \mathbb{F}_{ \scrN^R_{ m \delta } }
  \subseteq 
  \mathbb{F}_{ \scrN^R_{ v \delta } }
  .
\end{equation}
Combining this and \cref{eq:item_i_first_summand_measurability_1} 
with \cref{eq:item_i_first_summand_before_martingale_estimate} 
shows for all 
$ R, \delta \in (0,\infty) $, 
$ \varepsilon \in (0, ( 2 c )^{ - 1 } ) $, $ m, n \in \N_0 $ 
with $ m < n $ that 
\begin{equation}
\begin{split}
&
  \E\bigl[ 
    \| \Theta_{ \scrN^R_{ n \delta } } \|^2
    \mathbbm{1}_{
      \{ 
        \gamma_{ \scrN^R_{ m \delta } + 1 } \leq \varepsilon
      \}
      \cap 
      \{ \tau_R > N_{ n \delta } \}
    }
  \bigr]
\\ &
\leq
  3
  \exp\bigl(
    - 2 c ( n - m ) \delta 
  \bigr)
  \,
  \E\bigl[ 
    \| \Theta_{ N_{ m \delta } } \|^2
    \mathbbm{1}_{
      \{ 
        \gamma_{ N_{ m \delta } + 1 } \leq \varepsilon
      \}
      \cap 
      \{ \tau_R > N_{ m \delta } \}
    }
  \bigr]
\\ &
  +
  5
  \,
  \E\!\!\left[ 
  \sum_{ k = \scrN^R_{ m \delta } + 1 }^{ \scrN^R_{ n \delta } }
  \int_{ \ft_{ k - 1 } }^{ \ft_k }
  \exp\bigl(
    - 2 c ( \ft_{ \scrN^R_{ n \delta } } - s )
  \bigr)
  \gamma_k 
    \bigl[  
      \scrG_R
      +
      \scrD _R
    \bigr] 
    \mathbbm{1}_{ 
      \{ \gamma_{ \scrN^R_{ m \delta } + 1 } \leq \varepsilon \} 
      \cap 
      \{ \tau_R > N_{ n \delta } \}
    }
    \, \d s
  \right]
\\ &
  +
  6
  \sum_{ v = m }^{ n - 1 }
    \exp\bigl(
      - 2 c ( n - v ) \delta 
    \bigr)
  \Biggl(
  \E\!\Biggl[ 
\textstyle 
  \sum\limits_{ k = \scrN^R_{ v \delta } + 1 }^{ \scrN^R_{ (v + 1) \delta } }
  \biggl|
  \biggl[ 
    \prod\limits_{ l = \scrN^R_{ v \delta } + 1 }^k
    \bigl( 
      1 
      - 
      2 c \gamma_l 
      \mathbbm{1}_{ 
        \{ \gamma_{ \scrN^R_{ m \delta } + 1 } \leq \varepsilon \} 
      }
    \bigr)^{ - 1 }
  \biggr]
\\ &
\cdot 
    \gamma_k
    \left< 
      \Theta_{ k - 1 } 
      ,
      D_k
    \right>
      \mathbbm{1}_{ 
        \{ \gamma_{ \scrN^R_{ m \delta } + 1 } \leq \varepsilon \} 
        \cap 
        \{ \tau_R > N_{ m \delta } \}
      }
    \biggr|^2
  \Biggr]
  \Biggr)^{ \!\! \nicefrac{ 1 }{ 2 } }
  .
\end{split}
\end{equation}
The Cauchy-Schwarz inequality \hence shows for all 
$ R, \delta \in (0,\infty) $, 
$ \varepsilon \in (0, ( 2 c )^{ - 1 } ) $, $ m, n \in \N_0 $ 
with $ m < n $ that 
\begin{equation}
\begin{split}
&
  \E\bigl[ 
    \| \Theta_{ \scrN^R_{ n \delta } } \|^2
    \mathbbm{1}_{
      \{ 
        \gamma_{ \scrN^R_{ m \delta } + 1 } \leq \varepsilon
      \}
      \cap 
      \{ \tau_R > N_{ n \delta } \}
    }
  \bigr]
\\ &
\leq
  3
  \exp\bigl(
    - 2 c ( n - m ) \delta 
  \bigr)
  \,
  \E\bigl[ 
    \| \Theta_{ N_{ m \delta } } \|^2
    \mathbbm{1}_{
      \{ 
        \gamma_{ N_{ m \delta } + 1 } \leq \varepsilon
      \}
      \cap 
      \{ \tau_R > N_{ m \delta } \}
    }
  \bigr]
\\ &
  +
  5
  \,
  \E\!\!\left[ 
  \sum_{ k = \scrN^R_{ m \delta } + 1 }^{ \scrN^R_{ n \delta } }
  \int_{ \ft_{ k - 1 } }^{ \ft_k }
  \exp\bigl(
    - 2 c ( \ft_{ \scrN^R_{ n \delta } } - s )
  \bigr)
    \bigl[  
      \scrG_R
      +
      \scrD _R
    \bigr] 
  \gamma_{ \scrN^R_{ m \delta } + 1 } 
    \mathbbm{1}_{ 
      \{ \gamma_{ \scrN^R_{ m \delta } + 1 } \leq \varepsilon \} 
      \cap 
      \{ \tau_R > N_{ n \delta } \}
    }
    \, \d s
  \right]
\\ &
  +
  6
  \sum_{ v = m }^{ n - 1 }
    \exp\bigl(
      - 2 c ( n - v ) \delta 
    \bigr)
  \Biggl(
  \E\!\Biggl[ 
\textstyle 
  \sum\limits_{ k = \scrN^R_{ v \delta } + 1 }^{ \scrN^R_{ (v + 1) \delta } }
  \biggl[ 
    \prod\limits_{ l = \scrN^R_{ v \delta } + 1 }^k
    \bigl( 
      1 
      - 
      2 c \gamma_l 
      \mathbbm{1}_{ 
        \{ \gamma_{ \scrN^R_{ m \delta } + 1 } \leq \varepsilon \} 
      }
    \bigr)^{ - 2 }
  \biggr]
\\ &
\cdot 
    ( \gamma_k )^2
    \|
      \Theta_{ k - 1 } 
    \|^2
    \|
      D_k
    \|^2
      \mathbbm{1}_{ 
        \{ \gamma_{ \scrN^R_{ m \delta } + 1 } \leq \varepsilon \} 
        \cap 
        \{ \tau_R > N_{ m \delta } \}
      }
  \Biggr]
  \Biggr)^{ \!\! \nicefrac{ 1 }{ 2 } }
  .
\end{split}
\end{equation}
The assumption that 
$   
  \cup_{ n = 1 }^{ \infty }
  \{ 
    \gamma_{ n + 1 } \leq \gamma_n
  \} 
  = \Omega 
$ 
and \cref{cor:exp_bound}
\hence show 
for all 
$ R, \delta \in (0,\infty) $, 
$ m, n \in \N_0 $, 
$ 
  \varepsilon \in ( 0, \infty ) 
$ 
with 
$ m < n $
and 
$
  2 c \varepsilon \leq \exp(-1) < 1
$
that
\begin{equation}
\begin{split}
&
  \E\bigl[ 
    \| \Theta_{ \scrN^R_{ n \delta } } \|^2
    \mathbbm{1}_{
      \{ 
        \gamma_{ \scrN^R_{ m \delta } + 1 } \leq \varepsilon
      \}
      \cap 
      \{ \tau_R > N_{ n \delta } \}
    }
  \bigr]
\\ &
\leq
  3
  \exp\bigl(
    - 2 c ( n - m ) \delta 
  \bigr)
  \,
  \E\bigl[ 
    \| \Theta_{ N_{ m \delta } } \|^2
    \mathbbm{1}_{
      \{ 
        \gamma_{ N_{ m \delta } + 1 } \leq \varepsilon
      \}
      \cap 
      \{ \tau_R > N_{ m \delta } \}
    }
  \bigr]
\\ &
  +
  5
  \,
  \E\!\biggl[ 
  \int_{ \ft_{ \scrN^R_{ m \delta } } }^{ \ft_{ \scrN^R_{ n \delta } } }
  \exp\bigl(
    - 2 c ( \ft_{ \scrN^R_{ n \delta } } - s )
  \bigr)
    \bigl[  
      \scrG_R
      +
      \scrD _R
    \bigr] 
    \gamma_{ \scrN^R_{ m \delta } + 1 }
    \mathbbm{1}_{ 
      \{ \gamma_{ \scrN^R_{ m \delta } + 1 } \leq \varepsilon \} 
      \cap 
      \{ \tau_R > N_{ n \delta } \}
    }
    \, \d s
  \biggr]
\\ &
  +
  6 R^{ \nicefrac{ 1 }{ 2 } }
  \sum_{ v = m }^{ n - 1 }
    \exp\bigl(
      - 2 c ( n - v ) \delta 
    \bigr)
  \Biggl(
  \E\!\Biggl[ 
\textstyle 
  \sum\limits_{ k = \scrN^R_{ v \delta } + 1 }^{ \scrN^R_{ (v + 1) \delta } }
  \biggl[ 
    \prod\limits_{ l = \scrN^R_{ v \delta } + 1 }^k
    \Bigl[
      \exp\bigl(
        \exp(1)
        2 c \gamma_l 
        \mathbbm{1}_{ 
          \{ \gamma_{ \scrN^R_{ m \delta } + 1 } \leq \varepsilon \} 
        }
      \bigr)
    \Bigr]^2
  \biggr]
\\ &
\cdot 
    \scrD _R
    ( \gamma_k )^2
      \mathbbm{1}_{ 
        \{ \gamma_{ \scrN^R_{ m \delta } + 1 } \leq \varepsilon \} 
        \cap 
        \{ \tau_R > N_{ m \delta } \}
      }
  \Biggr]
  \Biggr)^{ \!\! \nicefrac{ 1 }{ 2 } }
  .
\end{split}
\end{equation}
\Hence 
for all 
$ R, \delta \in (0,\infty) $, 
$ m, n \in \N_0 $, 
$ \varepsilon \in ( 0, ( 2 c \exp(1) )^{ - 1 } ] $
with 
$ m < n $
that
\begin{equation}
\begin{split}
&
  \E\bigl[ 
    \| \Theta_{ \scrN^R_{ n \delta } } \|^2
    \mathbbm{1}_{
      \{ 
        \gamma_{ \scrN^R_{ m \delta } + 1 } \leq \varepsilon
      \}
      \cap 
      \{ \tau_R > N_{ n \delta } \}
    }
  \bigr]
\\ &
\leq
  3
  \exp\bigl(
    - 2 c ( n - m ) \delta 
  \bigr)
  \,
  \E\bigl[ 
    \| \Theta_{ N_{ m \delta } } \|^2
    \mathbbm{1}_{
      \{ 
        \gamma_{ N_{ m \delta } + 1 } \leq \varepsilon
      \}
      \cap 
      \{ \tau_R > N_{ m \delta } \}
    }
  \bigr]
\\ &
  +
  5
  \,
  \E\!\!\left[ 
  \int_{ 0 }^{ \ft_{ \scrN^R_{ n \delta } } - \ft_{ \scrN^R_{ m \delta } } }
    \exp\bigl(
      - 2 c s
    \bigr)
    \bigl[  
      \scrG_R
      +
      \scrD _R
    \bigr] 
    \gamma_{ \scrN^R_{ m \delta } + 1 }
    \mathbbm{1}_{ 
      \{ \gamma_{ \scrN^R_{ m \delta } + 1 } \leq \varepsilon \} 
      \cap 
      \{ \tau_R > N_{ n \delta } \}
    }
    \, \d s
  \right]
\\ &
  +
  6 R^{ \nicefrac{ 1 }{ 2 } }
  \sum_{ v = m }^{ n - 1 }
    \exp\bigl(
      - 2 c ( n - v ) \delta 
    \bigr)
\\ &
  \cdot 
  \Biggl(
  \E\!\Biggl[ 
\textstyle 
  \sum\limits_{ k = \scrN^R_{ v \delta } + 1 }^{ \scrN^R_{ (v + 1) \delta } }
    \exp\biggl(
      \exp(1)
      4 c 
      \biggl[
        \sum\limits_{ l = \scrN^R_{ v \delta } + 1 }^k 
        \gamma_l 
      \biggr]
    \biggr)
    \scrD _R
    ( \gamma_k )^2
      \mathbbm{1}_{ 
        \{ \gamma_{ \scrN^R_{ m \delta } + 1 } \leq \varepsilon \} 
        \cap 
        \{ \tau_R > N_{ m \delta } \}
      }
  \Biggr]
  \Biggr)^{ \!\! \nicefrac{ 1 }{ 2 } }
  .
\end{split}
\end{equation}
\Hence 
for all 
$ R, \delta \in (0,\infty) $, 
$ m, n \in \N_0 $, 
$ \varepsilon \in ( 0, ( 2 c \exp(1) )^{ - 1 } ] $
with 
$ m < n $
that
\begin{equation}
\begin{split}
&
  \E\bigl[ 
    \| \Theta_{ \scrN^R_{ n \delta } } \|^2
    \mathbbm{1}_{
      \{ 
        \gamma_{ \scrN^R_{ m \delta } + 1 } \leq \varepsilon
      \}
      \cap 
      \{ \tau_R > N_{ n \delta } \}
    }
  \bigr]
\leq
  3
  \exp\bigl(
    - 2 c ( n - m ) \delta 
  \bigr)
  \,
  \E\bigl[ 
    \| \Theta_{ N_{ m \delta } } \|^2
    \mathbbm{1}_{
      \{ 
        \gamma_{ N_{ m \delta } + 1 } \leq \varepsilon
      \}
      \cap 
      \{ \tau_R > N_{ m \delta } \}
    }
  \bigr]
\\ &
  +
  5
  \,
  \E\!\!\left[ 
  \int_0^{ \infty }
    \exp(
      - 2 c s
    )
    \bigl[  
      \scrG_R
      +
      \scrD _R
    \bigr] 
    \gamma_{ \scrN^R_{ m \delta } + 1 }
    \mathbbm{1}_{ 
      \{ \gamma_{ \scrN^R_{ m \delta } + 1 } \leq \varepsilon \} 
      \cap 
      \{ \tau_R > N_{ n \delta } \}
    }
    \, \d s
  \right]
  +
  6 R^{ \nicefrac{ 1 }{ 2 } }
  \sum_{ v = m }^{ n - 1 }
    \exp\bigl(
      - 2 c ( n - v ) \delta 
    \bigr)
\\ & 
  \cdot 
  \Biggl(
  \E\!\Biggl[ 
\textstyle 
  \sum\limits_{ k = \scrN^R_{ v \delta } + 1 }^{ \scrN^R_{ (v + 1) \delta } }
    \exp\biggl(
      \exp(1)
      4 c \bigl[ \ft_k - \ft_{ \scrN^R_{ v \delta } } \bigr]
    \biggr)
    \scrD _R
    ( \gamma_k )^2
      \mathbbm{1}_{ 
        \{ \gamma_{ \scrN^R_{ m \delta } + 1 } \leq \varepsilon \} 
        \cap 
        \{ \tau_R > N_{ m \delta } \}
      }
  \Biggr]
  \Biggr)^{ \!\! \nicefrac{ 1 }{ 2 } }
  .
\end{split}
\end{equation}
\Hence 
for all 
$ R, \delta \in (0,\infty) $, 
$ m, n \in \N_0 $, 
$ \varepsilon \in ( 0, ( 2 c \exp(1) )^{ - 1 } ] $
with 
$ m < n $
that
\begin{equation}
\begin{split}
&
  \E\bigl[ 
    \| \Theta_{ \scrN^R_{ n \delta } } \|^2
    \mathbbm{1}_{
      \{ 
        \gamma_{ \scrN^R_{ m \delta } + 1 } \leq \varepsilon
      \}
      \cap 
      \{ \tau_R > N_{ n \delta } \}
    }
  \bigr]
\leq
  3
  \exp\bigl(
    - 2 c ( n - m ) \delta 
  \bigr)
  \,
  \E\bigl[ 
    \| \Theta_{ N_{ m \delta } } \|^2
    \mathbbm{1}_{
      \{ 
        \gamma_{ N_{ m \delta } + 1 } \leq \varepsilon
      \}
      \cap 
      \{ \tau_R > N_{ m \delta } \}
    }
  \bigr]
\\ &
  +
  \frac{ 5 }{ 2 c }
  \,
  \E\bigl[ 
    (  
      \scrG_R
      +
      \scrD _R
    ) 
    \gamma_{ \scrN^R_{ m \delta } + 1 }
    \mathbbm{1}_{ 
      \{ \gamma_{ \scrN^R_{ m \delta } + 1 } \leq \varepsilon \} 
      \cap 
      \{ \tau_R > N_{ n \delta } \}
    }
  \bigr]
  +
  6 R^{ \nicefrac{ 1 }{ 2 } }
  \sum_{ v = m }^{ n - 1 }
    \exp\bigl(
      - 2 c ( n - v ) \delta 
    \bigr)
\\ &
  \cdot 
  \Biggl(
  \E\!\Biggl[ 
\textstyle 
  \sum\limits_{ k = \scrN^R_{ v \delta } + 1 }^{ \scrN^R_{ (v + 1) \delta } }
    \exp\Bigl(
      \exp(1)
      4 c \bigl[ \ft_{ \scrN^R_{ ( v + 1 ) \delta } } - \ft_{ \scrN^R_{ v \delta } } \bigr]
    \Bigr)
    \scrD _R
    ( \gamma_k )^2
      \mathbbm{1}_{ 
        \{ \gamma_{ \scrN^R_{ m \delta } + 1 } \leq \varepsilon \} 
        \cap 
        \{ \tau_R > N_{ m \delta } \}
      }
  \Biggr]
  \Biggr)^{ \!\! \nicefrac{ 1 }{ 2 } }
  .
\end{split}
\end{equation}
\Hence 
for all 
$ R, \delta \in (0,\infty) $, 
$ m, n \in \N_0 $, 
$ \varepsilon \in ( 0, ( 2 c \exp(1) )^{ - 1 } ] $
with 
$ m < n $
that
\begin{equation}
\begin{split}
&
  \E\bigl[ 
    \| \Theta_{ \scrN^R_{ n \delta } } \|^2
    \mathbbm{1}_{
      \{ 
        \gamma_{ \scrN^R_{ m \delta } + 1 } \leq \varepsilon
      \}
      \cap 
      \{ \tau_R > N_{ n \delta } \}
    }
  \bigr]
\\ &
\leq
  3
  \exp\bigl(
    - 2 c ( n - m ) \delta 
  \bigr)
  \,
  \E\bigl[ 
    \| \Theta_{ N_{ m \delta } } \|^2
    \mathbbm{1}_{
      \{ 
        \gamma_{ N_{ m \delta } + 1 } \leq \varepsilon
      \}
      \cap 
      \{ \tau_R > N_{ m \delta } \}
    }
  \bigr]
\\ &
  +
  3 c^{ - 1 }
  \E\bigl[ 
    (  
      \scrG_R
      +
      \scrD _R
    ) 
    \gamma_{ \scrN^R_{ m \delta } + 1 }
    \mathbbm{1}_{ 
      \{ \gamma_{ \scrN^R_{ m \delta } + 1 } \leq \varepsilon \} 
      \cap 
      \{ \tau_R > N_{ n \delta } \}
    }
  \bigr]
  +
  6 R^{ \nicefrac{ 1 }{ 2 } }
  \sum_{ v = m }^{ n - 1 }
    \exp\bigl(
      - 2 c ( n - v ) \delta 
    \bigr)
\\ & \cdot 
  \Biggl(
  \E\!\Biggl[ 
\textstyle 
  \sum\limits_{ k = \scrN^R_{ v \delta } + 1 }^{ \scrN^R_{ (v + 1) \delta } }
    \exp\bigl(
      \exp(1)
      4 c [ \delta + \varepsilon ] 
    \bigr)
    \scrD _R
      \gamma_k 
      \gamma_{ \scrN^R_{ v \delta } + 1 }
      \mathbbm{1}_{ 
        \{ \gamma_{ \scrN^R_{ m \delta } + 1 } \leq \varepsilon \} 
        \cap 
        \{ \tau_R > N_{ m \delta } \}
      }
  \Biggr]
  \Biggr)^{ \!\! \nicefrac{ 1 }{ 2 } }
  .
\end{split}
\end{equation}
\Hence 
for all 
$ R, \delta \in (0,\infty) $, 
$ m, n \in \N_0 $, 
$ \varepsilon \in ( 0, ( 2 c \exp(1) )^{ - 1 } ] $
with 
$ m < n $
that
\begin{equation}
\begin{split}
&
  \E\bigl[ 
    \| \Theta_{ \scrN^R_{ n \delta } } \|^2
    \mathbbm{1}_{
      \{ 
        \gamma_{ \scrN^R_{ m \delta } + 1 } \leq \varepsilon
      \}
      \cap 
      \{ \tau_R > N_{ n \delta } \}
    }
  \bigr]
\leq
  3
  \exp\bigl(
    - 2 c ( n - m ) \delta 
  \bigr)
  \,
  \E\bigl[ 
    \| \Theta_{ N_{ m \delta } } \|^2
    \mathbbm{1}_{
      \{ 
        \gamma_{ N_{ m \delta } + 1 } \leq \varepsilon
      \}
      \cap 
      \{ \tau_R > N_{ m \delta } \}
    }
  \bigr]
\\ &
  +
  3 c^{ - 1 }
  \E\bigl[ 
    (  
      \scrG_R
      +
      \scrD _R 
    ) 
    \gamma_{ \scrN^R_{ m \delta } + 1 }
    \mathbbm{1}_{ 
      \{ \gamma_{ \scrN^R_{ m \delta } + 1 } \leq \varepsilon \} 
      \cap 
      \{ \tau_R > N_{ n \delta } \}
    }
  \bigr]
\\ &
  +
  6 R^{ \nicefrac{ 1 }{ 2 } }
  \sum_{ v = m }^{ n - 1 }
    \exp\bigl(
      \exp(1)
      2 c [ \delta + \varepsilon ] 
      - 2 c ( n - v ) \delta 
    \bigr)
\\ &
\cdot 
  \Bigl(
  \E\Bigl[ 
    \bigl( 
      \ft_{ \scrN^R_{ (v + 1) \delta } }
      -
      \ft_{ \scrN^R_{ v \delta } + 1 } 
    \bigr)
    \scrD _R
    \gamma_{ \scrN^R_{ v \delta } + 1 }
    \mathbbm{1}_{ 
      \{ \gamma_{ \scrN^R_{ m \delta } + 1 } \leq \varepsilon \} 
      \cap 
      \{ \tau_R > N_{ m \delta } \}
    }
  \Bigr]
  \Bigr)^{ \! \nicefrac{ 1 }{ 2 } }
  .
\end{split}
\end{equation}
\Hence 
for all 
$ R, \delta \in (0,\infty) $, 
$ m, n \in \N_0 $, 
$ \varepsilon \in ( 0, ( 2 c \exp(1) )^{ - 1 } ] $
with 
$ m < n $
that
\begin{equation}
\begin{split}
&
  \E\bigl[ 
    \| \Theta_{ \scrN^R_{ n \delta } } \|^2
    \mathbbm{1}_{
      \{ 
        \gamma_{ \scrN^R_{ m \delta } + 1 } \leq \varepsilon
      \}
      \cap 
      \{ \tau_R > N_{ n \delta } \}
    }
  \bigr]
\\ &
\leq
  3
  \exp\bigl(
    - 2 c ( n - m ) \delta 
  \bigr)
  \,
  \E\bigl[ 
    \| \Theta_{ N_{ m \delta } } \|^2
    \mathbbm{1}_{
      \{ 
        \gamma_{ N_{ m \delta } + 1 } \leq \varepsilon
      \}
      \cap 
      \{ \tau_R > N_{ m \delta } \}
    }
  \bigr]
\\ &
  +
  3 c^{ - 1 }
  \E\bigl[ 
    (  
      \scrG_R
      +
      \scrD _R
    ) 
    \gamma_{ \scrN^R_{ m \delta } + 1 }
    \mathbbm{1}_{ 
      \{ \gamma_{ \scrN^R_{ m \delta } + 1 } \leq \varepsilon \} 
      \cap 
      \{ \tau_R > N_{ m \delta } \}
    }
  \bigr]
\\ &
  +
  6 R^{ \nicefrac{ 1 }{ 2 } }
  \sum_{ v = m }^{ n - 1 }
    \exp\bigl(
      2 c \delta \exp(1) + 2 c \varepsilon \exp(1) 
      - 2 c ( n - v ) \delta 
    \bigr)
\\ &
  \cdot 
    [ \delta + \varepsilon ]^{ \nicefrac{ 1 }{ 2 } }
  \bigl(
  \E\bigl[ 
    \scrD _R
    \gamma_{ \scrN^R_{ v \delta } + 1 }
    \mathbbm{1}_{ 
      \{ \gamma_{ \scrN^R_{ m \delta } + 1 } \leq \varepsilon \} 
      \cap 
      \{ \tau_R > N_{ m \delta } \}
    }
  \bigr]
  \bigr)^{ \! \nicefrac{ 1 }{ 2 } }
  .
\end{split}
\end{equation}
\Hence 
for all 
$ R, \delta \in (0,\infty) $, 
$ m, n \in \N_0 $, 
$ \varepsilon \in ( 0, ( 2 c \exp(1) )^{ - 1 } ] $
with 
$ m < n $
that
\begin{equation}
\begin{split}
&
  \E\bigl[ 
    \| \Theta_{ \scrN^R_{ n \delta } } \|^2
    \mathbbm{1}_{
      \{ 
        \gamma_{ \scrN^R_{ m \delta } + 1 } \leq \varepsilon
      \}
      \cap 
      \{ \tau_R > N_{ n \delta } \}
    }
  \bigr]
\\ &
\leq
  3
  \exp\bigl(
    - 2 c ( n - m ) \delta 
  \bigr)
  \,
  \E\bigl[ 
    \| \Theta_{ N_{ m \delta } } \|^2
    \mathbbm{1}_{
      \{ 
        \gamma_{ N_{ m \delta } + 1 } \leq \varepsilon
      \}
      \cap 
      \{ \tau_R > N_{ m \delta } \}
    }
  \bigr]
\\ &
  +
  3 c^{ - 1 }
  \E\bigl[ 
    (  
      \scrG_R
      +
      \scrD _R
    ) 
    \gamma_{ N_{ m \delta } + 1 }
    \mathbbm{1}_{ 
      \{ \gamma_{ N_{ m \delta } + 1 } \leq \varepsilon \} 
      \cap 
      \{ \tau_R > N_{ m \delta } \}
    }
  \bigr]
\\ &
  +
  6 R^{ \nicefrac{ 1 }{ 2 } }
  \sum_{ v = 1 }^{ n - m }
    \exp\bigl(
      2 c \delta \exp(1) + 1 
      - 2 c v \delta 
    \bigr)
    [ \delta + \varepsilon ]^{ \nicefrac{ 1 }{ 2 } }
  \bigl(
  \E\bigl[ 
    \scrD _R
    \gamma_{ N_{ m \delta } + 1 }
    \mathbbm{1}_{ 
      \{ \gamma_{ N_{ m \delta } + 1 } \leq \varepsilon \} 
      \cap 
      \{ \tau_R > N_{ m \delta } \}
    }
  \bigr]
  \bigr)^{ \! \nicefrac{ 1 }{ 2 } }
  .
\end{split}
\end{equation}
\Hence 
for all 
$ R, \delta \in (0,\infty) $, 
$ m, n \in \N_0 $, 
$ \varepsilon \in ( 0, ( 2 c \exp(1) )^{ - 1 } ] $
with 
$ m < n $
that
\begin{equation}
\begin{split}
&
  \E\bigl[ 
    \| \Theta_{ \scrN^R_{ n \delta } } \|^2
    \mathbbm{1}_{
      \{ 
        \gamma_{ \scrN^R_{ m \delta } + 1 } \leq \varepsilon
      \}
      \cap 
      \{ \tau_R > N_{ n \delta } \}
    }
  \bigr]
\leq
  3
  \exp\bigl(
    - 2 c ( n - m ) \delta 
  \bigr)
  \,
  \E\bigl[ 
    \| \Theta_{ N_{ m \delta } } \|^2
    \mathbbm{1}_{
      \{ \tau_R > N_{ m \delta } \}
    }
  \bigr]
\\ &
  +
  3 c^{ - 1 }
  \E\bigl[ 
    (  
      \scrG_R
      +
      \scrD _R
    ) 
    \gamma_{ N_{ m \delta } + 1 }
    \mathbbm{1}_{ 
      \{ \gamma_{ N_{ m \delta } + 1 } \leq \varepsilon \} 
    }
  \bigr]
\\ &
  +
  6 R^{ \nicefrac{ 1 }{ 2 } }
  \exp\bigl(
    1 +
    2 c \delta \exp(1) 
  \bigr)
  [ \delta + \varepsilon ]^{ \nicefrac{ 1 }{ 2 } }
  \biggl[
\textstyle 
    \sum\limits_{ v = 1 }^{ \infty }
    \exp\bigl(
      - 2 c \delta v 
    \bigr)
  \biggr]
  \bigl(
  \E\bigl[ 
    \scrD _R
    \gamma_{ N_{ m \delta } + 1 }
    \mathbbm{1}_{ 
      \{ \gamma_{ N_{ m \delta } + 1 } \leq \varepsilon \} 
    }
  \bigr]
  \bigr)^{ \nicefrac{ 1 }{ 2 } }
  .
\end{split}
\end{equation}
\Hence 
for all 
$ R, \delta \in (0,\infty) $, 
$ m, n \in \N_0 $, 
$ \varepsilon \in ( 0, ( 2 c \exp(1) )^{ - 1 } ] $
with 
$ m < n $
that
\begin{equation}
\begin{split}
&
  \E\bigl[ 
    \| \Theta_{ \scrN^R_{ n \delta } } \|^2
    \mathbbm{1}_{
      \{ 
        \gamma_{ N_{ m \delta } + 1 } \leq \varepsilon
      \}
      \cap 
      \{ \tau_R > N_{ n \delta } \}
    }
  \bigr]
\leq
  3
  \exp\bigl(
    - 2 c ( n - m ) \delta 
  \bigr)
  \,
  \E\bigl[ 
    \| \Theta_{ N_{ m \delta } } \|^2
    \mathbbm{1}_{
      \{ \tau_R > N_{ m \delta } \}
    }
  \bigr]
\\ &
  +
  3 c^{ - 1 }
  \E\bigl[ 
    (  
      \scrG_R
      +
      \scrD _R
    ) 
    \gamma_{ N_{ m \delta } + 1 }
    \mathbbm{1}_{ 
      \{ \gamma_{ N_{ m \delta } + 1 } \leq \varepsilon \} 
    }
  \bigr]
\\ &
  +
  6 R^{ \nicefrac{ 1 }{ 2 } }
  \exp\bigl(
    1 +
    2 c \delta \exp(1) 
  \bigr)
  [ \delta + \varepsilon ]^{ \nicefrac{ 1 }{ 2 } }
  \bigl[
    \exp( 2 c \delta ) - 1
  \bigr]^{ - 1 }
  \bigl(
  \E\bigl[ 
    \scrD _R
    \gamma_{ N_{ m \delta } + 1 }
    \mathbbm{1}_{ 
      \{ \gamma_{ N_{ m \delta } + 1 } \leq \varepsilon \} 
    }
  \bigr]
  \bigr)^{ \nicefrac{ 1 }{ 2 } }
  .
\end{split}
\end{equation}
The fact that for all $ x \in \R $ it holds that 
$
  \exp( x ) - 1 \geq x 
$
\hence 
shows that for all 
$ R, \delta \in (0,\infty) $, 
$ m, n \in \N_0 $, 
$ \varepsilon \in ( 0, ( 2 c \exp(1) )^{ - 1 } ] $
with 
$ m < n $
it holds that
\begin{equation}
\begin{split}
&
  \E\bigl[ 
    \| \Theta_{ \scrN^R_{ n \delta } } \|^2
    \mathbbm{1}_{
      \{ 
        \gamma_{ N_{ m \delta } + 1 } \leq \varepsilon
      \}
      \cap 
      \{ \tau_R > N_{ n \delta } \}
    }
  \bigr]
\leq
  3
  \exp\bigl(
    - 2 c ( n - m ) \delta 
  \bigr)
  \,
  \E\bigl[ 
    \| \Theta_{ N_{ m \delta } } \|^2
    \mathbbm{1}_{
      \{ \tau_R > N_{ m \delta } \}
    }
  \bigr]
\\ &
  +
  3 c^{ - 1 }
  \E\bigl[ 
    (  
      \scrG_R
      +
      \scrD _R
    ) 
    \gamma_{ N_{ m \delta } + 1 }
    \mathbbm{1}_{ 
      \{ \gamma_{ N_{ m \delta } + 1 } \leq \varepsilon \} 
    }
  \bigr]
\\ &
  +
  6 R^{ \nicefrac{ 1 }{ 2 } }
  \exp\bigl(
    1 +
    2 c \delta \exp(1) 
  \bigr)
  [ \delta + \varepsilon ]^{ \nicefrac{ 1 }{ 2 } }
  [
    2 c \delta 
  ]^{ - 1 }
  \bigl(
  \E\bigl[ 
    \scrD _R
    \gamma_{ N_{ m \delta } + 1 }
    \mathbbm{1}_{ 
      \{ \gamma_{ N_{ m \delta } + 1 } \leq \varepsilon \} 
    }
  \bigr]
  \bigr)^{ \nicefrac{ 1 }{ 2 } }
  .
\end{split}
\end{equation}
\Hence 
for all 
$ R, \delta \in (0,\infty) $, 
$ m, n \in \N_0 $, 
$ \varepsilon \in ( 0, ( 2 c \exp(1) )^{ - 1 } ] $
with 
$ m < n $
that
\begin{equation}
\begin{split}
&
  \E\bigl[ 
    \| \Theta_{ \scrN^R_{ n \delta } } \|^2
    \mathbbm{1}_{
      \{ 
        \gamma_{ N_{ m \delta } + 1 } \leq \varepsilon
      \}
      \cap 
      \{ \tau_R > N_{ n \delta } \}
    }
  \bigr]
\leq
  3 R
  \exp\bigl(
    - 2 c ( n - m ) \delta 
  \bigr)
\\ &
\textstyle 
  +
  \frac{ 3 }{ c }
  \,
  \E\bigl[ 
    (  
      \scrG_R
      +
      \scrD _R
    ) 
    \gamma_{ N_{ m \delta } + 1 }
    \mathbbm{1}_{ 
      \{ \gamma_{ N_{ m \delta } + 1 } \leq \varepsilon \} 
    }
  \bigr]
  +
  \frac{ 9 \sqrt{R} }{ c \delta } 
  \exp(
    6 c \delta  
  )
  [ \delta + \varepsilon ]^{ \nicefrac{ 1 }{ 2 } }
  \bigl(
  \E\bigl[ 
    \scrD _R
    \gamma_{ N_{ m \delta } + 1 }
    \mathbbm{1}_{ 
      \{ \gamma_{ N_{ m \delta } + 1 } \leq \varepsilon \} 
    }
  \bigr]
  \bigr)^{ \nicefrac{ 1 }{ 2 } }
  .
\end{split}
\end{equation}
\Hence 
for all 
$ R, \delta \in (0,\infty) $, 
$ m, n \in \N_0 $, 
$ \varepsilon \in ( 0, ( 2 c \exp(1) )^{ - 1 } ] $
with 
$ m \leq n $
that
\begin{equation}
\begin{split}
&
  \E\bigl[ 
    \| \Theta_{ \scrN^R_{ n \delta } } \|^2
    \mathbbm{1}_{
      \{ 
        \gamma_{ N_{ m \delta } + 1 } \leq \varepsilon
      \}
      \cap 
      \{ \tau_R > N_{ n \delta } \}
    }
  \bigr]
\leq
  3 R
  \exp\bigl(
    - 2 c ( n - m ) \delta 
  \bigr)
\\ &
\textstyle 
  +
  \frac{ 3 }{ c }
  \,
  \E\bigl[ 
    (  
      \scrG_R
      +
      \scrD _R
    ) 
    \min\{ \varepsilon, \gamma_{ N_{ m \delta } + 1 } \}
  \bigr]
  +
  \frac{ 9 \sqrt{R} }{ c \delta } 
  \exp(
    6 c \delta  
  )
  [ \delta + \varepsilon ]^{ \nicefrac{ 1 }{ 2 } }
  \bigl(
  \E\bigl[ 
    \scrD _R
    \min\{ \varepsilon, \gamma_{ N_{ m \delta } + 1 } \}
  \bigr]
  \bigr)^{ \nicefrac{ 1 }{ 2 } }
  .
\end{split}
\end{equation}
Combining this
and \cref{eq:almost_sure_convergence_gamma_N_m_delta} 
with Lebesgue's theorem of dominated convergence 
demonstrates for all $ R, \delta \in (0,\infty) $, 
$ \varepsilon \in ( 0, ( 2 c \exp(1) )^{ - 1 } ] $
that 
\begin{equation}
\label{eq:item_i_first_summand}
  \limsup_{ m \to \infty }
  \limsup_{ n \to \infty }
  \E\bigl[ 
    \| \Theta_{ \scrN^R_{ n \delta } } \|^2
    \mathbbm{1}_{
      \{ 
        \gamma_{ N_{ m \delta } + 1 } \leq \varepsilon
      \}
      \cap 
      \{ \tau_R > N_{ n \delta } \}
    }
  \bigr]  
  = 0
  .
\end{equation}
This proves that the first summand on the right hand side of 
\cref{eq:fundamental_estimate_for_item_i} vanishes. 
We can thus combine \cref{eq:fundamental_estimate_for_item_i}
with \cref{eq:item_i_first_summand,eq:item_i_second_summand,eq:item_i_third_summand} 
to obtain that 
for all $ \rho, \delta \in (0,\infty) $
it holds that 
\begin{equation}
\label{eq:convergence_in_probability_item_i}
  \lim_{ n \to \infty }
  \P\bigl(
    \| \Theta_{ N_{ n \delta } } \| \mathbbm{1}_{ \cE }
    > \rho 
  \bigr)
  = 0
  .
\end{equation}
Combining this with \cref{lem:convergence_in_probability_charac} establishes \cref{item:i}. 
We now prove \cref{item:ii}. 
We show \cref{item:ii} through 
an application of \cref{lem:convergence_in_probability_charac}. 
We thus need to verify that for all 
for all $ \rho, \delta \in (0,\infty) $
it holds that 
$
  \lim_{ n \to \infty }
  \P\bigl(
    (
      \max_{ N_{ n \delta } \leq j \leq N_{ (n+1) \delta } } 
      \| \Theta_{ N_{ n \delta } } \| 
    ) 
    \mathbbm{1}_{ \cE } 
    > \rho 
  \bigr)
  = 0
$
(see \cref{eq:convergence_in_probability_item_i} below). 
For this \nobs that 
the fact that for all $ A, B \in \cF $
it holds that 
$
  \P( A ) = \P( A \cap B ) + \P( A \cap ( \Omega \backslash B ) )
$
implies that 
for all 
$ R, \rho, \delta \in (0,\infty) $, 
$ n \in \N_0 $
it holds that
\begin{equation}
\begin{split}
&
  \P\bigl(
    \bigl(
      \max\nolimits_{ 
        N_{ n \delta } 
        \leq 
        j 
        \leq 
        N_{ (n + 1) \delta } 
      }
      \| \Theta_j \|
    \bigr)
    \mathbbm{1}_{
      \cE
    }
    >
    \rho 
  \bigr)
=
  \P\bigl(
    \bigl\{ 
        \max\nolimits_{ 
          N_{ n \delta } 
          \leq 
          j 
          \leq 
          N_{ (n + 1) \delta } 
        }
        \| \Theta_j \|
      >
      \rho 
    \bigr\}
    \cap 
    \cE
  \bigr)
\\ &
\leq 
  \P\bigl(
    \bigl\{ 
        \max\nolimits_{ 
          N_{ n \delta } 
          \leq 
          j 
          \leq 
          N_{ (n + 1) \delta } 
        }
        \| \Theta_j \|
      >
      \rho 
    \bigr\}
    \cap 
    \cE
    \cap 
    \{  
      \tau_R = \infty 
    \}
  \bigr)
  +
  \P\bigl(
    \cE
    \cap 
    \{  
      \tau_R < \infty 
    \}
  \bigr)
\\ & \leq 
  \P\bigl(
    \bigl\{ 
        \max\nolimits_{ 
          N_{ n \delta } 
          \leq 
          j 
          \leq 
          N_{ (n + 1) \delta } 
        }
        \| \Theta_j \|
      >
      \rho 
    \bigr\}
    \cap 
    \{  
      \tau_R = \infty 
    \}
  \bigr)
  +
  \P\bigl(
    \cE
    \cap 
    \{  
      \tau_R < \infty 
    \}
  \bigr)
\\ & \leq 
  \P\bigl(
    \bigl\{ 
        \max\nolimits_{ 
          N_{ n \delta } 
          \leq 
          j 
          \leq 
          N_{ (n + 1) \delta } 
        }
        \| \Theta_j \|
      >
      \rho 
    \bigr\}
    \cap 
    \{  
      \tau_R = \infty 
    \}
    \cap 
    \{ 
      \gamma_{ N_{ n \delta } + 1 } \leq \varepsilon
    \}
  \bigr)
\\ &
\quad 
  +
  \P\bigl(
    \gamma_{ N_{ n \delta } + 1 } > \varepsilon
  \bigr)
  +
  \P\bigl(
    \cE
    \cap 
    \{  
      \tau_R < \infty 
    \}
  \bigr)
\\ & \leq 
  \P\bigl(
    \bigl\{ 
        \max\nolimits_{ 
          N_{ n \delta } 
          \leq 
          j 
          \leq 
          N_{ (n + 1) \delta } 
        }
        \| \Theta_j \|
      >
      \rho 
    \bigr\}
    \cap 
    \{  
      \tau_R = \infty 
    \}
  \bigr)
  +
  \P\bigl(
    \cE
    \cap 
    \{  
      \tau_R < \infty 
    \}
  \bigr)
  .
\end{split}
\end{equation}
Combining this and again the fact that for all $ A, B \in \cF $
it holds that 
$
  \P( A ) = \P( A \cap B ) + \P( A \cap ( \Omega \backslash B ) )
$
with the Markov inequality
shows that for all 
$ R, \rho, \delta \in (0,\infty) $, 
$ n \in \N_0 $
it holds that
\begin{equation}
\begin{split}
&
  \P\bigl(
    \bigl(
      \max\nolimits_{ 
        N_{ n \delta } 
        \leq 
        j 
        \leq 
        N_{ (n + 1) \delta } 
      }
      \| \Theta_j \|
    \bigr)
    \mathbbm{1}_{
      \cE
    }
    >
    \rho 
  \bigr)
\\ & \leq 
  \P\bigl(
    \bigl\{ 
        \max\nolimits_{ 
          N_{ n \delta } 
          \leq 
          j 
          \leq 
          N_{ (n + 1) \delta } 
        }
        \| \Theta_j \|
      >
      \rho 
    \bigr\}
    \cap 
    \{  
      \tau_R = \infty 
    \}
    \cap 
    \{ 
      \gamma_{ N_{ n \delta } + 1 } \leq \varepsilon
    \}
  \bigr)
\\ &
\quad 
  +
  \P\bigl(
    \gamma_{ N_{ n \delta } + 1 } > \varepsilon
  \bigr)
  +
  \P\bigl(
    \cE
    \cap 
    \{  
      \tau_R < \infty 
    \}
  \bigr)
\\ & =
  \P\bigl(
    \bigl(
      \max\nolimits_{ 
        N_{ n \delta } 
        \leq 
        j 
        \leq 
        N_{ (n + 1) \delta } 
      }
      \| \Theta_j \|
    \bigr)
      \mathbbm{1}_{
        \{  
          \tau_R = \infty 
        \}
        \cap 
        \{ 
          \gamma_{ N_{ n \delta } + 1 } \leq \varepsilon
        \}
      }
    >
    \rho 
  \bigr)
\\ &
\quad 
  +
  \P\bigl(
    \gamma_{ N_{ n \delta } + 1 } > \varepsilon
  \bigr)
  +
  \P\bigl(
    \cE
    \cap 
    \{  
      \tau_R < \infty 
    \}
  \bigr)
\\ & 
\leq 
  \rho^{ - 2 }
  \E\bigl[
    \bigl(
      \max\nolimits_{ 
        N_{ n \delta } 
        \leq 
        j 
        \leq 
        N_{ (n + 1) \delta } 
      }
      \| \Theta_j \|^2
    \bigr)
      \mathbbm{1}_{
        \{ 
          \gamma_{ N_{ n \delta } + 1 } \leq \varepsilon
        \}
        \cap 
        \{  
          \tau_R = \infty 
        \}
      }
  \bigr]
\\ &
\quad 
  +
  \P\bigl(
    \gamma_{ N_{ n \delta } + 1 } > \varepsilon
  \bigr)
  +
  \P\bigl(
    \cE
    \cap 
    \{  
      \tau_R < \infty 
    \}
  \bigr)
  .
\end{split}
\end{equation}
This, \cref{eq:item_i_second_summand}, and \cref{eq:item_i_third_summand} ensure 
for all $ \delta \in (0,\infty) $ that 
\begin{equation}
\label{eq:fundamental_estimate_for_item_ii}
\begin{split}
&
  \limsup_{ n \to \infty }
  \P\bigl(
    \bigl(
      \max\nolimits_{ 
        N_{ n \delta } 
        \leq 
        j 
        \leq 
        N_{ (n + 1) \delta } 
      }
      \| \Theta_j \|
    \bigr)
    \mathbbm{1}_{
      \cE
    }
    >
    \rho 
  \bigr)
\\ &
=
  \limsup_{ R \to \infty }
  \limsup_{ n \to \infty }
  \P\bigl(
    \bigl(
      \max\nolimits_{ 
        N_{ n \delta } 
        \leq 
        j 
        \leq 
        N_{ (n + 1) \delta } 
      }
      \| \Theta_j \|
    \bigr)
    \mathbbm{1}_{
      \cE
    }
    >
    \rho 
  \bigr)
\\ & 
\leq 
  \rho^{ - 2 }
  \biggl[
  \limsup_{ n \to \infty }
  \E\bigl[
    \bigl(
      \max\nolimits_{ 
        N_{ n \delta } 
        \leq 
        j 
        \leq 
        N_{ (n + 1) \delta } 
      }
      \| \Theta_j \|^2
    \bigr)
      \mathbbm{1}_{
        \{ 
          \gamma_{ N_{ n \delta } + 1 } \leq \varepsilon
        \}
        \cap 
        \{  
          \tau_R = \infty 
        \}
      }
  \bigr]
  \biggr]
\\ &
\quad 
  +
  \limsup_{ n \to \infty }
  \P\bigl(
    \gamma_{ N_{ n \delta } + 1 } > \varepsilon
  \bigr)
  +
  \limsup_{ R \to \infty }
  \P\bigl(
    \cE
    \cap 
    \{  
      \tau_R < \infty 
    \}
  \bigr)
\\ & 
=
  \rho^{ - 2 }
  \biggl[
  \limsup_{ n \to \infty }
  \E\bigl[
    \bigl(
      \max\nolimits_{ 
        N_{ n \delta } 
        \leq 
        j 
        \leq 
        N_{ (n + 1) \delta } 
      }
      \| \Theta_j \|^2
    \bigr)
      \mathbbm{1}_{
        \{ 
          \gamma_{ N_{ n \delta } + 1 } \leq \varepsilon
        \}
        \cap 
        \{  
          \tau_R = \infty 
        \}
      }
  \bigr]
  \biggr]
  .
\end{split}
\end{equation}
To establish \cref{item:ii}, 
it thus remains to prove that 
the right hand side of \cref{eq:fundamental_estimate_for_item_ii} vanishes
and, thereafter, to combine this with 
\cref{eq:fundamental_estimate_for_item_ii} 
and \cref{lem:convergence_in_probability_charac} 
(see \cref{eq:item_ii_almost_done} below). 
We thus now verify that 
the right hand side of \cref{eq:fundamental_estimate_for_item_ii} vanishes. 
For this \nobs that 
\cref{eq:Theta_n_estimate_recursion} shows 
that 
for all $ v, m, n \in \N_0 $ with $ v \leq m \leq n $ 
it holds that 
\begin{equation}
\begin{split}
&
  \| \Theta_n \|^2
  \mathbbm{1}_{
    [ 0, ( 2 c )^{ - 1 } )
  }(
    \gamma_{ v + 1 }
  )
\leq
  \left\| \Theta_m \right\|^2
  \mathbbm{1}_{
    [ 0, ( 2 c )^{ - 1 } )
  }(
    \gamma_{ v + 1 }
  )
\\ &
  +
  \frac{ 3 }{ 2 }
  \left[ 
  \sum_{ k = m + 1 }^n
  \left( \gamma_k \right)^2
    \bigl[  
      \| g( \Theta_{ k - 1 } ) \|^2
      +
      \| D_k \|^2
    \bigr] 
    \mathbbm{1}_{
      [ 0, ( 2 c )^{ - 1 } )
    }(
      \gamma_{ v + 1 }
    )
  \right]
\\ &
  +
  2
  \left|
  \sum_{ k = m + 1 }^n
  \left[ 
    \prod_{ l = m + 1 }^k
    \left( 
      1 
      - 
      2 c \gamma_l 
      \mathbbm{1}_{
        [ 0, ( 2 c )^{ - 1 } )
      }(
        \gamma_{ v + 1 }
      )
    \right)^{ - 1 }
  \right]
  \gamma_k 
  \left< 
    \Theta_{ k - 1 } 
    ,
    D_k
  \right>
  \mathbbm{1}_{
    [ 0, ( 2 c )^{ - 1 } )
  }(
    \gamma_{ v + 1 }
  )
  \right|
  .
\end{split}
\end{equation}
\Hence 
for all $ R \in (0,\infty) $, $ v, m, n \in \N_0 $ with $ v \leq m \leq n $ 
it holds that 
\begin{equation}
\begin{split}
&
  \| \Theta_n \|^2
  \mathbbm{1}_{
    [ 0, ( 2 c )^{ - 1 } )
  }(
    \gamma_{ v + 1 }
  )
  \mathbbm{1}_{
    \{ 
      \tau_R = \infty 
    \} 
  }
\leq
  \left\| \Theta_m \right\|^2
  \mathbbm{1}_{
    [ 0, ( 2 c )^{ - 1 } )
  }(
    \gamma_{ v + 1 }
  )
  \mathbbm{1}_{
    \{ 
      \tau_R = \infty 
    \} 
  }
\\ &
  +
  \frac{ 3 }{ 2 }
  \left[ 
  \sum_{ k = m + 1 }^n
  \left( \gamma_k \right)^2
    \bigl[  
      \| g( \Theta_{ k - 1 } ) \|^2
      +
      \norm{D_k} ^2
    \bigr] 
    \mathbbm{1}_{
      [ 0, ( 2 c )^{ - 1 } )
    }(
      \gamma_{ v + 1 }
    )
    \mathbbm{1}_{
      \{ 
        \tau_R = \infty 
      \} 
    }
  \right]
\\ &
  +
  2
  \left|
  \sum_{ k = m + 1 }^n
  \left[ 
    \prod_{ l = m + 1 }^k
    \left( 
      1 
      - 
      2 c \gamma_l 
      \mathbbm{1}_{
        [ 0, ( 2 c )^{ - 1 } )
      }(
        \gamma_{ v + 1 }
      )
    \right)^{ - 1 }
  \right]
  \gamma_k 
  \left< 
    \Theta_{ k - 1 } 
    ,
    D_k
  \right>
  \mathbbm{1}_{
    [ 0, ( 2 c )^{ - 1 } )
  }(
    \gamma_{ v + 1 }
  )
  \right|
  \mathbbm{1}_{
    \{ 
      \tau_R = \infty 
    \} 
  }
  .
\end{split}
\end{equation}
\cref{eq:def_GR_DR} \hence proves
for all $ R \in (0,\infty) $, $ v, m, n \in \N_0 $ with $ v \leq m \leq n $ 
that 
\begin{equation}
\begin{split}
&
  \| \Theta_n \|^2
  \mathbbm{1}_{
    [ 0, ( 2 c )^{ - 1 } )
  }(
    \gamma_{ v + 1 }
  )
  \mathbbm{1}_{
    \{ 
      \tau_R = \infty 
    \} 
  }
\leq
  \left\| \Theta_m \right\|^2
  \mathbbm{1}_{
    [ 0, ( 2 c )^{ - 1 } )
  }(
    \gamma_{ v + 1 }
  )
  \mathbbm{1}_{
    \{ 
      \tau_R = \infty 
    \} 
  }
\\ &
  +
  \frac{ 3 }{ 2 }
  \bigl[ 
    ( \ft_n - \ft_m )
    \gamma_{ m + 1 }
    \bigl[  
      \scrG_R
      +
      \scrD_R
    \bigr] 
    \mathbbm{1}_{
      [ 0, ( 2 c )^{ - 1 } )
    }(
      \gamma_{ v + 1 }
    )
    \mathbbm{1}_{
      \{ 
        \tau_R = \infty 
      \} 
    }
  \bigr]
\\ &
  +
  2
  \left|
  \sum_{ k = m + 1 }^n
  \left[ 
    \prod_{ l = m + 1 }^k
    \left( 
      1 
      - 
      2 c \gamma_l 
      \mathbbm{1}_{
        [ 0, ( 2 c )^{ - 1 } )
      }(
        \gamma_{ v + 1 }
      )
    \right)^{ - 1 }
  \right]
  \gamma_k 
  \left< 
    \Theta_{ k - 1 } 
    ,
    D_k
  \right>
  \mathbbm{1}_{
    [ 0, ( 2 c )^{ - 1 } )
  }(
    \gamma_{ v + 1 }
  )
  \right|
  \mathbbm{1}_{
    \{ 
      \tau_R = \infty 
    \} 
  }
  .
\end{split}
\end{equation}
\Hence 
for all $ R \in (0,\infty) $, $ v, m, n \in \N_0 $ with $ v \leq m \leq n $ 
that 
\begin{equation}
\begin{split}
&
  \| \Theta_n \|^2
  \mathbbm{1}_{
    [ 0, ( 2 c )^{ - 1 } )
  }(
    \gamma_{ v + 1 }
  )
  \mathbbm{1}_{
    \{ 
      \tau_R = \infty 
    \} 
  }
\leq
  \min\{ R, \| \Theta_m \|^2 \}
  \mathbbm{1}_{
    [ 0, ( 2 c )^{ - 1 } )
  }(
    \gamma_{ v + 1 }
  )
  \mathbbm{1}_{
    \{ 
      \tau_R = \infty 
    \} 
  }
\\ &
  +
  \tfrac{ 3 }{ 2 }
    ( \ft_n - \ft_m )
    \gamma_{ m + 1 }
    (
      \scrG_R
      +
      \scrD _R
    ) 
    \mathbbm{1}_{
      [ 0, ( 2 c )^{ - 1 } )
    }(
      \gamma_{ v + 1 }
    )
\\ &
  +
  2
  \left|
  \sum_{ k = m + 1 }^n
  \left[ 
    \prod_{ l = m + 1 }^k
    \left( 
      1 
      - 
      2 c \gamma_l 
      \mathbbm{1}_{
        [ 0, ( 2 c )^{ - 1 } )
      }(
        \gamma_{ v + 1 }
      )
    \right)^{ - 1 }
  \right]
  \gamma_k 
  \left< 
    \Theta_{ k - 1 } 
    ,
    D_k
  \right>
  \mathbbm{1}_{
    [ 0, ( 2 c )^{ - 1 } )
  }(
    \gamma_{ v + 1 }
  )
  \right|
  \mathbbm{1}_{
    \{ 
      \tau_R = \infty 
    \} 
  }
  .
\end{split}
\end{equation}
\Hence 
for all $ R, \delta \in (0,\infty) $, $ n \in \N_0 $, 
$ \varepsilon \in ( 0, ( 2 c )^{ - 1 } ) $
that 
\begin{equation}
\begin{split}
&
  \bigl(
    \max\nolimits_{ 
      \scrN^R_{ n \delta } 
      \leq 
      j 
      \leq 
      \scrN^R_{ (n + 1) \delta } 
    }
    \| \Theta_j \|^2
  \bigr)
  \mathbbm{1}_{
    [ 0, \varepsilon ]
  }(
    \gamma_{ \scrN^R_{ n \delta } + 1 }
  )
  \mathbbm{1}_{
    \{ 
      \tau_R = \infty 
    \} 
  }
\\ &
\leq
  \min\bigl\{ 
    R, 
    \| \Theta_{ \scrN^R_{ n \delta } } \|^2 
  \bigr\}
  \mathbbm{1}_{
    [ 0, \varepsilon ]
  }(
    \gamma_{ \scrN^R_{ n \delta } + 1 }
  )
  \mathbbm{1}_{
    \{ 
      \tau_R = \infty 
    \} 
  }
\\ &
  +
  \tfrac{ 3 }{ 2 }
    \bigl(
      \ft_{ 
        \scrN^R_{ (n + 1) \delta }
      } 
      - 
      \ft_{
        \scrN^R_{ n \delta }
      }
    \bigr)
    \gamma_{ 
      \scrN^R_{ n \delta }
      + 1 
    }
    (
      \scrG_R
      +
      \scrD _R
    ) 
    \mathbbm{1}_{
    [ 0, \varepsilon ]
    }(
      \gamma_{ 
        \scrN^R_{ n \delta }
        + 1 
      }
    )
    \mathbbm{1}_{
      \{ 
        \tau_R = \infty 
      \} 
    }
\\ &
\textstyle 
  +
  2
  \Biggl[
  \max\limits_{ 
    \scrN^R_{ n \delta } \leq j \leq \scrN^R_{ (n + 1) \delta }
  }
  \biggl|
  \sum\limits_{ k = \scrN^R_{ n \delta } + 1 }^{ j }
  \biggl[ 
    \prod\limits_{ l = \scrN^R_{ n \delta } + 1 }^k
    \bigl( 
      1 
      - 
      2 c \gamma_l 
      \mathbbm{1}_{
        [ 0, \varepsilon ]
      }(
        \gamma_{ \scrN^R_{ n \delta } + 1 }
      )
    \bigr)^{ - 1 }
  \biggr]
\\ &
\cdot 
    \gamma_k 
    \left< 
      \Theta_{ k - 1 } 
      ,
      D_k
    \right>
    \mathbbm{1}_{
      [ 0, \varepsilon ]
    }(
      \gamma_{ \scrN^R_{ n \delta } + 1 }
    )
  \biggr|
  \Biggr]
  \mathbbm{1}_{
    \{ 
      \tau_R = \infty 
    \} 
  }
  .
\end{split}
\end{equation}
\Hence 
for all $ R, \delta \in (0,\infty) $, $ n \in \N_0 $, 
$
  \varepsilon \in ( 0, ( 2 c )^{ - 1 } )
$
that 
\begin{equation}
\begin{split}
&
  \E\bigl[
    \bigl(
      \max\nolimits_{ 
        N_{ n \delta } 
        \leq 
        j 
        \leq 
        N_{ (n + 1) \delta } 
      }
      \| \Theta_j \|^2
    \bigr)
    \mathbbm{1}_{
      \{ 
        \gamma_{ \scrN^R_{ n \delta } + 1 }
        \leq 
        \varepsilon
      \} 
      \cap 
      \{ 
        \tau_R = \infty 
      \} 
    }
  \bigr]
\\ &
\leq
  \E\bigl[ 
  \min\bigl\{ 
    R, 
    \| \Theta_{ N_{ n \delta } } \|^2 
  \bigr\}
  \mathbbm{1}_{
    \{ 
      \tau_R = \infty 
    \} 
  }
  \bigr]
\\ &
  +
  \tfrac{ 3 }{ 2 }
  \,
  \E\bigl[
    \bigl(
      \ft_{ 
        \scrN^R_{ (n + 1) \delta }
      } 
      - 
      \ft_{
        \scrN^R_{ n \delta }
      }
    \bigr)
    \gamma_{ 
      \scrN^R_{ n \delta }
      + 1 
    }
    (
      \scrG_R
      +
      \scrD _R
    ) 
    \mathbbm{1}_{
      \{ 
        \gamma_{ \scrN^R_{ n \delta } + 1 }
        \leq 
        \varepsilon
      \} 
      \cap 
      \{ 
        \tau_R = \infty 
      \} 
    }
  \bigr]
\\ &
\textstyle 
  +
  2 
  \Biggl(
  \E\Biggl[
  \max\limits_{ 
    \scrN^R_{ n \delta } \leq j \leq \scrN^R_{ (n + 1) \delta }
  }
  \biggl|
  \sum\limits_{ k = \scrN^R_{ n \delta } + 1 }^{ j }
  \biggl[ 
    \prod\limits_{ l = \scrN^R_{ n \delta } + 1 }^k
    \bigl( 
      1 
      - 
      2 c \gamma_l 
      \mathbbm{1}_{
        [ 0, \varepsilon ]
      }(
        \gamma_{ \scrN^R_{ n \delta } + 1 }
      )
    \bigr)^{ - 1 }
  \biggr]
\\ &
\cdot 
    \gamma_k 
    \left< 
      \Theta_{ k - 1 } 
      ,
      D_k
    \right>
    \mathbbm{1}_{
      \{ 
        \gamma_{ \scrN^R_{ n \delta } + 1 }
        \leq 
        \varepsilon
      \} 
      \cap 
      \{ 
        \tau_R = \infty 
      \} 
    }
  \biggr|^2
  \Biggr]
  \Biggr)^{ \!\! \nicefrac{ 1 }{ 2 } }
  .
\end{split}
\end{equation}
\Hence 
for all $ R, \delta \in (0,\infty) $, $ n \in \N_0 $, 
$ \varepsilon \in ( 0, ( 2 c )^{ - 1 } ) $
that 
\begin{equation}
\begin{split}
&
  \E\bigl[
    \bigl(
      \max\nolimits_{ 
        N_{ n \delta } 
        \leq 
        j 
        \leq 
        N_{ (n + 1) \delta } 
      }
      \| \Theta_j \|^2
    \bigr)
    \mathbbm{1}_{
      \{ 
        \gamma_{ \scrN^R_{ n \delta } + 1 }
        \leq 
        \varepsilon
      \} 
      \cap 
      \{ 
        \tau_R = \infty 
      \} 
    }
  \bigr]
\\ &
\leq
  \E\bigl[ 
    \min\bigl\{ 
      R, 
      \| \Theta_{ N_{ n \delta } } \|^2 
    \bigr\}
    \mathbbm{1}_{
      \{ 
        \sup_{ n \in \N_0 } \| \Theta_n \| < \infty 
      \} 
    }
  \bigr]
\\ &
  +
  \tfrac{ 3 }{ 2 }
  \,
  \E\bigl[
    \bigl(
      \ft_{ 
        N_{ (n + 1) \delta }
      } 
      - 
      \ft_{
        N_{ n \delta }
      }
    \bigr)
    \gamma_{ 
      N_{ n \delta }
      + 1 
    }
    (
      \scrG_R
      +
      \scrD _R
    ) 
    \mathbbm{1}_{
      \{ 
        \gamma_{ \scrN^R_{ n \delta } + 1 }
        \leq 
        \varepsilon
      \} 
      \cap 
      \{ 
        \tau_R = \infty 
      \} 
    }
  \bigr]
\\ &
\textstyle 
  +
  2 
  \Biggl(
  \E\Biggl[
  \max\limits_{ 
    \scrN^R_{ n \delta } \leq j \leq \scrN^R_{ (n + 1) \delta }
  }
  \biggl|
  \sum\limits_{ k = \scrN^R_{ n \delta } + 1 }^{ j }
  \biggl[ 
    \prod\limits_{ l = \scrN^R_{ n \delta } + 1 }^k
    \bigl( 
      1 
      - 
      2 c \gamma_l 
      \mathbbm{1}_{
        [ 0, \varepsilon ]
      }(
        \gamma_{ \scrN^R_{ n \delta } + 1 }
      )
    \bigr)^{ - 1 }
  \biggr]
\\ &
\cdot 
    \gamma_k 
    \left< 
      \Theta_{ k - 1 } 
      ,
      D_k
    \right>
    \mathbbm{1}_{
      \{ 
        \gamma_{ \scrN^R_{ n \delta } + 1 }
        \leq 
        \varepsilon
      \} 
      \cap 
      \{ 
        \tau_R > N_{ n \delta }
      \} 
    }
  \biggr|^2
  \Biggr]
  \Biggr)^{ \!\! \nicefrac{ 1 }{ 2 } }
  .
\end{split}
\end{equation}
\Hence 
for all $ R, \delta \in (0,\infty) $, $ n \in \N_0 $, 
$ \varepsilon \in ( 0, ( 2 c )^{ - 1 } ) $
that 
\begin{equation}
\begin{split}
&
  \E\bigl[
    \bigl(
      \max\nolimits_{ 
        N_{ n \delta } 
        \leq 
        j 
        \leq 
        N_{ (n + 1) \delta } 
      }
      \| \Theta_j \|^2
    \bigr)
    \mathbbm{1}_{
      \{ 
        \gamma_{ \scrN^R_{ n \delta } + 1 }
        \leq 
        \varepsilon
      \} 
      \cap 
      \{ 
        \tau_R = \infty 
      \} 
    }
  \bigr]
\\ &
\leq
  \E\bigl[ 
    \min\bigl\{ 
      R, 
      \| \Theta_{ N_{ n \delta } } \|^2 
    \bigr\}
    \mathbbm{1}_{
      \{ 
        \sup_{ n \in \N } \| \Theta_n \| < \infty 
      \} 
    }
  \bigr]
\\ &
  +
  \tfrac{ 3 }{ 2 }
  \,
  \E\bigl[
    \bigl(
      \ft_{ 
        N_{ (n + 1) \delta }
      } 
      - 
      \ft_{
        N_{ n \delta }
      }
    \bigr)
    \min\{ 
      \varepsilon,
      \gamma_{ 
        N_{ n \delta }
        + 1 
      }
    \}
    (
      \scrG_R
      +
      \scrD _R
    ) 
    \mathbbm{1}_{
      \{ 
        \gamma_{ \scrN^R_{ n \delta } + 1 }
        \leq 
        \varepsilon
      \} 
    }
  \bigr]
\\ &
\textstyle 
  +
  2 \kappa
  \Biggl(
  \E\Biggl[
    \sum\limits_{ 
      k = \scrN^R_{ n \delta } + 1 
    }^{ 
      \scrN^R_{ (n + 1) \delta }
    }
    \biggl|
    \biggl[ 
      \prod\limits_{ l = \scrN^R_{ n \delta } + 1 }^k
      \bigl( 
        1 
        - 
        2 c \gamma_l 
        \mathbbm{1}_{
          [ 0, \varepsilon ]
        }(
          \gamma_{ \scrN^R_{ n \delta } + 1 }
        )
      \bigr)^{ - 1 }
    \biggr]
\\ &
\cdot 
      \gamma_k 
      \left< 
        \Theta_{ k - 1 } 
        ,
        D_k
      \right>
      \mathbbm{1}_{
        \{ 
          \gamma_{ \scrN^R_{ n \delta } + 1 }
          \leq 
          \varepsilon
        \} 
        \cap 
        \{ 
          \tau_R > N_{ n \delta }
        \} 
      }
    \biggr|^2
  \Biggr]
  \Biggr)^{ \!\! \nicefrac{ 1 }{ 2 } }
  .
\end{split}
\end{equation}
The Cauchy-Schwarz inequality \hence shows 
for all 
$ R, \delta \in (0,\infty) $, $ n \in \N_0 $, 
$ \varepsilon \in ( 0, ( 2 c )^{ - 1 } ) $
that 
\begin{equation}
\begin{split}
&
  \E\bigl[
    \bigl(
      \max\nolimits_{ 
        N_{ n \delta } 
        \leq 
        j 
        \leq 
        N_{ (n + 1) \delta } 
      }
      \| \Theta_j \|^2
    \bigr)
    \mathbbm{1}_{
      \{ 
        \gamma_{ \scrN^R_{ n \delta } + 1 }
        \leq 
        \varepsilon
      \} 
      \cap 
      \{ 
        \tau_R = \infty 
      \} 
    }
  \bigr]
\\ &
\leq
  \E\bigl[ 
    \min\bigl\{ 
      R, 
      \| \Theta_{ N_{ n \delta } } \|^2 
    \bigr\}
    \mathbbm{1}_{ \cE }
  \bigr]
  +
  \tfrac{ 3 }{ 2 }
  \,
  \E\bigl[
    (
      \delta + \varepsilon 
    )
    \min\{ 
      \varepsilon ,
      \gamma_{ 
        N_{ n \delta }
        + 1 
      }
    \}
    (
      \scrG_R
      +
      \scrD _R
    ) 
    \mathbbm{1}_{
      \{ 
        \gamma_{ \scrN^R_{ n \delta } + 1 }
        \leq 
        \varepsilon
      \} 
    }
  \bigr]
\\ &
\textstyle 
  +
  2 \kappa
  \Biggl(
  \E\Biggl[
    \sum\limits_{ 
      k = \scrN^R_{ n \delta } + 1 
    }^{ 
      \scrN^R_{ (n + 1) \delta }
    }
    \biggl[ 
    \prod\limits_{ l = \scrN^R_{ n \delta } + 1 }^k
    \bigl( 
      1 
      - 
      2 c \gamma_l 
      \mathbbm{1}_{
        [ 0, \varepsilon ]
      }(
        \gamma_{ \scrN^R_{ n \delta } + 1 }
      )
    \bigr)^{ - 2 }
  \biggr]
\\ &
\cdot 
    ( \gamma_k )^2
    \| \Theta_{ k - 1 } \|^2
    \| D_k \|^2
    \mathbbm{1}_{
      \{ 
        \gamma_{ \scrN^R_{ n \delta } + 1 }
        \leq 
        \varepsilon
      \} 
      \cap 
      \{ 
        \tau_R > N_{ n \delta }
      \} 
    }
  \Biggr]
  \Biggr)^{ \!\! \nicefrac{ 1 }{ 2 } }
  .
\end{split}
\end{equation}
\cref{cor:exp_bound} \hence shows for all 
$ R, \delta \in (0,\infty) $, $ n \in \N_0 $, 
$ \varepsilon \in ( 0, \infty ) $
with 
$
  2 c \varepsilon \leq \exp( - 1 ) < 1
$
that 
\begin{equation}
\begin{split}
&
  \E\bigl[
    \bigl(
      \max\nolimits_{ 
        N_{ n \delta } 
        \leq 
        j 
        \leq 
        N_{ (n + 1) \delta } 
      }
      \| \Theta_j \|^2
    \bigr)
    \mathbbm{1}_{
      \{ 
        \gamma_{ \scrN^R_{ n \delta } + 1 }
        \leq 
        \varepsilon
      \} 
      \cap 
      \{ 
        \tau_R = \infty 
      \} 
    }
  \bigr]
\\ &
\leq
  \E\bigl[ 
    \min\bigl\{ 
      R, 
      \| \Theta_{ N_{ n \delta } } \|^2 
      \mathbbm{1}_{ \cE }
    \bigr\}
  \bigr]
  +
  \tfrac{ 3 ( \delta + \varepsilon ) }{ 2 }
  \,
  \E\bigl[
    \min\{ 
      \varepsilon ,
      \gamma_{ 
        N_{ n \delta }
        + 1 
      }
    \}
    (
      \scrG_R
      +
      \scrD _R
    ) 
    \mathbbm{1}_{
      \{ 
        \gamma_{ \scrN^R_{ n \delta } + 1 }
        \leq 
        \varepsilon
      \} 
    }
  \bigr]
\\ &
\textstyle 
  +
  2 \kappa
  \Biggl(
  \E\Biggl[
    \sum\limits_{ 
      k = \scrN^R_{ n \delta } + 1 
    }^{ 
      \scrN^R_{ (n + 1) \delta }
    }
    \biggl[ 
    \prod\limits_{ l = \scrN^R_{ n \delta } + 1 }^k
    \bigl[ 
      \exp\bigl(
        \exp(1)
        2 c \gamma_l 
        \mathbbm{1}_{
          [ 0, \varepsilon ]
        }(
          \gamma_{ \scrN^R_{ n \delta } + 1 }
        )
      \bigr)
    \bigr]^{ 2 }
  \biggr]
    ( \gamma_k )^2
    R
    \scrD _R
    \mathbbm{1}_{
      \{ 
        \gamma_{ \scrN^R_{ n \delta } + 1 }
        \leq 
        \varepsilon
      \} 
      \cap 
      \{ 
        \tau_R > N_{ n \delta }
      \} 
    }
  \Biggr]
  \Biggr)^{ \!\! \nicefrac{ 1 }{ 2 } }
  .
\end{split}
\end{equation}
\Hence for all 
$ R, \delta \in (0,\infty) $, $ n \in \N_0 $, 
$ \varepsilon \in ( 0, ( 2 c \exp(1) )^{ - 1 } ] $
that 
\begin{equation}
\begin{split}
&
  \E\bigl[
    \bigl(
      \max\nolimits_{ 
        N_{ n \delta } 
        \leq 
        j 
        \leq 
        N_{ (n + 1) \delta } 
      }
      \| \Theta_j \|^2
    \bigr)
    \mathbbm{1}_{
      \{ 
        \gamma_{ \scrN^R_{ n \delta } + 1 }
        \leq 
        \varepsilon
      \} 
      \cap 
      \{ 
        \tau_R = \infty 
      \} 
    }
  \bigr]
\\ &
\leq
  \E\bigl[ 
    \min\bigl\{ 
      R, 
      \| \Theta_{ N_{ n \delta } } \|^2 
      \mathbbm{1}_{ \cE }
    \bigr\}
  \bigr]
  +
  \tfrac{ 3 ( \delta + \varepsilon ) }{ 2 }
  \,
  \E\bigl[
    \min\{ 
      \varepsilon ,
      \gamma_{ 
        N_{ n \delta }
        + 1 
      }
    \}
    (
      \scrG_R
      +
      \scrD _R
    ) 
  \bigr]
\\ &
\textstyle 
  +
  2 \kappa
  \Biggl(
  \E\Biggl[
    \sum\limits_{ 
      k = \scrN^R_{ n \delta } + 1 
    }^{ 
      \scrN^R_{ (n + 1) \delta }
    }
    \exp\biggl(
      \exp(1)
      4 c 
      \biggl[
        \sum\limits_{ l = \scrN^R_{ n \delta } + 1 }^k
        \gamma_l 
      \biggr]
    \biggr)
    ( \gamma_k )^2
    R
    \scrD _R
    \mathbbm{1}_{
      \{ 
        \gamma_{ \scrN^R_{ n \delta } + 1 }
        \leq 
        \varepsilon
      \} 
      \cap 
      \{ 
        \tau_R > N_{ n \delta }
      \} 
    }
  \Biggr]
  \Biggr)^{ \!\! \nicefrac{ 1 }{ 2 } }
  .
\end{split}
\end{equation}
\Hence for all 
$ R, \delta \in (0,\infty) $, $ n \in \N_0 $, 
$ \varepsilon \in ( 0, ( 2 c \exp(1) )^{ - 1 } ] $
that 
\begin{equation}
\begin{split}
&
  \E\bigl[
    \bigl(
      \max\nolimits_{ 
        N_{ n \delta } 
        \leq 
        j 
        \leq 
        N_{ (n + 1) \delta } 
      }
      \| \Theta_j \|^2
    \bigr)
    \mathbbm{1}_{
      \{ 
        \gamma_{ \scrN^R_{ n \delta } + 1 }
        \leq 
        \varepsilon
      \} 
      \cap 
      \{ 
        \tau_R = \infty 
      \} 
    }
  \bigr]
\\ &
\leq
  \E\bigl[ 
    \min\bigl\{ 
      R, 
      \| \Theta_{ N_{ n \delta } } \|^2 
      \mathbbm{1}_{ \cE }
    \bigr\}
  \bigr]
  +
  \tfrac{ 3 ( \delta + \varepsilon ) }{ 2 }
  \,
  \E\bigl[
    \min\{ 
      \varepsilon ,
      \gamma_{ 
        N_{ n \delta }
        + 1 
      }
    \}
    (
      \scrG_R
      +
      \scrD _R
    ) 
  \bigr]
\\ &
\textstyle 
  +
  2 \kappa
  \Biggl(
  \E\Biggl[
    \sum\limits_{ 
      k = \scrN^R_{ n \delta } + 1 
    }^{ 
      \scrN^R_{ (n + 1) \delta }
    }
    \exp\bigl(
      \exp(1)
      4 c 
      \bigl[
        \ft_k - \ft_{ \scrN^R_{ n \delta } }
      \bigr]
    \bigr)
    \gamma_k 
    \gamma_{ \scrN^R_{ n \delta } + 1 }
    R
    \scrD _R
    \mathbbm{1}_{
      \{ 
        \gamma_{ \scrN^R_{ n \delta } + 1 }
        \leq 
        \varepsilon
      \} 
      \cap 
      \{ 
        \tau_R > N_{ n \delta }
      \} 
    }
  \Biggr]
  \Biggr)^{ \!\! \nicefrac{ 1 }{ 2 } }
  .
\end{split}
\end{equation}
\Hence for all 
$ R, \delta \in (0,\infty) $, $ n \in \N_0 $, 
$ \varepsilon \in ( 0, ( 2 c \exp(1) )^{ - 1 } ] $
that 
\begin{equation}
\begin{split}
&
  \E\bigl[
    \bigl(
      \max\nolimits_{ 
        N_{ n \delta } 
        \leq 
        j 
        \leq 
        N_{ (n + 1) \delta } 
      }
      \| \Theta_j \|^2
    \bigr)
    \mathbbm{1}_{
      \{ 
        \gamma_{ \scrN^R_{ n \delta } + 1 }
        \leq 
        \varepsilon
      \} 
      \cap 
      \{ 
        \tau_R = \infty 
      \} 
    }
  \bigr]
\\ &
\leq
  \E\bigl[ 
    \min\bigl\{ 
      R, 
      \| \Theta_{ N_{ n \delta } } \|^2 
      \mathbbm{1}_{ \cE }
    \bigr\}
  \bigr]
  +
  \tfrac{ 3 ( \delta + \varepsilon ) }{ 2 }
  \,
  \E\bigl[
    \min\{ 
      \varepsilon ,
      \gamma_{ 
        N_{ n \delta }
        + 1 
      }
    \}
    (
      \scrG_R
      +
      \scrD _R
    ) 
  \bigr]
\\ &
\textstyle 
  +
  2 \kappa
  \Biggl(
  \E\Biggl[
    \sum\limits_{ 
      k = \scrN^R_{ n \delta } + 1 
    }^{ 
      \scrN^R_{ (n + 1) \delta }
    }
    \exp\bigl(
      \exp(1)
      4 c 
      \bigl[
        \ft_{ \scrN^R_{ (n + 1) \delta } } - \ft_{ \scrN^R_{ n \delta } }
      \bigr]
    \bigr)
    \gamma_k 
    \gamma_{ N_{ n \delta } + 1 }
    R
    \scrD _R
    \mathbbm{1}_{
      \{ 
        \gamma_{ N_{ n \delta } + 1 }
        \leq 
        \varepsilon
      \} 
      \cap 
      \{ 
        \tau_R > N_{ n \delta }
      \} 
    }
  \Biggr]
  \Biggr)^{ \!\! \nicefrac{ 1 }{ 2 } }
  .
\end{split}
\end{equation}
\Hence for all 
$ R, \delta \in (0,\infty) $, $ n \in \N_0 $, 
$ \varepsilon \in ( 0, ( 2 c \exp(1) )^{ - 1 } ] $
that 
\begin{equation}
\begin{split}
&
  \E\bigl[
    \bigl(
      \max\nolimits_{ 
        N_{ n \delta } 
        \leq 
        j 
        \leq 
        N_{ (n + 1) \delta } 
      }
      \| \Theta_j \|^2
    \bigr)
    \mathbbm{1}_{
      \{ 
        \gamma_{ \scrN^R_{ n \delta } + 1 }
        \leq 
        \varepsilon
      \} 
      \cap 
      \{ 
        \tau_R = \infty 
      \} 
    }
  \bigr]
\\ &
\leq
  \E\bigl[ 
    \min\bigl\{ 
      R, 
      \| \Theta_{ N_{ n \delta } } \|^2 
      \mathbbm{1}_{ \cE }
    \bigr\}
  \bigr]
  +
  \tfrac{ 3 ( \delta + \varepsilon ) }{ 2 }
  \,
  \E\bigl[
    \min\{ 
      \varepsilon ,
      \gamma_{ 
        N_{ n \delta }
        + 1 
      }
    \}
    (
      \scrG_R
      +
      \scrD _R
    ) 
  \bigr]
\\ &
\textstyle 
  +
  2 \kappa
  \Bigl(
  \E\Bigl[
    \exp\bigl(
      \exp(1)
      4 c 
      \bigl[
        \delta + \varepsilon
      \bigr]
    \bigr)
    \bigl[
      \ft_{ \scrN^R_{ (n + 1) \delta } } - \ft_{ \scrN^R_{ n \delta } }
    \bigr]
    \min\{ \varepsilon, \gamma_{ N_{ n \delta } + 1 } \}
    R
    \scrD _R
    \mathbbm{1}_{
      \{ 
        \gamma_{ N_{ n \delta } + 1 }
        \leq 
        \varepsilon
      \} 
      \cap 
      \{ 
        \tau_R > N_{ n \delta }
      \} 
    }
  \Bigr]
  \Bigr)^{ \! \nicefrac{ 1 }{ 2 } }
  .
\end{split}
\end{equation}
\Hence for all 
$ R, \delta \in (0,\infty) $, $ n \in \N_0 $, 
$ \varepsilon \in ( 0, ( 2 c \exp(1) )^{ - 1 } ] $
that 
\begin{equation}
\begin{split}
&
  \E\bigl[
    \bigl(
      \max\nolimits_{ 
        N_{ n \delta } 
        \leq 
        j 
        \leq 
        N_{ (n + 1) \delta } 
      }
      \| \Theta_j \|^2
    \bigr)
    \mathbbm{1}_{
      \{ 
        \gamma_{ \scrN^R_{ n \delta } + 1 }
        \leq 
        \varepsilon
      \} 
      \cap 
      \{ 
        \tau_R = \infty 
      \} 
    }
  \bigr]
\\ &
\leq
  \E\bigl[ 
    \min\bigl\{ 
      R, 
      \| \Theta_{ N_{ n \delta } } \|^2 
      \mathbbm{1}_{ \cE }
    \bigr\}
  \bigr]
  +
  \tfrac{ 3 ( \delta + \varepsilon ) }{ 2 }
  \,
  \E\bigl[
    \min\{ 
      \varepsilon ,
      \gamma_{ 
        N_{ n \delta }
        + 1 
      }
    \}
    (
      \scrG_R
      +
      \scrD _R
    ) 
  \bigr]
\\ &
\textstyle 
  +
  2 R^{ \nicefrac{ 1 }{ 2 } }
  \kappa
    \exp\bigl(
      2 c \exp(1)
      (
        \delta + \varepsilon
      )
    \bigr)
  \bigl(
  \E\bigl[
    (
      \delta + \varepsilon
    )
    \scrD _R
    \min\{ \varepsilon, \gamma_{ N_{ n \delta } + 1 } \}
    \mathbbm{1}_{
      \{ 
        \gamma_{ N_{ n \delta } + 1 }
        \leq 
        \varepsilon
      \} 
      \cap 
      \{ 
        \tau_R > N_{ n \delta }
      \} 
    }
  \bigr]
  \bigr)^{ \! \nicefrac{ 1 }{ 2 } }
  .
\end{split}
\end{equation}
\Hence for all 
$ R, \delta \in (0,\infty) $, $ n \in \N_0 $, 
$ \varepsilon \in ( 0, ( 2 c \exp(1) )^{ - 1 } ] $
that 
\begin{equation}
\label{eq:max_estimate_for_item_ii}
\begin{split}
&
  \E\bigl[
    \bigl(
      \max\nolimits_{ 
        N_{ n \delta } 
        \leq 
        j 
        \leq 
        N_{ (n + 1) \delta } 
      }
      \| \Theta_j \|^2
    \bigr)
    \mathbbm{1}_{
      \{ 
        \gamma_{ \scrN^R_{ n \delta } + 1 }
        \leq 
        \varepsilon
      \} 
      \cap 
      \{ 
        \tau_R = \infty 
      \} 
    }
  \bigr]
\\ &
\leq
  \E\bigl[ 
    \min\bigl\{ 
      R, 
      \| \Theta_{ N_{ n \delta } } \|^2 
      \mathbbm{1}_{ \cE }
    \bigr\}
  \bigr]
  +
  \tfrac{ 3 ( \delta + \varepsilon ) }{ 2 }
  \,
  \E\bigl[
    (
      \scrG_R
      +
      \scrD _R
    ) 
    \min\{ 
      \varepsilon ,
      \gamma_{ 
        N_{ n \delta }
        + 1 
      }
    \}
  \bigr]
\\ &
\textstyle 
  +
  2 R^{ \nicefrac{ 1 }{ 2 } }
  \kappa
    \exp\bigl(
      1 +
      2 c \delta \exp(1)
    \bigr)
    (
      \delta + \varepsilon
    )^{ \nicefrac{ 1 }{ 2 } }
  \bigl(
  \E\bigl[
    \scrD _R
    \min\{ \varepsilon, \gamma_{ N_{ n \delta } + 1 } \}
  \bigr]
  \bigr)^{ \nicefrac{ 1 }{ 2 } }
  .
\end{split}
\end{equation}
\Moreover Vitali's convergence theorem 
(cf., for example, \cite[Theorem~6.25]{Klenke2014})
and \cref{item:i} show that 
\begin{equation}
\label{eq:consequence_vitali}
  \lim \nolimits_{ n \to \infty }
  \E\bigl[ 
    \min\bigl\{ 
      R, 
      \| \Theta_{ N_{ n \delta } } \|^2 
      \mathbbm{1}_{ \cE }
    \bigr\}
  \bigr]
  = 0
  .
\end{equation}
\Moreover Lebesgue's theorem of dominated convergence 
and the fact that $\forall \, R \in (0 , \infty ) \colon \E[ \scrD _R ] < \infty $
imply 
for all 
$ R, \delta, \varepsilon \in (0,\infty) $
that 
\begin{equation}
  \lim_{ n \to \infty }
  \E\bigl[
    (
      \scrG_R
      +
      \scrD _R
    ) 
    \min\{ 
      \varepsilon ,
      \gamma_{ 
        N_{ n \delta }
        + 1 
      }
    \}
  \bigr]
  =
  \lim_{ n \to \infty }
  \E\bigl[
    \scrD _R
    \min\{ 
      \varepsilon ,
      \gamma_{ 
        N_{ n \delta }
        + 1 
      }
    \}
  \bigr]
  =
  0 .
\end{equation}
This, \cref{eq:consequence_vitali}, 
and \cref{eq:max_estimate_for_item_ii} 
ensure for all 
$ R, \delta \in (0,\infty) $, $ n \in \N_0 $, 
$ \varepsilon \in ( 0, ( 2 c \exp(1) )^{ - 1 } ] $
that 
\begin{equation}
\label{eq:item_ii_almost_done}
  \lim_{ n \to \infty }
  \E\bigl[
    \bigl(
      \max\nolimits_{ 
        N_{ n \delta } 
        \leq 
        j 
        \leq 
        N_{ (n + 1) \delta } 
      }
      \| \Theta_j \|^2
    \bigr)
    \mathbbm{1}_{
      \{ 
        \gamma_{ \scrN^R_{ n \delta } + 1 }
        \leq 
        \varepsilon
      \} 
      \cap 
      \{ 
        \tau_R = \infty 
      \} 
    }
  \bigr] 
  = 0 .
\end{equation}
Combining this with \cref{lem:convergence_in_probability_charac} 
establishes \cref{item:ii}. 
It thus remains to prove \cref{item:iii}. 
For this we introduce more notation. 
Specifically, for every $ \delta \in (0,\infty) $
let $ \lfloor \cdot \rfloor_{ \delta } \colon \R \to \R $ satisfy 
for all $ t \in \R $ that 
\begin{equation}
\label{eq:round_down}
  \lfloor t \rfloor_{ \delta }
  =
  \min\bigl( 
    ( - \infty, t ]
    \cap 
    \{ 0, \delta, - \delta, 2 \delta, - 2 \delta, \dots \}
  \bigr)
  .
\end{equation}
Note that \cref{eq:round_down} 
and the fact that for all $ s, t \in [0,\infty) $
with $ s \leq t $ 
it holds that 
\begin{equation}
  N_s \leq N_t 
\end{equation}
ensure that 
for all $ \delta \in (0,\infty) $, $ t \in [0,\infty) $ it holds that
\begin{equation}
\textstyle 
  \| \Theta_{ N_t } \|
  \leq 
  \max_{ 
    N_{ \lfloor t \rfloor_{ \delta } } \leq m \leq N_{ \lfloor t \rfloor_{ \delta } + \delta } 
  }
  \| \Theta_m \|
  .
\end{equation}
\Hence for all $ \delta \in (0,\infty) $, $ t \in [0,\infty) $ that 
\begin{equation}
\textstyle 
  \E\bigl[ 
    \min\{ 
      1 ,
      \| \Theta_{ N_t } \|
    \}
  \bigr]
\leq 
  \E\bigl[ 
    \max_{ 
      N_{ \lfloor t \rfloor_{ \delta } } \leq m \leq N_{ \lfloor t \rfloor_{ \delta } + \delta } 
    }
    \min\{ 
      1 ,
      \| \Theta_m \|
    \}
  \bigr]
\end{equation}
Let $ \tau \colon \N \to \R $
satisfy for all $ n \in \N $ that 
$ 0 \leq \tau(n) \leq \tau( n + 1 ) $. 
\Hence for all 
$ \delta \in (0,\infty) $, $ n \in \N $ 
and all 
non-decreasing 
$ a \colon \N \to [0,\infty) $
that
\begin{equation}
\label{eq:small_a_estimate}
\textstyle 
  \| \Theta_{ N_{ a(n) } } \|
  \leq 
  \max_{ 
    N_{ \lfloor a(n) \rfloor_{ \delta } } \leq m \leq N_{ \lfloor a(n) \rfloor_{ \delta } + \delta } 
  }
  \| \Theta_m \|
  .
\end{equation}
\Moreover \cref{item:ii} ensures 
for all $ \delta \in (0,\infty) $ that 
\begin{equation}
  \lim_{ n \to \infty }
  \E\biggl[
    \max_{ 
      \scrN^R_{ n \delta } \leq m \leq N_{ n \delta + \delta }
    }
    \min\bigl\{ 
      1 ,
      \| 
        \Theta_m
      \|
    \bigr\} 
  \biggr]
  = 0
  .
\end{equation}
This demonstrates that for all $ \delta \in (0,\infty) $ 
and all 
non-decreasing $ A \colon \N \to \{ 0, \delta, 2 \delta, \dots \} $ 
with 
$ 
  \lim_{ n \to \infty } A(n) = \infty 
$
it holds that
\begin{equation}
  \lim_{ n \to \infty }
  \E\biggl[
    \max_{ 
      N_{ A(n) } \leq m \leq N_{ A(n) + \delta }
    }
    \min\bigl\{ 
      1 ,
      \| 
        \Theta_m
      \|
    \bigr\} 
  \biggr]
  = 0
  .
\end{equation}
Combining this with \cref{eq:small_a_estimate} shows that 
for all $ \delta \in (0,\infty) $
and all 
non-decreasing $ a \colon \N \to [0,\infty) $ 
with 
$ 
  \lim_{ n \to \infty } a(n) = \infty 
$
it holds that
\begin{equation}
  \limsup_{ n \to \infty }
  \E\bigl[
    \min\bigl\{ 
      1 ,
      \| 
        \Theta_{ N_{ a(n) } }
      \|
    \bigr\} 
  \bigr]
  \leq 
  \limsup_{ n \to \infty }
  \E\biggl[
    \max_{ 
      N_{ \lfloor a(n) \rfloor_{ \delta } } \leq m \leq N_{ \lfloor a(n) \rfloor_{ \delta } + \delta } 
    }
    \min\bigl\{ 
      1 ,
      \| \Theta_m \|
    \bigr\}
  \biggr]
  = 0
  .
\end{equation}
\Hence for all non-decreasing $ a \colon \N \to [0,\infty) $ 
with 
$ 
  \lim_{ n \to \infty } a(n) = \infty 
$
that
\begin{equation}
  \limsup_{ n \to \infty }
  \E\bigl[
    \min\bigl\{ 
      1 ,
      \| 
        \Theta_{ N_{ a(n) } }
      \|
    \bigr\} 
  \bigr]
  = 0
  .
\end{equation}
This establishes \cref{item:iii}. 
\end{cproof}

\subsection{Convergence in probability for general attractors}

Our next goal is to establish in \cref{cor:convergence_in_probability} a more general version of \cref{prop:convergence_in_probability} where the limit point of the \SGD\ sequence can be an arbitrary $\vartheta \in \R^\fd$ and where the assumptions 
on the random learning rate sequence $ (\gamma_n)_{ n \in \N } $ 
only hold with probability one and not necessarily everywhere on $ \Omega $.
For this we first state in
\cref{lem:augmentation} and \cref{lem:conditional_completions} below 
a few elementary properties on suitable augmentations 
of measurable and probability spaces (with null sets). 
\cref{lem:augmentation} is well-known; see, e.g., \cite[Theorem 3.3.5]{Dudley2002}.

\begin{lemma}[Augmentation of probability spaces]
\label{lem:augmentation}
Let $ ( \Omega, \cH, \P ) $ be a probability space, 
let $ \cF \subseteq \cH $ be a sigma-algebra on $ \Omega $, 
let $ \cN \subseteq \cH $ be given by 
$
  \cN = \{ A \in \cH \colon \P( A ) = 0 \}
$,
and let $ \cG \subseteq \cH $ be given by 
\begin{equation}
\label{eq:definition_of_G_completeness_augmentation}
  \cG = 
  \bigl\{
    \cA \in \cH 
    \colon 
    \bigl[
      \exists \, 
      A \in \cF, 
      B, C \in \cN
      \colon
      \cA = ( A \backslash B ) \cup C
    \bigr]
  \bigr\}
  .
\end{equation}
Then $ \cG = \sigma( \cF \cup \cN ) $.
\end{lemma}


\begin{lemma}[Conditional expectations with respect to augmentations]
\label{lem:conditional_completions}
Let $ ( \Omega, \cH, \P ) $ be a probability space, 
let $ X \colon \Omega \to \R $ be $ \cH $-measurable, 
assume $ \E[ | X | ] < \infty $, 
let $ \cF \subseteq \cH $ be a sigma-algebra on $ \Omega $, 
and let $ \cG \subseteq \cH $ be given by
\begin{equation}
\label{eq:augmentation_sigma_algebra}
  \cG =
  \sigma\bigl(
    \cF \cup 
    \{ 
      A \in \cH \colon \P( A ) = 0
    \}
  \bigr)
  .
\end{equation}
Then it holds $ \P $-a.s.\ that 
\begin{equation}
\label{eq:conditional_expectation_identity}
  \E\bigl[ X | \cF \bigr]
  =
  \E\bigl[ X | \cG \bigr]
  .
\end{equation}
\end{lemma}

\begin{cproof}{lem:conditional_completions}
Throughout this proof let $ Y \colon \Omega \to \R $ be a $ \cF $-measurable function which satisfies 
for all $ A \in \cF $ that 
\begin{equation}
\label{eq:definition_conditional_expectation_w_r_t_cF}
  \E[ X \mathbbm{1}_A ] 
  =
  \E[ Y \mathbbm{1}_A ] 
  .
\end{equation}
\Nobs that 
for every 
$ \cF $-measurable $ Z \colon \Omega \to \R $ 
and every  
$ A \in \cF $, $ B, C \in \cH $ with $ \P( B ) = \P( C ) = 0 $ it holds that 
\begin{equation}
\begin{split}
  \E\bigl[ 
    Z \mathbbm{1}_{ ( A \backslash B ) \cup C }
  \bigr]
& 
  =
  \E\bigl[ 
    Z \mathbbm{1}_{ ( A \backslash ( B \cup C ) ) \cup C }
  \bigr]
  =
  \E\bigl[ 
    Z 
    \bigl( 
      \mathbbm{1}_{ A \backslash ( B \cup C ) }
      +
      \mathbbm{1}_C
    \bigr)
  \bigr]
\\ &  
  =
  \E\bigl[ 
    Z 
    \mathbbm{1}_{ A \cap ( \Omega \backslash ( B \cup C ) ) }
    +
    Z
    \mathbbm{1}_C
  \bigr]
  =
  \E\bigl[ 
    Z 
    \mathbbm{1}_A 
    \mathbbm{1}_{ \Omega \backslash ( B \cup C ) }
  \bigr]
  +
  \E\bigl[ 
    Z
    \mathbbm{1}_C
  \bigr]
\\ & =
  \E\bigl[ 
    Z 
    \mathbbm{1}_A 
    \mathbbm{1}_{ ( \Omega \backslash B ) \cap ( \Omega \backslash C ) }
  \bigr]
  =
  \E\bigl[ 
    Z 
    \mathbbm{1}_A 
    \mathbbm{1}_{ \Omega \backslash B }
    \mathbbm{1}_{ \Omega \backslash C }
  \bigr]
  =
  \E\bigl[ 
    Z 
    \mathbbm{1}_A 
  \bigr]
  .
\end{split}
\end{equation}
Combining this with \cref{eq:definition_conditional_expectation_w_r_t_cF} ensures 
for all $ A \in \cF $, $ B, C \in \{ \cB \in \cH \colon \P( \cB ) = 0 \} $ that 
\begin{equation}
\label{eq:conditional_expectation_augmentation}
  \E\bigl[ 
    X \mathbbm{1}_{ ( A \backslash B ) \cup C }
  \bigr]
  =
  \E\bigl[ 
    X \mathbbm{1}_{ A }
  \bigr]
  =
  \E\bigl[ 
    Y \mathbbm{1}_{ A }
  \bigr]
  =
  \E\bigl[ 
    Y \mathbbm{1}_{ ( A \backslash B ) \cup C }
  \bigr]
  .
\end{equation}
\Moreover \cref{lem:augmentation} and \cref{eq:augmentation_sigma_algebra} show that
\begin{equation}
  \cG =
  \bigl\{
    \cA \in \cH \colon
    \bigl(
      \exists \, 
      A \in \cF, 
      B, C \in 
      \{ 
        \cB \in \cH 
        \colon
        \P( \cB ) = 0
      \}
      \colon
      \cA = ( A \backslash B ) \cup C
    \bigr)
  \bigr\}
  .
\end{equation}
This and \cref{eq:conditional_expectation_augmentation} 
show for all $ \cA \in \cG $ that 
\begin{equation}
  \E[ 
    X \mathbbm{1}_{ \cA }
  ]
  =
  \E[ 
    Y \mathbbm{1}_{ \cA }
  ]
  .
\end{equation}
Combining this with \cref{eq:definition_conditional_expectation_w_r_t_cF}
establishes \cref{eq:conditional_expectation_identity}. 
\end{cproof}

\begin{cor}[Convergence in probability for general attractors]
\label{cor:convergence_in_probability}
Let $ ( \Omega, \mathcal{F}, ( \mathbb{F}_n )_{ n \in \N_0 }, \P ) $ 
be a filtered probability space, 
for every $ n \in \N $ let $ \gamma_n \colon \Omega \to \R $ 
be $ \mathbb{F}_{ n - 1 } $-measurable, 
assume 
$
  \P\bigl( 
    \lim_{ n \to \infty }
    \gamma_n 
    =
    0
  \bigr)
  = 
  \P\bigl( 
    \sum_{ n = 1 }^{ \infty } | \gamma_n |
    =
    \infty
  \bigr)
  = 1
$, 
assume for all $ n \in \N $ that 
$
  \P(
    \gamma_{ n + 1 } \leq \gamma_n
  )
  = 1
$,
let $ \fd \in \N $, $ \vartheta \in \R^{ \fd } $, 
let $ g \colon \R^{ \fd } \to \R^\fd $ be measurable and locally bounded, 
let $ c \in (0,\infty) $ satisfy for all $ \theta \in \R^{ \fd } $ that 
\begin{equation}
\label{eq:coercivity_general}
  \left< \theta - \vartheta , g( \theta ) \right> 
  \leq 
  - c \| \theta - \vartheta \|^2 
  ,
\end{equation}
let 
$ \Theta \colon \N_0 \times \Omega \to \R $
be an $ ( \mathbb{F}_n )_{ n \in \N_0 } $-adapted stochastic process, 
let 
$ D \colon \N \times \Omega \to \R^{ \fd } $ 
be an $ ( \mathbb{F}_n )_{ n \in \N } $-adapted stochastic process, 
assume 
$
  \E[ \sup_{ n \in \N } \| D_n \|^2 ] < \infty 
$, 
assume that for all $ n \in \N $ it holds $ \P $-a.s.\ that
\begin{equation}
\label{eq:definition_of_Theta_n_cor}
  \Theta_n 
  = \Theta_{ n - 1 } 
  + 
  \gamma_n 
  g( \Theta_{ n - 1 } )
  +
  \gamma_n 
  D_n 
\qqandqq
  \E[ 
    D_n | \mathbb{F}_{ n - 1 }
  ] 
  = 0
  ,
\end{equation}
for every $ t \in [0,\infty) $ 
let 
$
  N_t \colon \Omega \to \R
$
satisfy 
\begin{equation}
\label{eq:definition_of_Nt_in_cor}
\textstyle 
  N_t
  =
  \inf\bigl(
    \bigl\{
      n \in \N_0
      \colon 
      \sum_{ k = 1 }^n
      \gamma_k
      \geq 
      t
    \bigr\}
  \bigr)
  ,
\end{equation}
and let $ \cE \subseteq \Omega $ satisfy 
$
  \cE = \{ \sup_{ n \in \N } \| \Theta_n \| < \infty \} 
$. 
Then 
\begin{equation}
\label{eq:convergence_in_probability_cor}
  \lim_{ t \to \infty }
  \E\bigl[ 
    \min\bigl\{ 1,  
      \| \Theta_{ N_t } - \vartheta \|
      \mathbbm{1}_{ \cE }
    \bigr\}
  \bigr]
  = 0
  .
\end{equation}
\end{cor}

\begin{cproof}{cor:convergence_in_probability}
Throughout this proof 
let $ \varOmega \in \cF $ satisfy 
\begin{equation}
\label{eq:definition_of_varOmega}
\begin{split}
\textstyle 
  \varOmega 
&
\textstyle 
  = 
  \{ 
    \lim_{ n \to \infty }
    \gamma_n 
    =
    0
  \}
  \cap 
  \bigl(
    \cap_{ n = 1 }^{ \infty }
    \{ 
      \gamma_{ n + 1 } \leq \gamma_n
    \} 
  \bigr)
  \cap 
  \{ 
    \sum_{ n = 1 }^{ \infty } \gamma_n 
    =
    \infty
  \} 
\\ &
\textstyle 
\quad 
  \cap 
  \bigl(
    \cap_{ n = 1 }^{ \infty }
    \{ 
      \Theta_n 
      = \Theta_{ n - 1 } 
      + 
      \gamma_n 
      g( \Theta_{ n - 1 } )
      +
      \gamma_n 
      D_n 
    \}
  \bigr)
  ,
\end{split}
\end{equation}
for every $ n \in \N_0 $ let 
$ \mathbb{G}_n \subseteq \cF $
satisfy 
\begin{equation}
\label{eq:definition_of_filtration_Gn}
  \mathbb{G}_n = 
  \sigma(
    \mathbb{F}_n \cup 
    \{
      A \in \cF
      \colon 
      \P( A ) = 0
    \}
  )
  ,
\end{equation}
let $ G \colon \R^{ \fd } \to \R^{ \fd } $ satisfy 
for all $ \theta \in \R^{ \fd } $ that
\begin{equation}
\label{eq:definition_of_G}
  G( \theta ) = g( \theta + \vartheta )
  ,
\end{equation}
let $ \varTheta \colon \N_0 \times \Omega \to \R^{ \fd } $,  
$
  \cD \colon \N \times \Omega \to \R^{ \fd }
$, 
and 
$
  \lambda \colon \N \times \Omega \to [0,\infty)
$
satisfy for all $ n \in \N $ that
\begin{equation}
\label{eq:definition_of_varTheta}
  \varTheta_{ n - 1 } 
  = 
  ( \Theta_{ n - 1 } - \vartheta ) \mathbbm{1}_{ \varOmega }
  ,
\qquad 
  \cD_n
  = 
  D_n
  \mathbbm{1}_{ \varOmega }
  - 
  G( 0 )
  \mathbbm{1}_{ \Omega \backslash \varOmega }
  ,
\qqandqq 
  \lambda_n = 
  \gamma_n
  \mathbbm{1}_{ \varOmega }
  +
  n^{ - 1 }
  \mathbbm{1}_{ \Omega \backslash \varOmega }
  ,
\end{equation}
and for every $ t \in [0,\infty) $ let $ M_t \colon \Omega \to \R $ satisfy 
\begin{equation}
\label{eq:definition_of_Mt_in_cor}
\textstyle 
  M_t
  =
  \inf\bigl(
    \bigl\{
      n \in \N_0
      \colon 
      \sum_{ k = 1 }^n
      \lambda_k
      \geq 
      t
    \bigr\}
  \bigr)
  .
\end{equation}
\Nobs that \cref{eq:coercivity_general,eq:definition_of_G} show 
for all $ \theta \in \R^{ \fd } $ that 
\begin{equation}
\label{eq:G_coercivity}
\begin{split}
  \left< \theta, G( \theta ) \right>
& =
  \left< \theta, g( \theta + \vartheta ) \right>
=
  \left< ( \theta + \vartheta ) - \vartheta , g( \theta + \vartheta ) \right>
\leq 
  - c
  \| ( \theta + \vartheta ) - \vartheta \|^2 
=
  - c
  \| \theta \|^2 
  .
\end{split}
\end{equation}
\Moreover \cref{eq:definition_of_varOmega,eq:definition_of_varTheta}
show that
\begin{equation}
\label{eq:varTheta_n_recursion}
\begin{split}
  \varTheta_n 
& =
  ( \Theta_n - \vartheta ) \mathbbm{1}_{ \varOmega }
=  
  \bigl(
    \Theta_{ n - 1 } - \vartheta 
    +
    \gamma_n g( \Theta_{ n - 1 } )
    + 
    \gamma_n D_n
  \bigr)
  \mathbbm{1}_{ \varOmega }
\\ &
=
  (
    \Theta_{ n - 1 } - \vartheta 
  )
  \mathbbm{1}_{ \varOmega }
  +
  \lambda_n g( \Theta_{ n - 1 } )
  \mathbbm{1}_{ \varOmega }
  + 
  \lambda_n D_n
  \mathbbm{1}_{ \varOmega }
\\ & =
  \varTheta_{ n - 1 } 
  +
  \lambda_n g( \Theta_{ n - 1 } - \vartheta + \vartheta )
  \mathbbm{1}_{ \varOmega }
  + 
  \lambda_n 
  \bigl( 
    \cD_n
    +
    G(0)
    \mathbbm{1}_{ \Omega \backslash \varOmega }
  \bigr)
\\ & =
  \varTheta_{ n - 1 } 
  +
  \lambda_n 
  g( 
    ( \Theta_{ n - 1 } - \vartheta ) \mathbbm{1}_{ \varOmega } 
    + \vartheta 
  )
  \mathbbm{1}_{ \varOmega }
  + 
  \lambda_n 
  \cD_n
  +
  \lambda_n 
  G(0)
  \mathbbm{1}_{ \Omega \backslash \varOmega }
\\ & =
  \varTheta_{ n - 1 } 
  +
  \lambda_n 
  g( 
    \varTheta_{ n - 1 } 
    + \vartheta 
  )
  \mathbbm{1}_{ \varOmega }
  + 
  \lambda_n 
  \cD_n
  +
  \lambda_n 
  G(0)
  \mathbbm{1}_{ \Omega \backslash \varOmega }
\\ & =
  \varTheta_{ n - 1 } 
  +
  \lambda_n 
  G( 
    \varTheta_{ n - 1 } 
  )
  \mathbbm{1}_{ \varOmega }
  +
  \lambda_n 
  G(0)
  \mathbbm{1}_{ \Omega \backslash \varOmega }
  + 
  \lambda_n 
  \cD_n
\\ &
  =
  \varTheta_{ n - 1 } 
  +
  \lambda_n 
  \bigl[
    G( 
      \varTheta_{ n - 1 } 
    )
    \mathbbm{1}_{ \varOmega }
    +
    G(0)
    \mathbbm{1}_{ \Omega \backslash \varOmega }
  \bigr]
  + 
  \lambda_n 
  \cD_n
\\ &
  =
  \varTheta_{ n - 1 } 
  +
  \lambda_n 
  \bigl[
    G( 
      \varTheta_{ n - 1 } 
    )
    \mathbbm{1}_{ \varOmega }
    +
    G( \varTheta_{ n - 1 } )
    \mathbbm{1}_{ \Omega \backslash \varOmega }
  \bigr]
  + 
  \lambda_n 
  \cD_n
=
  \varTheta_{ n - 1 } 
  +
  \lambda_n 
  G( \varTheta_{ n - 1 } )
  + 
  \lambda_n 
  \cD_n
  .
\end{split}
\end{equation}
\Moreover \cref{eq:definition_of_Theta_n_cor} proves that 
$ \P( \varOmega ) = 1 $. 
This and \cref{eq:definition_of_filtration_Gn} show that 
for all $ n \in \N_0 $ it holds that 
\begin{equation}
\label{eq:varOmega_has_full_measure_and_is_measurable}
  \varOmega \in \mathbb{G}_n 
\qqandqq
  \P( \varOmega ) = 1
  .
\end{equation}
Combining this with \cref{eq:definition_of_varTheta} 
and the fact that $ \Theta $ is a 
$ ( \mathbb{G}_n )_{ n \in \N_0 } $-adapted stochastic process 
shows that 
$ \varTheta $ is a 
$ ( \mathbb{G}_n )_{ n \in \N_0 } $-adapted stochastic process. 
\Moreover 
\cref{eq:definition_of_varTheta}, 
\cref{eq:varOmega_has_full_measure_and_is_measurable}, 
and the fact that $ D $ is a 
$ ( \mathbb{G}_n )_{ n \in \N } $-adapted stochastic process 
ensure that $ \cD $ is an $ ( \mathbb{G}_n )_{ n \in \N } $-adapted stochastic process. 
\cref{lem:conditional_completions}, 
\cref{eq:definition_of_Theta_n_cor}, 
\cref{eq:definition_of_filtration_Gn}, 
\cref{eq:definition_of_varTheta}, 
\cref{eq:varOmega_has_full_measure_and_is_measurable}, 
and the fact that for all $ n \in \N $ that 
it holds that 
$ \E[ \| D_n \| ] < \infty $ 
\hence 
show that for all $ n \in \N $ it holds $ \P $-a.s.\ that 
\begin{equation}
\label{eq:Dn_conditional_expectation_is_zero}
\begin{split}
  \E[
    \cD_n | \mathbb{G}_{ n - 1 }
  ]
& =
  \E\bigl[
    (
      D_n \mathbbm{1}_{ \varOmega }
      -
      G(0) \mathbbm{1}_{ \Omega \backslash \varOmega }
    )
    | \mathbb{G}_{ n - 1 }
  \bigr]
  =
  \E\bigl[
    D_n \mathbbm{1}_{ \varOmega }
    | \mathbb{G}_{ n - 1 }
  \bigr]
\\ & 
  =
  \E[
    D_n | \mathbb{G}_{ n - 1 }
  ]
  =
  \E[
    D_n | \cG_{ n - 1 }
  ]
  =
  0
  .
\end{split}
\end{equation}
Combining 
\cref{prop:convergence_in_probability} 
(applied with 
$ 
  ( \Omega, \cF, ( \mathbb{F}_n )_{ n \in \N_0 }, \P ) 
  \with 
  ( \Omega, \cF, ( \mathbb{G}_n )_{ n \in \N_0 }, \P ) 
$, 
$
  ( \gamma_n )_{ n \in \N }
  \with 
  ( \lambda_n )_{ n \in \N }
$, 
$
  \fd \with \fd 
$,
$
  g \with G
$, 
$
  c \with c
$, 
$
  D \with \cD
$, 
$
  \Theta \with \varTheta
$, 
$
  ( N_t )_{ t \in [0,\infty) }
  \with 
  ( M_t )_{ t \in [0,\infty) }
$
in the notation of \cref{prop:convergence_in_probability}), 
\cref{eq:definition_of_varOmega}, 
\cref{eq:definition_of_varTheta}, 
\cref{eq:definition_of_Mt_in_cor}, 
\cref{eq:G_coercivity}, 
\cref{eq:varTheta_n_recursion}, 
\cref{eq:Dn_conditional_expectation_is_zero}, 
and 
the fact that 
$ \varTheta $ is a 
$ ( \mathbb{G}_n )_{ n \in \N_0 } $-adapted stochastic process
\hence proves that 
\begin{equation}
  \lim_{ t \to \infty }
  \E\bigl[ 
    \min\bigl\{ 1,  
      \| \varTheta_{ M_t } \| 
      \mathbbm{1}_{ 
        \{ 
          \sup_{ n \in \N } \| \varTheta_n \| < \infty 
        \} 
      }
    \bigr\}
  \bigr]
  = 0
  .
\end{equation}
This, 
\cref{eq:definition_of_varOmega}, 
\cref{eq:definition_of_varTheta}, 
\cref{eq:definition_of_Mt_in_cor}, 
and 
\cref{eq:varOmega_has_full_measure_and_is_measurable} 
show that
\begin{equation}
\begin{split}
&
  \limsup_{ t \to \infty }
  \E\bigl[ 
    \min\bigl\{ 1,  
      \| \Theta_{ N_t } - \vartheta \| 
      \mathbbm{1}_{ 
        \cE
      }
    \bigr\}
  \bigr]
  =
  \limsup_{ t \to \infty }
  \E\bigl[ 
    \min\bigl\{ 1,  
      \| \Theta_{ N_t } - \vartheta \| 
      \mathbbm{1}_{ 
        \cE
      }
    \bigr\}
    \mathbbm{1}_{ \varOmega }
  \bigr]
\\ 
&
  =
  \limsup_{ t \to \infty }
  \E\bigl[ 
    \min\bigl\{ 1,  
      \| \Theta_{ N_t } - \vartheta \| 
      \mathbbm{1}_{ 
        \{ 
          \sup_{ n \in \N } \| \Theta_n \| < \infty 
        \} 
      }
    \bigr\}
    \mathbbm{1}_{ \varOmega }
  \bigr]
\\ 
&
  =
  \limsup_{ t \to \infty }
  \E\bigl[ 
    \min\bigl\{ 1,  
      \| \Theta_{ N_t } - \vartheta \| 
      \mathbbm{1}_{ 
        \{ 
          \sup_{ n \in \N } \| \Theta_n - \vartheta \| < \infty 
        \} 
      }
    \bigr\}
    \mathbbm{1}_{ \varOmega }
  \bigr]
\\ &
  =
  \limsup_{ t \to \infty }
  \E\bigl[ 
    \min\bigl\{ 1,  
      \| \varTheta_{ M_t } \| 
      \mathbbm{1}_{ 
        \{ 
          \sup_{ n \in \N } \| \varTheta_n \| < \infty 
        \} 
      }
    \bigr\}
    \mathbbm{1}_{ \varOmega }
  \bigr]
\\ &
  =
  \limsup_{ t \to \infty }
  \E\bigl[ 
    \min\bigl\{ 1,  
      \| \varTheta_{ M_t } \| 
      \mathbbm{1}_{ 
        \{ 
          \sup_{ n \in \N } \| \varTheta_n \| < \infty 
        \} 
      }
    \bigr\}
  \bigr]
  = 0
  .
\end{split}
\end{equation}
This establishes \cref{eq:convergence_in_probability_cor}.
\end{cproof}

\subsection{Convergence in probability for quadratic optimization problems}

In this section we specialize \cref{cor:convergence_in_probability} to the case where the considered objective function is a quadratic function of the form $\fl ( \theta , \vartheta ) = \param \norm{\theta - \vartheta } ^2$.
In this case it turns out the \SGD\ sequence is automatically bounded, so the restriction to the event $\cE$ in \cref{cor:convergence_in_probability} can be removed.
We first collect in \cref{lem:boundedness0},
\cref{lem:boundedness}, and \cref{cor:boundedness2} a few well-known results regarding almost sure boundedness of random variables.
The proofs of these results are omitted, they can be found in the \href{https://arxiv.org/abs/2406.14340}{arXiv-version} of this article.

\begin{lemma}[Almost sure boundedness of measurable functions]
\label{lem:boundedness0}
Let $ ( \Omega, \cF, \mu ) $ be a measure space 
and let $ X \colon \Omega \to \R $ be measurable.  
Then 
\begin{equation}
\label{eq:lem_boundedness0}
\begin{split}
&
\textstyle 
  \sup\bigl\{ 
    c \in \R 
    \colon
    \mu\bigl( | X | \geq c \bigr) > 0
  \bigr\} 
=
  \sup\bigl\{ 
    c \in \R 
    \colon
    \mu\bigl( | X | > c \bigr) > 0
  \bigr\} 
\\ &
\textstyle 
=
  \inf\bigl\{ 
    c \in [0,\infty] 
    \colon
    \mu\bigl( | X | \geq c \bigr) = 0
  \big\}
=
  \inf\bigl\{ 
    c \in [0,\infty] 
    \colon
    \mu\bigl( | X | > c \bigr) = 0 
  \big\}
  .
\end{split}
\end{equation}
\end{lemma}

\begin{lemma}[Almost sure boundedness of i.i.d.\ random variables]
\label{lem:boundedness}
Let $ ( \Omega, \cF, \P ) $ be a probability space
and let $ X_n \colon \Omega \to \R $, $ n \in \N $, be i.i.d.\ random variables.  
Then 
\begin{equation}
\label{eq:L_infty_identity_for_X_n}
\begin{split}
&
\textstyle 
  \E\bigl[
    \sup_{ n \in \N } | X_n |
  \bigr]
=
  \inf\bigl\{ 
    c \in [0,\infty] 
    \colon
    \P\bigl( | X_1 | \leq c \bigr) = 1 
  \big\}
  .
\end{split}
\end{equation}
\end{lemma}

\begin{cor}[Almost sure boundedness of i.i.d.\ random variables]
\label{cor:boundedness2}
Let $ \dd \in \N $, 
let $ ( \Omega, \cF, \P ) $ be a probability space, 
and let $ X_n \colon \Omega \to \R^{ \dd } $, $ n \in \N $, be i.i.d.\ random variables. 
Then the following two statements are equivalent:
\begin{enumerate}[label = (\roman*)]
\item 
\label{item1:equivalence_bounded_RVs}
It holds that 
$
  \E\bigl[
    \sup_{ n \in \N } \| X_n \| 
  \bigr] < \infty 
$.
\item 
\label{item2:equivalence_bounded_RVs}
There exists $ c \in \R $ such that 
$
  \P( \| X_1 \| \leq c ) = 1
$. 
\end{enumerate}
\end{cor}
\begin{cproof}{cor:boundedness2}
\Nobs that \cref{lem:boundedness} establishes 
(\cref{item1:equivalence_bounded_RVs}
$
  \leftrightarrow 
$
\cref{item2:equivalence_bounded_RVs}).
\end{cproof}

We now establish in \cref{cor:convergence_in_probability_ex} the final result of this section,
which proves convergence of the \SGD\ method with decreasing random learning rates in the case of a quadratic objective function.
The proof of \cref{cor:convergence_in_probability_ex} relies on an application of \cref{cor:convergence_in_probability}, after verifying the condition on the stochastic error terms $D_n$.

\begin{cor}[Strong mean square convergence for quadratic optimization problems]
\label{cor:convergence_in_probability_ex}
Let 
$ \param \in (0,\infty) $, 
$ \fd \in \N $, 
let 
$ 
  \fl = ( \fl( \theta, x ) )_{ ( \theta, x ) \in \R^{ \fd } \times \R^{ \fd } } \colon \R^{ \fd } \times \R^{ \fd } \to \R 
$ 
satisfy for all $ \theta, x \in \R^{ \fd } $ that 
\begin{equation}
  \fl( \theta, x ) 
  =
  \tfrac{\param}{2} \| \theta - x \|^2
  ,
\end{equation}
let $ ( \Omega, \mathcal{F}, ( \mathbb{F}_n )_{ n \in \N_0 }, \P ) $ 
be a filtered probability space, 
for every $ n \in \N $ 
let $ \gamma_n \colon \Omega \to \R $
be $ \mathbb{F}_{ n - 1 } $-measurable, 
assume 
$
  \P\bigl( 
    \lim_{ n \to \infty }
    \gamma_n 
    =
    0
  \bigr)
  = 
  \P\bigl( 
    \sum_{ n = 1 }^{ \infty } | \gamma_n |
    =
    \infty
  \bigr)
  = 1
$, 
assume for all $ n \in \N $ that 
$
  \P(
    \gamma_{ n + 1 } \leq \gamma_n \leq \param^{ - 1 }
  )
  = 1
$,
let $ M \in \N $, 
let 
$ X_{ n, m } \colon \Omega \to \R^{ \fd } $, $ (n,m) \in \N \times \{ 1, 2, \dots, M \} $, 
be i.i.d.\ random variables, 
let 
$ \Theta \colon \N_0 \times \Omega \to \R^{ \fd } $ 
be an $ ( \mathbb{F}_n )_{ n \in \N_0 } $-adapted stochastic process 
which satisfies for all $ n \in \N $ that 
\begin{equation}
\label{eq:in_corollary_quadratic_definition_of_Theta_n}
  \Theta_n 
  = \Theta_{ n - 1 } 
  - 
  \frac{ \gamma_n }{ M }
  \left[ 
    \sum_{ m = 1 }^M
    ( \nabla_{ \theta } \fl )( \Theta_{ n - 1 }, X_{ n, m } ) 
  \right]
  ,
\end{equation}
assume for every $ n \in \N $, $ m \in \{ 1, 2, \dots, M \} $ 
that $ X_{ n, m } $
is $ \mathbb{F}_n $-measurable, 
assume for all $ n \in \N $ that 
$ \sigma( ( X_{ n, m } )_{ m \in \{ 1, 2, \dots, M \} } ) $ and $ \mathbb{F}_{ n - 1 } $ are independent, 
assume 
$ 
  \E\bigl[ 
    \| \Theta_0 \|^2 +
    \sup_{ n \in \N }
    \| X_{ n, 1 } \| 
  \bigr] < \infty 
$, 
and for every $ t \in [0,\infty) $ 
let 
$
  N_t \colon \Omega \to \R
$
satisfy 
\begin{equation}
\textstyle 
  N_t
  =
  \inf\bigl(
    \bigl\{
      n \in \N_0
      \colon 
      \sum_{ k = 1 }^n
      \gamma_k
      \geq 
      t
    \bigr\}
  \bigr)
  .
\end{equation}
Then 
\begin{equation}
\label{eq:convergence_in_probability_ex}
  \lim_{ t \to \infty }
  \E\bigl[ 
    \| \Theta_{ N_t } - \E[ X_{ 1, 1 } ] \|^2 
  \bigr]
  = 0
  .
\end{equation}
\end{cor}

\begin{cproof}{cor:convergence_in_probability_ex}
Throughout this proof let 
$
  g \colon \R^{ \fd } \to \R^{ \fd } 
$
satisfy for all $ \theta \in \R^{ \fd } $ that 
\begin{equation}
  g( \theta ) = 
  -
  \nabla_{ \theta } 
  \bigl(
    \E\bigl[
      \fl( \theta, X_{ 1, 1 } )
    \bigr]
  \bigr)
  =
  -
  \E\bigl[
    ( \nabla_{ \theta } \fl )( \theta, X_{ 1, 1 } )
  \bigr]
  =
  - \param
  \bigl(
    \theta 
    -
    \E[
      X_{ 1, 1 }
    ]
  \bigr)
\end{equation}
(cf.\ \cref{cor:boundedness2}) 
and let $ D \colon \N \times \Omega \to \R^{ \fd } $ satisfy
for all $ n \in \N $ that 
\begin{equation}
\label{eq:definition_of_D_n}
  D_n = 
  -
  g( \Theta_{ n - 1 } ) 
  -
  \frac{ 1 }{ M }
  \left[ 
    \sum_{ m = 1 }^M
    ( \nabla_{ \theta } \fl )( \Theta_{ n - 1 }, X_{ n, m } ) 
  \right]
  .
\end{equation}
\Nobs that \cref{eq:definition_of_D_n} demonstrates that 
\begin{equation}
\label{eq:definition_of_D_n_2}
  D_n = 
\textstyle 
  \param
  \bigl( 
    \Theta_{ n - 1 }
    -
    \E[ X_{ 1, 1 } ]
  \bigr)
  -
  \frac{ \param }{ M }
  \!
  \left[ 
    \sum\limits_{ m = 1 }^M
    \bigl(
      \Theta_{ n - 1 } - X_{ n, m }
    \bigr)
  \right]
  =
  \frac{ \param }{ M }
  \!
  \left[ 
    \sum\limits_{ m = 1 }^M
    \bigl(
      X_{ n, m }
      -
      \E[ X_{ n, m } ]
    \bigr)
  \right]
  .
\end{equation}
The assumption that 
$ \Theta \colon \N_0 \times \Omega \to \R^{ \fd } $ 
is an $ ( \mathbb{F}_n )_{ n \in \N_0 } $-adapted stochastic process 
and the assumption that
for all $ n \in \N $, $ m \in \{ 1, 2, \dots, M \} $ 
it holds that $ X_{ n, m } $
is $ \mathbb{F}_n $-measurable \hence 
prove that 
$ D \colon \N \times \Omega \to \R^{ \fd } $ 
is an $ ( \mathbb{F}_n )_{ n \in \N } $-adapted stochastic process. 
\Moreover 
\cref{cor:boundedness2}, 
the fact that 
$
  \E[ \sup_{ n \in \N } \| X_{ n, 1 } \| ] < \infty 
$, 
and 
the assumption that 
$
  X_{ n, m } \colon \Omega \to \R^{ \fd }
$, 
$
  (n,m) \in \{ 1, 2, \dots, M \} \times \N
$,
are i.i.d.\ random variables
imply that there exists $ \kappa \in \R $ which satisfies 
\begin{equation}
\label{eq:X_bound_by_c}
  \P\bigl(
    \| X_{ 1, 1 } \| \leq \kappa
  \bigr) = 1
  .
\end{equation}
\Nobs that 
the assumption that 
$
  X_{ n, m } \colon \Omega \to \R^{ \fd }
$, 
$
  (n,m) \in \N \times \{ 1, 2, \dots, M \}
$,
are i.i.d.\ random variables
and \cref{eq:X_bound_by_c} show that
\begin{equation}
\label{eq:X_bound_by_c_2}
\textstyle 
  \P\bigl(
    \sup_{ (n, m) \in \N \times \{ 1, 2, \dots, M \} }
    \| X_{ n, m } \| \leq \kappa 
  \bigr) = 1
  .
\end{equation}
This and \cref{eq:definition_of_D_n_2} demonstrate that 
$
  \P(
    \sup_{ n \in \N }
    \| D_n \|
    \leq 
    2 \kappa \param
  ) = 1 
$. 
\Hence that 
\begin{equation}
\label{eq:D_n_square_integrability_condition}
\textstyle 
  \E\bigl[
    \sup_{ n \in \N }
    \| D_n \|^2
  \bigr]
  \leq 
  4 \kappa^2 \param^2 
  < \infty
  .
\end{equation}
Combining \cref{eq:definition_of_D_n_2} 
and the assumption that for all $ n \in \N $ it holds that 
$ \sigma( ( X_{ n, m } )_{ m \in \{ 1, 2, \dots, M \} } ) $ and $ \mathbb{F}_{ n - 1 } $ are independent 
\hence ensures that for all $ n \in \N $ it holds $ \P $-a.s.\ that 
\begin{equation}
\label{eq:conditional_expectation_D_n_on_F_n_1}
  \E[ 
    D_n | \mathbb{F}_{ n - 1 }
  ]
  = 0 
  .
\end{equation}
\Moreover \cref{eq:X_bound_by_c_2} proves that 
\begin{equation}
\label{eq:X_bound_by_c_3}
\textstyle 
  \P\bigl(
    \sup_{ n \in \N }
    \bigl\| 
      \frac{ 1 }{ M }
      \bigl[ 
        \sum_{ m = 1 }^M
        X_{ n, m }
      \bigr]
    \bigr\| 
    \leq \kappa 
  \bigr) = 1
  .
\end{equation}
\Moreover \cref{eq:in_corollary_quadratic_definition_of_Theta_n} and \cref{eq:definition_of_D_n} show 
for all $ n \in \N $ that 
\begin{equation}
\label{eq:Theta_n_g_plus_D_n}
\begin{split}
  \Theta_n 
& = 
  \Theta_{ n - 1 } + \gamma_n g( \Theta_{ n - 1 } ) + \gamma_n D_n 
\\ &
\textstyle 
  = 
  \Theta_{ n - 1 } 
  - 
  \param \gamma_n 
  ( \Theta_{ n - 1 } - \E[ X_{ 1, 1 } ] ) 
  + 
  \displaystyle 
  \frac{ \param \gamma_n }{ M } 
  \textstyle 
  \left[ 
    \sum\limits_{ m = 1 }^M
    \bigl(
      X_{ n, m } - \E[ X_{ n, m } ] 
    \bigr)
  \right] 
\\ &
\textstyle 
  = 
  \Theta_{ n - 1 } 
  - 
  \param \gamma_n 
  \biggl( 
    \Theta_{ n - 1 } 
    -
    \displaystyle 
    \frac{ 1 }{ M } 
    \textstyle 
    \biggl[
      \sum\limits_{ m = 1 }^M
      X_{ n, m } 
    \biggr]
  \biggr) 
  .
\end{split}
\end{equation}
The fact that 
$ \E[ \| \Theta_0 \|^2 ] < \infty $, 
the fact that 
$
  \P( 
    \sup_{ n \in \N } | \param \gamma_n | \leq 1 
  ) = 1
$, 
induction, 
and 
\cref{eq:X_bound_by_c_3}
\hence prove that 
\begin{equation}
\label{eq:a_priori_bound_in_proof}
\textstyle 
  \E\bigl[ \sup_{ n \in \N_0 } \| \Theta_n \|^2 \bigr] < \infty 
  .
\end{equation}
The fact that for all $ v, w \in \R^{ \fd } $ it holds that 
$
  \| v - w \|^2 \leq 2 \| v \|^2 + 2 \| w \|^2
$
\hence shows that
\begin{equation}
\label{eq:a_priori_bound_in_proof_2}
\begin{split}
\textstyle 
  \E\bigl[ 
    \sup_{ t \in [0,\infty) }
    \| \Theta_{ N_t } - \E[ X_{ 1, 1 } ] \|^2 
  \bigr]
& 
\textstyle 
  \leq 
  \E\bigl[ 
    \sup_{ n \in \N_0 }
    \| \Theta_n - \E[ X_{ 1, 1 } ] \|^2 
  \bigr]
\\ &
\textstyle 
\leq 
  2 \,
  \E\bigl[ 
    \sup_{ n \in \N_0 }
    \| \Theta_n \|^2
  \bigr]
  +
  2 \, \bigl| \E\bigl[ \| X_{ 1, 1 } \| \bigr] \bigr|^2 
< 
  \infty 
  .
\end{split}
\end{equation}
\Moreover \cref{eq:a_priori_bound_in_proof} ensures that 
$
  \P\bigl(
    \sup_{ n \in \N_0 } \| \Theta_n \|
    < \infty
  \bigr) 
  =
  \P\bigl(
    \sup_{ n \in \N_0 } \| \Theta_n \|^2
    < \infty
  \bigr) 
  = 1
$.
Combining this, 
\cref{eq:D_n_square_integrability_condition}, 
\cref{eq:conditional_expectation_D_n_on_F_n_1}, 
\cref{eq:Theta_n_g_plus_D_n}, 
\cref{cor:convergence_in_probability}, 
and the fact that 
$ D \colon \N \times \Omega \to \R^{ \fd } $ 
is an $ ( \mathbb{F}_n )_{ n \in \N } $-adapted stochastic process 
shows that 
\begin{equation}
\begin{split}
& 
  \limsup_{ t \to \infty }
  \E\bigl[ 
    \min\bigl\{ 
      1 ,  
      \| \Theta_{ N_t } - \E[ X_{ 1, 1 } ] \|^2 
    \bigr\}
  \bigr]
\\ &
  =
  \limsup_{ t \to \infty }
  \E\bigl[ 
    \min\bigl\{ 
      1 ,  
      \| \Theta_{ N_t } - \E[ X_{ 1, 1 } ] \|^2 
      \mathbbm{1}_{ 
        \{ 
          \sup_{ n \in \N } \| \Theta_n \|
          < \infty 
        \} 
      }
    \bigr\}
  \bigr]
  =
  0
  .
\end{split}
\end{equation}
\cref{lem:convergence_in_probability_charac} \hence demonstrates for all $ \varepsilon \in (0,\infty) $ that 
\begin{equation}
  \limsup_{ t \to \infty }
  \P\bigl(
    \| \Theta_{ N_t } - \E[ X_{ 1, 1 } ] \|^2 
    > \varepsilon
  \bigr)
  =
  0
  .
\end{equation}
Combining this and \cref{eq:a_priori_bound_in_proof_2} with 
Lebesgue's theorem of dominated convergence 
(cf., for example, \cite[Definition 6.2 (i) and Corollary~6.26]{Klenke2014}) 
establishes \cref{eq:convergence_in_probability_ex}.
\end{cproof}

\section{Analysis of the test loss for SGD with constant learning rate}
\label{sec:test_loss_strictly_increases}

The next goal is to apply \cref{cor:convergence_in_probability_ex} to a simplified version of the methodology described in \cref{sec:algorithmic_description} where the learning rate is decreased whenever the loss value increased in the last training step.
To do this, we need to verify the assumptions of \cref{cor:convergence_in_probability_ex}.
The main step is to ensure convergence of the learning rates to zero, for which we need to verify that the loss value increases almost surely infinitely many times during the training process.
This is the content of the next two sections.

\subsection{A priori bounds for certain stochastic processes}

We first show in \cref{lem:a_priori} suitable elementary a priori bounds for the \SGD\ iterates 
in the case of a simple class of quadratic target functions.
Our proof of \cref{lem:a_priori} uses the elementary estimate in \cref{lem:elementary}.

\begin{lemma}[Elementary real estimate]
\label{lem:elementary}
It holds for all $ x \in (0,2) $ that
\begin{equation}
  \frac{ x }{ 
    ( 1 - | 1 - x | )
  }
= 
  \frac{ 
    1 + | 1 - x | 
  }{ 
    ( 2 - x ) 
  }
\leq 
  \frac{ 2 }{ ( 2 - x ) }
  .
\end{equation}
\end{lemma}

\begin{cproof}{lem:elementary}
\Nobs that the fact that for all 
$ y \in \R $ it holds that 
$
  ( 1 - y ) ( 1 + y ) = 1 - y^2
$
ensures that 
for all $ x \in \R \backslash \{ 0 \} $ it holds that 
\begin{equation}
\begin{split}
  \frac{ x }{ 
    ( 1 - | 1 - x | )
  }
& =
  \frac{ 
    x 
    ( 1 + | 1 - x | )
  }{ 
    ( 1 - | 1 - x | )
    ( 1 + | 1 - x | )
  }
=
  \frac{ 
    x 
    ( 1 + | 1 - x | )
  }{ 
    ( 
      1 - ( 1 - x )^2 
    )
  }
=
  \frac{ 
    x 
    ( 1 + | 1 - x | )
  }{ 
    ( 
      1 - ( 1 - 2 x + x^2 ) 
    )
  }
\\ &
=
  \frac{ 
    x 
    ( 1 + | 1 - x | )
  }{ 
    ( 
      2 x - x^2  
    )
  }
=
  \frac{ 
    x 
    ( 1 + | 1 - x | )
  }{ 
    x
    ( 
      2 - x  
    )
  }
=
  \frac{ 
    ( 1 + | 1 - x | )
  }{ 
    ( 
      2 - x  
    )
  }
  .
\end{split}
\end{equation}
\end{cproof}

\begin{lemma}[A priori bounds]
\label{lem:a_priori}
Let 
$ \fd \in \N $, 
$ \param \in (0,\infty) $, 
let $ \fl \colon \R^{ \fd } \times \R^{ \fd } \to \R $ satisfy 
for all $ \theta, \vartheta \in \R^{ \fd } $ that 
\begin{equation}
\label{eq:definition_of_fl_lem_a_priori}
  \fl( \theta, \vartheta ) 
  =
\textstyle 
  \frac{ \param }{ 2 } \| \theta - \vartheta \|^2
  ,
\end{equation}
let $ ( \Omega, \mathcal{F}, \P ) $ be a probability space, 
let $ X_{ n, m } \colon \Omega \to \R^{ \fd } $, $ n, m \in \Z $, be i.i.d.\ random variables, 
let $ M, \fM \in \N $, $ \gamma \in ( 0, \nicefrac{ 2 }{ \param } ) $, 
and let 
$ \Theta \colon \N_0 \times \Omega \to \R^{ \fd } $
satisfy for all $ n \in \N $ that 
\begin{equation}
\label{eq:example_definition_GD_lem_a_priori}
  \Theta_n 
  = \Theta_{ n - 1 } 
  - 
  \frac{ \gamma }{ M }
  \left[ 
    \sum_{ m = 1 }^M
    ( \nabla_{ \theta } \fl )( \Theta_{ n - 1 }, X_{ n, m } ) 
  \right]
  .
\end{equation}
Then 
\begin{enumerate}[label = (\roman*)]
\item 
\label{item:i_a_priori_bounds}
it holds that 
$
  | 1 - \param \gamma | < 1
$,
\item 
\label{item:ii_a_priori_bounds}
it holds for all $ n \in \N_0 $ that 
\begin{equation}
  \Theta_n 
  =
  ( 1 - \param \gamma )^n \Theta_0 
  +
  \sum_{ k = 1 }^n
  \left[
  \frac{ 
    ( 1 - \param \gamma )^{ n - k }
    \param \gamma 
  }{ M } 
  \left[ 
    \sum\limits_{ m = 1 }^M X_{ k, m }
  \right] 
  \right] 
  ,
\end{equation}
and 
\item 
\label{item:iii_a_priori_bounds}
it holds for all $ n \in \N_0 $ that 
\begin{equation}
\begin{split}
\textstyle 
  \| \Theta_n \|
&
\textstyle 
  \leq 
  | 1 - \param \gamma |^n 
  \| \Theta_0 \|
  +
  \param \gamma 
  \biggl[
    \sum\limits_{ k = 0 }^{ n - 1 }
    | 1 - \param \gamma |^k
  \biggr]
  \displaystyle 
  \biggl[ 
    \sup_{ v, w \in \N_0 }
    \| 
      X_{ v, w }
    \|
  \biggr]  
\\ &
\textstyle 
\leq 
  \| \Theta_0 \|
  +
  \biggl[
  \displaystyle
    \frac{ 2 }{ 2 - \param \gamma }
  \biggr]
  \biggl[ 
    \sup_{ v, w \in \N_0 }
    \| 
      X_{ v, w }
    \|
  \biggr] 
  .
\end{split}
\end{equation}
\end{enumerate}
\end{lemma}

\begin{cproof}{lem:a_priori}
Throughout this proof for every $ n \in \N $ 
let $ \cX_n \colon \Omega \to \R^{ \fd } $ 
satisfy 
\begin{equation}
\label{eq:definition_of_cX_n}
  \cX_n 
  =
  \frac{ \param \gamma }{ M } 
  \left[ 
    \sum_{ m = 1 }^M X_{ n, m }
  \right] 
  .
\end{equation}
\Nobs that the fact that $ 0 < \gamma < \nicefrac{ 2 }{ \param } $ 
ensures that 
$
  0 < \param \gamma < 2
$. 
\Hence that
$
  -1 < \param \gamma - 1 < 1
$. 
\Hence that 
$
  1 > 1 - \param \gamma > -1
$.
This establishes \cref{item:i_a_priori_bounds}. 
\Nobs that \cref{eq:example_definition_GD_lem_a_priori,eq:definition_of_cX_n} 
ensure that for all $ n \in \N $ it holds that 
\begin{equation}
  \Theta_n 
  =
  ( 1 - \param \gamma ) 
  \Theta_{ n - 1 } 
  +
  \cX_n 
  .
\end{equation}
\Hence for all $ n \in \N_0 $ that 
\begin{equation}
\label{eq:Theta_n_recursion}
  \Theta_n 
  =
  ( 1 - \param \gamma )^n \Theta_0 
  +
  \sum_{ k = 1 }^n
  ( 1 - \param \gamma )^{ n - k }
  \cX_k
  .
\end{equation}
Combining this and \cref{eq:definition_of_cX_n} 
establishes \cref{item:ii_a_priori_bounds}. 
\Nobs that 
\cref{eq:definition_of_cX_n}
and 
\cref{eq:Theta_n_recursion} 
ensure for all $ n \in \N_0 $ that 
\begin{equation}
\begin{split}
  \| \Theta_n \|
& 
  \leq 
  | 1 - \param \gamma |^n 
  \| \Theta_0 \|
  +
  \left[
    \sum_{ k = 1 }^n
    | 1 - \param \gamma |^{ n - k }
  \right] 
  \left[ 
    \sup_{ k \in \N_0 }
    \| \cX_k \|
  \right] 
\\ &
  =
  | 1 - \param \gamma |^n 
  \| \Theta_0 \|
  +
  \left[
    \sum_{ k = 1 }^n
    | 1 - \param \gamma |^k
  \right] 
  \left[ 
    \sup_{ k \in \N_0 }
    \| \cX_k \|
  \right] 
\\ &
  \leq 
  | 1 - \param \gamma |^n 
  \| \Theta_0 \|
  +
  \param \gamma 
  \left[
    \sum_{ k = 1 }^n
    | 1 - \param \gamma |^k
  \right] 
  \left[ 
    \sup_{ v, w \in \N_0 }
    \| X_{ v, w } \|
  \right] 
\\ &
  \leq 
  | 1 - \param \gamma |^n 
  \| \Theta_0 \|
  +
  \param \gamma 
  \left[
    \sum_{ k = 1 }^{ \infty }
    | 1 - \param \gamma |^k
  \right] 
  \left[ 
    \sup_{ v, w \in \N_0 }
    \| X_{ v, w } \|
  \right] 
\\ &
  =
  | 1 - \param \gamma |^n 
  \| \Theta_0 \|
  +
  \left[
    \frac{ 
      \param \gamma 
    }{
      (
        1 - | 1 - \param \gamma |
      )
    }
  \right]
  \left[ 
    \sup_{ v, w \in \N_0 }
    \| X_{ v, w } \|
  \right] 
  .
\end{split}
\end{equation}
\cref{lem:elementary} \hence shows for all $ n \in \N_0 $ that 
\begin{equation}
\begin{split}
  \| \Theta_n \|
& \leq 
  | 1 - \param \gamma |^n 
  \| \Theta_0 \|
  +
  \param \gamma 
  \left[
    \sum_{ k = 1 }^n
    | 1 - \param \gamma |^k
  \right] 
  \left[ 
    \sup_{ v, w \in \N_0 }
    \| X_{ v, w } \|
  \right] 
\\ & \leq 
  | 1 - \param \gamma |^n 
  \| \Theta_0 \|
  +
  \left[
    \frac{
      2 
    }{
      2 - \param \gamma 
    }
  \right]
  \left[ 
    \sup_{ v, w \in \N_0 }
    \| X_{ v, w } \|
  \right] 
  .
\end{split}
\end{equation}
This establishes \cref{item:iii_a_priori_bounds}. 
\end{cproof}

\subsection{Properties of certain stationary stochastic processes}

In this subsection we establish two elementary results on 
stationary SGD processes for the considered class of 
quadratic optimization problems. 
More specifically, we show in \cref{lem:invariant_2} below 
that suitable \SGD\ sequences obtained from a linear recursion are stationary stochastic processes 
provided that the initial distribution is chosen accordingly.
Thereafter, we establish in \cref{lem:invariant} a version of the boundedness result in \cref{lem:a_priori} which holds in the case of this initial distribution.

\begin{lemma}[Stationary processes]
\label{lem:invariant_2}
Let 
$ \fd \in \N $, 
let $ ( \Omega, \mathcal{F}, \P ) $ be a probability space, 
let $ \cX_n \colon \Omega \to \R^{ \fd } $, $ n \in \Z $, be i.i.d.\ random variables, 
assume 
$
  \{ 
    \sup_{ n \in \N } \| \cX_{ - n } \| < \infty 
  \} = \Omega 
$, 
let $ \Gamma \in (0,1) $, 
let 
$ \Theta \colon \N_0 \times \Omega \to \R^{ \fd } $
satisfy for all $ n \in \N $ that 
$
  \Theta_0 
  = 
  \sum_{ n = 0 }^{ \infty }
  \Gamma^n \cX_{ - n } 
$
and 
\begin{equation}
\label{eq:example_definition_GD_lem_invariant_2}
  \Theta_n 
  = 
  \Gamma
  \Theta_{ n - 1 }
  +
  \cX_n 
  ,
\end{equation}
and let $ N \in \N_0 $. 
Then 
it holds for all $ n, m \in \N_0 $ that 
\begin{equation}
\label{eq:identity_stationary_process_distribution}
  ( \Theta_n, \Theta_{ n + 1 }, \dots, \Theta_{ n + N } )( \P ) 
  = 
  ( \Theta_m, \Theta_{ m + 1 }, \dots, \Theta_{ m + N } )( \P ) 
  .
\end{equation}
\end{lemma}

\begin{cproof}{lem:invariant_2}
Throughout this proof let 
$
  \psi \colon \R^{ \fd } \times ( \R^{ \fd } )^N \to ( \R^{ \fd } )^{ N + 1 }
$
satisfy for all $ \theta, x_1, \allowbreak x_2, \allowbreak \dots, x_N \in \R^{ \fd } $ that
\begin{equation}
\label{eq:definition_of_psi}
\textstyle 
  \psi( \theta, x_1, x_2, \dots,  x_N ) 
  = 
  \bigl(
    \theta , 
    \Gamma \theta + x_1 ,
    \Gamma^2 \theta + \Gamma x_1 + x_2 ,
    \dots ,
    \Gamma^N \theta 
    + 
    \sum_{ k = 1 }^N \Gamma^{ N - k } x_k 
  \bigr)
  .
\end{equation}
\Nobs that \cref{eq:example_definition_GD_lem_invariant_2} shows for all 
$ n, m \in \N_0 $ with $ m \leq n $ that
\begin{equation}
  \Theta_n 
  =
  \Gamma^{ n - m }
  \Theta_m
  +
  \sum_{ k = m + 1 }^n
  \Gamma^{ n - k }
  \cX_k
  .
\end{equation}
This and \cref{eq:definition_of_psi} 
ensure for all $ n \in \N_0 $ that
\begin{equation}
\label{eq:psi_representation}
  ( \Theta_n , \Theta_{ n + 1 }, \dots, \Theta_{ n + N } )
  =
  \psi( \Theta_n, \cX_{ n + 1 }, \cX_{ n + 2 }, \dots, \cX_{ n + N } )
  .
\end{equation}
\Moreover \cref{eq:example_definition_GD_lem_invariant_2} 
and the assumption that 
$
  \Theta_0 
  = 
  \sum_{ n = 0 }^{ \infty }
  \Gamma^n \cX_{ - n } 
$
assure that 
\begin{equation}
\label{eq:representation_Theta_n_history}
\begin{split}
  \Theta_n 
& 
  = 
  \Gamma^n \Theta_0 
  + 
  \left[ 
    \sum_{ k = 1 }^n \Gamma^{ n - k } \cX_k
  \right] 
  =
  \Gamma^n
  \left[
    \sum_{ k = 0 }^{ \infty } \Gamma^k \cX_{ - k }
  \right] 
  + 
  \left[ 
    \sum_{ k = 1 }^n \Gamma^{ n - k } \cX_k
  \right] 
\\ &
  =
  \left[
    \sum_{ k = 0 }^{ \infty } \Gamma^{ n + k } \cX_{ - k }
  \right] 
  + 
  \left[
    \sum_{ k = 1 }^n \Gamma^{ n - k } \cX_k
  \right]
  =
  \left[
    \sum_{ k = n }^{ \infty } \Gamma^{ k } \cX_{ n - k }
  \right] 
  + 
  \left[
    \sum_{ k = 0 }^{ n - 1 } \Gamma^k \cX_{ n - k }
  \right]
  =
  \sum_{ k = 0 }^{ \infty } \Gamma^k \cX_{ n - k }
  .
\end{split}
\end{equation}
The assumption that 
$ \cX_n \colon \Omega \to \R^{ \fd } $, $ n \in \Z $, 
are i.i.d.\ random variables \hence shows for all $ n, m \in \N_0 $ that
\begin{equation}
\label{eq:stationary_process_identity}
  \Theta_n( \P ) 
  = 
  \Theta_m( \P ) 
  .
\end{equation}
\Moreover 
\cref{eq:representation_Theta_n_history}
and the assumption that 
$ \cX_n \colon \Omega \to \R^{ \fd } $, $ n \in \Z $, 
are i.i.d.\ random variables 
show that for all 
$ n \in \N_0 $ it holds that 
%
\begin{equation}
\label{eq:stationary_process_independence_property}
  \Theta_n, 
  \qquad 
  \cX_{ n + 1 },
  \qquad 
  \cX_{ n + 2 },
  \qquad 
  ,
  \dots 
  ,
  \qquad 
  \cX_{ n + N - 1 }, 
  \qquad 
  \text{and}
  \qquad 
  \cX_{ n + N }
\end{equation}
are independent. 
Combining this 
and the change-of-variables formula for pushforward measures 
with \cref{eq:psi_representation} demonstrates that for all 
$ n \in \N_0 $, 
$ 
  \varphi \in 
  C( 
    ( \R^{ \fd } )^N, \R 
  ) 
$ 
with $ \sup_{ x \in ( \R^{ \fd } )^N } | \varphi(x) | < \infty $ 
it holds that 
\begin{equation}
\begin{split}
&
  \E\bigl[
    \varphi( \Theta_n , \Theta_{ n + 1 }, \dots, \Theta_{ n + N } )
  \bigr]
  =
  \E\bigl[
    \varphi\bigl( \psi( \Theta_n, \cX_{ n + 1 }, \cX_{ n + 2 }, \dots, \cX_{ n + N } ) \bigr)
  \bigr]
\\ &
=
  \int_{ ( \R^{ \fd } )^{ N + 1 } }
  \varphi( \psi( x ) )
  \,
  \bigl( 
    ( \Theta_n, \cX_{ n + 1 }, \cX_{ n + 2 }, \dots, \cX_{ n + N } )( \P )
  \bigr)( \d x )
\\ & 
=
  \int_{ ( \R^{ \fd } )^{ N + 1 } }
  \varphi( \psi( x ) )
  \,
  \bigl( 
    ( \Theta_n( \P ) ) 
    \otimes 
    ( \cX_{ n + 1 }( \P ) )
    \otimes 
    \dots 
    \otimes 
    ( \cX_{ n + N }( \P ) )
  \bigr)( \d x )
\end{split}
\end{equation}
The assumption that 
$ \cX_n \colon \Omega \to \R^{ \fd } $, $ n \in \Z $, 
are identically distributed random variables, 
\cref{eq:stationary_process_identity}, 
and 
\cref{eq:stationary_process_independence_property}
\hence show that for all 
$ n, m \in \N_0 $, 
$ 
  \varphi \in 
  C( 
    ( \R^{ \fd } )^N, \R 
  ) 
$ 
with $ \sup_{ x \in ( \R^{ \fd } )^N } | \varphi(x) | < \infty $ 
it holds that 
\begin{equation}
\begin{split}
&
  \E\bigl[
    \varphi( \Theta_n , \Theta_{ n + 1 }, \dots, \Theta_{ n + N } )
  \bigr]
\\ & 
=
  \int_{ ( \R^{ \fd } )^{ N + 1 } }
  \varphi( \psi( x ) )
  \,
  \bigl( 
    ( \Theta_n( \P ) ) 
    \otimes 
    ( \cX_{ m + 1 }( \P ) )
    \otimes 
    \dots 
    \otimes 
    ( \cX_{ m + N }( \P ) )
  \bigr)( \d x )
\\ & 
=
  \int_{ ( \R^{ \fd } )^{ N + 1 } }
  \varphi( \psi( x ) )
  \,
  \bigl( 
    ( \Theta_m( \P ) ) 
    \otimes 
    ( \cX_{ m + 1 }( \P ) )
    \otimes 
    \dots 
    \otimes 
    ( \cX_{ m + N }( \P ) )
  \bigr)( \d x )
\\ & 
=
  \int_{ ( \R^{ \fd } )^{ N + 1 } }
  \varphi( \psi( x ) )
  \,
  \bigl( 
    ( \Theta_m, \cX_{ m + 1 }, \cX_{ m + 2 }, \dots, \cX_{ m + N } )( \P )
  \bigr)( \d x )
  .
\end{split}
\end{equation}
The change-of-variables formula for pushforward measures 
\hence shows for all $ n, m \in \N_0 $ that 
\begin{equation}
\begin{split}
  \E\bigl[
    \varphi( \Theta_n , \Theta_{ n + 1 }, \dots, \Theta_{ n + N } )
  \bigr]
& =
  \E\bigl[
    \varphi\bigl( \psi( \Theta_m, \cX_{ m + 1 }, \cX_{ m + 2 }, \dots, \cX_{ m + N } ) \bigr)
  \bigr]
\\ & =
  \E\bigl[
    \varphi( \Theta_m , \Theta_{ m + 1 }, \dots, \Theta_{ m + N } )
  \bigr]
  .
\end{split}
\end{equation}
This establishes \cref{eq:identity_stationary_process_distribution}. 
\end{cproof}

\begin{lemma}[A priori bounds for stationary processes]
\label{lem:invariant}
Let 
$ \fd \in \N $, 
$ \param \in (0,\infty) $, 
let $ \fl \colon \R^{ \fd } \times \R^{ \fd } \to \R $ satisfy 
for all $ \theta, \vartheta \in \R^{ \fd } $ that 
\begin{equation}
\label{eq:definition_of_fl_lem_invariant}
  \fl( \theta, \vartheta ) 
  =
\textstyle 
  \frac{ \param }{ 2 } \| \theta - \vartheta \|^2
  ,
\end{equation}
let $ ( \Omega, \mathcal{F}, \P ) $ be a probability space, 
let $ X_{ n, m } \colon \Omega \to \R^{ \fd } $, $ n, m \in \Z $, be i.i.d.\ random variables, 
let 
$ M \in \N $, 
$ \gamma \in ( 0, \nicefrac{ 2 }{ \param } ) $, 
let 
$ \Theta \colon \N_0 \times \Omega \to \R^{ \fd } $
satisfy for all $ n \in \N $ that 
\begin{equation}
\label{eq:example_definition_GD_lem_invariant}
  \Theta_n 
  = \Theta_{ n - 1 } 
  - 
  \frac{ \gamma }{ M }
  \left[ 
    \sum_{ m = 1 }^M
    ( \nabla_{ \theta } \fl )( \Theta_{ n - 1 }, X_{ n, m } ) 
  \right]
  ,
\end{equation}
assume 
$ 
  \{ \sup_{ n, m \in \Z } \| X_{ n, m } \| \} = \Omega 
$, 
and let $ \kappa \in \Z $ satisfy
\begin{equation}
\label{eq:definition_of_chi_lem_invariant}
  \Theta_0 = 
  \sum_{ n = 0 }^{ \infty }
  \left(
    \frac{ 
      \param \gamma ( 1 - \param \gamma )^n
    }{ M }
    \left[
      \sum_{ m = \kappa + 1 }^{ \kappa + M }
      X_{ -n, m }
    \right]
  \right)
  .
\end{equation}
Then 
\begin{equation}
\label{eq:a_priori_invariant_measure}
\begin{split}
  \sup_{ n \in \N_0 }
  \| \Theta_n \|
&
\textstyle 
\leq 
  \biggl[
  \displaystyle
    \frac{ 2 }{ 2 - \param \gamma }
  \biggr]
  \biggl[ 
    \sup_{ v, w \in \Z }
    \| 
      X_{ v, w }
    \|
  \biggr] 
  .
\end{split}
\end{equation}
\end{lemma}

\begin{cproof}{lem:invariant}
\Nobs that \cref{eq:definition_of_chi_lem_invariant} 
ensures for all $ n \in \N_0 $ that 
\begin{equation}
\begin{split}
  \| \Theta_0 \| 
& 
\leq 
  \sum_{ n = 0 }^{ \infty }
  \left\|
    \frac{ 
      \param \gamma ( 1 - \param \gamma )^n
    }{ M }
    \left[
      \textstyle 
      \sum\limits_{ m = \kappa + 1 }^{ \kappa + M }
      \displaystyle 
      X_{ -n, m }
    \right]
  \right\|
\leq 
  \sum_{ n = 0 }^{ \infty }
  \biggl(
    \frac{ 
      \param \gamma | 1 - \param \gamma |^n
    }{ M }
    \left[
      \textstyle 
      \sum\limits_{ m = \kappa + 1 }^{ \kappa + M }
      \displaystyle 
      \|
        X_{ -n, m }
      \|
    \right]
  \biggr)
\\ &
\leq 
  \sum_{ n = 0 }^{ \infty }
  \biggl(
    \param \gamma | 1 - \param \gamma |^n
    \biggl[
      \textstyle 
      \sup\limits_{ v, w \in \Z }
      \displaystyle 
      \|
        X_{ v, w }
      \|
    \biggr]
  \biggr)
=
  \param \gamma 
  \biggl[
    \textstyle 
    \sum\limits_{ k = 0 }^{ \infty }
    | 1 - \param \gamma |^k
  \biggr]
    \biggl[
      \textstyle 
      \sup\limits_{ v, w \in \Z }
      \displaystyle 
      \|
        X_{ v, w }
      \|
    \biggr]
  .
\end{split}
\end{equation}
\cref{lem:a_priori} \hence implies for all $ n \in \N_0 $ that
\begin{equation}
\begin{split}
  \| \Theta_n \|
& \leq 
  | 1 - \param \gamma |^n \| \Theta_0 \|
  +
  \param \gamma 
  \biggl[
    \textstyle 
    \sum\limits_{ k = 0 }^{ n - 1 } 
    \displaystyle 
    | 1 - \param \gamma |^k 
  \biggr]
  \biggl[
    \sup_{ v, w \in \N_0 }
    \| X_{ v, w } \|
  \biggr]
\\ & \leq 
  \param \gamma 
  | 1 - \param \gamma |^n 
  \biggl[
    \textstyle 
    \sum\limits_{ k = 0 }^{ \infty }
    \displaystyle 
    | 1 - \param \gamma |^k 
  \biggr]
  \biggl[
    \sup_{ v, w \in \Z }
    \| X_{ v, w } \|
  \biggr]
  +
  \param \gamma 
  \biggl[
    \textstyle 
    \sum\limits_{ k = 0 }^{ n - 1 } 
    \displaystyle 
    | 1 - \param \gamma |^k 
  \biggr]
  \biggl[
    \sup_{ v, w \in \Z }
    \| X_{ v, w } \|
  \biggr]
\\ & =
  \param \gamma 
  \biggl[
    \textstyle 
    \sum\limits_{ k = n }^{ \infty }
    \displaystyle 
    | 1 - \param \gamma |^k 
  \biggr]
  \biggl[
    \sup_{ v, w \in \Z }
    \| X_{ v, w } \|
  \biggr]
  +
  \param \gamma 
  \biggl[
    \textstyle 
    \sum\limits_{ k = 0 }^{ n - 1 } 
    \displaystyle 
    | 1 - \param \gamma |^k 
  \biggr]
  \biggl[
    \sup_{ v, w \in \Z }
    \| X_{ v, w } \|
  \biggr]
\\ & =
  \param \gamma 
  \biggl[
    \textstyle 
    \sum\limits_{ k = 0 }^{ \infty }
    \displaystyle 
    | 1 - \param \gamma |^k 
  \biggr]
  \biggl[
    \sup_{ v, w \in \Z }
    \| X_{ v, w } \|
  \biggr]
  .
\end{split}
\end{equation}
The fact that $ | 1 - \param \gamma | < 1 $ 
and \cref{lem:elementary} \hence show for all $ n \in \N_0 $ that
\begin{equation}
\begin{split}
  \| \Theta_n \| 
& 
\leq 
  \param \gamma 
  \biggl[
    \frac{ 1 }{ 
      1 - | 1 - \param \gamma |
    }
  \biggr]
  \biggl[
    \textstyle 
    \sup\limits_{ v, w \in \Z }
    \displaystyle 
    \|
      X_{ v, w }
    \|
  \biggr]
\leq 
  \biggl[
    \frac{ 2 }{ 
      2 - \param \gamma 
    }
  \biggr]
  \biggl[
    \textstyle 
    \sup\limits_{ v, w \in \Z }
    \displaystyle 
    \|
      X_{ v, w }
    \|
  \biggr]
  .
\end{split}
\end{equation}
This establishes \cref{eq:a_priori_invariant_measure}. 
\end{cproof}

\subsection{Lower bounds for the probability that the test loss strictly increases}

The next result, \cref{prop:SGD_constant}, is the main result of this section.
It shows for the considered quadratic optimization problem that the test loss, computed with respect to the random variables $X_{n, -m}$, $m \in \cu{1, 2, \ldots, \fM }$, which are independent of the data used for the training process,
increases infinitely many times during the training with probability $1$ 
(cf.~\cref{eq:to_prove_prop:SGD_constant} below).
For this we assume that the probability that 
the test loss increases in a single step starting in the invariant measure is positive 
(cf.~\cref{eq:assumption_invariant_measure} below).
This assumption will be verified in \cref{sec:invariant_measure}.

\begin{prop}[\SGD\ with constant learning rates]
\label{prop:SGD_constant}
Let 
$ \fd \in \N $, 
$ \param \in (0,\infty) $, 
let $ \fl \colon \R^{ \fd } \times \R^{ \fd } \to \R $ satisfy 
for all $ \theta, \vartheta \in \R^{ \fd } $ that 
\begin{equation}
\label{eq:definition_of_fl}
  \fl( \theta, \vartheta ) 
  =
\textstyle 
  \frac{ \param }{ 2 } \| \theta - \vartheta \|^2
  ,
\end{equation}
let $ \xi \in \R^{ \fd } $, $ \gamma \in ( 0, \nicefrac{ 2 }{ \param } ) $, 
let $ ( \Omega, \mathcal{F}, \P ) $ be a probability space, 
let $ X_{ n, m } \colon \Omega \to \R^{ \fd } $, $ n, m \in \Z $, be i.i.d.\ random variables, 
let $ M, \fM \in \N $, 
let 
$ \Theta \colon \N_0 \times \Omega \to \R^{ \fd } $
satisfy for all $ n \in \N $ that 
\begin{equation}
\label{eq:example_definition_GD}
  \Theta_n 
  = \Theta_{ n - 1 } 
  - 
  \frac{ \gamma }{ M }
  \left[ 
    \sum_{ m = 1 }^M
    ( \nabla_{ \theta } \fl )( \Theta_{ n - 1 }, X_{ n, m } ) 
  \right]
  ,
\end{equation}
assume 
$ 
  \{ \sup_{ n, m \in \Z } \| X_{ n, m } \| < \infty \} = \Omega 
$, 
assume that 
$
  ( X_{ n, m } )_{ (n, m) \in \Z^2 }
$
and 
$ \Theta_0 $
are independent, 
let $ \chi \colon \Omega \to \R^{ \fd } $ satisfy 
\begin{equation}
\label{eq:definition_of_chi}
  \chi = 
  \sum_{ n = 0 }^{ \infty }
  \left(
    \frac{ 
      \param \gamma ( 1 - \param \gamma )^n
    }{ M }
    \left[
      \sum_{ m = 1 }^M
      X_{ -n, m }
    \right]
  \right)
\end{equation}
assume 
\begin{equation}
\label{eq:assumption_invariant_measure}
\textstyle 
  \P\biggl(
  \,
    \biggl\| 
      ( 1 - \param \gamma ) \chi 
      + 
      \displaystyle 
      \frac{ \param \gamma }{ M } 
      \textstyle 
      \biggl[ 
        \sum\limits_{ m = 1 }^{ M } X_{ 1, m }
      \biggr]
      - 
      \displaystyle 
      \frac{ 1 }{ \fM } 
      \textstyle 
      \biggl[ 
        \sum\limits_{ m = 1 }^{ \fM } X_{ 2, m }
      \biggr]
    \biggr\|
    > 
    \biggl\| 
      \chi 
      - 
      \displaystyle 
      \frac{ 1 }{ \fM } 
      \textstyle 
      \biggl[ 
        \sum\limits_{ m = 1 }^{ \fM } X_{ 2, m }
      \biggr]
    \biggr\|
    \,
  \biggr) > 0
  ,
\end{equation}
and let $ N \in \N $. 
Then 
\begin{equation}
\label{eq:to_prove_prop:SGD_constant}
\textstyle 
  \P\biggl(
    \bigcup\limits_{ n = N }^{ \infty }
    \biggl\{ 
      \sum\limits_{ m = 1 }^{ \fM }
      \fl ( \Theta_n, X_{ n, - m } ) 
      >
      \sum\limits_{ m = 1 }^{ \fM }
      \fl ( \Theta_{ n - 1 }, X_{ n, - m } ) 
    \biggr\} 
  \biggr)
  = 1
  .
\end{equation}
\end{prop}

Regarding the proof of \cref{prop:SGD_constant},
we note that the distribution of the random variable $ \chi $ (cf.~\cref{eq:definition_of_chi}) 
is the invariant measure of the \SGD\ process.
In our proof of \cref{prop:SGD_constant} 
we establish exponential convergence towards this measure in a suitable sense (see \cref{eq:error_estimate_theta_minus_phi} below),
by constructing a coupling with certain stationary processes $\phi^V$ satisfying a similar recursion
(see \cref{eq:definition_of_phi} and \cref{eq:definition_of_phi_2}).
It is therefore enough to verify the increase property for these stationary processes (see \cref{eq:probability_varepsilon_strictly_bigger_than_zero}).

Similar results regarding the invariant distribution of certain \SGD\ type stochastic processes have been obtained in Dieuleveut et al.~\cite{DieuleveutBach2020}.

\begin{cproof}{prop:SGD_constant}
Throughout this proof 
let 
$
  \varTheta \colon \N \times \Omega 
  \to
  \R^{ \fd } \times \R^{ \fd } \times ( \R^{ \fd } )^{ \fM } 
$
satisfy for all $ n \in \N $ that 
\begin{equation}
\label{eq:definition_of_varTheta_process}
  \varTheta_n 
  = 
  ( 
    \Theta_n, \Theta_{ n - 1 }, 
    X_{ n, -1 }, X_{ n, -2 }, \dots, X_{ n, -\fM } 
  )
  ,
\end{equation}
let 
$ f \colon \R^{ \fd } \times \R^{ \fd } \times ( \R^{ \fd } )^{ \fM } \to \R $
satisfy 
for all 
$ \theta, \vartheta \in \R^{ \fd } $, 
$ x_1, x_2, \dots, x_{ \fM } \in \R^{ \fd } $
that 
\begin{equation}
\label{eq:definition_of_f}
\textstyle 
  f( \theta, \vartheta, x_1, \dots, x_{ \fM } ) 
  = 
    \biggl[
      \sum\limits_{ m = 1 }^{ \fM }
      \fl ( \theta, x_m ) 
    \biggr]
    -
    \biggl[
      \sum\limits_{ m = 1 }^{ \fM }
      \fl ( \vartheta, x_m ) 
    \biggr]
  ,
\end{equation}
for every 
$ \Va \in \Z $
let 
$ 
  \Phi^{ \Va } = 
  ( \phi^{ \Va }, \varphi^{ \Va }, \psi^{ \Va, 1 }, \dots, \psi^{ \Va, \fM } )
  \colon 
  \N_0 \times \Omega \to 
  \R^{ \fd } \times \R^{ \fd } \times ( \R^{ \fd } )^{ \fM } 
$ 
satisfy for all 
$ \Va \in \Z $, 
$ n \in \N $, 
$ k \in \{ 1, 2, \dots, \fM \} $ 
that 
\begin{equation}
\label{eq:definition_of_phi}
  \phi_0^{ \Va } 
  = 
  \sum_{ w = 0 }^{ \infty }
  \left(
    \frac{ 
      \param \gamma ( 1 - \param \gamma )^w
    }{ M }
    \left[
      \sum_{ m = V M + 1 }^{ V M + M }
      X_{ -w, m }
    \right]
  \right)
,
\qquad 
  \varphi_n^{ \Va }
  =
  \phi_{ n - 1 }^{ \Va }
,
\qquad 
  \psi_n^{ \Va, k }
  =
  X_{ n, - k }
,
\end{equation}
\begin{equation}
\label{eq:definition_of_phi_2}
  \text{and}
\qquad 
  \phi_n^{ \Va }
  =
  ( 1 - \param \gamma )
  \phi_{ n - 1 }^{ \Va }
  +
  \frac{ \param \gamma }{ M } 
  \Biggl[ 
    \sum\limits_{ 
      m = M \max\{ \Va - n , 0 \} + 1 
    }^{ 
      M \max\{ \Va - n, 0 \} + M
    } 
    X_{ n, m }
  \Biggr] 
  ,
\end{equation}
and let 
$
  \Xi \colon \Omega \to \R
$
satisfy 
\begin{equation}
\label{eq:definition_of_Xi}
  \Xi
  =
  4 \param \fM 
  \biggl[
    \max\biggl\{ 
      1 ,
      \| \Theta_0 \|
      +
      \biggl[ 
        \displaystyle 
        \frac{ 4 }{ 2 - \param \gamma }
        \textstyle 
      \biggr]
      \biggl[ 
        \sup\limits_{ v, w \in \Z } \| X_{ v, w } \|
      \biggr]
    \biggr\}
  \biggr]^2
  .
\end{equation}
\Nobs that 
\cref{eq:definition_of_phi}, 
\cref{eq:definition_of_phi_2}, 
and the fact that 
$
  Z_{ n, m } \colon \Omega \to \R^{ \fd }
$, 
$
  (n,m) \in \Z^2
$, 
are i.i.d.\ random variables prove that 
for all $ K \in \N \backslash \{ 1 \} $, 
$ \Va_1, \Va_2, \dots, \Va_{ K + 1 } \in \N_0 $ 
with 
$
  \Va_1 < \Va_2 < \ldots < \Va_{ K + 1 }
$
it holds that 
\begin{equation}
  ( \phi^{ \Va_1 }_k )_{ 
    k \in \N_0 \cap [0, \Va_2) 
  }
  ,
  \quad 
  ( \phi^{ \Va_2 }_k )_{ 
    k \in \N_0 \cap [0, \Va_3) 
  }
  ,
  \quad 
  \dots 
  ,
  \quad 
  ( 
    \phi^{ \Va_{K-1} }_k 
  )_{ 
    k \in \N_0 \cap [0, \Va_{ K } )
  }
  ,
\qandq 
  ( 
    \phi^{ \Va_K }_k 
  )_{ 
    k \in \N_0 \cap [0, \Va_{ K + 1 } )
  }
\end{equation}
are independent. 
\Hence that
for all $ K \in \N \backslash \{ 1 \} $, 
$ \Va_1, \Va_2, \dots, \Va_{ K + 1 } \in \N $ 
with 
$
  \Va_1 < \Va_2 < \ldots < \Va_{ K + 1 }
$
it holds that 
\begin{equation}
  ( 
    \phi^{ V_1 }_{ V_2 - 1 } ,
    \phi^{ V_1 }_{ V_2 - 2 }
  )
  , 
  \qquad 
  (
    \phi^{ V_2 }_{ V_3 - 1 } ,
    \phi^{ V_2 }_{ V_3 - 2 } 
  )
  , 
  \qquad 
  \dots 
  ,
\qqandqq 
  (
    \phi^{ V_K }_{ V_{ K + 1 } - 1 } ,
    \phi^{ V_K }_{ V_{ K + 1 } - 2 }
  )
\end{equation}
are independent. 
Combining this with \cref{eq:definition_of_phi} shows that 
for all $ K \in \N \backslash \{ 1 \} $, 
$ \Va_1, \Va_2, \dots, \Va_{ K + 1 } \in \N $ 
with 
$
  \Va_1 < \Va_2 < \ldots < \Va_{ K + 1 }
$
it holds that 
\begin{equation}
  ( 
    \phi^{ V_1 }_{ V_2 - 1 } ,
    \varphi^{ V_1 }_{ V_2 - 1 }
  )
  , 
  \qquad 
  (
    \phi^{ V_2 }_{ V_3 - 1 } ,
    \varphi^{ V_2 }_{ V_3 - 1 } 
  )
  , 
  \qquad 
  \dots 
  ,
\qqandqq 
  (
    \phi^{ V_K }_{ V_{ K + 1 } - 1 } ,
    \varphi^{ V_K }_{ V_{ K + 1 } - 1 }
  )
\end{equation}
are independent. 
This, \cref{eq:definition_of_phi}, 
and \cref{eq:definition_of_phi_2} 
ensure that 
for all $ K \in \N \backslash \{ 1 \} $, 
$ \Va_1, \Va_2, \dots, \Va_{ K + 1 } \in \N $ 
with 
$
  \Va_1 < \Va_2 < \ldots < \Va_{ K + 1 }
$
it holds that 
$
  ( 
    \phi^{ V_1 }_{ V_2 - 1 } ,
\allowbreak 
    \varphi^{ V_1 }_{ V_2 - 1 } ,
\allowbreak 
    \psi^{ V_1, 1 }_{ V_2 - 1 } , 
\allowbreak 
    \dots ,
\allowbreak 
    \psi^{ V_1, \fM }_{ V_2 - 1 } 
  )
$, 
$
  (
    \phi^{ V_2 }_{ V_3 - 1 } ,
\allowbreak 
    \varphi^{ V_2 }_{ V_3 - 1 } ,
\allowbreak 
    \psi^{ V_2, 1 }_{ V_3 - 1 } , 
\allowbreak 
    \dots , 
\allowbreak 
    \psi^{ V_2, \fM }_{ V_3 - 1 } 
  )
$, 
$
  \dots 
$, 
$
  (
    \phi^{ V_K }_{ V_{ K + 1 } - 1 } ,
\allowbreak 
    \varphi^{ V_K }_{ V_{ K + 1 } - 1 } ,
\allowbreak 
    \psi^{ V_K, 1 }_{ V_{ K + 1 } - 1 } ,
\allowbreak 
    \dots , 
\allowbreak 
    \psi^{ V_K, \fM }_{ V_{ K + 1 } - 1 }
  )
$
are independent. 
The fact that for all 
$ k \in \{ 1, 2, \dots, K \} $
it holds that 
\begin{equation}
  \Phi^{ V_k }_{ V_{ k + 1 } - 1 }
  =
  ( 
    \phi^{ V_k }_{ V_k - 1 } ,
    \varphi^{ V_k }_{ V_k - 1 } ,
    \psi^{ V_k, 1 }_{ V_k - 1 } , 
    \dots ,
    \psi^{ V_k, \fM }_{ V_k - 1 } 
  )
\end{equation}
\hence shows that 
for all $ K \in \N \backslash \{ 1 \} $, 
$ 
  \Va_1, \Va_2, \dots, 
\allowbreak 
  \Va_{ K + 1 } 
\allowbreak 
  \in \N 
$ 
with 
$
  \Va_1 < \Va_2 < \ldots < \Va_{ K + 1 }
$
it holds that 
\begin{equation}
\label{eq:fundamental_indepence}
  \Phi^{ V_1 }_{ V_2 - 1 }
  , 
  \qquad 
  \Phi^{ V_2 }_{ V_3 - 1 }
  , 
  \qquad 
  \dots 
  ,
\qqandqq 
  \Phi^{ V_K }_{ V_{ K + 1 } - 1 }
\end{equation}
are independent. 
\Moreover 
\cref{eq:definition_of_phi}, \cref{eq:definition_of_phi_2}, 
and the fact that 
$ X_{ n, m } \colon \Omega \to \R^{ \fd } $, $ (n, m) \in \Z^2 $, are i.i.d.\ random variables 
show that for all 
$ \Va, \Vb \in \Z $
it holds that 
\begin{equation}
\label{eq:stationary_process}
  \bigl( 
    ( \Phi^{ \Va }_n )_{ n \in \N } 
  \bigr)( \P )
  =
  \bigl( 
    ( \Phi^{ \Vb }_n )_{ n \in \N } 
  \bigr)( \P )
  .
\end{equation}
\Moreover 
\cref{lem:invariant_2}, 
\cref{eq:definition_of_phi}, \cref{eq:definition_of_phi_2}, 
and the fact that 
$ X_{ n, m } \colon \Omega \to \R^{ \fd } $, 
$ ( n, m ) \in \Z^2 $, 
are i.i.d.\ random variables 
show that 
for all 
$ \Va \in \Z $, $ n, m \in \N_0 $
it holds that 
\begin{equation}
  ( \phi^V_n, \phi^V_{ n + 1 } )( \P )
  =
  ( \phi^V_m, \phi^V_{ m + 1 } )( \P )
  .
\end{equation}
\Hence for all 
$ \Va \in \Z $, $ n, m \in \N $
that 
\begin{equation}
  ( \phi^V_n, \varphi^V_n )( \P )
  =
  ( \phi^V_n, \phi^V_{ n - 1 } )( \P )
  =
  ( \phi^V_m, \phi^V_{ m - 1 } )( \P )
  =
  ( \phi^V_m, \varphi^V_m )( \P )
  .
\end{equation}
Combining this, 
\cref{eq:definition_of_phi}, and \cref{eq:definition_of_phi_2} 
with the fact that 
$ X_{ n, m } \colon \Omega \to \R^{ \fd } $, 
$ ( n, m ) \in \Z^2 $, 
are i.i.d.\ random variables 
shows that 
for all 
$ \Va, n, m \in \N $
it holds that 
\begin{equation}
\label{eq:stationary_process_2}
\begin{split}
&
  \bigl( 
    \Phi^{ \Va }_n 
  \bigr)( \P )
  =
  ( 
    \phi^V_n, \varphi^V_n, \psi^{ V, 1 }_n, \dots, \psi^{ V, \fM }_n  
  )( \P )
  =
  \bigl(
    ( 
      \phi^V_n, \varphi^V_n
    )( \P )
  \bigr)
  \otimes 
  \bigl(
    ( 
      \psi^{ V, 1 }_n, \dots, \psi^{ V, \fM }_n  
    )( \P )
  \bigr)
\\ &
  =
  \bigl(
    ( 
      \phi^V_m, \varphi^V_m
    )( \P )
  \bigr)
  \otimes 
  \bigl(
    ( 
      \psi^{ V, 1 }_m, \dots, \psi^{ V, \fM }_m
    )( \P )
  \bigr)
  =
  ( 
    \phi^V_m, \varphi^V_m, \psi^{ V, 1 }_m, \dots, \psi^{ V, \fM }_m  
  )( \P )
  =
  \bigl( 
    \Phi^{ \Va }_m 
  \bigr)( \P )
  .
\end{split}
\end{equation}
\Moreover \cref{eq:definition_of_fl,eq:example_definition_GD} 
ensure for all $ n \in \N $ that
\begin{equation}
\label{eq:recursion_of_Theta_n}
  \Theta_n 
  = 
  ( 1 - \param \gamma )
  \Theta_{ n - 1 } 
  +
  \frac{ \param \gamma }{ M } 
  \left[ 
    \sum_{ m = 1 }^M
    X_{ n, m }  
  \right]
  .
\end{equation}
Combining this with 
\cref{eq:definition_of_phi_2}
shows for all $ \Va \in \Z $, $ n \in \N $ 
with $ n \geq \Va $ that 
\begin{equation}
\begin{split}
\textstyle 
  \| \Theta_n - \phi^{ \Va }_n \|
&
\textstyle 
=
  \left\|
    \left[ 
      ( 1 - \param \gamma )
      \Theta_{ n - 1 } 
      +
      \frac{ \param \gamma }{ M } 
      \left[ 
        \sum\limits_{ m = 1 }^M
        X_{ n, m }  
      \right]
    \right] 
    -
    \left[ 
      ( 1 - \param \gamma )
      \phi^V_{ n - 1 } 
      +
      \frac{ \param \gamma }{ M } 
      \left[ 
        \sum\limits_{ m = 1 }^M
        X_{ n, m }  
      \right]
    \right] 
  \right\|
\\ &
=
  \|
    ( 1 - \param \gamma )
    \Theta_{ n - 1 } 
    -
    ( 1 - \param \gamma )
    \phi^V_{ n - 1 } 
  \|
=
  | 1 - \param \gamma |
  \,
  \| \Theta_{ n - 1 } - \phi_{ n - 1 }^{ \Va } \|
  .
\end{split}
\end{equation}
\Hence for all $ \Va \in \Z $, $ n, m \in \N_0 $ 
with $ n \geq m \geq \Va - 1 $ that 
\begin{equation}
\label{eq:error_estimate_theta_minus_phi}
\begin{split}
&
\textstyle 
  \| \Theta_n - \phi_n^{ \Va } \|
=
  | 1 - \param \gamma |^{ n - m }
  \| \Theta_m - \phi_m^{ \Va } \|
  .
\end{split}
\end{equation}
\Moreover 
\cref{eq:definition_of_fl}
assures 
for all $ \theta, x_1, x_2, \dots, x_{ \fM }, \fe \in \R^{ \fd } $
with 
$
  \fe = \frac{ 1 }{ \fM } [ \sum_{ m = 1 }^{ \fM } x_m ]
$
that 
\begin{equation}
\begin{split}
&
  \frac{ 1 }{ \fM }
  \left[ 
    \sum_{ m = 1 }^{ \fM }
    \fl( \theta, x_m )
  \right] 
  =
  \frac{ 1 }{ \fM }
  \left[ 
    \sum_{ m = 1 }^{ \fM }
    \fl( \theta - \fe , x_m - \fe )
  \right] 
\\ &
  =
  \frac{ 1 }{ \fM }
  \left[ 
    \sum_{ m = 1 }^{ \fM }
    \bigl[
      \fl(
        \theta, \fe 
      )
      +
      \tfrac{ \param }{ 2 } 
      \langle \theta - \fe , x_m - \fe \rangle 
      +
      \fl(
        x_m, \fe 
      )
    \bigr]
  \right] 
=
  \fl(
    \theta, \fe 
  )
  +
  \frac{ 1 }{ \fM }
  \left[ 
    \sum_{ m = 1 }^{ \fM }
    \fl(
      x_m, \fe 
    )
  \right] 
\\ &
  =
  \fl\biggl(
    \theta 
    , 
    \frac{ 1 }{ \fM }
  \textstyle 
    \biggl[
      \sum\limits_{ m = 1 }^{ \fM }
      x_m 
    \biggr]
  \biggr)
  \displaystyle 
  +
  \frac{ 1 }{ \fM }
  \left[ 
    \sum_{ m = 1 }^{ \fM }
    \fl\biggl( 
      x_m ,
      \frac{ 1 }{ \fM }
      \textstyle 
      \biggl[
        \sum\limits_{ n = 1 }^{ \fM }
        x_n
      \biggr]
    \biggr)
  \right] 
  .
\end{split}
\end{equation}
\Hence 
for all $ \theta, \vartheta, x_1, x_2, \dots, x_{ \fM } \in \R^{ \fd } $
that 
\begin{equation}
\begin{split}
&
  \frac{ 1 }{ \fM }
  \left[ 
    \sum_{ m = 1 }^{ \fM }
    \fl( \theta, x_m )
  \right] 
  -
  \frac{ 1 }{ \fM }
  \left[ 
    \sum_{ m = 1 }^{ \fM }
    \fl( \vartheta, x_m )
  \right] 
  =
  \fl\biggl(
    \theta 
    , 
    \frac{ 1 }{ \fM }
  \textstyle 
    \biggl[
      \sum\limits_{ m = 1 }^{ \fM }
      x_m 
    \biggr]
  \biggr)
  \displaystyle 
  -
  \fl\biggl(
    \vartheta 
    , 
    \frac{ 1 }{ \fM }
  \textstyle 
    \biggl[
      \sum\limits_{ m = 1 }^{ \fM }
      x_m 
    \biggr]
  \biggr)
  \displaystyle 
  .
\end{split}
\end{equation}
Combining this,  
\cref{eq:definition_of_f}, 
\cref{eq:definition_of_phi}, 
and 
\cref{eq:definition_of_phi_2} 
ensures that
\begin{equation}
\begin{split}
&  
  \lim_{ \varepsilon \searrow 0 }
  \P\bigl(
    f( \Phi^0_1 )
    > \varepsilon
  \bigr)
=
  \P\bigl(
  \cup_{ \varepsilon \in (0,\infty) }
  \bigl\{ 
    f( \Phi^0_1 )
    > \varepsilon
  \bigr\}
  \bigr)
=
  \P\bigl(
    f( \Phi^0_1 )
    > 0 
  \bigr)
\\
& 
\textstyle 
  =
  \P\biggl(
    \biggl[
      \sum\limits_{ m = 1 }^{ \fM }
      \fl\bigl( 
        \phi_1^0
        , 
        \psi^{ 0, (m) }_1
      \bigr) 
    \biggr]
    - 
    \biggl[
      \sum\limits_{ m = 1 }^{ \fM }
      \fl\bigl( 
        \varphi_1^0, 
        \psi^{ 0, (m) }_1
      \bigr)
    \biggr]
    > 0
  \biggr)
\\
& 
\textstyle 
  =
  \P\biggl(
    \biggl[
      \sum\limits_{ m = 1 }^{ \fM }
      \fl\biggl( 
        ( 1 - \param \gamma ) 
        \phi_0^0
        + 
        \displaystyle 
        \frac{ \param \gamma }{ M } 
        \textstyle 
        \biggl[ 
          \sum\limits_{ n = 1 }^{ M } X_{ 1, n }
        \biggr]
        , 
        X_{ 1, - m }
      \biggr) 
    \biggr]
    - 
    \biggl[
      \sum\limits_{ m = 1 }^{ \fM }
      \fl\bigl( 
        \phi_0^0, 
        X_{ 1, - m }
      \bigr)
    \biggr]
    > 0
  \biggr)
\\
& 
\textstyle 
  =
  \P\biggl(
    \displaystyle
    \frac{ 1 }{ \fM }
    \textstyle 
    \biggl[
      \sum\limits_{ m = 1 }^{ \fM }
      \fl\biggl( 
        ( 1 - \param \gamma ) 
        \phi_0^0
        + 
        \displaystyle 
        \frac{ \param \gamma }{ M } 
        \textstyle 
        \biggl[ 
          \sum\limits_{ n = 1 }^{ M } X_{ 1, n }
        \biggr]
        , 
        X_{ 1, -m }
      \biggr) 
    \biggr]
    - 
    \displaystyle
    \frac{ 1 }{ \fM }
    \textstyle 
    \biggl[
      \sum\limits_{ m = 1 }^{ \fM }
      \fl\bigl( 
        \phi_0^0,
        X_{ 1, - m }
      \bigr)
    \biggr]
    > 0
  \biggr)
\\
& =
  \P\biggl(
    \fl\biggl( 
      ( 1 - \param \gamma ) 
      \phi_0^0
      + 
      \displaystyle 
      \frac{ \param \gamma }{ M } 
      \textstyle 
      \biggl[ 
        \sum\limits_{ n = 1 }^{ M } X_{ 1, n }
      \biggr]
      , 
      \displaystyle 
      \frac{ 1 }{ \fM } 
      \textstyle 
      \biggl[ 
        \sum\limits_{ m = 1 }^{ \fM } X_{ 1, - m }
      \biggr]
    \biggr) 
    - 
    \fl\biggl( 
      \phi_0^0, 
      \displaystyle 
      \frac{ 1 }{ \fM } 
      \textstyle 
      \biggl[ 
        \sum\limits_{ m = 1 }^{ \fM } X_{ 1, - m }
      \biggr]      
    \biggr)
    > 0
  \biggr)
  .
\end{split}
\end{equation}
The fact that 
$ X_{ n, m } \colon \Omega \to \R^{ \fd } $, $ (n, m) \in \Z^2 $, 
are i.i.d.\ random variables,
\cref{eq:definition_of_fl}, 
\cref{eq:definition_of_chi}, 
\cref{eq:assumption_invariant_measure}, 
and 
\cref{eq:definition_of_phi} 
\hence show that 
\begin{equation}
\begin{split}
&  
  \lim_{ \varepsilon \searrow 0 }
  \P\bigl(
    f( \Phi^0_1 )
    > \varepsilon
  \bigr)
\\
& =
  \P\biggl(
    \fl\biggl( 
      ( 1 - \param \gamma ) 
      \phi_0^0
      + 
      \displaystyle 
      \frac{ \param \gamma }{ M } 
      \textstyle 
      \biggl[ 
        \sum\limits_{ n = 1 }^{ M } X_{ 1, n }
      \biggr]
      , 
      \displaystyle 
      \frac{ 1 }{ \fM } 
      \textstyle 
      \biggl[ 
        \sum\limits_{ m = 1 }^{ \fM } X_{ 2, m }
      \biggr]
    \biggr) 
    - 
    \fl\biggl( 
      \phi_0^0, 
      \displaystyle 
      \frac{ 1 }{ \fM } 
      \textstyle 
      \biggl[ 
        \sum\limits_{ m = 1 }^{ \fM } X_{ 2, m }
      \biggr]      
    \biggr)
    > 0
  \biggr)
\\
& =
  \P\biggl(
  \,
    \biggl\| 
      ( 1 - \param \gamma ) 
      \chi
      + 
      \displaystyle 
      \frac{ \param \gamma }{ M } 
      \textstyle 
      \biggl[ 
        \sum\limits_{ m = 1 }^{ M } X_{ 1, m }
      \biggr]
      - 
      \displaystyle 
      \frac{ 1 }{ \fM } 
      \textstyle 
      \biggl[ 
        \sum\limits_{ m = 1 }^{ \fM } X_{ 2, m }
      \biggr]
    \biggr\|^2
    -
    \biggl\| 
      \chi
      - 
      \displaystyle 
      \frac{ 1 }{ \fM } 
      \textstyle 
      \biggl[ 
        \sum\limits_{ m = 1 }^{ \fM } X_{ 2, m }
      \biggr]
    \biggr\|^2
    > 0
    \,
  \biggr) 
\\
& =
  \P\biggl(
  \,
    \biggl\| 
      ( 1 - \param \gamma ) 
      \chi
      + 
      \displaystyle 
      \frac{ \param \gamma }{ M } 
      \textstyle 
      \biggl[ 
        \sum\limits_{ m = 1 }^{ M } X_{ 1, m }
      \biggr]
      - 
      \displaystyle 
      \frac{ 1 }{ \fM } 
      \textstyle 
      \biggl[ 
        \sum\limits_{ m = 1 }^{ \fM } X_{ 2, m }
      \biggr]
    \biggr\|
    >
    \biggl\| 
      \chi
      - 
      \displaystyle 
      \frac{ 1 }{ \fM } 
      \textstyle 
      \biggl[ 
        \sum\limits_{ m = 1 }^{ \fM } X_{ 2, m }
      \biggr]
    \biggr\|
    \,
  \biggr) 
  > 0
  .
\end{split}
\end{equation}
\Hence that 
there exists $ \varepsilon \in (0,\infty) $ 
which satisfies
\begin{equation}
\label{eq:probability_varepsilon_strictly_bigger_than_zero}
  \P\bigl(
    f( \Phi_1^0 ) > \varepsilon
  \bigr)
  > 0
  .
\end{equation}
\Nobs that 
\cref{eq:stationary_process}, 
\cref{eq:stationary_process_2},  
and 
\cref{eq:probability_varepsilon_strictly_bigger_than_zero}
prove that for all $ \Va, n \in \N $ it holds that 
\begin{equation}
\label{eq:probability_smaller_than_1}
  \P\bigl(
    f( \Phi_n^{ \Va } ) \leq \varepsilon
  \bigr)
=
  \P\bigl(
    f( \Phi_1^{ \Va } ) \leq \varepsilon
  \bigr)
=
  \P\bigl(
    f( \Phi_1^0 ) \leq \varepsilon
  \bigr)
=
  1 - 
  \P\bigl(
    f( \Phi_1^0 ) > \varepsilon
  \bigr)
< 
  1 .
\end{equation}
\Moreover 
the triangle inequality ensures 
for all 
$ k \in \N $, 
$ 
  v = ( v_1, \dots, v_k ) 
$, 
$
  w = ( w_1, \dots, w_k ) 
  \in \R^k 
$
that
\begin{equation}
  \| v \| - \| w \|
  =
  \| v - w + w \| - \| w \|
  \leq 
  \| v - w \| + \| w \| - \| w \|
  \leq 
  \| v - w \|
  .
\end{equation}
\Hence 
for all 
$ k \in \N $, 
$ 
  v = ( v_1, \dots, v_k ) 
$, 
$
  w = ( w_1, \dots, w_k ) 
  \in \R^k 
$
that
\begin{equation}
  \bigl| 
    \| v \| - \| w \|
  \bigr|
  \leq 
  \| v - w \|
  .
\end{equation}
\Hence 
for all 
$ \theta_1, \theta_2, x_1, x_2, \dots, x_m \in \R^{ \fd } $ 
that 
\begin{equation}
\begin{split}
&
\textstyle 
  \biggl|
    \biggl[ 
      \sum\limits_{ m = 1 }^{ \fM } 
      \| \theta_1 - x_m \|^2 
    \biggr]
    - 
    \biggl[ 
      \sum\limits_{ m = 1 }^{ \fM } 
      \| \theta_2 - x_m \|^2
    \biggr]
  \biggr|
= 
  \biggl|
    \sum\limits_{ m = 1 }^{ \fM } 
    \bigl( 
      \| \theta_1 - x_m \|^2
      -
      \| \theta_2 - x_m \|^2
    \bigr) 
  \biggr|
\\ &
\textstyle 
= 
  \biggl|
    \sum\limits_{ m = 1 }^{ \fM } 
    \bigl(
      \| \theta_1 - x_m \|
      -
      \| \theta_2 - x_m \|
    \bigr)
    \bigl(
      \| \theta_1 - x_m \|
      +
      \| \theta_2 - x_m \|
    \bigr)
  \biggr|
\\ &
\textstyle 
\leq 
    \sum\limits_{ m = 1 }^{ \fM } 
    \| \theta_1 - \theta_2 \|
    \bigl(
      \| \theta_1 - x_m \|
      +
      \| \theta_2 - x_m \|
    \bigr)
\leq 
  \| \theta_1 - \theta_2 \|
  \biggl[
    \sum\limits_{ m = 1 }^{ \fM } 
    \bigl(
      \| \theta_1 - x_m \|
      +
      \| \theta_2 - x_m \|
    \bigr)
  \biggr]
\\ & 
\textstyle 
\leq 
  \| \theta_1 - \theta_2 \|
  \biggl[
    \sum\limits_{ m = 1 }^{ \fM } 
    \bigl(
      \| \theta_1 \|
      +
      \| \theta_2 \|
      +
      2
      \| x_m \|
    \bigr)
  \biggr]
\leq 
  \| \theta_1 - \theta_2 \|
  \biggl[
    \fM
    \| \theta_1 \|
    +
    \fM
    \| \theta_2 \|
    +
    2
    \biggl(
      \sum\limits_{ m = 1 }^{ \fM } 
      \| x_m \|
    \biggr)
  \biggr]
  .
\end{split}
\end{equation}
\Hence 
for all 
$ \theta_1, \theta_2, \vartheta_1, \vartheta_2, x_1, x_2, \dots, x_m \in \R^{ \fd } $ 
that 
\begin{equation}
\begin{split}
&
\textstyle 
  \biggl|
  \biggl(
    \biggl[ 
      \sum\limits_{ m = 1 }^{ \fM } 
      \| \theta_1 - x_m \|^2 
    \biggr]
    - 
    \biggl[ 
      \sum\limits_{ m = 1 }^{ \fM } 
      \| \vartheta_1 - x_m \|^2
    \biggr]
  \biggr)
\\ &
\textstyle 
  -
  \biggl(
    \biggl[ 
      \sum\limits_{ m = 1 }^{ \fM } 
      \| \theta_2 - x_m \|^2 
    \biggr]
    - 
    \biggl[ 
      \sum\limits_{ m = 1 }^{ \fM } 
      \| \vartheta_2 - x_m \|^2
    \biggr]
  \biggr)
  \biggr|
\\ & 
\textstyle 
\leq 
  \| \theta_1 - \theta_2 \|
  \biggl[
    \fM
    \| \theta_1 \|
    +
    \fM
    \| \theta_2 \|
    +
    2
    \biggl(
      \sum\limits_{ m = 1 }^{ \fM } 
      \| x_m \|
    \biggr)
  \biggr]
\\ &
\textstyle 
  +
  \| \vartheta_1 - \vartheta_2 \|
  \biggl[
    \fM
    \| \vartheta_1 \|
    +
    \fM
    \| \vartheta_2 \|
    +
    2
    \biggl(
      \sum\limits_{ m = 1 }^{ \fM } 
      \| x_m \|
    \biggr)
  \biggr]
\\ &
\textstyle 
\leq 
  \bigl(
    \| \theta_1 - \theta_2 \|
    +
    \| \vartheta_1 - \vartheta_2 \|
  \bigr)
  \biggl[
    \fM
    \bigl(
      \max\{ \| \theta_1 \|, \| \vartheta_1 \| \}
      +
      \max\{ \| \theta_2 \|, \| \vartheta_2 \| \}
    \bigr)
    +
    2
    \biggl(
      \sum\limits_{ m = 1 }^{ \fM } 
      \| x_m \|
    \biggr)
  \biggr]
  .
\end{split}
\end{equation}
Combining this with \cref{eq:definition_of_fl} and \cref{eq:definition_of_f} shows
for all 
$ \theta_1, \theta_2, \vartheta_1, \vartheta_2, x_1, x_2, \dots, x_m \in \R^{ \fd } $ 
that 
\begin{equation}
\begin{split}
&
  \bigl|
    f( \theta_1, \vartheta_1, x_1, x_2, \dots, x_{ \fM } )
    -
    f( \theta_2, \vartheta_2, x_1, x_2, \dots, x_{ \fM } )
  \bigr|
\\ &
\textstyle 
\leq
  \tfrac{ \param }{ 2 }
  \bigl(
    \| \theta_1 - \theta_2 \|
    +
    \| \vartheta_1 - \vartheta_2 \|
  \bigr)
  \biggl[
    \fM
    \biggl(
      \sum\limits_{ i = 1 }^2
      \max\{ \| \theta_i \|, \| \vartheta_i \| \}
    \biggr)
    +
    2
    \biggl(
      \sum\limits_{ m = 1 }^{ \fM } 
      \| x_m \|
    \biggr)
  \biggr]
  .
\end{split}
\end{equation}
\Hence 
for all $ \Va \in \Z $, $ n \in \N $ that 
\begin{equation}
\begin{split}
&
\textstyle 
  \bigl|
    f( \varTheta_n )
    -
    f( \Phi^{ \Va }_n )
  \bigr|
\leq
  \bigl(
    \| \Theta_n - \phi^{ \Va }_n \|
    +
    \| \Theta_{ n - 1 } - \phi^{ \Va }_{ n - 1 } \|
  \bigr)
\\ &
\textstyle
  \cdot 
  \biggl[
    \frac{ \param \fM }{ 2 } 
    \bigl(
      \max\{ 
        \| \Theta_n \|, \| \Theta_{ n - 1 } \|
      \}
      +
      \max\{ 
        \| \phi^{ \Va }_n \|, \| \phi^{ \Va }_{ n - 1 } \|
      \}
    \bigr)
    +
    \param
    \biggl(
      \sum\limits_{ m = 1 }^{ \fM } 
      \| X_{ n, - m } \|
    \biggr)
  \biggr]
  .
\end{split}
\end{equation}
Combining this and \cref{eq:error_estimate_theta_minus_phi} 
with the fact that $ | 1 - \param \gamma | \leq 1 $ 
shows that for all $ n, \cN \in \N_0 $ 
with $ \cN < n $ it holds that 
\begin{equation}
\label{eq:f_varTheta_minus_f_Phi_estimate}
\begin{split}
&
\textstyle 
  \bigl|
    f( \varTheta_n )
    -
    f( \Phi^{ n - \cN }_n )
  \bigr|
\\ &
\textstyle 
\leq 
  \bigl(
    | 1 - \param \gamma |^{ \cN + 1 }
    +
    | 1 - \param \gamma |^{ \cN }
  \bigr)
  \|
    \Theta_{ n - \cN - 1}
    -
    \phi^{ n - \cN }_{ n - \cN - 1 }
  \|
\\ &
\textstyle 
  \cdot 
  \biggl[
    \frac{ \param \fM }{ 2 } 
    \bigl(
      \max\{ 
        \| \Theta_n \|, \| \Theta_{ n - 1 } \|
      \}
      +
      \max\{ 
        \| \phi_n^{ n - \cN } \|, \| \phi_{ n - 1 }^{ n - \cN } \|
      \} 
    \bigr)
    +
    \param
    \biggl(
      \sum\limits_{ m = 1 }^{ \fM } 
      \| X_{ n, - m } \|
    \biggr)
  \biggr]
\\ &
\textstyle 
\leq 
  | 1 - \param \gamma |^{ \cN }
  \|
    \Theta_{ n - \cN - 1}
    -
    \phi^{ n - \cN }_{ n - \cN - 1 }
  \|
\\ &
\textstyle 
  \cdot 
  \biggl[
    \param \fM 
    \bigl(
      \max\{ 
        \| \Theta_n \|, \| \Theta_{ n - 1 } \|
      \}
      +
      \max\{ 
        \| \phi_n^{ n - \cN } \|, \| \phi_{ n - 1 }^{ n - \cN } \|
      \} 
    \bigr)
    +
    2 \param
    \biggl(
      \sum\limits_{ m = 1 }^{ \fM } 
      \| X_{ n, - m } \|
    \biggr)
  \biggr]
  .
\end{split}
\end{equation}
\Moreover \cref{item:iii_a_priori_bounds} in \cref{lem:a_priori} implies 
that for all $ n \in \N_0 $ it holds that 
\begin{equation}
\label{eq:Theta_n_xi_a_priori_bound}
  \| \Theta_n \|
  \leq 
  \| \Theta_0 \|
  +
  \left[ 
    \frac{ 2 }{ 2 - \param \gamma }
  \right] 
  \left[ 
    \sup_{ v, w \in \N_0 } \| X_{ v, w } \|
  \right]
\leq 
  \| \Theta_0 \|
  +
  \left[ 
    \frac{ 2 }{ 2 - \param \gamma }
  \right] 
  \left[ 
    \sup_{ v, w \in \Z } \| X_{ v, w } \|
  \right]
  .
\end{equation}
\Moreover \cref{eq:definition_of_phi_2} 
and the fact that for all 
$ \theta, x \in \R^{ \fd } $
it holds that 
$
  ( \nabla_{ \theta } \fl )( \theta, x ) 
  =
  \param ( \theta - x ) 
$
prove that for all $ \Va \in \Z $, $ n \in \N $ it holds that 
\begin{equation}
\begin{split}
  \phi^{ \Va }_n 
& 
  =
  \phi^{ \Va }_{ n - 1 }
  -
  \frac{ \param \gamma }{ M }
  \left[ 
    \sum_{ 
      m = M \max\{ \Va - n, 0 \} + 1 
    }^{
      M \max\{ \Va - n, 0 \} + M
    }
    \bigl( 
      \phi^{ \Va }_{ n - 1 }
      -
      X_{ n, m } 
    \bigr)
  \right] 
\\ &
  =
  \phi^{ \Va }_{ n - 1 }
  -
  \frac{ \gamma }{ M }
  \left[ 
    \sum_{ 
      m = M \max\{ \Va - n, 0 \} + 1 
    }^{
      M \max\{ \Va - n, 0 \} + M
    }
    ( \nabla_{ \theta } \fl )( \phi^{ \Va }_{ n - 1 }, X_{ n, m } )
  \right] 
\\ &
  =
  \phi^{ \Va }_{ n - 1 }
  -
  \frac{ \gamma }{ M }
  \left[ 
    \sum_{ 
      m = 1 
    }^{
      M
    }
    ( 
      \nabla_{ \theta } \fl 
    )( \phi^{ \Va }_{ n - 1 }, X_{ n, m + M \max\{ \Va - n, 0 \} } )
  \right] 
  .
\end{split}
\end{equation}
\cref{lem:invariant} and \cref{eq:definition_of_phi} \hence show 
for all $ \Va \in \Z $ that
\begin{equation}
  \sup_{ n \in \N_0 }
  \| \phi^{ \Va }_n \|
  \leq 
  \left[ 
    \frac{ 2 }{ 2 - \param \gamma }
  \right]
  \left[ 
    \sup_{ v, w \in \Z } \| X_{ v, w } \|
  \right]
  .
\end{equation}
Combining this and \cref{eq:Theta_n_xi_a_priori_bound} 
with \cref{eq:f_varTheta_minus_f_Phi_estimate}
demonstrates that
for all $ n, \cN \in \N_0 $ 
with $ \cN < n $ it holds that 
\begin{equation}
\begin{split}
&
\textstyle 
  \bigl|
    f( \varTheta_n )
    -
    f( \Phi^{ n - \cN }_n )
  \bigr|
\\ &
\textstyle 
\leq 
  | 1 - \param \gamma |^{ \cN }
  \|
    \Theta_{ n - \cN - 1}
    -
    \phi^{ n - \cN }_{ n - \cN - 1 }
  \|
\\ &
\textstyle 
  \cdot 
  \biggl[
    \param \fM 
    \biggl(
      \| \Theta_0 \|
      +
      \biggl[ 
        \displaystyle 
        \frac{ 4 }{ 2 - \param \gamma }
        \textstyle 
      \biggr]
      \biggl[ 
        \sup\limits_{ v, w \in \Z } \| X_{ v, w } \|
      \biggr]
    \biggr)
    +
    2 \param
    \biggl(
      \sum\limits_{ m = 1 }^{ \fM } 
      \| X_{ n, - m } \|
    \biggr)
  \biggr]
\\ &
\textstyle 
\leq 
  | 1 - \param \gamma |^{ \cN }
  \biggl(
    \| \Theta_0 \|
    +
    \biggl[ 
      \displaystyle 
      \frac{ 4 }{ 2 - \param \gamma }
      \textstyle 
    \biggr]
    \biggl[ 
      \sup\limits_{ v, w \in \Z } \| X_{ v, w } \|
    \biggr]
  \biggr)
\\ &
\textstyle 
  \cdot 
  \biggl[
    \param \fM 
    \biggl(
      \| \Theta_0 \|
      +
      \biggl[ 
        \displaystyle 
        \frac{ 4 }{ 2 - \param \gamma }
        \textstyle 
      \biggr]
      \biggl[ 
        \sup\limits_{ v, w \in \Z } \| X_{ v, w } \|
      \biggr]
    \biggr)
    +
    2 \param
    \biggl(
      \sum\limits_{ m = 1 }^{ \fM } 
      \| X_{ n, - m } \|
    \biggr)
  \biggr]
  .
\end{split}
\end{equation}
This and \cref{eq:definition_of_Xi} 
show for all $ \cN \in \N_0 $, $ n \in \N \cap ( \cN, \infty ) $ 
that 
\begin{equation}
\label{eq:error_Chi_estimate}
\begin{split}
&
\textstyle 
  \bigl|
    f( \varTheta_n )
    -
    f( \Phi^{ n - \cN }_n )
  \bigr|
\\ &
\textstyle 
\leq 
  | 1 - \param \gamma |^{ \cN }
  \biggl(
    \| \Theta_0 \|
    +
    \biggl[ 
      \displaystyle 
      \frac{ 4 }{ 2 - \param \gamma }
      \textstyle 
    \biggr]
    \biggl[ 
      \sup\limits_{ v, w \in \Z } \| X_{ v, w } \|
    \biggr]
  \biggr)
\\ &
\textstyle 
  \cdot 
  \biggl[
    \param \fM 
    \biggl(
      \| \Theta_0 \|
      +
      \biggl[ 
        \displaystyle 
        \frac{ 4 }{ 2 - \param \gamma }
        \textstyle 
      \biggr]
      \biggl[ 
        \sup\limits_{ v, w \in \Z } \| X_{ v, w } \|
      \biggr]
    \biggr)
    +
    2 \param \fM
    \biggl[
      \sup\limits_{ v, w \in \Z }
      \| X_{ v, w } \|
    \biggr]
  \biggr]
\\ &
\textstyle 
\leq 
  | 1 - \param \gamma |^{ \cN }
  2
  \param \fM 
  \biggl[
    \max\biggl\{ 
      1 ,
      \| \Theta_0 \|
      +
      \biggl[ 
        \displaystyle 
        \frac{ 4 }{ 2 - \param \gamma }
        \textstyle 
      \biggr]
      \biggl[ 
        \sup\limits_{ v, w \in \Z } \| X_{ v, w } \|
      \biggr]
    \biggr\}
  \biggr]^2
=
\displaystyle 
  \frac{ 
    \Xi
    | 1 - \param \gamma |^{ \cN }
  }{ 2 }
  .
\end{split}
\end{equation}
\Moreover the triangle inequality assures that 
for all $ \cN \in \N \cap ( N, \infty ) $ it holds that 
\begin{equation}
\begin{split}
&
\textstyle 
  \P\biggl(
    \bigcup\limits_{ n = N }^{ \infty }
    \bigl\{ 
      f( \varTheta_n ) > \frac{ \varepsilon }{ 2 } 
    \bigr\} 
  \biggr)
=
  \P\biggl(
    \bigcup\limits_{ n = N + 1 }^{ \infty }
    \bigl\{ 
      f( \varTheta_{ n - 1 } ) > \frac{ \varepsilon }{ 2 } 
    \bigr\} 
  \biggr)
\geq 
  \P\biggl(
    \bigcup\limits_{ n = \cN + 1 }^{ \infty }
    \bigl\{ 
      f( \varTheta_{ n - 1 } ) > \frac{ \varepsilon }{ 2 } 
    \bigr\} 
  \biggr)
\\ &
\textstyle 
\geq 
  \P\biggl(
    \bigcup\limits_{ n = \cN + 1 }^{ \infty }
    \Bigl[
      \bigl\{ 
        f( \Phi_{ n - 1 }^{ n - \cN } ) > \varepsilon
      \bigr\} 
      \cap 
      \bigl\{ 
        | 
          f( \varTheta_{ n - 1 } ) - f( \Phi_{ n - 1 }^{ n - \cN } ) 
        | < \frac{ \varepsilon }{ 2 } 
      \bigr\} 
    \Bigr]
  \biggr)
\\ &
\textstyle 
\geq 
  \P\biggl(
    \bigcup\limits_{ n = \cN + 1 }^{ \infty }
    \Bigl[
      \bigl\{ 
        f( \Phi_{ n - 1 }^{ n - \cN } ) > \varepsilon
      \bigr\} 
      \cap 
      \Bigl(
        \cap_{ m = \cN + 1 }^{ \infty }
        \bigl\{ 
          | 
            f( \varTheta_{ m - 1 } ) 
            - 
            f( \Phi_{ m - 1 }^{ m - \cN } ) 
          | < \frac{ \varepsilon }{ 2 } 
        \bigr\} 
      \Bigr)
    \Bigr]
  \biggr)
\\ &
\textstyle 
=
  \P\Bigl(
    \bigl(
      \cup_{ n = \cN + 1 }^{ \infty }
      \bigl\{ 
        f( \Phi_{ n - 1 }^{ n - \cN } ) > \varepsilon
      \bigr\} 
    \bigr)
    \cap 
    \bigl(
      \cap_{ n = \cN + 1 }^{ \infty }
      \bigl\{ 
        | 
          f( \varTheta_{ n - 1 } ) - f( \Phi_{ n - 1 }^{ n - \cN } ) 
        | < \frac{ \varepsilon }{ 2 } 
      \bigr\} 
    \bigr)
  \Bigr)
\\ &
\textstyle 
=
  \P\Bigl(
    \bigl(
      \cup_{ n = \cN + 1 }^{ \infty }
      \bigl\{ 
        f( \Phi_{ n - 1 }^{ n - \cN } ) > \varepsilon
      \bigr\} 
    \bigr)
    \backslash 
    \bigl(
      \cup_{ n = \cN + 1 }^{ \infty }
      \bigl\{ 
        | 
          f( \varTheta_{ n - 1 } ) 
          - 
          f( \Phi_{ n - 1 }^{ n - \cN } ) 
        | 
        \geq \frac{ \varepsilon }{ 2 } 
      \bigr\} 
    \bigr)
  \Bigr)
  .
\end{split}
\end{equation}
The fact that for all $ A, B \in \cF $ it holds that 
$ \P( A \backslash B ) \geq \P( A ) - \P( B ) $ 
and \cref{eq:error_Chi_estimate}
\hence show that 
for all $ \cN \in \N \cap ( N, \infty ) $, 
$ \cM \in \N \cap ( \cN, \infty ) $
it holds that 
\begin{equation}
\begin{split}
&
\textstyle 
  \P\biggl(
    \bigcup\limits_{ n = N }^{ \infty }
    \bigl\{ 
      f( \varTheta_n ) > \frac{ \varepsilon }{ 2 } 
    \bigr\} 
  \biggr)
\\ & 
\textstyle 
\geq 
  \P\biggl(
    \bigcup\limits_{ n = \cN + 1 }^{ \infty }
    \bigl\{ 
      f( \Phi_{ n - 1 }^{ n - \cN } ) > \varepsilon
    \bigr\} 
  \biggr)
  -
  \P\biggl(
    \bigcup\limits_{ n = \cN + 1 }^{ \infty }
    \bigl\{ 
      | 
        f( \varTheta_{ n - 1 } ) 
        - 
        f( \Phi_{ n - 1 }^{ n - \cN } ) 
      | \geq \frac{ \varepsilon }{ 2 } 
    \bigr\} 
  \biggr)
\\ & 
\textstyle 
\geq 
  \P\biggl(
    \bigcup\limits_{ n = \cN + 1 }^{ \infty }
    \bigl\{ 
      f( \Phi_{ n \cN - 1 }^{ n \cN - \cN } ) > \varepsilon
    \bigr\} 
  \biggr)
  -
  \P\biggl(
    \bigcup\limits_{ n = \cN }^{ \infty }
    \bigl\{ 
      | 
        f( \varTheta_n ) 
        - 
        f( \Phi_n^{ n - ( \cN - 1 ) } ) 
      | \geq \frac{ \varepsilon }{ 2 } 
    \bigr\} 
  \biggr)
\\ & 
\textstyle 
\geq 
  \P\biggl(
    \bigcup\limits_{ n = \cN + 1 }^{ \cM }
    \bigl\{ 
      f( \Phi_{ n \cN - 1 }^{ n \cN - \cN } ) > \varepsilon
    \bigr\} 
  \biggr)
  -
  \P\biggl(
    \bigcup\limits_{ n = \cN }^{ \infty }
    \bigl\{  
      2^{ - 1 } \Xi | 1 - \param \gamma |^{ \cN - 1 } 
      \geq \frac{ \varepsilon }{ 2 } 
    \bigr\} 
  \biggr)
  .
\end{split}
\end{equation}
Combining this with \cref{eq:fundamental_indepence}, 
\cref{eq:stationary_process}, and \cref{eq:stationary_process_2}
proves that 
for all $ \cN \in \N \cap ( N, \infty ) $, 
$ \cM \in \N \cap ( \cN, \infty ) $
it holds that 
\begin{equation}
\begin{split}
&
\textstyle 
  \P\biggl(
    \bigcup\limits_{ n = N }^{ \infty }
    \bigl\{ 
      f( \varTheta_n ) > \frac{ \varepsilon }{ 2 } 
    \bigr\} 
  \biggr)
\\ & 
\textstyle 
\geq 
  1
  -
  \P\biggl(
    \bigcap\limits_{ n = \cN + 1 }^{ \cM }
    \bigl\{ 
      f( \Phi_{ n \cN - 1 }^{ n \cN - \cN } ) \leq \varepsilon
    \bigr\} 
  \biggr)
  -
  \P\biggl(
    \bigcup\limits_{ n = \cN }^{ \infty }
    \bigl\{ 
      \Xi 
      | 1 - \param \gamma |^{ \cN - 1 }
      \geq 
      \varepsilon
    \bigr\} 
  \biggr)
\\ & 
\textstyle 
=
  1
  -
  \biggl[
    \prod\limits_{ n = \cN + 1 }^{ \cM }
    \P\bigl(
      f( \Phi_{ n \cN - 1 }^{ n \cN - \cN } ) 
      \leq \varepsilon
    \bigr)
  \biggr]
  -
  \P\bigl(
    \Xi 
    | 1 - \param \gamma |^{ \cN - 1 }
    \geq 
    \varepsilon
  \bigr)
\\ & 
\textstyle 
=
  1
  -
  \biggl[
    \prod\limits_{ n = \cN + 1 }^{ \cM }
    \P\bigl(
      f( \Phi_1^{ n \cN - \cN } ) 
      \leq \varepsilon
    \bigr)
  \biggr]
  -
  \P\bigl(
    \Xi 
    | 1 - \param \gamma |^{ \cN - 1 }
    \geq 
    \varepsilon
  \bigr)
\\ & 
\textstyle 
=
  1
  -
  \biggl[
    \prod\limits_{ n = \cN + 1 }^{ \cM }
    \P\bigl(
      f( \Phi_1^0 ) 
      \leq \varepsilon
    \bigr)
  \biggr]
  -
  \P\bigl(
    \Xi 
    | 1 - \param \gamma |^{ \cN - 1 }
    \geq 
    \varepsilon
  \bigr)
\\ & 
\textstyle 
=
  1
  -
  \bigl[
    \P\bigl(
      f( \Phi_1^0 ) \leq \varepsilon
    \bigr)
  \bigr]^{ \cM - \cN }
  -
  \P\bigl(
    \Xi 
    | 1 - \param \gamma |^{ \cN - 1 }
    \geq 
    \varepsilon
  \bigr)
  .
\end{split}
\end{equation}
This and \cref{eq:probability_smaller_than_1} imply that 
\begin{equation}
\begin{split}
&
\textstyle 
  \P\biggl(
    \bigcup\limits_{ n = N }^{ \infty }
    \bigl\{ 
      f( \Theta_n ) > \frac{ \varepsilon }{ 2 } 
    \bigr\} 
  \biggr)
\\ & 
\textstyle 
\geq 
  \limsup\limits_{ \cN \to \infty }
  \limsup\limits_{ \cM \to \infty }
  \Bigl[
    1
    -
    \bigl[
      \P\bigl(
        f( \Phi_1^0 ) \leq \varepsilon
      \bigr)
    \bigr]^{ \cM - \cN }
    -
    \P\bigl(
      \Xi 
      | 1 - \param \gamma |^{ \cN - 1 }
      \geq 
      \varepsilon
    \bigr)
  \Bigr]
\\ &
\textstyle 
=
  \limsup\limits_{ \cN \to \infty }
  \Bigl[
    1
    -
    \P\bigl(
      \Xi 
      | 1 - \param \gamma |^{ \cN - 1 }
      \geq 
      \varepsilon
    \bigr)
  \Bigr]
=
  1
  -
  \liminf\limits_{ \cN \to \infty }
  \P\bigl(
    \Xi 
    | 1 - \param \gamma |^{ \cN - 1 }
    \geq 
    \varepsilon
  \bigr)
\\ &
\textstyle 
=
  1
  -
  \lim\limits_{ \cN \to \infty }
  \P\bigl(
    \Xi 
    | 1 - \param \gamma |^{ \cN - 1 }
    \geq 
    \varepsilon
  \bigr)
=
  1
  -
  \P\biggl(
    \bigcap\limits_{ \cN = 1 }^{ \infty }
    \{
      \Xi 
      | 1 - \param \gamma |^{ \cN - 1 }
      \geq 
      \varepsilon
    \} 
  \biggr)
\\ &
\textstyle 
=
  1
  -
  \P\Bigl(
    \forall \, \cN \in \N \colon 
    \Xi 
      | 1 - \param \gamma |^{ \cN - 1 }
      \geq 
      \varepsilon
  \Bigr)
=
  1
  -
  \P\Bigl(
    \inf_{ \cN \in \N }
    \bigl(
      \Xi 
      | 1 - \param \gamma |^{ \cN - 1 }
    \bigr)
    \geq 
    \varepsilon
  \Bigr)
= 1
  .
\end{split}
\end{equation}
\Hence that
\begin{equation}
\begin{split}
&
\textstyle 
  \P\biggl(
    \bigcup\limits_{ n = N }^{ \infty }
    \bigl\{ 
      f( \varTheta_n ) > 0
    \bigr\} 
  \biggr)
\geq 
  \P\biggl(
    \bigcup\limits_{ n = N }^{ \infty }
    \bigl\{ 
      f( \varTheta_n ) > \frac{ \varepsilon }{ 2 }
    \bigr\} 
  \biggr)
  = 1 .
\end{split}
\end{equation}
Combining this with 
\cref{eq:definition_of_varTheta_process,eq:definition_of_f} 
establishes \cref{eq:to_prove_prop:SGD_constant}. 
\end{cproof}

\section{Analysis of invariant measures of SGD processes}
\label{sec:invariant_measure}
	
The main goal of this section is to establish \cref{prop:positive:prob},
which shows that the assumption in \cref{eq:assumption_invariant_measure} in \cref{prop:SGD_constant} 
(increase of the test loss with positive probability for the invariant distribution)
holds under fairly general assumptions on the training data.

\subsection{Supports of random variables}

We first collect a few well-known definitions and properties regarding the support of measures.
In particular, \cref{def:supp_mu} can be found, for example, in \cite[Chapter 5]{Ambrosio2008_Book}.
The proofs of the following results are elementary and thus omitted, they can be found in the \href{https://arxiv.org/abs/2406.14340}{arXiv-version} of this article.

\begin{definition}[Support of a measure]
\label{def:supp_mu}
Let $ ( E, \cE ) $ be a topological space 
and let $ \mu \colon \cB( E ) \to [0,\infty] $ be a measure. 
Then we denote by $ \operatorname{supp}( \mu ) \subseteq E $ the set given by
\begin{equation}
\label{eq:definition_of_support}
  \operatorname{supp}( \mu )
  =
  \bigl\{ 
    x \in E 
    \colon 
    \bigl(
      \forall \, A \in \{ B \in \cE \colon x \in B \} 
      \colon 
      \mu( A ) > 0 
    \bigr)
  \bigr\}
  .
\end{equation}
\end{definition}
\cfclear 
\begin{lemma}[Elementary properties of the support]
\label{lem:closedness_of_support}
Let $ ( E, \cE ) $ be a topological space 
and let $ \mu \colon \cB( E ) \to [0,\infty] $ be a measure. 
Then 
\begin{enumerate}[label=(\roman*)]
\item it holds that $ \operatorname{supp}( \mu ) $ is closed, 
\item it holds that 
$
  E \backslash \operatorname{supp}( \mu ) 
  = 
  \cup_{ A \in \{ B \in \cE \colon \mu( B ) = 0 \} }
  A
$, 
and 
\item 
it holds 
for every $ A \in \cE $
that 
$
  \mu( A ) > 0 
$
if and only if 
$
  ( A \cap \operatorname{supp}( \mu ) ) \neq \emptyset  
$.
\end{enumerate}
\cfadd{def:supp_mu}\cfload.
\end{lemma}
%
%

\cfclear 
\begin{lemma}[Full measures of supports]
\label{lem:full_measure_support}
Let $ d \in \N $ and let $ \mu \colon \cB( \R^d ) \to [0,\infty] $ be a measure. 
Then \allowbreak 
$ \mu( \R^d \backslash \operatorname{supp}( \mu ) ) = 0 $
\cfadd{def:supp_mu}\cfload.
\end{lemma}

The support of a random variable is defined 
as the support of the probability distribution 
of the random variable 
(as the support of the push foward measure of the random variable, 
as the support of the law of the random variable). 
This is the subject of the next notion.

\cfclear
\begin{definition}[Support of a random variable]
\label{def:Supp_rand}
Let $ ( \Omega, \cF, \P ) $ be a probability space, 
let $ ( E, \cE ) $ be a topological space, 
and let $ X \colon \Omega \to E $ be a random variable. 
Then we denote by $ \spp( X ) \subseteq E $ the set given by
\begin{equation}
\label{eq:def_sppX}
  \spp( X )
  =
  \operatorname{supp}( 
    \P_X
  )
\end{equation}
\cfadd{def:supp_mu}\cfload.
\end{definition}

\cfclear
\begin{cor}[Support of a random variable]
\label{cor:Supp_rand}
Let $ ( \Omega, \cF, \P ) $ be a probability space, 
let $ d \in \N $, 
and let $ X \colon \Omega \to \R^d $ be a random variable. 
Then 
$
  \P(
    X \in 
    \spp( X )
  )
  = 1 
$
\cfadd{def:Supp_rand}\cfload.
\end{cor}
%
%

\subsection{Non-degenerate random variables}

We now introduce in \cref{def:supp_rand} a suitable condition for random variables that will be used in the following analysis.
\begin{definition}[Non-degenerate random variables]
\label{def:supp_rand}
Let $ ( \Omega, \cF, \P ) $ be a probability space, 
let $ ( E, \cE ) $ be a topological space, 
and let $ X \colon \Omega \to E $ be a random variable. 
Then we say that $ X $ is non-degenerate if and only if 
the interior of $ \spp( X ) $ is non-empty 
\cfadd{def:Supp_rand}\cfload.
\end{definition} 
\cfclear

For every dimension $ d \in \N $ we have that every random variable 
with a strictly positive density on some non-empty open subset of $ \R^d $ 
is non-degenerate. 
This is the subject of the next elementary lemma.

\begin{lemma}[Random variables with absolutely continuous distributions]
\label{lem:random_variable_with_density}
Let $ (\Omega, \cF, \P) $ be a probability space,
let $d \in \N$, let $X \colon \Omega \to \R^d$ be a random variable, 
let $ U \subseteq \R^d $ be a non-empty open set, 
let $p \colon U \to (0 , \infty ) $ be measurable, 
and assume for every open $ A \subseteq U $ that $ \P( X \in A ) = \int_A p(x) \, \d x $.
Then 
\begin{enumerate}[label=(\roman*)]
\item 
it holds that $ U \subseteq \spp( X ) $ and 
\item 
it holds that $ X $ is non-degenerate
\end{enumerate}
\cfadd{def:Supp_rand,def:supp_rand}\cfload.
\end{lemma}
\begin{cproof}{lem:random_variable_with_density}
Throughout this proof let $ \lambda \colon \cB( \R^d ) \to [0,\infty] $ be the Lebesgue--Borel measure on $ \R^d $. 
\Nobs that the fact that for all $ x \in U $ it holds that $ p(x) > 0 $ implies 
that for every measurable $ A \subseteq U $ with 
$
  \int_A p( x ) \, \d x = 0 
$
it holds that 
\begin{equation}
  \lambda( A ) 
  =
  \lambda(
    \{ x \in A \colon p(x) \neq 0 \}
  ) = 0
  .
\end{equation}
\Hence for every measurable $ A \subseteq U $ with 
$ \lambda( A ) > 0 $
that
$
  \int_A p(x) \, \d x > 0
$. 
\Hence for every non-empty open $ A \subseteq U $ 
that
$
  \P_X( A )
  =
  \P( X \in A ) = \int_A p(x) \, \d x > 0
$. 
Combining this with \cref{eq:definition_of_support} establishes 
that 
$
  U \subseteq \operatorname{supp}( \P_X ) = \spp( X )
$.
\end{cproof}

In the remainder of this subsection we present several elementary and well-known properties 
for supports of (independent) random variables and linear combinations of sets.
For the detailed proofs we again refer the reader to the \href{https://arxiv.org/abs/2406.14340}{arXiv-version}.

\cfclear 
\begin{lemma}[Support of the sum of two random variables]
\label{lem:support:sum}
Let $ (\Omega, \cF, \P) $ be a probability space, let $ \dd \in \N $,
and let $ X \colon \Omega \to \R^\dd $ 
and $ Y \colon \Omega \to \R^{ \dd } $ be independent random variables.
Then\footnote{Note that for every $ d \in \N $, every set $ A \subseteq \R^d $, and every 
set $ B \subseteq \R^d $ it holds that 
$ 
  A + B = 
  \{ 
    c \in \R^d \colon 
    (
      \exists \, a \in A, b \in B \colon c = a + b
    ) 
  \}
$
(\href{https://en.wikipedia.org/wiki/Minkowski_addition}{Minkowski sum} of $ A $ and $ B $). \Moreover 
for every $ d \in \N $, $ r \in \R $ and every set $ A \subseteq \R^d $
it holds that 
$
  r A = \{ x \in \R^d \colon ( \exists \, a \in A \colon x = r a ) \}
$.} 
\begin{equation}
  ( \spp( X ) + \spp( Y ) )
  \subseteq 
  \spp( X + Y ) 
\end{equation}
\cfadd{def:Supp_rand}\cfload.
\end{lemma}
%
%

\cfclear 
\begin{lemma}[Support of the scalar multiplication of a random variable]
\label{lem:support:multiplication}
Let $ (\Omega, \cF, \P) $ be a probability space, let $ \dd \in \N $, 
$ \gamma \in \R $, 
and let $ X \colon \Omega \to \R^\dd $ 
be a random variable.
Then
\begin{equation}
  \spp( \gamma X ) 
  =
  ( \gamma \spp( X ) )
\end{equation}
\cfadd{def:Supp_rand}\cfload.
\end{lemma}
%
%

\cfclear
\begin{cor}[Support of a linear combination of random variables]
\label{cor:support_sum_RV}
\cfadd{def:Supp_rand}
Let $ \dd, N \in \N $, 
$ \gamma_1, \gamma_2, \allowbreak \dots, \allowbreak \gamma_N \in \R $, 
let $ (\Omega, \cF, \P) $ be a probability space, 
and let $ X_n \colon \Omega \to \R^d $, $ n \in \{ 1, 2, \dots, N \} $, 
be independent random variables.
Then 
\begin{equation}
\label{eq:linear_combination_support}
\textstyle 
  \bigl( 
    \sum_{ n = 1 }^N \gamma_n \spp( X_n ) 
  \bigr)
  \subseteq 
  \spp\bigl( 
    \sum_{ n = 1 }^N \gamma_n X_n 
  \bigr)
\end{equation}
\cfload.
\end{cor}
\begin{cproof}{cor:sum_nondeg}
\Nobs that \cref{lem:support:sum}, \cref{lem:support:multiplication}, and induction establish \cref{eq:linear_combination_support}. 
\end{cproof}

\cfclear
\begin{lemma}[Minkowski sums for open sets]
\label{lem:open_Minkowski_sums}
Let $ d \in \N $ and let $ A \subseteq \R^d $ be open. 
Then 
\begin{enumerate}[label=(\roman*)]
\item it holds for all $ \gamma \in \R \backslash \{ 0 \} $ that $ \gamma A $ is open and 
\item it holds for all $ B \subseteq \R^d $ that $ A + B $ is open. 
\end{enumerate}
\end{lemma}


\cfclear
\begin{cor}[Linear combinations of non-degenerate random variables]
\label{cor:sum_nondeg}
\cfadd{def:supp_rand}
Let $ \dd, N \in \N $, 
$ \gamma_1, \allowbreak \gamma_2, \allowbreak \dots, \allowbreak \gamma_N \in \R $, 
let $ ( \Omega, \cF, \P ) $ be a probability space, 
let $ X_n \colon \Omega \to \R^d $, $ n \in \{ 1, 2, \dots, N \} $, be independent random variables, 
let $ m \in \{ 1, 2, \dots, N \} $ satisfy $ \gamma_m \neq 0 $, 
and assume that $ X_m $ is non-degenerate. 
Then $ \sum_{ n = 1 }^N \gamma_n X_n $ is non-degenerate
\cfload.
\end{cor}

\subsection{Supports of series of random variables}

We now apply \cref{cor:support_sum_RV} to show a property of the support of infinite sums of random variables.
This will be useful in the following since the invariant distribution of the \SGD\ method for quadratic target functions is precisely of this form (see \cref{prop:SGD_constant} above).

\cfclear
\begin{lemma}[Support of a series of random variables]
\label{lem:series:support}
\cfadd{def:Supp_rand}
Let $(\Omega, \cF, \P)$ be a probability space,
let $\dd \in \N$, $\eta \in (0 , 1)$,
let $ X_n \colon \Omega \to \R^\dd $, $ n \in \N $, be i.i.d.~random variables,
assume $ \E[ \norm{X_1 } ] < \infty $, 
let $ A \subseteq \Omega $ satisfy 
$
  A = \cu{ \sum_{ n = 1 }^{ \infty } \eta ( 1 - \eta )^{ n - 1 } \norm{ X_n } < \infty } 
$, 
and let 
$
  \Ysum \colon \Omega \to \R ^\dd
$ be a random variable 
which satisfies for all $ \omega \in A $ that 
$
  \Ysum(\omega ) = \sum_{ n = 1 }^\infty \eta ( 1 - \eta )^{ n - 1 } X_n( \omega )
$. 
Then
\begin{enumerate}[label=(\roman*)]
\item \label{lem:series:support:item1}
it holds that $A \in \cF$ and $\P ( A ) = 1$ and
\item \label{lem:series:support:item2}
it holds that $ \spp( X_1 ) \subseteq \spp( \Ysum ) $
\end{enumerate}
\cfload.
\end{lemma}
\begin{cproof}{lem:series:support}
Throughout this proof 
for every $ x \in \spp( X_1 ) $, $ \varepsilon \in (0, \infty) $ 
let $ \fN_{ x, \varepsilon } \in \N $ satisfy 
\begin{equation}
\label{lem:series:support:eq_3}
  (1 - \eta )^{ \fN_{ x, \varepsilon } }
  \bigl( 
    \norm{x} + \varepsilon^{-1} 
    \E[ \norm{X_1} ] 
  \bigr) 
  < \varepsilon .
\end{equation}
\Nobs that the fact that 
$ 
  \E\bigl[ 
    \sum_{ n = 1 }^\infty \eta ( 1 - \eta )^{ n - 1 } \norm{ X_n } 
  \bigr] < \infty 
$ 
establishes \cref{lem:series:support:item1}. 
\Nobs that \cref{cor:support_sum_RV} assures for all $ N \in \N $ that
\begin{equation}
\label{lem:series:support:eq_1}
\begin{split}
&
  \spp\bigl( 
    \ssum_{ n = 1 }^N \eta ( 1 - \eta )^{ n - 1 } X_n 
  \bigr)
\\ &
\supseteq 
  \bigl(
    \ssum_{ n = 1 }^N \br*{ \eta ( 1 - \eta )^{ n - 1 } \spp( X_n ) }
  \bigr)
  = 
  \bigl(
    \ssum_{ n = 1 }^N \br*{ \eta ( 1 - \eta )^{ n - 1 } \spp( X_1 ) } 
  \bigr)
\\
  & 
  \supseteq 
  \bigl(
    \bigl[ 
      \ssum_{ n = 1 }^N \eta ( 1 - \eta )^{ n - 1 } 
    \bigr]
    \spp( X_1 ) 
  \bigr)
   = ( 1 - ( 1 - \eta )^N ) \spp( X_1 ) .
\end{split}
\end{equation}
\Hence 
for all 
$ N \in \N $, 
$ x \in \spp( X_1 ) $, $ \varepsilon \in (0,\infty) $ 
that
$
  \P\bigl( 
    \norm{ 
      (1 - (1 - \eta )^N ) x 
      -
      \ssum_{ n = 1 }^{ N } 
      \eta ( 1 - \eta )^{ n - 1 } X_n 
    } < \varepsilon 
  \bigr) > 0 
$. 
\Hence for all 
$ x \in \spp( X_1 ) $, $ \varepsilon \in (0,1) $ that
\begin{equation}
\label{eq:support_series_probability_positive_1}
\begin{split}
&
  \P\bigl(
    \|
      x 
      -
      \ssum_{ n = 1 }^{ \fN_{ x, \varepsilon } } 
      \eta ( 1 - \eta )^{ n - 1 } X_n 
    \|
    < 2 \varepsilon
  \bigr)
\\ &
=
  \P\bigl(
    \|
      ( 1 - \eta )^{ \fN_{ x, \varepsilon } } x
      +
      ( 
        1 - ( 1 - \eta )^{ \fN_{ x, \varepsilon } } 
      )
      x 
      -
      \ssum_{ n = 1 }^{ \fN_{ x, \varepsilon } } 
      \eta ( 1 - \eta )^{ n - 1 } X_n 
    \|
    < 2 \varepsilon
  \bigr)
\\ &
\geq 
  \P\bigl(
    \bigl\{
      \|
        ( 1 - \eta )^{ \fN_{ x, \varepsilon } }
        x
      \|
      < 
      \varepsilon
    \bigr\}
    \cap 
    \bigl\{ 
      \norm{
        ( 
          1 
          - 
          ( 1 - \eta )^{ \fN_{ x, \varepsilon } } 
        ) x
        -
        \ssum_{ n = 1 }^{ \fN_{ x, \varepsilon } } 
        \eta ( 1 - \eta )^{ n - 1 } X_n 
      }
      < \varepsilon
    \bigr\}
  \bigr)
\\ &
  =
  \P\bigl(
      \norm{
        ( 
          1 
          - 
          ( 1 - \eta )^{ \fN_{ x, \varepsilon } } 
        ) x
        -
        \ssum_{ n = 1 }^{ \fN_{ x, \varepsilon } } 
        \eta ( 1 - \eta )^{ n - 1 } X_n 
      }
      < \varepsilon
  \bigr)
  > 0 
  .
\end{split}
\end{equation}
\Moreover the Markov inequality implies 
for all $ N \in \N $, $ R \in (0, \infty ) $ that
\begin{equation}
\label{lem:series:support:eq_2}
\begin{split}
&
  \P\bigl(
    \ssum_{ n = {N+1} }^{ \infty } 
    \eta ( 1 - \eta )^{ n - 1 } \norm{ X_n } \ge R 
  \bigr)
\le 
  R^{-1} 
  \E\br*{
    \ssum_{ n = {N+1} }^{ \infty } 
    \eta ( 1 - \eta )^{ n - 1 } 
    \norm{ X_n }
  } 
\\
& 
=
  R^{-1} 
  \bigl( 
    \ssum_{ n = {N+1} }^{ \infty } 
    \eta ( 1 - \eta )^{ n - 1 } \E [ \norm{ X_n } ]
  \bigr)
  = R^{ - 1 } (1 - \eta )^N \E[ \norm{X_1} ] < \infty 
  .
\end{split}
\end{equation}
Combining this with \cref{lem:series:support:eq_3} shows 
for all $ x \in \spp( X_1 ) $, $ \varepsilon \in (0,1) $ that 
\begin{equation}
\begin{split}
&
  \P\bigl( 
    \ssum_{ n = { \fN_{ x, \varepsilon } + 1 } }^{ \infty } 
    \eta ( 1 - \eta )^{ n - 1 } \norm{ X_n }
    < \varepsilon
  \bigr) 
=
  1 -
  \P\bigl( 
    \ssum_{ n = { \fN_{ x, \varepsilon } + 1 } }^{ \infty } 
    \eta ( 1 - \eta )^{ n - 1 } \norm{ X_n }
    \geq \varepsilon
  \bigr) 
\\ &
\geq 
  1 - \varepsilon^{ - 1 } 
  ( 1 - \eta )^{ \fN_{ x, \varepsilon } } \E[ \norm{X_1} ] 
\geq
  1 - 
  ( 1 - \eta )^{ \fN_{ x, \varepsilon } } 
  \bigl(
    \| x \|
    +
    \varepsilon^{ - 1 } \E[ \norm{X_1} ] 
  \bigr)
>
  1 - \varepsilon
> 
  0
  .
\end{split}
\end{equation}
The triangle inequality 
and \cref{eq:support_series_probability_positive_1}
\hence show 
for all $ x \in \spp( X_1 ) $, $ \varepsilon \in (0,1) $ 
that
\begin{equation}
\begin{split}
& 
  \P\bigl( \norm{ x - \Ysum } < 3 \varepsilon \bigr)
\\ &
\ge 
  \P\bigl(
    \bigl\{ 
      \norm{ 
        x
        -
        \ssum_{ n = 1 }^{ \fN_{ x, \varepsilon } } \eta ( 1 - \eta )^{ n - 1 } X_n  
      } < 2 \varepsilon 
    \bigr\}
    \cap 
    \bigl\{
      \ssum_{ n = { \fN_{ x, \varepsilon } + 1} }^\infty 
      \eta ( 1 - \eta )^{ n - 1 } \norm{ X_n }
      < \varepsilon 
    \bigr\}
  \bigr)
\\
&
  = 
  \P\bigl( 
    \norm{ 
      x
      -
      \ssum_{ n = 1 }^{ \fN_{ x, \varepsilon } } 
      \eta ( 1 - \eta )^{ n - 1 } X_n 
    } 
    < 2 \varepsilon 
  \bigr)
  \,
  \P\bigl(
    \ssum_{ n = { \fN_{ x, \varepsilon } + 1 } }^{ \infty } 
    \eta ( 1 - \eta )^{ n - 1 }   
    \norm{ X_n } 
    < \varepsilon 
  \bigr) 
  > 0
  .
\end{split}
\end{equation}
\Hence for all $ x \in \spp( X_1 ) $ that $ x \in \spp( \Ysum ) $. 
This establishes \cref{lem:series:support:item2}.
\end{cproof}

\subsection{Analysis of the test loss under the invariant measure}

In this part we prove in \cref{prop:positive:prob} below 
that, roughly speaking, the probability that the test loss after one \SGD\ step strictly increases 
is strictly positive.
In our mathematical analyses of \SGD\ processes 
in \cref{cor:SGD_constant_nondeg} and \cref{statement_adaptive_convergence} below 
we apply \cref{prop:positive:prob} in the situation 
where the random variables $ Y_k \colon \Omega \to \R^{ \dd } $, 
$ k \in \N_0 $, satisfy for all $ n \in \N $ that 
$ Y_n = M^{ - 1 } \sum_{ m = 1 }^M X_{ n, m } $ 
where $ M \in \N $ is the batch-size for the considered \SGD\ method 
and where $ ( X_{ n, m } )_{ (n, m) \in \N^2 }  $ are i.i.d.\ samples drawn from some distribution 
(see \cref{cor:SGD_constant_nondeg} below for details). 
In this context we note that \cref{cor:sum_nondeg} ensures that we have that 
if the random variables $ X_{ n, m } \colon \Omega \to \R^{ \dd } $, $ (n, m) \in \N^2 $, 
are non-degenerate, then it holds that the random variables 
$ ( Y_n )_{ n \in \N } $ are non-degenerate as well.

\cfclear
\begin{prop}[Test loss strictly increases with strictly positive probability]
\label{prop:positive:prob}
\cfadd{def:supp_rand}
Let $ \dd \in \N $, $ \eta \in (0, 1) $, 
let $ (\Omega, \cF, \P) $ be a probability space,
let $ Y_k \colon \Omega \to \R^\dd $, $ k \in \N_0 $, be non-degenerate i.i.d.\ random variables, 
assume $ \E[ \norm{ Y_0 } ] < \infty $, 
let $ \X \colon \Omega \to \R^\dd $ be a non-degenerate random variable,
assume that $ ( Y_k )_{ k \in \N_0 } $ and $ \X $ are independent, 
let $ A \subseteq \Omega $ satisfy 
$
  A = \{ \sum_{ k = 1 }^{ \infty } \eta ( 1 - \eta )^{ k - 1 } \| Y_k \| < \infty \}
$, 
and let $ \Ysum \colon \Omega \to \R^{ \dd } $ be a random variable 
which satisfies for all $ \omega \in A $ that
$
  \Ysum( \omega ) = \sum_{ k = 1 }^\infty \eta ( 1 - \eta ) ^{ k - 1 } Y_k( \omega )
$ \cfload. 
Then
\begin{equation}
\label{eq:probability_strictly_positive}
  \P\bigl( 
    \norm{ (1 - \eta ) \Ysum + \eta  Y_0 - \X } > \norm{\Ysum - \X } 
  \bigr) > 0 .
\end{equation}
\end{prop}

Regarding the proof of \cref{prop:positive:prob},
due to the non-degeneracy of $Y_0$ and $Z$ we are able to find a direction from $\chi$ in which, roughly speaking, the distance to $Z$ increases with positive probability and such that the random vector $Y_0$ points approximately into that direction with positive probability
(see \cref{positive:prob:proof_claim_1} below for details).

\begin{cproof}{prop:positive:prob}
Throughout this proof for every $ x \in \R^d $, $ r \in (0,\infty) $ 
let $ B_r( x ) \subseteq \R^d $ satisfy 
\begin{equation}
  B_r( x ) = \{ y \in \R^d \colon \| x - y \| < r \}
  .
\end{equation}
\Nobs the assumption that $ Y_0 $ is non-degenerate 
assures that there exist $ x \in \R^\dd \backslash \{ 0 \} $, $ r \in (0, \infty) $
which satisfy 
\begin{equation}
\label{eq:B_r_x_in_spp_Y_0}
  B_r( x ) \subseteq \spp( Y_0 ) .
\end{equation}
In the following we assume without loss of generality that 
$ \norm{x} = 1 $ (otherwise we consider the random variables 
$ \norm{x}^{ - 1 } Y_k $, $ k \in \N_0 $, 
and 
$
  \norm{x}^{ - 1 } Z 
$).
\Moreover \cref{lem:series:support} shows that 
\begin{equation}
\label{eq:spp_Y_0_in_spp_chi}
  \spp( Y_0 ) \subseteq \spp( \chi ) .
\end{equation}
\Moreover the assumption that $ \X $ is non-degenerate implies 
that the interior of $ \spp( \X ) $ is non-empty. 
Combining this with the fact that 
the interior of 
\begin{equation}
  \{ z \in \R^d \colon \spro{ x, z } = \norm{x}^2 \} 
  =
  \{ z \in \R^d \colon \spro{ x, z - x } = 0 \} 
\end{equation}
is the empty set proves that 
it does not hold that 
$
  \spp( \X )
  \subseteq 
  \{ z \in \R^d \colon \spro{ x, z - x } = 0 \} 
$.
This implies that 
\begin{equation}
  (
    \spp( \X ) 
    \cap 
    \{ z \in \R^d \colon \spro{ x, z - x } \neq 0 \} 
  )
  =
  (
    \spp( \X ) 
    \cap 
    ( \R^d \backslash \{ z \in \R^d \colon \spro{ x, z - x } = 0 \} )
  )
  \neq \emptyset 
  .
\end{equation}
\cref{lem:closedness_of_support} 
and the fact that 
$
  \{ z \in \R^d \colon \spro{ x, z - x } \neq 0 \} 
$
is open \hence show that 
$
  \P_Z(
    \{ z \in \R^d \colon \spro{ x, z - x } \neq 0 \} 
  )
  > 0 
$.
\Hence that 
\begin{equation}
\begin{split}
&
  \lim_{ Q \nearrow \infty }
  \lim_{ \delta \searrow 0 }
  \P\bigl(
    \bigl\{ 
      | \spro{ x, Z - x } |
      > \delta 
    \bigr\} 
    \cap 
    \{ 
      \| Z \| < Q
    \} 
  \bigr)
\\ &
  =
  \P\bigl(
    \cup_{ \delta, Q \in (0,\infty) \cap \Q }
    \bigl[
      \bigl\{ 
        | \spro{ x, Z - x } |
        > \delta 
      \bigr\} 
      \cap 
      \{ 
        \| Z \| < Q
      \} 
    \bigr]
  \bigr)
\\ & 
  =
  \P\bigl(
    \cup_{ \delta \in (0,\infty) \cap \Q }
    \bigl[
      \{ 
        | \spro{ x, Z - x } |
        > \delta 
      \} 
      \cap 
      \{ 
        \| Z \| < \infty 
      \}
    \bigr]
  \bigr)
  =
  \P\bigl(
    \cup_{ \delta \in (0,\infty) \cap \Q }
    \bigl\{ 
      | \spro{ x, Z - x } |
      > \delta 
    \bigr\} 
  \bigr)
\\ &
  =
  \P\bigl(
    | \spro{ x, Z - x } |
    > 0 
  \bigr)
  =
  \P\bigl(
    \spro{ x, Z - x }
    \neq 0 
  \bigr)
  > 0 
  .
\end{split}
\end{equation}
\Hence 
that there exist 
$
  \delta \in ( 0, \frac{ r \eta }{ 2 } )
$, 
$ Q \in ( 1, \infty ) $ 
which satisfy 
\begin{equation}
\label{eq:property_delta_Q}
  \P\bigl( 
    \bigl\{ 
      | \spro{ x, \X - x } | > \delta 
    \bigr\}
    \cap 
    \bigl\{ 
      \norm{\X} < Q 
    \bigr\}
  \bigr)
  > 0
  .
\end{equation}
\Nobs that \cref{eq:property_delta_Q} shows that 
\begin{equation}
\begin{split}
&
    \P\bigl( 
      \bigl\{ 
        \spro{ x, \X - x } > \delta 
      \bigr\}
      \cap 
      \bigl\{ 
        \norm{\X} < Q 
      \bigr\}
    \bigr)
    +
    \P\bigl( 
      \bigl\{ - \spro{ x, \X - x } > \delta \bigr\} 
      \cap 
      \bigl\{ 
        \norm{\X} < Q 
      \bigr\} 
    \bigr) 
\\ & 
=
    \P\bigl( 
      \bigl\{ 
        \spro{ x, \X - x } > \delta 
      \bigr\}
      \cap 
      \bigl\{ 
        \norm{\X} < Q 
      \bigr\}
    \bigr)
    +
    \P\bigl( 
      \bigl\{ \spro{ x, \X - x } < - \delta \bigr\} 
      \cap 
      \bigl\{ 
        \norm{\X} < Q 
      \bigr\} 
    \bigr) 
\\ & 
  =
  \P\bigl( 
    \bigl\{ 
      | \spro{ x, \X - x } | > \delta 
    \bigr\}
    \cap 
    \bigl\{ 
      \norm{\X} < Q 
    \bigr\}
  \bigr)
  > 0
  .
\end{split}
\end{equation}
\Hence there exists $ \sparam \in \{ -1, 1 \} $ which satisfies 
\begin{equation}
\label{eq:OBDA}
  \P\bigl( 
    \bigl\{ 
      \sparam \spro{ x, \X - x } > \delta 
    \bigr\}
    \cap 
    \bigl\{ 
      \norm{\X} < Q 
    \bigr\}
  \bigr)
  > 0 
  .
\end{equation}
In the following let $ \varepsilon \in \R $, $ z_1, z_2 \in \R^d $ satisfy 
for all $ k \in \{ 1, 2 \} $ that
\begin{equation}
\label{positive:prob:proof_eq_defe}
  \varepsilon = 
  \min\cu*{ 
    \frac{ \delta }{ 2 } , 
    \frac{ 
      \delta^2 
    }{ 
      9 Q 
    }
  } 
\qqandqq
  z_k = 
  \bigl(
    1 + \varXi ( - 1 )^k ( \delta - \varepsilon ) 
  \bigr) x .
\end{equation}
In the next step of our proof of \cref{eq:probability_strictly_positive} 
we show that
\begin{equation}
\label{positive:prob:proof_claim_1} 
  \P\bigl( 
    \bigl\{ 
      ( 1 - \eta ) \Ysum + \eta Y_0 
      \in 
      B_{ \varepsilon }( z_1 ) 
    \bigr\} 
    \cap 
    \bigl\{ 
      \Ysum \in B_\varepsilon ( z_2 ) 
    \bigr\} 
  \bigr) > 0
  .
\end{equation}
We will use \cref{positive:prob:proof_claim_1} together with \cref{eq:OBDA} 
to establish \cref{eq:probability_strictly_positive} 
(see 
\cref{positive:prob:proof_claim_2}--\cref{eq:probability_strictly_bigger_than_zero_final_step_of_the_proof} 
below). 
To establish \cref{positive:prob:proof_claim_1} we note that 
the fact that $ \delta < \frac{ r \eta }{ 2 } $ 
and the assumption that $ \eta < 1 $
demonstrate that 
$ 2 \eta^{ - 1 } - 1 > 1 $
and 
$
  \delta 
  < 
  2^{ - 1 } \eta r 
  =
  ( 2 \eta^{ - 1 } )^{ - 1 } r
  \leq ( 2 \eta^{ - 1 } - 1 )^{ - 1 } r
$. 
\Hence that 
\begin{equation}
  0 < (\delta - \varepsilon)(2 \eta^{-1} - 1 ) < r - \varepsilon ( 2 \eta^{-1} - 1  ) < r - \varepsilon
  .
\end{equation}
The triangle inequality and the fact that 
\begin{equation}
\begin{split}
&
  \eta^{ - 1 } ( z_1 - (1 - \eta) z_2 ) 
  =
  \eta^{ - 1 } ( z_1 - z_2 + \eta z_2 ) 
  =
  \eta^{ - 1 } ( z_1 - z_2 ) + z_2 
\\ &
  = 
  \eta^{ - 1 } 
  \bigl[ 
    2 \varXi ( \varepsilon - \delta ) x 
  \bigr]
  +
  ( 1 + \varXi ( \delta - \varepsilon ) ) x
  = 
  2 \varXi \eta^{ - 1 } 
  ( \varepsilon - \delta ) x 
  + ( 1 + \varXi ( \delta - \varepsilon ) ) x
\\ &
  =
  \bigl[ 
    -
    \varXi 
    ( \delta - \varepsilon )
    2 \eta^{ - 1 } 
    +
    1 
    + 
    \varXi ( \delta - \varepsilon )
  \bigr]
  x
  = 
  \bigl[
    1 + 
    \varXi 
    (\delta - \varepsilon ) ( 1 - 2 \eta^{ - 1 } ) 
  \bigr]
  x
\end{split}
\end{equation}
therefore show that 
for all 
$
  y \in 
  B_{ \varepsilon }(
    \eta^{ - 1 } 
    ( z_1 - (1 - \eta) z_2 ) 
  ) 
$
it holds that 
\begin{equation}
  \| x - y \|
\leq 
  \| 
    x 
    -
    \eta^{ - 1 } 
    ( z_1 - (1 - \eta) z_2 ) 
  \|
  +
  \|
    \eta^{ - 1 } 
    ( z_1 - (1 - \eta) z_2 ) 
    -
    y
  \|
<
  ( \delta - \varepsilon ) ( 2 \eta^{ - 1 } - 1 )
  +
  \varepsilon
<
  r
  .
\end{equation}
\Hence that 
$
  B_{ \varepsilon }(
    \eta^{ - 1 } 
    ( z_1 - (1 - \eta) z_2 ) 
  ) 
  \subseteq B_r ( x )
$. 
This and \cref{eq:B_r_x_in_spp_Y_0} show that 
$
  B_{ \varepsilon }(
    \eta^{ - 1 } 
    ( z_1 - (1 - \eta) z_2 ) 
  ) 
  \subseteq 
  \spp( Y_0 )
$. 
\Hence that 
\begin{equation}
\label{eq:ball_intersect_no_emptyset}
  \bigl(
    B_{ \varepsilon }(
      \eta^{ - 1 } 
      ( z_1 - (1 - \eta) z_2 ) 
    ) 
    \cap 
    \spp( Y_0 )
  \bigr)
  \neq 
  \emptyset 
  .
\end{equation}
\Moreover 
\cref{positive:prob:proof_eq_defe} and 
the assumption that $ \eta < 1 $ 
show that
for all 
$ 
  y \in B_{ \varepsilon }( z_2 )
$
it holds that 
\begin{equation}
  \| x - y \|
  \leq 
  \| x - z_2 \|
  +
  \| z_2 - y \|
  <
  \| 
    x - 
    ( 1 + \varXi ( \delta - \varepsilon ) ) x
  \|
  +
  \varepsilon
  =
  ( \delta - \varepsilon )
  +
  \varepsilon
  =
  \delta 
  <
\textstyle 
  \frac{ r \eta }{ 2 }
  <
  r 
  .
\end{equation}
\Hence that 
$
  B_{ \varepsilon }( z_2 )
  \subseteq B_r( x )
$. 
Combining this with \cref{eq:B_r_x_in_spp_Y_0,eq:spp_Y_0_in_spp_chi} 
demonstrates that 
$
  B_{ \varepsilon }( z_2 ) 
  \subseteq B_r( x )
  \subseteq \spp( Y_0 )
  \subseteq \spp( \chi )
$.
\Hence that 
$
  ( 
    B_{ \varepsilon }( z_2 ) \cap \spp( \chi )
  )
  \neq 
  \emptyset 
$. 
The fact that $ \Ysum $ and $ Y_0 $ are independent random variables, 
\cref{eq:ball_intersect_no_emptyset}, 
and \cref{lem:closedness_of_support} 
\hence show that 
\begin{equation}
\label{eq:strictly_positive_probability_preparation}
\begin{split}
&
  \P\Bigl( 
    \bigl\{ 
      \Ysum \in B_\varepsilon ( z_2 ) 
    \bigr\}  
    \cap 
    \bigl\{ 
    \bigl[
      ( 1 - \eta ) \Ysum + \eta Y_0 
    \bigr]
    \in 
    \bigl[
      ( 1 - \eta ) B_{ \varepsilon }( z_2 ) 
      +
      \eta 
      B_{ \varepsilon }( \eta^{ - 1 } ( z_1 - (1 - \eta) z_2 ) ) 
    \bigr]
    \bigr\} 
  \Bigr) 
\\ &
  \geq
  \P\bigl( 
    \bigl\{ 
      \Ysum \in B_\varepsilon ( z_2 ) 
    \bigr\}  
    \cap 
    \bigl\{  
      Y_0 \in B_\varepsilon ( \eta^{-1} ( z_1 - (1 - \eta) z_2 ) ) 
    \bigr\} 
  \bigr) 
\\ &
  =
  \P\bigl( 
    \Ysum \in B_\varepsilon ( z_2 ) 
  \bigr)
  \,
  \P\bigl(
    Y_0 \in B_\varepsilon ( \eta^{-1} ( z_1 - (1 - \eta) z_2 ) ) 
  \bigr) 
  > 0 .
\end{split}
\end{equation}
\Moreover 
the triangle inequality proves that for all 
$ 
  v \in B_{ ( 1 - \eta ) \varepsilon }( (1 - \eta ) z_2 ) 
$, 
$ 
  w \in B_{ \eta \varepsilon }( z_1 - ( 1 - \eta ) z_2 )
$
it holds that 
\begin{equation}
\begin{split}
  \|
    z_1
    -
    ( v + w )
  \|
& =
  \bigl\|
    \bigl[
      z_1
      -
      ( 1 - \eta ) z_2
      -
      w 
    \bigr]
    +
    \bigl[ 
      ( 1 - \eta ) z_2 - v 
    \bigr]
  \bigr\|
\\ &
\leq 
  \|
    z_1
    -
    ( 1 - \eta ) z_2
    -
    w 
  \|
  +
  \|
    ( 1 - \eta ) z_2 - v 
  \|
<
  \eta \varepsilon
  +
  ( 1 - \eta ) \varepsilon
=
  \varepsilon
  .
\end{split}
\end{equation}
\Hence that 
\begin{equation}
\begin{split}
&
  (1 - \eta ) B_\varepsilon ( z_2 ) + \eta B_\varepsilon (\eta^{-1} ( z_1 - (1 - \eta) z_2 ) )
\\ &
  = B_{(1 - \eta ) \varepsilon} ( (1 - \eta ) z_2 )
  + B_{\eta \varepsilon} ( z_1 - (1 - \eta) z_2 ) 
  \subseteq B_\varepsilon ( z_1)
  .
\end{split}
\end{equation}
Combining this with \cref{eq:strictly_positive_probability_preparation} 
establishes \cref{positive:prob:proof_claim_1}. 
In the next step of our proof of \cref{eq:probability_strictly_positive} 
we show that 
for all 
$ 
  u \in 
  \{ 
    v \in \R^\dd \colon 
    (
      \varXi \spro{ x, v - x } > \delta 
    )
    \wedge 
    ( 
      \norm{v} < Q 
    ) 
  \} 
$, 
$
  y_1 \in B_{ \varepsilon }( z_1 ) 
$, 
$
  y_2 \in B_{ \varepsilon }( z_2 ) 
$
it holds that 
\begin{equation}
\label{positive:prob:proof_claim_2} 
  \norm{ u - y_1 } > \norm{ u - y_2 }
  .
\end{equation}
We will use \cref{positive:prob:proof_claim_2} together with 
\cref{eq:OBDA} and \cref{positive:prob:proof_claim_1}
to prove \cref{eq:probability_strictly_positive} 
(see \cref{eq:preparation_for_combination,eq:probability_strictly_bigger_than_zero_final_step_of_the_proof} 
below). 
To establish \cref{positive:prob:proof_claim_2} we note that 
the fact that $ \norm{x} = 1 $ 
shows that 
for all $ u, p \in \R^\dd $, $ \lambda \in \R $
with 
$ p = \spro{ x, u } x $
it holds that 
$
  \spro{ \lambda x, u - p } 
  = 
  \lambda \spro{ x, u - \spro{x,u} x } 
  =
  \lambda \spro{ x, u } - 
  \lambda \spro{x,u} \spro{ x, x }
  = 0
$. 
\Hence 
for all $ u, p \in \R^\dd $
with 
$
  \norm{u} < Q
$
and 
$ p = \spro{ x, u } x $
that 
\begin{equation}
  \spro{ u - p, p } 
  =
  \spro{ u - p, p - z_1 } 
  = 
  \spro{ u - p , p - z_2 } 
  = 0 
\end{equation}
and 
$
  \norm{ u - p }^2 
  \le 
  \norm{ u - p }^2 + \norm{ p }^2 
  = 
  \norm{ u }^2 < Q^2 
$.
\Hence 
for all $ u, p \in \R^{ \dd } $
with 
$
  \varXi \spro{ x, u - x } > \delta
$, 
$
  \norm{u} < Q
$, 
and
$ p = \spro{ x, u } x $
that 
\begin{equation}
\begin{split}
&
  \norm{ u - z_1 } - \norm{ u - z_2 }
  =
  \sqrt{
    \| ( u - p ) + ( p - z_1 ) \|^2
  }
  -
  \sqrt{
    \| ( u - p ) + ( p - z_2 ) \|^2
  }
\\ &
  = 
  \sqrt{ \norm{ u - p }^2 + \norm{ p - z_1 }^2 }
  - 
  \sqrt{ \norm{u - p } ^2 + \norm{p - z_2 } ^2 }
\\ 
&
  = 
  \frac{ 
    (
      \norm{ u - p }^2
      +
      \norm{ p - z_1 }^2 
    )
    - 
    (
      \norm{ u - p }^2
      +
      \norm{ p - z_2 }^2 
    )
  }{ 
    \sqrt{
      \norm{ u - p }^2 + \norm{ p - z_1 }^2 
    }
    + 
    \sqrt{ 
      \norm{ u - p }^2 + \norm{ p - z_2 }^2 
    }
  } 
\\ 
&
  = 
  \frac{ 
    \norm{ p - z_1 }^2 
    - 
    \norm{ p - z_2 }^2 
  }{ 
    \sqrt{
      \norm{ u - p }^2 + \norm{ p - z_1 }^2 
    }
    + 
    \sqrt{ 
      \norm{ u - p }^2 + \norm{ p - z_2 }^2 
    }
  } 
\\ &
  >
  \frac{ 
    \norm{ \spro{ x, u } x - ( 1 - \varXi ( \delta - \varepsilon ) ) x }^2 
    - 
    \norm{ \spro{ x, u } x - ( 1 + \varXi ( \delta - \varepsilon ) ) x }^2 
  }{ 
    \sqrt{
      Q^2 + \norm{ \spro{ x, u } x - ( 1 - \varXi ( \delta - \varepsilon ) ) x }^2 
    }
    + 
    \sqrt{ 
      Q^2 + \norm{ \spro{ x, u } x - ( 1 + \varXi ( \delta - \varepsilon ) ) x }^2 
    }
  } 
\\
&
  = 
  \frac{
    | \spro{ x, u - x } + \varXi ( \delta - \varepsilon ) |^2 
    - 
    | \spro{ x, u - x } - \varXi ( \delta - \varepsilon ) |^2 
  }{
    \sqrt{ 
      Q^2 
      + 
      | \spro{ x, u - x } + \varXi ( \delta - \varepsilon ) |^2 
    } 
    + 
    \sqrt{ 
      Q^2 
      + 
      | \spro{ x, u - x } - \varXi ( \delta - \varepsilon ) |^2 
    }
  } 
\\
&
  = 
  \frac{
    | \varXi \spro{ x, u - x } + ( \delta - \varepsilon ) |^2 
    - 
    | \varXi \spro{ x, u - x } - ( \delta - \varepsilon ) |^2 
  }{
    \sqrt{ 
      Q^2 
      + 
      | \varXi \spro{ x, u - x } + ( \delta - \varepsilon ) |^2 
    } 
    + 
    \sqrt{ 
      Q^2 
      + 
      | \varXi \spro{ x, u - x } - ( \delta - \varepsilon ) |^2 
    }
  } 
\\
& 
  =
  \frac{ 
    2 ( \delta - \varepsilon )
    (
      2 \varXi \spro{ x, u - x }
    )
  }{ 
    \sqrt{ 
      Q^2 
      + 
      | \varXi \spro{ x, u - x } + ( \delta - \varepsilon ) |^2 
    } 
    + 
    \sqrt{ 
      Q^2 
      + 
      | \varXi \spro{ x, u - x } - ( \delta - \varepsilon ) |^2 
    }
  } 
\\
& 
  \ge 
  \frac{ 
    4 
    ( \delta - \varepsilon ) 
    ( \varXi \spro{ x, u - x } )
  }{ 
    2 
    \sqrt{ 
      Q^2 
      + 
      | 
        \varXi \spro{ x, u - x } + ( \delta - \varepsilon ) 
      |^2 
    } 
  } 
  >
  \frac{ 
    2 \delta 
    ( \delta - \varepsilon) 
  }{
    \sqrt{ 
      Q^2 
      + 
      | 
        \varXi \spro{ x, u - x } + \delta 
      |^2 
    } 
  }
\\ &
> 
  \frac{ 
    2 \delta ( \delta - ( \nicefrac{ \delta }{ 2 } ) ) 
  }{ 
    \sqrt{ Q^2 + ( 2 \varXi \spro{ x, u - x } )^2 } 
  } 
=
  \frac{ 
    \delta^2  
  }{ 
    \sqrt{ Q^2 + 4 ( \spro{ x, u } - 1 )^2 } 
  } 
\geq 
  \frac{ 
    \delta^2 ) 
  }{ \sqrt{ Q^2 + 4 ( \| x \| \| u \| + 1 )^2 } } 
\\ & \geq
  \frac{ 
    \delta^2 
  }{ \sqrt{ Q^2 + 4 ( 2 Q )^2 } } 
=
  \frac{ 
    \delta^2 
  }{ \sqrt{ 17 Q^2 } } 
\geq 
  \frac{ 
    \delta^2 
  }{ Q \sqrt{ 17 } } 
  .
\end{split}
\end{equation}
Combining this with \cref{positive:prob:proof_eq_defe} ensures 
for all 
$ 
  u \in 
  \{ 
    v \in \R^\dd \colon 
    (
      \spro{ x, v - x } > \delta 
    )
    \wedge 
    ( 
      \norm{v} < Q 
    ) 
  \} 
$ 
that
\begin{equation}
  \norm{u - z_1 } - \norm{u - z_2 }
  > 
  \frac{ \delta^2 }{ 
    Q \sqrt{ 17 }
  }
  =
  2 
  \left[ 
    \frac{ \delta^2 }{ 2 Q \sqrt{17} }
  \right] 
  =
  2 
  \left[ 
    \frac{ \delta^2 }{ Q \sqrt{68} }
  \right] 
  \geq 
  2 
  \left[ 
    \frac{ \delta^2 }{ 9 Q }
  \right] 
  \geq 
  2 \varepsilon 
  .
\end{equation}
The triangle inequality \hence demonstrates 
for all 
$ 
  u \in 
  \{ 
    v \in \R^\dd \colon 
    (
      \varXi \spro{ x, v - x } > \delta 
    )
    \wedge 
    ( 
      \norm{v} < Q 
    ) 
  \} 
$, 
$ y_1 \in B_{ \varepsilon }( z_1) $, 
$ y_2 \in B_{ \varepsilon }( z_2 ) $
that
\begin{equation}
\begin{split}
&
  \norm{ u - y_1 } - \norm{ u - y_2 }
=
  \| ( u - z_1 ) - ( y_1 - z_1 ) \|
  -
  \| ( u - z_2 ) + ( z_2 - y_2 ) \|
\\ &
\geq  
  \bigl(
    \| u - z_1 \| 
    -
    \| z_1 - y_1 \|
  \bigr)
  -
  \bigl(
    \| u - z_2 \|
    +
    \| z_2 - y_2 \|
  \bigr)
\\ &
  =
  \bigl(
    \| u - z_1 \| - \| u - z_2 \|
  \bigr)
  -
  \bigl(
    \| y_1 - z_1 \|
    +
    \| y_2 - z_2 \|
  \bigr)
\\ &
  >
  2 \varepsilon 
  -
  \bigl(
    \| y_1 - z_1 \|
    +
    \| y_2 - z_2 \|
  \bigr)  
  > 
  2 \varepsilon - 2 \varepsilon = 0 .
\end{split}
\end{equation}
This establishes \cref{positive:prob:proof_claim_2}. 
We now combine \cref{positive:prob:proof_claim_1} and \cref{positive:prob:proof_claim_2} with 
\cref{eq:OBDA} to establish \cref{eq:probability_strictly_positive}. 
\Nobs that \cref{positive:prob:proof_claim_2} shows that 
\begin{equation}
\label{eq:preparation_for_combination}
\begin{split}
&
  \bigl(
    \bigl\{ \sparam \spro{ x, \X - x } > \delta \bigr\}  
    \cap 
    \bigl\{ \norm{ \X } < Q \bigr\} 
    \cap 
    \bigl\{ (1 - \eta ) \Ysum + \eta  Y_0 \in B_\varepsilon ( z_1 ) \bigr\}
    \cap  
    \bigl\{ \Ysum \in B_{ \varepsilon }( z_2 ) \bigr\}
  \bigr)
\\ &
\subseteq 
  \bigl\{ 
    \bigl\| \X - \bigl[ (1 - \eta ) \Ysum + \eta  Y_0 \bigr] \bigr\|
    > 
    \bigl\| \X - \Ysum \bigr\|
  \bigr\} 
\end{split}
\end{equation}
Combining this 
and the fact that $ ( Y_0, \chi ) $ and $ Z $ are independent random variables 
with \cref{eq:OBDA} 
and \cref{positive:prob:proof_claim_1} 
demontrates that
\begin{equation}
\label{eq:probability_strictly_bigger_than_zero_final_step_of_the_proof}
\begin{split}
& 
  \P\bigl( 
    \norm{ (1 - \eta ) \Ysum + \eta  Y_0 - \X } 
    > 
    \norm{ \Ysum - \X } 
  \bigr)
\\
& 
  \ge 
  \P\bigl( 
    \bigl\{ (1 - \eta ) \Ysum + \eta  Y_0 \in B_\varepsilon ( z_1 ) \bigr\}
    \cap  
    \bigl\{ \Ysum \in B_\varepsilon ( z_2 ) \bigr\}
    \cap 
    \bigl\{ \varXi \spro{ x, \X - x } > \delta \bigr\}  
    \cap 
    \bigl\{ \norm{ \X } < Q \bigr\} 
  \bigr)
\\
&
  = 
  \P\bigl( 
    \bigl\{ 
      (1 - \eta ) \Ysum + \eta  Y_0 \in B_{ \varepsilon }( z_1 ) 
    \bigr\} 
    \cap 
    \bigl\{ 
      \Ysum \in B_\varepsilon ( z_2 ) 
    \bigr\} 
  \bigr)
  \,
  \P\bigl( 
    \bigl\{ \varXi \spro{ x, \X - x } > \delta \bigr\} 
    \cap 
    \bigl\{ \norm{ \X} < Q \bigr\} 
  \bigr)
\\ &
  > 0 .
\end{split}
\end{equation}
\end{cproof}

\cfclear
\begin{cor}[Test loss strictly increases with strictly positive probability]
\label{cor:SGD_constant_nondeg}
\cfadd{def:supp_rand}
Let $ \fd \in \N $, $ \param \in (0,\infty) $, 
let $ \fl \colon \R^{ \fd } \times \R^{ \fd } \to \R $ satisfy 
for all $ \theta, \vartheta \in \R^{ \fd } $ that 
\begin{equation}
  \fl( \theta, \vartheta ) 
  =
\textstyle 
  \frac{ \param }{ 2 } \| \theta - \vartheta \|^2
  ,
\end{equation}
let $ \gamma \in ( 0, \param^{ - 1 } ) $, 
let $ ( \Omega, \mathcal{F}, \P ) $ be a probability space, 
let $ X_{ n, m } \colon \Omega \to \R^{ \fd } $, $ (n, m) \in \Z^2 $, be non-degenerate i.i.d.\ random variables,  
let $ M, \fM \in \N $, 
let 
$ \Theta \colon \N_0 \times \Omega \to \R^{ \fd } $
satisfy for all $ n \in \N $ that 
\begin{equation}
			\Theta_n 
			= \Theta_{ n - 1 } 
			- 
			\frac{ \gamma }{ M }
			\left[ 
			\sum_{ m = 1 }^M
			( \nabla_{ \theta } \fl )( \Theta_{ n - 1 }, X_{ n, m } ) 
			\right]
			,
\end{equation}
assume 
$ 
  \{ \sup_{ n, m \in \Z } \| X_{ n, m } \| < \infty \} = \Omega 
$, 
assume that 
$
  ( X_{ n, m } )_{ (n, m) \in \Z^2 }
$
and 
$ \Theta_0 $
are independent, 
and let $ N \in \N $ \cfload. 
Then 
\begin{equation}
			\label{eq:to_prove_cor:SGD_constant}
			\textstyle 
			\P\biggl(
			\bigcup\limits_{ n = N }^{ \infty }
			\biggl\{ 
			\sum\limits_{ m = 1 }^{ \fM }
			\fl ( \Theta_n, X_{ n, - m } ) 
			>
			\sum\limits_{ m = 1 }^{ \fM }
			\fl ( \Theta_{ n - 1 }, X_{ n, - m } ) 
			\biggr\} 
			\biggr)
			= 1
			.
\end{equation}
\end{cor}
\begin{cproof}{cor:SGD_constant_nondeg}
Throughout this proof for every $ n \in \N_0 $ 
let $ Y_n \colon \Omega \to \R^\fd $ satisfy
$
  Y_n = \frac{ 1 }{ M } ( \sum_{m=1}^M X_{1 - n , m } )
$
and let 
$ \Ysum \colon \Omega \to \R^{ \fd } $ 
and 
$ \X \colon \Omega \to \R^{ \fd } $ 
satisfy 
\begin{equation}
  \Ysum = 
  \ssuml_{ n = 1 }^{ \infty } \param \gamma ( 1 - \param \gamma )^{n-1}
  Y_n
  \qqandqq 
  \X = \displaystyle\frac{ 1 }{ \fM } 
  \left[ 
    \ssuml_{ m = 1 }^{ \fM } X_{ 2, m }
  \right] 
  .
\end{equation}
\Nobs that \cref{cor:sum_nondeg} proves 
that $ Y_n $, $ n \in \N_0 $, and $ \X $ are non-degenerate random variables. 
\Moreover the assumption that $ X_{ n, m } $, $ (n, m) \in \Z^2 $, are i.i.d.\ random variables 
ensures that $ Y_n $, $n \in \N_0$, are i.i.d.\ random variables. 
\Moreover the assumption that $ X_{ n, m } $, $ (n, m) \in \Z^2 $, 
are i.i.d.\ random variables 
shows that $ ( Y_n )_{ n \in \N_0 } $ and $ \X $ are independent. 
\Moreover the assumption that 
$ \{ \sup_{ n, m \in \Z } \| X_{ n, m } \| < \infty \} = \Omega $ 
implies that $ \E\bigl[ \norm{ Y_0 } \bigr] < \infty $.
\cref{prop:positive:prob} 
(applied with $ \eta \with \param \gamma $ in the notation of \cref{prop:positive:prob}) 
\hence demonstrates that
\begin{equation}
  \P\rbr[\big]{
    \norm*{(1 - \param \gamma ) \Ysum + \param \gamma Y_0 - \X } > \norm{\Ysum - \X } 
  } > 0 .
\end{equation}
Combining this with \cref{prop:SGD_constant} establishes \cref{eq:to_prove_cor:SGD_constant}. 
\end{cproof}

\section{Convergence analysis of SGD with adaptive learning rates} 
\label{sec:convergence_analysis}

We now prove in \cref{statement_adaptive_convergence} the main theoretical result of this article and, thereby, also establish \cref{thm:intro} from the introduction.

\cfclear
\begin{theorem}[\SGD\ with adaptive learning rates]
\label{statement_adaptive_convergence}
\cfadd{def:supp_rand}
Let 
$ \fd \in \N $, 
$ \param \in (0,\infty) $, 
let $ ( \Omega, \mathcal{F}, \P ) $ be a probability space, 
let 
$ 
  X_{ n, m } \colon \Omega \to \R^{ \fd } 
$, $ (n, m) \in \Z^2 $, be non-degenerate i.i.d.\ random variables, 
let 
$ 
  \fl = 
  ( \fl( \theta, x ) )_{ ( \theta, x ) \in \R^{ \fd } \times \R^{ \fd } } \colon \R^{ \fd } \times \R^{ \fd } 
  \to \R 
$ 
satisfy for all $ \theta, x \in \R^{ \fd } $ that 
\begin{equation}
  \fl ( \theta, x ) 
  =
  \param \| \theta - x \|^2
  ,
\end{equation}
let $ M, \fM \in \N $, 
let 
$ \Theta \colon \N_0 \times \Omega \to \R^{ \fd } $
and 
$ \gamma \colon \N \times \Omega \to \R $ 
be stochastic processes which satisfy for all $ n \in \N $ that 
\begin{equation}
\label{eq:definition_of_Theta_n_conjecture}
  \Theta_n 
  = \Theta_{ n - 1 } 
  - 
  \frac{ \gamma_n }{ M }
  \left[ 
    \sum_{ m = 1 }^M
    ( \nabla_{ \theta } \fl )( \Theta_{ n - 1 }, X_{ n, m } ) 
  \right]
  ,
\end{equation}
let $ ( \nu_n )_{ n \in \N } \subseteq \R $ be strictly decreasing,
assume 
$
  \lim_{ n \to \infty } \nu_n = 0
$
and 
$
  \sum_{ n = 1 }^{ \infty } \nu_n = \infty 
$,
assume for all $ n \in \N $ that 
$
  \gamma_1 = \nu_1 < \param^{ - 1 }
$
and 
\begin{equation}
\label{eq:definition_of_adaptive_gamma}
  \gamma_{ n + 1 } = 
  \begin{cases}
    \gamma_n
  &
    \colon
      \sum_{ m = 1 }^{ \fM }
      \fl ( \Theta_n, X_{ n, - m } ) 
    \leq 
      \sum_{ m = 1 }^{ \fM }
      \fl ( \Theta_{ n - 1 }, X_{ n, - m } ) 
  \\
    \max\bigl(
      ( 0, \gamma_n )
      \cap 
      [
        \cup_{
          k = 1
        }^{ \infty }
        \{ \nu_k \} 
      ]
    \bigr)
  &
    \colon
      \sum_{ m = 1 }^{ \fM }
      \fl ( \Theta_n, X_{ n, - m } ) 
    >
      \sum_{ m = 1 }^{ \fM }
      \fl ( \Theta_{ n - 1 }, X_{ n, - m } ) 
      ,
  \end{cases}
\end{equation}
for every $ t \in [0,\infty) $ let $ N_t \colon \Omega \to \R $ satisfy 
\begin{equation}
\textstyle 
  N_t = 
  \inf\bigl\{ 
    n \in \N_0 \colon 
    \sum_{ k = 1 }^n \gamma_k \geq t
  \bigr\} 
  ,
\end{equation}
assume that 
$ 
  ( X_{ n, m } )_{ ( n, m ) \in \Z^2 }
$ 
and 
$ \Theta_0 $ 
are independent, 
and assume 
$  
  \E\bigl[ 
    \| \Theta_0 \|^2 
    + 
    \sup_{ n \in \N } 
    \| X_{ n, 0 } \| 
  \bigr] < \infty 
$ 
\cfload. 
Then 
\begin{equation}
\label{eq:convergence_for_quadratic_optimization_problems_with_adaptive_learning_rates}
  \lim_{ t \to \infty }
  \E\bigl[
    \|  
      \Theta_{ N_t } - \E[ X_{ 1, 1 } ]
    \|^2
  \bigr]
  = 0
  .
\end{equation}
\end{theorem}

The proof of \cref{statement_adaptive_convergence} relies on an application of \cref{cor:convergence_in_probability_ex},
and thus we need to verify the required assumptions.
The main difficulty lies in establishing that $\gamma_n \to 0$ almost surely (see \cref{eq:gamma_n_converge_to_zero} below).
The idea to show this is as follows:
If $\gamma_n \to 0$ woud not hold almost surely then the learning rates would stop decreasing at some (random) time with positive probability, and hence the test loss would always decrease from that time on due to the definition of $\gamma_n$ in \cref{eq:definition_of_adaptive_gamma}. 
By the strong Markov property, this would contradict \cref{cor:SGD_constant_nondeg}.

\begin{cproof}{statement_adaptive_convergence}
Throughout this proof 
let $ \psi \colon (0,\infty) \to \R $ satisfy for all $ x \in (0,\infty) $ that 
\begin{equation}
\label{eq:definition_of_psi_round_down}
  \psi( x ) 
  =
  \max\bigl(
    ( 0, x )
    \cap 
    [
      \cup_{
        k = 1
      }^{ \infty }
      \{ \nu_k \} 
    ]
  \bigr)
  ,
\end{equation}
for every $ n \in \N_0 $ 
let $ \mathbb{F}_n \subseteq \cF $
satisfy 
\begin{equation}
\label{eq:filtration_introduction}
  \mathbb{F}_n = 
  \sigma\bigl( 
    \Theta_0 , 
    ( X_{ a, b } )_{ 
      (a, b) \in 
      ( \Z \cap [0,n] ) \times \Z
      \}
    }
  \bigr)
  ,
\end{equation}
for every $ \lambda \in \R $, $ k \in \N_0 $  
let 
$
  \varTheta^{ \lambda, k } \colon \N_0 \times \Omega \to \R^{ \fd }
$
satisfy for all $ n \in \N $ that
$
  \varTheta^{ \lambda, k }_0 = \Theta_k
$
and 
\begin{equation}
\label{eq:definition_of_varTheta_gamma_k}
  \varTheta^{ \lambda, k }_n  
  =
  \varTheta^{ \lambda, k }_{ n - 1 }
  -
  \frac{ \lambda }{ M }
  \left[ 
    \sum_{ m = 1 }^M
    ( 
      \nabla_{ \theta } \fl
    )( \varTheta^{ \lambda, k }_{ n - 1 }, X_{ { k + n }, m } ) 
  \right] 
  ,
\end{equation}
let 
$
  \fL \colon \Z \times \R^{ \fd } \times \Omega \to \R
$
satisfy for all $ n \in \Z $, $ \theta \in \R^{ \fd } $ that
\begin{equation}
\label{eq:definition_of_fracL}
  \fL( n, \theta )
  =
  \sum_{ m = 1 }^{ \fM }
  \fl ( \theta, X_{ n, - m } ) 
  ,
\end{equation}
and for every $ \lambda \in \R $, $ k \in \N_0 $ 
let $ \supertau_{ \lambda, k } \colon \Omega \to \N \cup \{ \infty \} $
satisfy 
\begin{equation}
\label{eq:definition_of_tau_lambda_k}
\textstyle 
  \supertau_{ \lambda, k } = 
  \inf\bigl(
    \bigl\{ 
      n \in \N \colon
      \fL( k + n, \varTheta^{ \lambda, k }_{ n } )
      >
      \fL( k + n, \varTheta^{ \lambda, k }_{ n - 1 } )
    \bigr\}
    \cup 
    \bigl\{ \infty \bigr\}
  \bigr)
  .
\end{equation}
\Nobs that 
\cref{eq:definition_of_Theta_n_conjecture,eq:definition_of_varTheta_gamma_k} ensure that 
for all 
$ 
  n \in \N
$, 
$ \omega \in \Omega $
it holds that 
\begin{equation}
  \Theta_n( \omega )
  =
  \varTheta^{ \gamma_n( \omega ), n - 1 }_1( \omega )
  .
\end{equation}
\Moreover 
\cref{eq:definition_of_adaptive_gamma,eq:definition_of_fracL} 
show that 
for all 
$ 
  n \in \N 
$
it holds that 
$ \gamma_1 = \nu_1 $ 
and 
\begin{equation}
\label{eq:definition_of_adaptive_gamma_using_fL}
  \gamma_{ n + 1 } = 
  \begin{cases}
    \gamma_{ n }
  &
    \colon
    \fL( n, \Theta_{ n } )
    \leq 
    \fL( n, \Theta_{ n - 1 } ) 
  \\
    \psi( \gamma_n )
  &
    \colon
    \fL( n, \Theta_{ n } )
    >
    \fL( n, \Theta_{ n - 1 } )
  \end{cases}
  .
\end{equation}
\Hence for all $ n \in \N $ that 
\begin{equation}
\label{eq:equal_learning_rates_conclusion}
  \{ 
    \gamma_{ n + 1 } = \gamma_{ n }
  \} 
  =
  \{ 
    \fL( n, \Theta_{ n } ) 
    \leq 
    \fL( n, \Theta_{ n - 1 } ) 
  \} 
  .
\end{equation}
We intend to prove \cref{statement_adaptive_convergence} 
through an application of \cref{cor:convergence_in_probability_ex}. 
We thus need to verify the assumptions of \cref{cor:convergence_in_probability_ex}. 
For this we note that \cref{eq:filtration_introduction} ensures 
that
\begin{equation}
\label{eq:filtered_probability_space}
  ( \Omega, \mathcal{F}, ( \mathbb{F}_n )_{ n \in \N_0 }, \P ) 
\end{equation}
is a filtered probability space. 
\Moreover \cref{eq:definition_of_adaptive_gamma_using_fL} ensures that for all $ n \in \N $ that 
\begin{equation}
\label{eq:definition_of_adaptive_gamma_using_fL_2}
  \gamma_{ n + 1 }
  =
  \gamma_n 
  \mathbbm{1}_{
    \{ 
      \fL( n, \Theta_{ n } )
      \leq 
      \fL( n, \Theta_{ n - 1 } )
    \}
  }
  +
  \psi( \gamma_n )
  \mathbbm{1}_{
    \Omega \backslash
    \{ 
      \fL( n, \Theta_{ n } )
      \leq 
      \fL( n, \Theta_{ n - 1 } )
    \}
  }
  .
\end{equation}
\Moreover \cref{eq:definition_of_psi_round_down} 
ensures for all $ x, y \in (0,\infty) $ with $ x \leq y $ that 
$
  \psi( x ) \leq \psi( y )
$. \Hence that $ \psi $ is measurable. 
Combining this with 
\cref{eq:definition_of_Theta_n_conjecture,eq:filtration_introduction,eq:definition_of_fracL,eq:definition_of_adaptive_gamma_using_fL_2} 
implies that 
for all $ n \in \N $ 
it holds that 
$
  \gamma_n 
$
is $ \mathbb{F}_{ n - 1 } $-measurable. 
\Moreover 
\cref{eq:definition_of_adaptive_gamma_using_fL_2}, 
the fact that for all $ x \in (0,\infty) $ it holds that $ \psi(x) < x $, 
and 
the fact that for all $ n \in \N $, $ x \in (0,\infty) $ with $ x \geq \nu_n $ 
it holds that 
$
  \psi( x ) \geq \psi( \nu_n ) = \nu_{ n + 1 }
$
demonstrate that for all $ n \in \N $ 
it holds that
\begin{equation}
\label{eq:gamma_monotonicity}
  \gamma_{ n + 1 } \leq \gamma_n 
\qquad 
\text{and}
\qquad 
  \gamma_n 
  \geq \nu_n 
  .
\end{equation}
The assumption that $ \sum_{ n = 1 }^{ \infty } \nu_n = \infty $
\hence shows that 
\begin{equation}
\label{eq:gamma_n_summability}
\textstyle 
    \sum_{ n = 1 }^{ \infty } \gamma_n = \infty 
  .
\end{equation}
In the next step we verify that
\begin{equation}
\label{eq:gamma_n_converge_to_zero}
\textstyle 
  \P\bigl(
    \lim_{ n \to \infty } \gamma_n = 0
  \bigr)
  = 1 .
\end{equation}
We show \cref{eq:gamma_n_converge_to_zero} by contradiction. 
We thus assume that 
\begin{equation}
\label{eq:contradiction}
\textstyle 
  \P\bigl(
    \limsup_{ n \to \infty } \gamma_n > 0
  \bigr)
  > 0
  .
\end{equation}
\Nobs that \cref{eq:definition_of_adaptive_gamma}, \cref{eq:gamma_monotonicity}, 
and the assumption that $ \lim_{ n \to \infty } \nu_n = 0 $ 
assure that 
\begin{equation}
\textstyle 
  \bigl\{ 
    \limsup_{ n \to \infty } \gamma_n = 0 
  \bigr\} 
  = 
  \bigl\{
    \omega \in \Omega 
    \colon 
    \#\bigl( 
      \bigl\{ 
        n \in \N 
        \colon
        \gamma_{ n + 1 }( \omega ) < \gamma_n( \omega )
      \bigr\} 
    \bigr)
    = \infty 
  \bigr\}
  .
\end{equation}
\Hence that 
\begin{equation}
\begin{split}
\textstyle 
  \bigl\{ 
    \limsup_{ n \to \infty } \gamma_n > 0 
  \bigr\} 
& =
  \bigl\{ 
    \omega \in \Omega 
    \colon
    \bigl(
      \exists \, n \in \N \colon 
      \forall \, m \in \N \cap [ n, \infty ) \colon
      \gamma_{ m + 1 }( \omega ) = \gamma_m( \omega )
    \bigr)
  \bigr\} 
\\ & =
  \bigl\{ 
    \omega \in \Omega 
    \colon
    \bigl(
      \exists \, n \in \N \colon 
      \forall \, m \in \N \cap [ n, \infty ) \colon
      \gamma_m( \omega ) = \gamma_n( \omega )
    \bigr)
  \bigr\} 
\\ &
  =
  \cup_{ n \in \N }
  \cap_{ m = n }^{ \infty }
  \{ 
    \gamma_m = \gamma_n
  \} 
  .
\end{split}
\end{equation}
The fact that the measure $ \P $ is continuous from below 
and \cref{eq:contradiction} \hence show that
\begin{equation}
\textstyle 
  0
<
  \P\bigl(
    \limsup_{ n \to \infty } \gamma_n > 0 
  \bigr)
=
  \P\bigl(
    \cup_{ n \in \N }
    \cap_{ m = n }^{ \infty }
    \{ 
      \gamma_m = \gamma_n
    \} 
  \bigr)
=
  \lim_{ n \to \infty }
  \P\bigl(
    \cap_{ m = n }^{ \infty }
    \{ 
      \gamma_m = \gamma_n
    \} 
  \bigr)
  .
\end{equation}
This 
implies that there exists $ \cN \in \N $ 
which satisfies
\begin{equation}
\label{eq:probability_bigger_than_zero}
  \P\bigl(
    \cap_{ m = \cN }^{ \infty }
    \{ 
      \gamma_m = \gamma_{ \cN }
    \} 
  \bigr)
  > 0 .
\end{equation}
\Nobs that \cref{eq:definition_of_varTheta_gamma_k} shows that for all 
$ 
  \omega \in 
  (
    \cap_{ m = \cN }^{ \infty }
    \{ 
      \gamma_m = \gamma_{ \cN }
    \} 
  )
$, 
$ m \in \N_0 $
it holds that
\begin{equation}
  \varTheta^{ \gamma_{ \cN }( \omega ), \, \cN }_m( \omega )
  =
  \Theta_{ \cN + m }( \omega )
  .
\end{equation}
\Hence that 
\begin{equation}
\label{eq:varTheta_and_Theta_coincide}
  (
    \cap_{ m = \cN }^{ \infty }
    \{ 
      \gamma_m = \gamma_{ \cN }
    \} 
  )
\subseteq
  \bigl(
    \cap_{ m = 0 }^{ \infty }
    \bigl\{
      \varTheta^{ \gamma_{ \cN }, \, \cN }_m
      =
      \Theta_{ \cN + m }
    \bigr\}
  \bigr)
  .
\end{equation}
\Moreover 
\cref{cor:SGD_constant_nondeg} 
and the assumption that $ X_{ n, m } \colon \Omega \to \R^{ \fd } $, $ (n,m) \in \Z^2 $, are non-degenerate i.i.d.\ random variables 
show that for all $ \lambda \in (0, \param^{-1}) $, $ k \in \N_0 $ it holds that 
\begin{equation}
  \P\bigl( 
    \supertau_{ \lambda, k } < \infty 
  \bigr) = 1
  .
\end{equation}
The fact that for all $ n \in \N $ it holds that $ \nu_n < \nu_1 < \param^{ - 1 } $
\hence demonstrates for all $ n \in \N $, $ k \in \N_0 $ that
\begin{equation}
  \P\bigl( 
    \supertau_{ \nu_n, k } < \infty 
  \bigr) = 1
  .
\end{equation}
\Hence for all $ k \in \N_0 $ that 
$
  \P\bigl( \cap_{ n = 1 }^{ \infty } \{ \supertau_{ \nu_n, k } < \infty \} \bigr) = 1
$. 
Combining this with \cref{eq:probability_bigger_than_zero} 
shows that 
\begin{equation}
  \P\bigl(
    \bigl(
      \cap_{ m = \cN }^{ \infty }
      \{ 
        \gamma_m = \gamma_{ \cN }
      \} 
    \bigr)
    \cap 
    \bigl( 
      \cap_{ n = 1 }^{ \infty } 
      \{ 
        \supertau_{ \nu_n, \, \cN } < \infty 
      \} 
    \bigr)
  \bigr) 
  > 0
  .
\end{equation}
\Hence that 
\begin{equation}
\label{eq:not_equal_to_empty_set}
\textstyle
    \bigl(
      \cap_{ m = \cN }^{ \infty }
      \{ 
        \gamma_m = \gamma_{ \cN }
      \} 
    \bigr)
    \cap 
    \bigl( 
      \cap_{ w = 1 }^{ \infty } 
      \{ 
        \supertau_{ \nu_w, \, \cN } < \infty 
      \} 
    \bigr)
  \neq 
  \emptyset 
  .
\end{equation}
In the following let $ \cT \colon \Omega \to \N \cup \{ \infty \} $ 
satisfy 
\begin{equation}
  \cT = \supertau_{ \gamma_{ \cN }, \, \cN } 
  .
\end{equation}
\Moreover \cref{eq:definition_of_tau_lambda_k} 
assures that for all $ \lambda \in \R $, $ k \in \N_0 $ 
it holds that 
\begin{equation}
  \{ 
    \supertau_{ \lambda, k } < \infty 
  \}  
  =
  \{ 
    \supertau_{ \lambda, k } < \infty 
  \}  
  \cap 
  \bigl\{ 
    \fL( 
      k + \supertau_{ \lambda, k }, \varTheta^{ \lambda, k }_{ \supertau_{ \lambda, k } } 
    )
    >
    \fL( 
      k + \supertau_{ \lambda, k }, \varTheta^{ \lambda, k }_{ \supertau_{ \lambda, k } - 1 } 
    )
  \bigr\} 
  .
\end{equation}
\Hence that 
\begin{equation}
\begin{split}
&
\textstyle 
  \cap_{ w = 1 }^{ \infty }
  \{ 
    \supertau_{ \nu_w, \, \cN } < \infty 
  \}  
\\ &
\textstyle 
=
  \cap_{ w = 1 }^{ \infty }
  \bigl(
  \{ 
    \supertau_{ \nu_w, \, \cN } < \infty 
  \}  
  \cap 
  \bigl\{ 
    \fL( 
      \cN + \supertau_{ \nu_w, \, \cN }, 
      \varTheta^{ \nu_w, \, \cN }_{ \supertau_{ \nu_w, \, \cN } } 
    )
    >
    \fL( 
      \cN + \supertau_{ \nu_w, \, \cN }, 
      \varTheta^{ \nu_w, \, \cN }_{ \supertau_{ \nu_w, \, \cN } - 1 } 
    )
  \bigr\} 
  \bigr)
\\ &
\textstyle 
= 
  \bigl(
    \cap_{ w = 1 }^{ \infty }
    \{ 
      \supertau_{ \nu_w, \, \cN } < \infty 
    \}
  \bigr)
  \cap 
  \bigl(
    \cap_{ w = 1 }^{ \infty }
    \bigl\{ 
      \fL( 
        \cN + \supertau_{ \nu_w, \, \cN }, 
        \varTheta^{ \nu_w, \, \cN }_{ \supertau_{ \nu_w, \, \cN } } 
      )
      >
      \fL( 
        \cN + \supertau_{ \nu_w, \, \cN }, 
        \varTheta^{ \nu_w, \, \cN }_{ \supertau_{ \nu_w, \, \cN } - 1 } 
      )
    \bigr\} 
  \bigr)
\\ &
\textstyle 
= 
  \bigl(
    \cap_{ w = 1 }^{ \infty }
    \{ 
      \supertau_{ \nu_w, \, \cN } < \infty 
    \}
  \bigr)
  \cap 
  \bigl(
    \cap_{ w = 1 }^{ \infty }
    \bigl\{ 
      \fL( 
        \cN + \supertau_{ \nu_w, \, \cN }, 
        \varTheta^{ \nu_w, \, \cN }_{ \supertau_{ \nu_w, \, \cN } } 
      )
      >
      \fL( 
        \cN + \supertau_{ \nu_w, \, \cN }, 
        \varTheta^{ \nu_w, \, \cN }_{ \supertau_{ \nu_w, \, \cN } - 1 } 
      )
    \bigr\} 
  \bigr)
\\ &
\quad 
  \cap
  \bigl\{ 
    \fL(
      \cN + \cT, 
      \varTheta^{ \gamma_{ \cN }, \, \cN }_{ \cT } 
    )
    >
    \fL( 
      \cN + \cT, 
      \varTheta^{ \gamma_{ \cN }, \, \cN }_{ \cT - 1 } 
    )
  \bigr\}
\\ &
\textstyle 
= 
  \bigl(
    \cap_{ w = 1 }^{ \infty }
    \{ 
      \supertau_{ \nu_w, \, \cN } < \infty 
    \}
  \bigr)
  \cap
  \bigl\{ 
    \fL( 
      \cN + \cT, 
      \varTheta^{ \gamma_{ \cN }, \, \cN }_{ \cT } 
    )
    >
    \fL( 
      \cN + \cT, 
      \varTheta^{ \gamma_{ \cN }, \, \cN }_{ \cT - 1 } 
    )
  \bigr\}
\\ &
\textstyle 
\subseteq
  \bigl\{ 
    \cT < \infty 
  \bigr\} 
  \cap 
  \bigl\{ 
    \fL( 
      \cN + \cT, 
      \varTheta^{ \gamma_{ \cN }, \, \cN }_{ \cT } 
    )
    >
    \fL( 
      \cN + \cT, 
      \varTheta^{ \gamma_{ \cN }, \, \cN }_{ \cT - 1 } 
    )
  \bigr\}
  .
\end{split}
\end{equation}
Combining this with 
\cref{eq:equal_learning_rates_conclusion}, 
\cref{eq:varTheta_and_Theta_coincide}, 
and \cref{eq:not_equal_to_empty_set}
shows that 
\begin{equation}
\begin{split}
\textstyle
  \emptyset 
& 
\neq 
  \bigl(
    \cap_{ m = \cN }^{ \infty }
    \{ 
      \gamma_m = \gamma_{ \cN }
    \} 
  \bigr)
  \cap 
  \bigl( 
    \cap_{ w = 1 }^{ \infty } 
    \{ 
      \supertau_{ \nu_w, \, \cN } < \infty 
    \} 
  \bigr)
\\ &
\subseteq 
  \bigl(
  \bigl\{ 
    \cT < \infty 
  \bigr\} 
  \cap 
  \bigl\{ 
    \fL( 
      \cN + \cT, 
      \varTheta^{ \gamma_{ \cN }, \, \cN }_{ \cT } 
    )
    >
    \fL( 
      \cN + \cT, 
      \varTheta^{ \gamma_{ \cN }, \, \cN }_{ \cT - 1 } 
    )
  \bigr\}
  \cap 
  \bigl(
    \cap_{ m = \cN }^{ \infty }
    \{ 
      \gamma_m = \gamma_{ \cN }
    \} 
  \bigr)
  \bigr)
\\ &
\subseteq 
  \bigl(
  \bigl\{ 
    \cT < \infty 
  \bigr\} 
  \cap 
  \bigl\{ 
    \fL( 
      \cN + \cT, 
      \varTheta^{ \gamma_{ \cN }, \, \cN }_{ \cT } 
    )
    >
    \fL( 
      \cN + \cT, 
      \varTheta^{ \gamma_{ \cN }, \, \cN }_{ \cT - 1 } 
    )
  \bigr\}
\\ &
\quad 
  \cap 
  \bigl[
  \textstyle
    \cap_{ m = 0 }^{ \infty }
    \bigl\{
      \varTheta^{ \gamma_{ \cN }, \, \cN }_m
      =
      \Theta_{ \cN + m }
    \bigr\}
  \bigr]
  \cap 
  \bigl[
    \cap_{ m = \cN }^{ \infty }
    \{ 
      \gamma_m = \gamma_{ \cN }
    \} 
  \bigr]
  \bigr)
\\ &
\subseteq 
  \bigl(
  \bigl\{ 
    \cT < \infty 
  \bigr\} 
  \cap 
  \bigl\{ 
    \fL( 
      \cN + \cT, 
      \Theta_{ \cN + \cT } 
    )
    >
    \fL( 
      \cN + \cT , 
      \Theta_{ \cN + \cT - 1 } 
    )
  \bigr\}
\\ &
\quad 
  \cap 
  \bigl[
  \textstyle
    \cap_{ m = 0 }^{ \infty }
    \bigl\{
      \varTheta^{ \gamma_{ \cN }, \, \cN }_m
      =
      \Theta_{ \cN + m }
    \bigr\}
  \bigr]
  \cap 
    \bigl\{ 
      \gamma_{ 
        \cN + \cT + 1
      }
      =
      \gamma_{ 
        \cN + \cT 
      }
    \bigr\} 
  \bigr)
\\ &
\subseteq 
  \bigl(
  \bigl\{ 
    \cT < \infty 
  \bigr\} 
  \cap 
  \bigl\{ 
    \fL( 
      \cN + \cT, 
      \Theta_{ \cN + \cT } 
    )
    >
    \fL( 
      \cN + \cT , 
      \Theta_{ \cN + \cT - 1 } 
    )
  \bigr\}
  \cap 
    \bigl\{ 
      \gamma_{ 
        \cN + \cT + 1
      }
      =
      \gamma_{ 
        \cN + \cT
      }
    \bigr\} 
  \bigr)
\\ &
=
  \bigl(
  \bigl\{ 
    \cT < \infty 
  \bigr\} 
  \cap 
  \bigl\{ 
    \fL( 
      \cN + \cT, 
      \Theta_{ \cN + \cT } 
    )
    >
    \fL( 
      \cN + \cT , 
      \Theta_{ \cN + \cT - 1 } 
    )
  \bigr\}
\\ &
\quad 
  \cap 
    \bigl\{ 
      \fL( 
        \cN + \cT, 
        \Theta_{ \cN + \cT } 
      )
      \leq 
      \fL( 
        \cN + \cT, 
        \Theta_{ \cN + \cT - 1 } 
      )
    \bigr\} 
  \bigr)
\\ & =
  \bigl\{ 
    \cT < \infty 
  \bigr\} 
  \cap 
  \bigl\{ 
    \fL( 
      \cN + \cT , 
      \Theta_{ \cN + \cT - 1 } 
    )
    <
    \fL( 
      \cN + \cT, 
      \Theta_{ \cN + \cT } 
    )
    \leq 
    \fL( 
      \cN + \cT , 
      \Theta_{ \cN + \cT - 1 } 
    )
  \bigr\}
\\ & 
=
  \emptyset 
  .
\end{split}
\end{equation}
This contradiction proves \cref{eq:gamma_n_converge_to_zero}. 
Combining 
\cref{eq:gamma_n_converge_to_zero}
with \cref{eq:filtered_probability_space}, 
\cref{eq:gamma_monotonicity}, 
\cref{eq:gamma_n_summability}, 
the fact that 
$ \Theta \colon \N_0 \times \Omega \to \R $ 
is an $ ( \mathbb{F}_n )_{ n \in \N_0 } $-adapted stochastic process, 
the fact that for all $ n \in \N $ it holds that 
$ \gamma_n $ is $ \mathbb{F}_{ n - 1 } $-measurable, 
and the fact that for all $ n \in \N $ 
it holds that 
$
  \sigma( ( X_{ n, m } )_{ m \in \{ 1, 2, \dots, M \} } )
$
and 
$
  \mathbb{F}_{ n - 1 }
$
are independent
enables us to apply \cref{cor:convergence_in_probability_ex}
to obtain that
\begin{equation}
  \limsup_{ t \to \infty }
  \E\bigl[ 
    \| \Theta_{ N_t } - \E[ X_{ 1, 1 } ] \|^2
  \bigr]
  = 0
  .
\end{equation}
This establishes \cref{eq:convergence_for_quadratic_optimization_problems_with_adaptive_learning_rates}. 
\end{cproof}

\subsection*{Acknowledgments}

This work has been partially funded by the European Union (ERC, MONTECARLO, 101045811). 
The views and the opinions expressed in this work are however those of the authors only 
and do not necessarily reflect those of the European Union 
or the European Research Council (ERC). Neither the European Union nor 
the granting authority can be held responsible for them. 
In addition, this work has been partially funded by the 
Deutsche Forschungsgemeinschaft (DFG, German Research Foundation) 
under Germany's Excellence Strategy EXC 2044-390685587, 
Mathematics Münster: Dynamics-Geometry-Structure.

\bibliographystyle{abbrvurl}
\bibliography{bibfile}

\end{document}